**The work of Robert Langlands**

**James G. Arthur**

Robert Langlands was awarded the Abel Prize on May 22, 2018. This summary of his work will appear in **The Abel Laureates 2018–2022**, edited by Helge Holden and Ragni Piene, and to be published by Springer in 2024.

# The work of Robert Langlands

James G. Arthur

## Foreword

A more accurate title might have been *On the Work of Robert Langlands in Representation Theory, Automorphic Forms, Number Theory and Arithmetic Geometry*. For I have left out a significant part of Langlands' work, his papers in percolation theory and in mathematical physics, published in the years 1988–2000. I have however included a brief description of his recent work on the geometric theory. It occurs near the end of §11, the last section otherwise devoted to Beyond Endoscopy.

There is more than enough to discuss in the subjects I will be considering. I hope to be able to communicate the remarkable continuity that runs throughout all of the work of Langlands, with its roots in several fundamental areas of mathematics. What is now known as the Langlands program represents a unification of some of the deepest parts of these areas.

For example, it has been suggested that the Langlands program is ultimately a theory of $L$-functions, with its roots in analytic number theory. With this interpretation it goes back to Euler. Others might see the Langlands program at its most striking in its discovery of the long sought reciprocity laws for nonabelian class field theory, a culmination of perhaps two centuries of study of algebraic number theory. This of course goes back to Gauss, and his law of quadratic reciprocity. Both interpretations are equally valid, and they come together in the work from the 1920s and 1930s of Emil Artin. But they are by no means the full story. The origins of the subject for Langlands, and the power he was able to bring to it right from the beginning, came from harmonic analysis and group representations, specifically the work of Harish-Chandra. He was the first mathematician to gain a deep understand-

J.G. Arthur
Department of Mathematics, University of Toronto,
40 St. George Street,
Toronto, Ontario,
Canada M5S 2E4,
e-mail: arthur@math.utoronto.ca

ing of the fundamental contributions of Harish-Chandra to the representation theory of semisimple/reductive Lie groups. This analytic side of the subject has continued to inform the Langlands program right up to the present time. I note, for no particular reason except perhaps for some further sense of unity, that Harish-Chandra was a graduate student of Paul Dirac at Cambridge. However, his early career as a physicist ended when he switched completely to mathematics a few years later.

My hope has been to bring the work of Langlands to a more general mathematical audience. In the attempt to emphasize the continuity of the work, I have tried to write the report as a narrative that evolves with time. This entails fewer statements of formal theorems, and more efforts to describe underlying ideas. It also includes more repetition than would a formal paper. The fundamental ideas of Langlands occur again and again, often in different guises. Seeing them appear in this way might give a broader sense of the symmetry of the subject.

For example, Section 2 on Langlands' fundamental manuscript on Eisenstein series is certainly among the more technical sides of the report. We then follow it in Section 3 with a relatively leisurely introduction to class field theory. In general, I hope that a nonspecialist reader will be encouraged by the more elementary parts, and initially at least, not feel the need to take the more difficult passages as seriously.

A reader might also find it helpful to consult Langlands' later commentary on his various papers, to be found in the different sections of his website. I have certainly gained insight from it in the preparation of this report. I should add that I know some of the papers better than others, and I apologize in advance for any misstatements in my attempts to make this deep and fundamental work as accessible as I can.

# Contents

# 1 Group representations and harmonic analysis

As we have noted in the Foreword above, Langlands' early work came from group representations and harmonic analysis. These areas have remained at the heart of much of what is now known as the Langlands program. The analytic power in their methods has been indispensable in many of Langlands' greatest discoveries.

The groups in question are (unimodular) locally compact groups $H$. A unitary representation of $H$ is a (weakly continuous) homomorphism

$$R \colon H \to \mathrm{U}(\mathscr{H})$$

from $H$ to the group of unitary operators on a Hilbert space $\mathscr{H}$. For example, one could take $H$ to be a product

$$H = G \times G$$

of a group $G$ with itself, and $\mathscr{H}$ to be the Hilbert space $L^2(G)$ of square integrable complex valued functions on $G$ with respect to a Haar measure. One then has the regular representation

$$(R_H(y_1, y_2)\phi)(x) = \phi(y_1^{-1} x y_2), \quad \phi \in L^2(G), \quad x, y_1, y_2 \in G,$$

of $G \times G$ on $\mathscr{H}$. Another broad example is associated with a discrete subgroup $\Gamma$ of $H = G$, in which one assumes that the quotient space $\Gamma \backslash H$ of right cosets has finite volume with respect to an $H$-invariant measure. In this case, one has the unitary representation

$$(R_\Gamma(y)\phi)(x) = \phi(xy), \quad \phi \in L^2(\Gamma \backslash H), \quad x \in \Gamma \backslash H, \quad y \in H,$$

of $H$ by right translation on $\mathscr{H} = L^2(\Gamma \backslash H)$.

A representation $\pi$ of $H$ on a Hilbert space $V$ is *irreducible* if $V$ has no closed, $\pi$-invariant subspaces other than $\{0\}$ and $V$. Recall also that two unitary representations $(\pi, V)$ and $(\pi', V')$ of $H$ are *(unitarily) equivalent* if

$$\pi'(y) = U\pi(y)U^{-1}, \quad y \in H,$$

for a unitary, linear, intertwining isomorphism $U$ from $V$ to $V'$. For any given $H$, one would like to classify $\Pi_{\mathrm{unit}}(H)$, the set of equivalence classes of irreducible unitary representations of $H$. This can be regarded as the fundamental problem in group representations.

The fundamental problem for harmonic analysis would apply to any natural unitary representation $R$ of $H$, such as $R = R_H$ or $R = R_\Gamma$ as above. It is to find an *explicit* decomposition of $R$ into irreducible representations. This presupposes a knowledge of $\Pi_{\mathrm{unit}}(H)$, or at least a subset of $\Pi_{\mathrm{unit}}(H)$ that is large enough to support the measure class of $\Pi_{\mathrm{unit}}(H)$ that governs the decomposition of $R$.

For example, in the special case of the additive group $H = \mathbb{R} \oplus \mathbb{R}$, the irreducible unitary representations are the one-dimensional representations

$$\pi(y_1, y_2) = e^{\lambda_1 y_1} e^{\lambda_2 y_2}, \quad (y_1, y_2) \in \mathbb{R} \oplus \mathbb{R},$$

where $(\lambda_1, \lambda_2)$ ranges over the points in the imaginary space $i\mathbb{R} \oplus i\mathbb{R}$. The decomposition of $R_H$ is just classical harmonic analysis, the *Fourier transform* $\phi \to \widehat{\phi}$ being an explicit unitary isomorphism from $L^2(\mathbb{R})$ onto $L^2(i\mathbb{R})$ such that

$$(R_H(y_1, y_2)\phi)^\wedge(\lambda) = e^{-\lambda_1 y_1} \widehat{\phi} e^{\lambda_2}, \quad \lambda \in i\mathbb{R}.$$

From the other broad example, consider the special case that $H = \mathbb{R}$ and $\Gamma = \mathbb{Z}$. In this case, the map that assigns Fourier coefficients to functions is an explicit unitary isomorphism from $L^2(\mathbb{Z} \setminus \mathbb{R})$ onto the Hilbert space $L^2(2\pi i\mathbb{Z})$ of functions on the subset of irreducible representations of $\mathbb{R}$ that occur in the decomposition of $\mathbb{R}$.

The mathematical area Langlands entered in 1960 was considerably more elaborate than these basic examples would suggest. In particular, the groups $H$ were nonabelian, which meant that the irreducible unitary representations $\pi \in \Pi_{\mathrm{unit}}(H)$ were typically infinite dimensional. One consequence of this for harmonic analysis was that the decompositions of representations $R_H$ and $R_\Gamma$ typically had both a continuous part, qualitatively like the theory of Fourier transforms, and a discrete part, like the theory of Fourier series. By 1960, the theory was already rich and sophisticated, thanks in large measure to the ongoing efforts of Harish-Chandra.

From the beginning, Harish-Chandra had limited his efforts to semisimple Lie groups, such as the special linear groups $\mathrm{SL}(n, \mathbb{R})$, the special orthogonal groups $\mathrm{SO}(p, q, \mathbb{R})$ and the symplectic groups $\mathrm{Sp}(2n, \mathbb{R})$. In contrast to abstract locally compact groups, semisimple Lie groups have a rich structure, which gives rise to an even richer structure for their representations. Harish-Chandra's goal was to establish the Plancherel formula for any such $G$. It includes the problem of the explicit decomposition of the regular representation $R_{G \times G}$. The problem is actually a little more precise. A solution would in fact include a natural measure, the Plancherel measure, within the measure class that gives the decomposition into irreducible representations. Harish-Chandra ultimately established the Plancherel formula around 1975, but by 1960, he was well on his way to constructing the *discrete series*[1] for $G$. These are irreducible representations $\pi \in \Pi_{\mathrm{unit}}(G)$ of $G$ such that the products

$$\pi^\vee \otimes \pi \to (y_1, y_2) = {}^t\pi(y_1^{-1}) \otimes \pi(y_2)$$

give the representations that occur discretely in the decomposition of $R_{G \times G}$. His construction of the discrete series was completed around 1965, and is among Harish-Chandra's greatest achievements.

In 1960, Harish-Chandra's papers were regarded by many as being simply too difficult for anyone to read. Nonetheless, Langlands set about doing a comprehensive study of Harish-Chandra's work. Within a couple of years, he had acquired a mastery of at least part of it, including the nascent discrete series. This was demonstrated widely in three remarkable contributions to the 1965 AMS Summer Sympo-

---

[1] Harish-Chandra no doubt appropriated this term from Valentine Bargmann, the physicist who classified the irreducible unitary representations of the group $\mathrm{SL}(2, \mathbb{R})$ in 1947

sium in Boulder, Colorado, a broad conference designed to assess the general state of the subject.

Before I discuss Langlands' Boulder contributions, let me mention a later paper, because it bears directly on the work of Harish-Chandra. There are a couple of points to be made first. Harish-Chandra had found that, paradoxically, it was more natural to study the set $\Pi(G)$ of equivalence classes of *all* irreducible representations of $G$, rather than just the unitary ones. To be sure, his Plancherel measure was to be supported on the subset $\Pi_{\text{temp}}(G)$ of tempered representations that he would introduce in 1966 [88], and these are all unitary. But much of his developing analytic power required a command of the full set $\Pi(G)$. Another curious fact is that there are interesting unitary representations $\pi \in \Pi_{\text{unit}}(G)$ that do not lie in $\Pi_{\text{temp}}(G)$. This contradicts our intuition from classical Fourier analysis, in which the irreducible representations of $\mathbb{R}$ are one-dimensional quasi-characters

$$x \to e^{-\lambda x}, \quad x \in \mathbb{R}, \lambda \in \mathbb{C},$$

while both the unitary and the tempered representations coincide with the set of characters on $\mathbb{R}$, in which $\lambda$ is purely imaginary.

In 1973, Langlands gave a classification of the full set of irreducible representations $\Pi(G)$ [151], modulo a knowledge of the subset $\Pi_{\text{temp}}(G)$ of tempered representations. By that time, the set $\Pi_{\text{temp}}(G)$ had been classified by Harish-Chandra up to a set of Plancherel measure 0. The remaining singular representations in $\Pi_{\text{temp}}(G)$ were classified soon afterwards by Knapp and Zuckerman [114] [115], which meant that the Langlands classification then gave all the irreducible representations. The setting in Langlands' paper was actually slightly different. He replaced Harish-Chandra's semisimple Lie groups by groups $G(\mathbb{R})$ of real points, in which $G$ here represents a *reductive algebraic group* over $\mathbb{R}$. (A reductive algebraic group is a finite quotient of the product of a semisimple algebraic group with an algebraic torus.) The Langlands classification was soon extended to $p$-adic groups, in the weaker sense that it classifies all irreducible representations $\pi \in \Pi(G(\mathbb{Q}_p))$ in terms of representations that are tempered. (The tempered representations $\Pi_{\text{temp}}(G(\mathbb{Q}_p))$ of a $p$-adic group are another story, which we will come to in due course.) It thus pertains to all the local ingredients of general automorphic representations. For obvious reasons, the Langlands classification has become very influential. We shall return to it in Section 10, as a foundation for Langlands' later theory of endoscopy.

The three earlier contributions of Langlands to the Boulder proceedings were all striking, but one of them, on the theory of Eisenstein series [134], dominates the others in depth and importance. We leave it until the next section, and describe the other two here.

All three of Langlands' Boulder articles pertain to the quotients $\Gamma \backslash H$, which come with the associated representations $R_\Gamma$ of $H$ on $L^2(\Gamma \backslash H)$. The paper [135] concerns the space $\Gamma \backslash H$ itself. It applies to the case

$$\Gamma = G(\mathbb{Z}) \subset G(\mathbb{R}) = H$$

for a split (Chevalley) algebraic group $G$ over $\mathbb{Q}$, such as for example the group $G = \mathrm{SL}(n)$. The problem arose from the work of Tamagawa and Weil on the volume of $\Gamma \setminus H$, with respect to a canonical measure obtained from what is known as the Tamagawa measure.

Langlands established an explicit formula for the volume of $\Gamma \setminus H$ in this case, in terms of special values of the Riemann zeta function. His short proof was an ingenious combination of an interesting contour integral motivated by his theory of Eisenstein series, a well known real variable integral formula of Gindikin and Karpelevich (which had been used by Harish-Chandra in a different context), and some standard properties of Chevalley groups. At the time there was also a famous conjecture for the volume attached to any $G$ [251], which Weil had formulated in adelic form for simply connected groups. Langlands' paper confirmed Weil's conjecture in the special case of split groups. This was extended to quasi-split groups in 1972 by K. F. Lai [130], using the method of Langlands. Then later, following the suggestion of Jacquet and Langlands on p. 525 of [103], Kottwitz [122] extended Lai's result to arbitrary $G$, using a simple form [22, Corollary 23.6] of the general trace formula.[2] More precisely, Kottwitz extended the result to any $G$ that satisfies the Hasse principle, which was known at the time for any group without factors of type $E_8$. The Hasse principle was later established for that last case by Chernousov [49]. Langlands' Boulder paper [135] thus became the foundation for a general proof of Weil's conjecture.

The second Boulder paper [133] to discuss here does concern the representation $R_\Gamma$. Langlands used it to introduce a version in higher rank of the Selberg trace formula for compact quotient.

Suppose that $H = G$ is a semisimple Lie group, and that $\Gamma$ is a discrete subgroup with $\Gamma \setminus G$ compact. Among other things, Selberg introduced a formula that could be applied to the representation $R_\Gamma$ in this case. It is an identity

$$\sum_{\{\gamma\}} \mathrm{vol}(\Gamma_\gamma \setminus G_\gamma) \int_{G_\gamma \setminus G} f(x^{-1}\gamma x)\,dx = \sum_{\{\pi\}} \mathrm{mult}_\Gamma(\pi)\,\mathrm{tr}(\pi(f)), \qquad (1)$$

in which $f \in C_c^\infty(G)$ is to be regarded as a general test function on $G$. On the left hand side, $\{\gamma\}$ stands for a set of representations of conjugacy classes in $\Gamma$, $\Gamma_\gamma$ is the centralizer of $\gamma$ in $\Gamma$, $G_\gamma$ is the centralizer of $\gamma$ in $G$, and $\mathrm{vol}(\Gamma_\gamma \setminus G_\gamma)$ is the volume of the quotient $\Gamma_\gamma \setminus G_\gamma$ with respect to the right invariant measure defined by a fixed Haar measure on $G_\gamma$. The integral over $G_\gamma \setminus G$ is taken with respect to the quotient of a fixed Haar measure on $G$ by the chosen measure on $G_\gamma$. On the right hand side $\{\pi\}$ ranges over $\Pi_{\mathrm{unit}}(G)$, $\mathrm{mult}_\Gamma(\pi)$ is the multiplicity (a finite nonnegative integer) with which $\pi$ occurs discretely in the irreducible decomposition of $L^2(\Gamma \setminus G)$, and $\mathrm{tr}(\pi(f))$ is the trace of the operator

$$\pi(f) = \int_G f(x)\pi(x)\,dx$$

---

[2] For simplicity we shall often restrict references for the general trace formula to this introductory survey. The reader can then consult the original articles listed there, as needed.

on the Hilbert space $V$ on which $\pi$ acts. The left hand side is often called the *geometric side*, since the objects $\{\gamma\}$ have natural geometric interpretations. The right hand side is called the *spectral side*, since the coefficients $\text{mult}_\Gamma(\pi)$ concern spectral data in the decomposition of $R_\Gamma$. We note that if $G$ is a finite group, the formula becomes the well known theorem of Frobenius reciprocity, or rather the special case that applies to the trivial one-dimensional representation of the subgroup $\Gamma \subset G$.

In [133], Langlands derived (1) from first principles. He followed the original argument of Selberg, but the form of (1) is in some sense new. It differs from that of Selberg in that it reflects the theory of group representations, and in particular, the work of Harish-Chandra. Langlands then posed the question of using (1) to derive an explicit formula for the multiplicity $\text{mult}_\Gamma(\pi)$ of $\pi$. One does not expect a closed formula for all $\pi$. Langlands was asking about the case that $\pi$ lies in the subset $\Pi_2(G)$ of discrete series. Incidentally, the subscript 2 here means "square integrable", in the sense that every matrix coefficient

$$x \to (\pi(x)\phi, \psi), \quad \phi, \psi \in V, \tag{2}$$

of $\pi$ is a square integrable function of $x \in G$. Harish-Chandra had earlier noted that $\pi$ belongs to the discrete series if and only if it is square integrable. Armed with the simple trace formula (1), Langlands proposed letting the test function $f$ be a matrix coefficient of $\pi$ as above (with $\phi, \psi \neq 0$). The problem with this, however, was that $f$ is not compactly supported. In particular, the integrals in (1) need not converge. Langlands added the condition that $\pi$ lie in the smaller subset $\Pi_1(G)$ of *integrable* discrete series. With this assumption, the integrals do converge, and he was able to use the work of Harish-Chandra to compute the terms in (1) explicitly. He thereby obtained a simple, explicit formula for $\text{mult}_\Gamma(\pi)$.

As a supplement, Langlands added a conjectural interpretation of his multiplicities $\text{mult}_\Gamma(\pi)$ in terms of the cohomology of complex vector bundles. It was a generalization of the Borel–Weil formula for compact groups, suggested perhaps by a later proof of Kostant. It was very appealing at the time, especially to mathematicians with a background in complex analysis, and for a few years was sometimes known as "The Langlands Conjecture". This was of course before the sweeping conjectures that evolved into the Langlands program. The original Langlands conjecture from [133] was established a few years later by Wilfried Schmid [201].

Langlands' formula for $\text{mult}_\Gamma(\gamma)$ is an explicit, finite linear combination of terms. These were attached to the values of the character of $\pi$ at regular, elliptic elements, namely points $\gamma$ at which the centralizer $G_\gamma$ is a compact abelian subgroup of $G$. We shall not pause here to describe the formula precisely, or to recall how Harish-Chandra was able to construct the character of the infinite dimensional representation $\pi$ as a locally integrable class function on $G$. We note only that the formula obtained by Langlands, simple as it may be, is part of the foundation of the ongoing comparison of spectral data in representation theory with arithmetic spectral data in Shimura varieties. To be sure, it is the adelic version of the formula with its enrichment by Hecke operators that is relevant today, as well as a further stabilization of the formula. But it is still interesting to think that this early result of Langlands would have an implicit role in the study of Shimura varieties he began ten years later. We shall discuss these matters in Section 8.

It is also interesting that we now have a proper chain

$$\Pi_1(G) \subset \Pi_2(G) \subset \Pi_{\text{temp}}(G) \subset \Pi_{\text{unit}}(G) \subset \Pi(G) \tag{3}$$

of sets of (equivalence classes of) irreducible representations of $G$. Each of these families has its own role to play in some aspect of the theory. It is not hard to describe $\Pi_1(G)$ explicitly as a proper subset of $\Pi_2(G)$ in terms of the parametrization Harish-Chandra gave to the discrete series, a problem that arose with the publication of Langlands' paper, and that was solved shortly thereafter.

Langlands' paper on multiplicities provides a good introduction to the trace formula (1) for compact quotient. There is now a trace formula for general arithmetic quotients $\Gamma \setminus G$, such as

$$\Gamma \setminus G = \text{SL}(n, \mathbb{Z}) \setminus \text{SL}(n, \mathbb{R}),$$

which are typically noncompact (see [22]). This is much more difficult to establish, owing to the rather severe singularities at the boundary (which manifest themselves analytically as badly divergent integrals). However, the trace formula has been central to the subject. It represents an essential part of the work of Langlands, as a force behind proofs of fundamental theorems such as those for inner twists of $\text{GL}(2)$ [103], base change for $\text{GL}(2)$ [149], $L$-indistinguishably for $\widehat{\text{SL}}(2)$ [127] and the cohomology of Shimura varieties [140], as well as foundation of broader theories, such as Endoscopy [150] and Beyond Endoscopy [155] that remain conjectural.

## 2 Eisenstein series

The remaining Boulder article [134] of Langlands was a concise summary (with some supplementary ideas that were later used [11] in the development of a general trace formula) of his unpublished 270-page manuscript *On the Functional Equations Satisfied by Eisenstein Series*. The manuscript was later published, with four supplementary Appendices, as the monograph [141]. It was well ahead of its time, and was described in the preface of [141] of having been "almost impenetrable". Harish-Chandra temporarily suspended his work on the Plancherel formula to study the manuscript. The result was a set of expository lecture notes [89], which were a little closer to his own perhaps more familiar style, but which did not contain the most difficult part of the manuscript, the final Chapter 7. Langlands' theory of Eisenstein series has gradually become more widely understood. There are now a number of expositions of varying length, the most comprehensive being the monograph [178] of Mœglin and Waldspurger. Our aim here is to discuss some of the background and the content of Langlands' theory.

Eisenstein series have to do with the spectral decomposition of a space $L^2(\Gamma \setminus G)$. In general, there is an orthogonal decomposition

$$L^2(\Gamma \setminus G) = L^2_{\mathrm{disc}}(\Gamma \setminus G) \oplus L^2_{\mathrm{cont}}(\Gamma \setminus G)$$

where $L^2_{\mathrm{cont}}(\Gamma \setminus G)$ is an $R_\Gamma$-invariant subspace of $L^2(\Gamma \setminus G)$ that decomposes into a continuous direct sum (sometimes known as a direct integral), and $L^2_{\mathrm{disc}}(\Gamma \setminus G)$ is a subspace that decomposes discretely. Langlands' manuscript provided an explicit description of the continuous part of $R_\Gamma$. More precisely, it gave a decomposition of the continuous spectrum $L^2_{\mathrm{cont}}(\Gamma \setminus G)$ in terms of discrete spectra $L^2_{\mathrm{disc}}(M \cap \Gamma \setminus M)$, for a finite set of proper subgroups $M$ of $G$. Before we recall how this works, we should review the general setting of Langlands' work, and the way it is usually formulated today.

Langlands took $G$ to be a semisimple Lie group and $\Gamma$ to be a discrete subgroup such that the quotient $\Gamma \setminus G$ satisfies some general axioms. By far the most important case is that of an arithmetic quotient $\Gamma(N) \setminus G$, where $G$ is now a *reductive* algebraic group over $\mathbb{Q}$, with a $\mathbb{Z}$-scheme structure (such as $\mathrm{SL}(n)$), and

$$\Gamma = \Gamma(N) = \{ \gamma \in G(\mathbb{Z}) : \gamma \equiv 1 \, (\mathrm{mod}\, N) \}$$

is a principal congruence subgroup of $G(\mathbb{Z})$. The right regular representation of $G(\mathbb{R})$ on $L^2(\Gamma(N) \setminus G(\mathbb{R}))$ is best studied in the modern adelic formulation, which we pause briefly to review. (The reader can also refer to Langlands' adelic reformulation of his results in Appendix II of [141].)

The real field $\mathbb{R}$ is the completion of $\mathbb{Q}$ with respect to the usual archimedean absolute value $|\cdot| = |\cdot|_\infty$. For every prime number $p$, there is also a $p$-adic absolute value $|\cdot|_p$ on $\mathbb{Q}$, defined by setting $|x|_p = p^{-r}$ if

$$x = (ab^{-1})p^r,$$

for integers $a$, $b$ and $r$ with $(a, p) = (b, p) = 1$, and $|x|_p = 0$ if $x = 0$. Like $\mathbb{R}$, its completion $\mathbb{Q}_p$ is a locally compact field in which $\mathbb{Q}$ embeds as a dense subfield, and to which $|\cdot|_p$ extends continuously. Unlike $\mathbb{R}$, $\mathbb{Q}_p$ has an open, compact subring

$$\mathbb{Z}_p = \{x_p \in \mathbb{Q}_p : |x|_p \leq 1\}.$$

The ring product of all these completions is no longer locally compact. Suppose however that

$$S \subset \{v\} = \{\infty\} \cup \{p \text{ prime}\}$$

indexes a finite subset of these valuations that contains the archimedean absolute value $|\cdot|_\infty$. The product

$$\widehat{\mathbb{Z}}^S = \prod_{p \notin S} \mathbb{Z}_p$$

is then a compact ring, while the larger product

$$\mathbb{A}_S = \left( \prod_{v \in S} \mathbb{Q}_v \right) \widehat{\mathbb{Z}}^S$$

is a locally compact ring. The topological direct limit

$$\mathbb{A} = \varinjlim_S \mathbb{A}_S$$

of adeles is therefore a locally compact ring. It contains the diagonal image of the field $\mathbb{Q}$ as a discrete, cocompact subring.

For the given algebraic group $G$ over $\mathbb{Q}$, we can form the group $G(\mathbb{A})$ of points with values in the adele ring $\mathbb{A} = \mathbb{A}_{\mathbb{Q}}$ of $\mathbb{Q}$. It is a locally compact group, which contains $G(\mathbb{Q})$ as a discrete subgroup. The pair

$$\Gamma = G(\mathbb{Q}) \subset H = G(\mathbb{A})$$

is thus an example of the kind of object we have been considering, and to which we can attach a Hilbert space $L^2(G(\mathbb{Q}) \setminus G(\mathbb{A}))$. This looks very different from the concrete Hilbert space $L^2(\Gamma(N) \setminus G(\mathbb{R}))$ above. It is not.

Suppose that $K = K^\infty$ is an open compact subgroup of the nonarchimedean part

$$G(\mathbb{A}^\infty) = \left\{ x = \prod_v x_v \in G(\mathbb{A}) : x_\infty = 1 \right\}$$

of $G(\mathbb{A})$. The product $G(\mathbb{R})K$ is then an open subgroup of $G(\mathbb{A})$. Under some natural conditions on $G$, the set of $(G(\mathbb{Q}), G(\mathbb{R})K)$ double cosets in $G(\mathbb{A})$ is then finite. Writing

$$G(\mathbb{A}) = \coprod_{i=1}^{n} \left( G(\mathbb{Q}) \cdot x^i \cdot G(\mathbb{R})K \right),$$

for elements $x^1 = 1, x^2, \ldots, x^n$ in $G(\mathbb{A}^\infty)$, we obtain a right $G(\mathbb{R})$-invariant decomposition

$$G(\mathbb{Q}) \setminus G(\mathbb{A})/K = \coprod_{i=1}^n G(\mathbb{Q}) \setminus \big(G(\mathbb{Q}) \cdot x^i \cdot G(\mathbb{R})K\big)/K$$

$$\cong \coprod_{i=1}^n \big(\Gamma^i \setminus G(\mathbb{R})\big),$$

for discrete subgroups

$$\Gamma^i = G(\mathbb{R}) \cap \big(G(\mathbb{Q}) \cdot x^i K (x^i)^{-1}\big)$$

of $G(\mathbb{R})$, and a $G(\mathbb{R})$-isomorphism of Hilbert spaces

$$L^2(G(\mathbb{Q}) \setminus G(\mathbb{A})/K) = \bigoplus_{i=1}^n \big(\Gamma^i \setminus G(\mathbb{R})\big). \tag{4}$$

Each discrete group $\Gamma^i$ is defined by congruence conditions, which are determined by the choice of $K$ and $x^i$. Each Hilbert space on the right hand side of (4) is thus a modest generalization of the space $L^2(\Gamma(N) \setminus G(\mathbb{R}))$ above. But we can see directly that it contains more information that just the regular representation of $G(\mathbb{R})$.

On one hand, the action of $G(\mathbb{R})$ on the spaces on either side of (4) corresponds to the action by right convolution on either space by functions in the algebra $C_c(G(\mathbb{R}))$. But there is a supplementary convolution algebra, the Hecke algebra $\mathcal{H}(G(\mathbb{A}^\infty), K)$ of compactly supported functions that are left and right invariant under translation by the group $K$. It also acts by right convolution on the left hand side of (4), in a way that clearly commutes with the action of $G(\mathbb{R})$. The corresponding action of $\mathcal{H}(G(\mathbb{A}^\infty), K)$ on the right hand side of (4) includes general analogues of the operators defined by Hecke on classical modular forms.

Hecke operators for general groups are at the heart of the theory. They contain the data that according to Langlands' later conjectures govern much of the arithmetic world. One builds them into the representation theory by setting

$$C_c^\infty(G(\mathbb{A})) = C_c(G(\mathbb{R})) \otimes C_c^\infty(G(\mathbb{A}^\infty))$$

where $C_c^\infty(G(\mathbb{A}^\infty))$ denotes the space of locally constant, complex valued functions of compact support on $G(\mathbb{A}^\infty)$. Any function in $C_c^\infty(G(\mathbb{A}))$ is bi-invariant under translation by an open compact subgroup $K = K^\infty$ of $G(\mathbb{A}^\infty)$, and therefore acts by right convolution on the corresponding space (4). When we vary $f$, and hence $K$, we are working with the understanding that the action of the convolution algebra $C_c^\infty(G(\mathbb{A}))$ on the Hilbert space $L^2(G(\mathbb{Q}) \setminus G(\mathbb{A}))$ is equivalent to the right regular representation $R_G$ of $G(\mathbb{A})$ on the space. This is the standard modern setting, which despite possible appearances, streamlines the joint study of the right regular representation of $G(\mathbb{R})$ and the underlying Hecke operators.

We would now like to describe Langlands' work on Eisenstein series. Following the statement in [22], we shall describe the main results in complete detail, making the rest of this section one of the more technical parts of our report. For a start, the results are formulated in terms of the basic structure of algebraic groups, which we will have to apply with only limited comments. We should first take care of the minor annoyance that $G(\mathbb{A})$ can have noncompact centre, which implies for trivial reasons that $L^2(G(\mathbb{Q}) \setminus G(\mathbb{A}))$ will have no discrete spectrum. One deals with it either by replacing $G(\mathbb{A})$ by the quotient $A_G(\mathbb{R})^0 \setminus G(\mathbb{A})$, in which $A_G$ is the $\mathbb{Q}$-split component of the centre $Z = Z(G)$ of $G$, or by the subgroup

$$G(\mathbb{A})^1 = \{x \in G(\mathbb{A}) : |\chi(x)| = 1, \quad \chi \in X(G)_{\mathbb{Q}}\},$$

in which $X(G)_{\mathbb{Q}}$ is the group of characters from $G$ to $\mathbb{G}_m = \mathrm{GL}(1)$ defined over $\mathbb{Q}$, and

$$|\chi(x)| = \prod_v |\chi(x_v)|_v$$

is the (adelic) absolute value on the group $\mathbb{A}^* = \mathrm{GL}(1, \mathbb{A})$ of ideles.

The two modifications are equivalent. It is clear that $G(\mathbb{Q})$ embeds as a discrete subgroup of either $G(\mathbb{A})^1$ or $A_G(\mathbb{R})^0 \setminus G(\mathbb{A})$. Moreover, there is a surjective homomorphism $H_G$ from the group $G(\mathbb{A})$ onto the real vector space

$$\mathfrak{a}_G = \mathrm{Hom}(X(G)_{\mathbb{Q}}, \mathbb{R}),$$

defined by setting

$$e^{\langle H_G(x), \chi \rangle} = |\chi(x)|, \quad x \in G(\mathbb{A}), \chi \in X(G)_{\mathbb{Q}},$$

whose kernel is $G(\mathbb{A})^1$, and whose restriction to the central subgroup $A_G(\mathbb{R})^0$ of $G(\mathbb{A})$ is an isomorphism onto $\mathfrak{a}_G$. It follows that

$$G(\mathbb{A}) = G(\mathbb{A})^1 \times A_G(\mathbb{R})^0.$$

In particular, there is a unitary $(G(\mathbb{A}) = G(\mathbb{A})^1 \times A_G(\mathbb{R})^0)$-equivariant isomorphism between the two Hilbert spaces

$$L^2(G(\mathbb{Q}) \setminus G(\mathbb{A})^1) \cong L^2(G(\mathbb{Q})A_G(\mathbb{R})^0 \setminus G(\mathbb{A})),$$

each of which will have nontrivial discrete spectra

$$\mathrm{L}^2_{\mathrm{disc}}(G(\mathbb{Q}) \setminus G(\mathbb{A})^1) \cong \mathrm{L}^2_{\mathrm{disc}}(G(\mathbb{Q})A_G(\mathbb{R})^0 \setminus G(\mathbb{A})), \tag{5}$$

(unless $G$ is a split torus over $\mathbb{Q}$). The spaces on the left are perhaps more natural. For among other things, we can identify the irreducible unitary representations $\pi \in \Pi_{\mathrm{unit}}(G)$ of $G(\mathbb{A})^1$ with the $i\mathfrak{a}_G^*$-orbits

$$\{\pi_\lambda(x) = \pi_0(x)e^{-\lambda(H_G(x))} : x \in G(\mathbb{A}), \lambda \in i\mathfrak{a}_G^*\}$$

of irreducible unitary representations of $G(\mathbb{A})$. The base point $\pi_0$ can be any irreducible unitary representation of $G(\mathbb{A})$ whose restriction of $G(\mathbb{A})^1$ is $\pi$. A similar convention holds if $\pi$ is replaced by any representation $R$ of $G(\mathbb{A})^1$, such as for example the representation $R_{G,\mathrm{disc}}$ of $G(\mathbb{A})^1$ on $\mathrm{L}^2_{\mathrm{disc}}(G(\mathbb{Q}) \setminus G(\mathbb{A})^1)$, and $\lambda$ is any point in the complex vector space $\mathfrak{a}^*_{G,\mathbb{C}} = \mathfrak{a}^*_G \otimes \mathbb{C}$, with the understanding that $R_0$ is a representation of $A_G(\mathbb{R})^0 \setminus G(\mathbb{A})$. These conventions are due to Harish-Chandra, and are quite natural, even if they might seem cumbersome at first.

We fix a minimal parabolic subgroup $P_0$ of $G$ over $\mathbb{Q}$, together with a Levi decomposition $P_0 = M_0 N_0$, where $M_0$ (resp $N_0$) is a reductive (resp. unipotent) subgroup of $P_0$ over $\mathbb{Q}$. We also fix a suitable maximal compact subgroup $K_0 = \prod_v K_v$ of $G(\mathbb{A})$ [12, p. 9] with $G(\mathbb{A}) = P_0(\mathbb{A}) K_0$. We shall then work in what remains of this section with the finite set of standard parabolic subgroups, namely the subgroups $P$ of $G$ which contain $P_0$. Any such $P$ has a unique Levi decomposition

$$P = M_P N_P,$$

in which the Levi component $M_P$ contains $M_0$. Since $P$ contains $P_0$, we then have a decomposition

$$G(\mathbb{A}) = P(\mathbb{A}) K_0 = N_P(\mathbb{A}) M_P(\mathbb{A}) K_0 = N_P(\mathbb{A}) M_P(\mathbb{A})^1 A_P(\mathbb{R})^0 K_0,$$

where $A_P = A_{M_P}$ in the notation above (with $M_P$ in place of $G$). This allows us to define a continuous mapping

$$H_P \colon G(\mathbb{A}) \to \mathfrak{a}_P, \quad \mathfrak{a}_P = \mathfrak{a}_{M_P},$$

by setting

$$H_P(nmk) = H_{M_P}(m), \quad n \in N_P(\mathbb{A}), m \in M_P(\mathbb{A}), k \in K,$$

where $H_{M_P}$ and $\mathfrak{a}_P = \mathfrak{a}_{M_P}$ are again as above. Finally, we fix Haar measures $dx$, $dn$, $dm$, $da$ and $dk$ on the groups $G(\mathbb{A})$, $N_P(\mathbb{A})$, $M_P(\mathbb{A})^1$ (or $M_P(\mathbb{A})$, depending on the context), $A_P(\mathbb{R})^0$ and $K$ such that for the decomposition above, we have

$$dx = e^{2\rho_P(H_P(a))} \, dn \, dm \, da \, dk.$$

Here $\rho_P$ is the familiar vector in $\mathfrak{a}_P$ such that $e^{2\rho_P(H_P(\cdot))}$ is the modular function on the (nonunimodular) group $P(\mathbb{A})$. The notation is again due to Harish-Chandra, and is convenient for working with representations of $G(\mathbb{A})$ induced from $P(\mathbb{A})$.

Suppose that $P$ is a standard parabolic subgroup of $G$, and that $\lambda$ lies in $\mathfrak{a}^*_{P,\mathbb{C}}$. We write

$$y \to \mathscr{I}_P(\lambda, y), \quad y \in G(\mathbb{A}),$$

for the induced representation

$$\mathrm{Ind}^{G(\mathbb{A})}_{P(\mathbb{A})}(I_{N_P(\mathbb{A})} \otimes R_{M_P,\mathrm{disc},\lambda})$$

of $G(\mathbb{A})$ obtained from $\lambda$ and the discrete spectrum of the reductive group $M_P$. This representation acts on the Hilbert space $\mathscr{H}_P$ of measurable functions

$$\phi\colon N_P(\mathbb{A})M_P(\mathbb{Q})A_P(\mathbb{R})^0 \setminus G(\mathbb{A}) \to \mathbb{C}$$

such that the function

$$\phi_x\colon m \to \phi(mx),$$

belongs to $L^2_{\mathrm{disc}}(M_P(\mathbb{Q}) \setminus M_P(\mathbb{A})^1)$ for almost all $x \in G(\mathbb{A})$, and such that

$$\|\phi\|^2 = \int_K \int_{M_P(\mathbb{Q}) \setminus M_P(\mathbb{A})^1} |\phi(mk)|^2 \, dm \, dk < \infty.$$

For any $y \in G(\mathbb{A})$, $\mathscr{I}_P(\lambda, y)$ maps a function $\phi \in \mathscr{H}_P$ to the function

$$(\mathscr{I}_P(\lambda, y)\phi)(x) = \phi(xy)e^{(\lambda+\rho_P)(H_P(xy))}e^{-(\lambda+\rho_P)(H_P(x))}.$$

We have put the twist by $\lambda$ into the operator $\mathscr{I}_P(\lambda, y)$, rather than the underlying Hilbert space $\mathscr{H}_P$, in order that $\mathscr{H}_P$ be independent of $\lambda$. The function $e^{\rho_P(H_P(\cdot))}$ is included in the definition in order that the representation $\mathscr{I}_P(\lambda)$ be unitary whenever the inducing representation is unitary, which is to say, whenever $\lambda$ belongs to the subset $i\mathfrak{a}_P^*$ of $\mathfrak{a}_{P,\mathbb{C}}^*$.

Suppose that

$$R_{M_P,\mathrm{disc}} \cong \bigoplus_\pi \pi \cong \bigoplus_\pi \left( \bigotimes_v \pi_v \right)$$

is the decomposition of $R_{M_P,\mathrm{disc}}$ into irreducible representations $\pi = \bigotimes_v \pi_v$ of $M_P(\mathbb{A})/A_P(\mathbb{R})^0$. The induced representation $\mathscr{I}_P(\lambda)$ then has a corresponding decomposition

$$\mathscr{I}_P(\lambda) \cong \bigoplus_\pi \mathscr{I}_P(\pi_\lambda) \cong \bigoplus_\pi \left( \bigotimes_v \mathscr{I}_P(\pi_{v,\lambda}) \right)$$

in terms of induced representations $\mathscr{I}_P(\pi_{v,\lambda})$ of local groups $G(\mathbb{Q}_v)$. This follows form the definition of induced representation , and the fact that

$$e^{\lambda(H_{M_P}(m))} = \prod_v e^{\lambda(H_{M_P}(m_v))},$$

for any point $m = \prod_v m_v$ in $M_P(\mathbb{A})$. If $\lambda \in i\mathfrak{a}_P^*$ is in general position, all of the induced representations $\mathscr{I}_P(\pi_{v,\lambda})$ are irreducible. Thus, if we understand the decomposition of the discrete spectrum of $M_P$ into irreducible representations of the local groups $M_P(\mathbb{Q}_v)$, we understand the decomposition of the generic induced representations $\mathscr{I}_P(\lambda)$ into irreducible representations of the local groups $G(\mathbb{Q}_v)$.

The aim of the theory of Eisenstein series is to construct intertwining operators between the induced representations $\mathscr{I}_P(\lambda)$ and the continuous part of the regular representation $R$ of $G(\mathbb{A})$. The problem includes being able to construct intertwining operators among the representations $\mathscr{I}_P(\lambda)$, as $P$ and $\lambda$ vary. The symmetries among pairs $(P, \lambda)$ are given by the Weyl sets $W(\mathfrak{a}_P, \mathfrak{a}_{P'})$ of Langlands. For a given pair $P$ and $P'$ of standard parabolic subgroups, $W(\mathfrak{a}_P, \mathfrak{a}_{P'})$ is defined as the set of distinct linear isomorphisms from $\mathfrak{a}_P \subset \mathfrak{a}_0$ onto $\mathfrak{a}_{P'} \subset \mathfrak{a}_0$ obtained by restrictions of elements in the (restricted) Weyl group

$$W_0 = \mathrm{Norm}(G, A_0) / \mathrm{Cent}(G, A_0).$$

Suppose for example that $G = \mathrm{GL}(n)$. If $P$ and $P'$ correspond to the partitions $(n_1, \ldots n_p)$ and $(n'_1, \ldots n'_{p'})$ of $n$, the set $W(\mathfrak{a}_P, \mathfrak{a}_{P'})$ is empty unless $p = p'$, in which case

$$W(\mathfrak{a}_P, \mathfrak{a}_{P'}) \cong \{s \in S_p : n'_i = n_{s(i)}, \ 1 \leq s \leq p\}.$$

In general, we say that $P$ and $P'$ are *associated* if the set $W(\mathfrak{a}_P, \mathfrak{a}_{P'})$ is nonempty. We would expect a pair of induced representations $\mathscr{I}_P(\lambda)$ and $\mathscr{I}_{P'}(\lambda')$ to be equivalent if $P$ and $P'$ belong to the same associated class, and $\lambda' = s\lambda$ for some element $s \in W(\mathfrak{a}_P, \mathfrak{a}_{P'})$.

The formal definitions apply to any elements $x \in G(\mathbb{A})$, $\phi \in \mathscr{H}_P$ and $\lambda \in \mathfrak{a}^*_{M,\mathbb{C}}$. The associated Eisenstein series is

$$E(x, \phi, \lambda) = \sum_{\delta \in P(\mathbb{Q}) \backslash G(\mathbb{Q})} \phi(\delta x) e^{(\lambda + \rho_P)(H_P(\delta x))}. \tag{6}$$

If $s$ belongs to $W(\mathfrak{a}_P, \mathfrak{a}_{P'})$, the operator

$$M(s, \lambda) \colon \mathscr{H}_P \to \mathscr{H}_{P'}$$

that intertwines $\mathscr{I}_P(\lambda)$ with $\mathscr{I}_{P'}(s\lambda)$ is defined by

$$M(s, \lambda)(x) = \int \phi(w_s^{-1} n x) e^{(\lambda + \rho_P)(H_P(w_s^{-1} n x))} e^{(-s\lambda + \rho_{P'})(H_{P'}(x))} \, dn, \tag{7}$$

where the integral is taken over the quotient

$$N_{P'(\mathbb{A})} \cap w_s N_P(\mathbb{A}) w_s^{-1} \backslash N_{P'}(\mathbb{A}),$$

and $w_s$ is any representative of $s$ in $G(\mathbb{Q})$. A reader so inclined could motivate both definitions in terms of finite group theory. Each definition is a formal analogue of a general construction of Mackey [169] for the space of intertwining operators between two induced representations $\mathrm{Ind}_{H_1}^H(\rho_1)$ and $\mathrm{Ind}_{H_2}^H(\rho_2)$ of a finite group $H$.

It follows formally from the definitions that

$$E(x, \mathscr{I}_P(\lambda, y)\phi, \lambda) = E(xy, \phi, \lambda)$$

and

$$M(s, \lambda)\mathscr{I}_P(\lambda, y) = \mathscr{I}_{P'}(s\lambda, y)M(s, \lambda).$$

These are the desired intertwining properties. However, (6) and (7) are defined by sums and integrals over noncompact spaces. They do not generally converge. It is this fact that makes the theory of Eisenstein series so difficult.

Let $\mathscr{H}_P^0$ be the subspace of vectors $\phi \in \mathscr{H}_P$ that are $K_0$-finite, in the sense that the subset

$$\{\mathscr{I}_P(\lambda, k)\phi : k \in K_0\}$$

of $\mathscr{H}_P$ spans a finite dimensional space, and that lie in a finite sum of irreducible subspaces of $\mathscr{H}_P$ under the action of $\mathscr{I}_P(\lambda)$ of $G(\mathbb{A})$. The two conditions do not depend on the choice of $\lambda$. Taken together, they are equivalent to the requirement that the function

$$\phi(x_\infty x^\infty), \quad x_\infty \in G(\mathbb{R}), x^\infty \in G(\mathbb{A}^\infty)$$

be locally constant in $x^\infty$, and smooth, $K_\mathbb{R}$-finite and $\mathscr{Z}_\infty$-finite in $x_\infty$, where $\mathscr{Z}_\infty = \mathscr{Z}_{G,\infty}$ denotes the algebra of bi-invariant differential operators on $G(\mathbb{R})$. The space $\mathscr{H}_P^0$ is dense in $\mathscr{H}_P$.

For any $P$, we can form the chamber

$$(\mathfrak{a}_P^*)^+ = \{\Lambda \in \mathfrak{a}_P^* : \Lambda(\alpha^\vee) > 0, \alpha \in \Delta_P\}$$

in $\mathfrak{a}_P^*$, in which $\Delta_P$ denotes the set of simple parabolic roots of $(P, A_P)$.

**Lemma  (Langlands)** *Suppose that $\phi \in \mathscr{H}_P^0$ and that $\lambda$ lies in the open subset*

$$\{\lambda \in \mathfrak{a}_{P,\mathbb{C}}^* : \mathrm{Re}(\lambda) \in \rho_P + (\mathfrak{a}_P^*)^+\}$$

*of $\mathfrak{a}_{P,\mathbb{C}}^*$. Then the sum* (6) *and the integral* (7) *that define $E(x, \phi, \lambda)$ and $(M(s, \lambda)\phi)(x)$ both converge absolutely to analytic functions of $\lambda$.*

For spectral theory, one is interested in points $\lambda$ such that $\mathscr{I}_P(\lambda)$ is unitary, which is to say that $\lambda$ belongs to the real subspace $i\mathfrak{a}_P^*$ of $\mathfrak{a}_{P,\mathbb{C}}^*$. This is outside the domain of absolute convergence for (6) and (7). The problem is to show that the functions $E(x, \phi, \lambda)$ and $M(s, \lambda)\phi$ have analytic continuation to this space. The following theorem summarizes Langlands' main results in Eisenstein series.

**Main Theorem  (Langlands)**

(a)    *Suppose that $\phi \in \mathscr{H}_P^0$. Then $E(x, \phi, \lambda)$ and $M(s, \lambda)\phi$ can be analytically continued to meromorphic functions of $\lambda \in \mathfrak{a}_{P,\mathbb{C}}^*$ that satisfy the functional equations*

$$E(x, M(s, \lambda)\phi, s\lambda) = E(x, \phi, \lambda) \tag{8}$$

*and*
$$M(ts, \lambda) = M(t, s\lambda)M(s, \lambda), \quad t \in W(\mathfrak{a}_P, \mathfrak{a}_{P'}). \tag{9}$$

*If $\lambda \in i\mathfrak{a}_P^*$, both $E(x, \phi, \lambda)$ and $M(s, \lambda)$ are analytic, and $M(s, \lambda)$ extends to a unitary operator from $\mathscr{H}_P$ to $\mathscr{H}_{P'}$.*

(b)    *Given an associated class $\mathscr{P} = \{P\}$, define $\widehat{L}_{\mathscr{P}}$ to be the Hilbert space of families of measurable functions*

$$F = \{F_P \colon i\mathfrak{a}_P^* \to \mathscr{H}_P, \quad P \in \mathscr{P}\}$$

*that satisfy the symmetry condition*

$$F_{P'}(s\lambda) = M(s,\lambda)F_P(\lambda), \quad s \in W(\mathfrak{a}_P, \mathfrak{a}_{P'}),$$

*and the finiteness condition*

$$\|F\|^2 = \sum_{P \in \mathscr{P}} n_P^{-1} \int_{i\mathfrak{a}_P^*} \|F_P(\lambda)\|^2 \, d\lambda < \infty,$$

*where*

$$n_P = \sum_{P' \in \mathscr{P}} |W(\mathfrak{a}_P, \mathfrak{a}_{P'})|$$

*for any $P \in \mathscr{P}$. Then the mapping that sends $F$ to the function*

$$\sum_{P \in \mathscr{P}} n_P^{-1} \int_{i\mathfrak{a}_P^*} E(x, F_P(\lambda), \lambda) \, d\lambda, \quad x \in G(\mathbb{A}),$$

*defined whenever $F_P(\lambda)$ is a smooth, compactly supported function of $\lambda$ with values in a finite dimensional subspace of $\mathscr{H}_P^0$, extends to a unitary mapping from $\widehat{L}_{\mathscr{P}}$ onto a closed $G(\mathbb{A})$-invariant subspace $L_{\mathscr{P}}^2(G(\mathbb{Q}) \setminus G(\mathbb{A}))$ of $L^2(G(\mathbb{Q}) \setminus G(\mathbb{A}))$. Moreover, the original space $L^2(G(\mathbb{Q}) \setminus G(\mathbb{A}))$ has an orthogonal direct sum decomposition*

$$L^2(G(\mathbb{Q}) \setminus G(\mathbb{A})) = \bigoplus_{\mathscr{P}} L_{\mathscr{P}}^2(G(\mathbb{Q}) \setminus G(\mathbb{A})). \tag{10}$$

The theorem gives a qualitative description of the decomposition of $R$. It provides a finite decomposition

$$R = \bigoplus_{\mathscr{P}} R_{\mathscr{P}},$$

where $R_{\mathscr{P}}$ is the restriction of $R$ to the invariant subspace $L_{\mathscr{P}}^2(G(\mathbb{Q}) \setminus G(\mathbb{A}))$ of $L^2(G(\mathbb{Q}) \setminus G(\mathbb{A}))$. It also provides a unitary intertwining operator from $R_{\mathscr{P}}$ onto the representation $\widehat{R}_{\mathscr{P}}$ of $G(\mathbb{A})$ on $\widehat{L}_{\mathscr{P}}$ defined by

$$(\widehat{R}_P(y)F)_P(\lambda) = \mathscr{I}_P(\lambda, y)F_P(\lambda), \quad F \in \widehat{L}_{\mathscr{P}}^2, P \in \mathscr{P}.$$

The theorem is thus compatible with the general intuition we retain from the theory of Fourier series and Fourier transforms.

In summary, let us say that while it might seem a little overwhelming at first, the theorem is quite comprehensive, and ultimately, remarkably simple. It is exactly what one might hope for in terms of an explicit decomposition of the continuous spectrum $L_{\text{cont}}^2(G(\mathbb{Q}) \setminus G(\mathbb{A}))$ into irreducible representations. On the other hand, the proof of the theorem is long and complex, to an extent that is hard to quantify in a few words. Langlands had to overcome many obstacles, the most severe being the analytic problems treated in Chapter 7 of [141]. This last chapter consists of a sophisticated residue scheme, designed to construct inaccessible constituents of discrete spectra $\mathrm{L}^2_{\mathrm{disc}}(M(\mathbb{Q}) \setminus M(\mathbb{A})^1)$ and their Eisenstein series from residues of cuspidal Eisenstein series.

We define a locally integrable function $\phi$ on $G(\mathbb{Q}) \setminus G(\mathbb{A})$ to be *cuspidal* if for every standard parabolic subgroup $P = MN$ distinct from $G$, the integral

$$\int_{N(\mathbb{Q}) \setminus N(\mathbb{A})} \phi(nx)\, dn, \quad x \in G(\mathbb{A}),$$

vanishes for almost all $x$. The main property of these functions is that they lie in the relative discrete spectrum [10], [141]. In other words, the $G(\mathbb{A})$-invariant subspace $\mathrm{L}^2_{\mathrm{cusp}}(G(\mathbb{Q}) \setminus G(\mathbb{A}))$ of cuspidal functions in $L^2(G(\mathbb{Q}) \setminus G(\mathbb{A}))$ lies in the summand $\mathrm{L}^2_G(G(\mathbb{Q}) \setminus G(\mathbb{A}))$ of the decomposition (10) attached to $\mathscr{P} = \{G\}$, which is the space of functions $\phi$ whose restriction to $G(\mathbb{A})^1$ lies in $\mathrm{L}^2_{\mathrm{disc}}(G(\mathbb{Q}) \setminus G(\mathbb{A})^1)$. For any $P$, let $\mathscr{H}^0_{P,\mathrm{cusp}}$ be the subspace of functions $\phi$ in $\mathscr{H}^0_P$ such that the function $\phi_x(m) = \phi(mx)$ defined above is cuspidal for almost all $x$. A *cuspidal Eisenstein series* is then a series (6) in which $\phi$ lies in the subspace $\mathscr{H}^0_{P,\mathrm{cusp}}$ of $\mathscr{H}^0_P$.

Langlands studied cuspidal Eisenstein series in the first six chapters of the volume. These objects were difficult enough, but there were some available techniques of Selberg, especially for the case of rank 1 (in which $\dim(A_P) = 1$). Langlands used these techniques, and others that he created. By the end of Chapter 6, he had established the analytic continuation and functional equations from Part (b) of the theorem, for cuspidal Eisenstein series.

Noncuspidal Eisenstein series, however, were a different matter. We now have a classification for the noncuspidal discrete spectrum for general linear groups [177], as well as a conjectural classification [18] in general, but it does not help us with their Eisenstein series. It was Langlands' indirect residue scheme that ultimately led to the required properties. His resolution was a very delicate interplay between the desired analytic continuation and the required spectral properties, all within an extended induction argument.

In general, the irreducible representations (typically induced) in the spectral decomposition of each space $\mathrm{L}^2_{\mathscr{P}}(G(\mathbb{Q}) \setminus G(\mathbb{A}))$ are called *automorphic representations*, (while the associated Eisenstein series are called *automorphic forms*). Let us write $\Pi_{\mathrm{temp}}(G)$ for the set of such representations, where the subscript stands for "globally tempered", in whatever sense this term would have if there were a global Schwartz space on $G(\mathbb{Q}) \setminus G(\mathbb{A})$. (It does *not* mean that the local constituents of representations are tempered). We can also write $\Pi_2(G)$ and $\Pi_1(G)$ for the subset of automorphic representations in the decompositions of the spaces $\mathrm{L}^2_G(G(\mathbb{Q}) \setminus G(\mathbb{A}))$ and $\mathrm{L}^2_{\mathrm{cusp}}(G(\mathbb{Q}) \setminus G(\mathbb{A}))$. Finally, we should say that the formal definition of automorphic representations [35, 145] gives a wider class of irreducible representations $\Pi(G)$ of $G(\mathbb{A})$, which of course includes the subset $\Pi_{\mathrm{unit}}(G)$ of unitary automorphic representations. We obtain a proper chain

$$\Pi_1(G) \subset \Pi_2(G) \subset \Pi_{\mathrm{temp}}(G) \subset \Pi_{\mathrm{unit}}(G) \subset \Pi(G) \tag{11}$$

of sets of irreducible automorphic representations of $G(\mathbb{A})$, with obvious analogy to the chain (3) for real groups. This global notation leaves us free to write $\Pi_{\mathrm{cusp}}(G)$ for the set of *all* cuspidal automorphic representations of $G(\mathbb{A})$, thereby allowing for nonunitary central characters. We then have

$$\Pi_1(G) = \Pi_{\mathrm{cusp}}(G) \cap \Pi_2(G) = \Pi_{\mathrm{cusp}}(G) \cap \Pi_{\mathrm{temp}}(G) = \Pi_{\mathrm{cusp}}(G) \cap \Pi_{\mathrm{unit}}(G).$$

The analogy between the local and global chains (3) and (10) is a little fanciful, but suggests analogies between the two kinds of representations. The local chain (3) was actually taken for a *semisimple* Lie group $G$, in which the centre is finite. At this point, we would take $G$ to be a Lie group that has been implicitly identified with the group of real points $G(\mathbb{R})$ of a *reductive* group over $\mathbb{R}$. The symbol $\Pi_2(G)$ in (3) would then stand for the *relative discrete series*, which is to say, tempered representations $\pi$ of $G(\mathbb{R})$ whose restrictions to $G(\mathbb{R})^1$ are in the discrete series. The equivalent term *square integrable* would then be understood to mean square integrable modulo the centre of $G$. This slight generalization was in fact Harish-Chandra's original definition.

We will not try to formulate an analogue of the chain (10) for spaces $\mathscr{A}(G)$ of automorphic forms. However, this is a good opportunity to say a few words about these objects, even if we do not recall their precise definition from [35]. We are of course speaking here of the modern day generalizations of classical modular forms on the upper half plane from which the subject as a whole now takes its name. (See [31].)

Roughly speaking, an automorphic form is a function on $G(\mathbb{Q}) \backslash G(\mathbb{A})$, while an (irreducible) automorphic representation is a representation $\pi$ of $G(\mathbb{A})$ by right translation on a space of automorphic forms. We usually think of $\pi$ as a representation (often unitary) on a Hilbert space. However, an automorphic form on $G(\mathbb{A})$ is required to satisfy finiteness conditions akin to those on the pre-Hilbert space $\mathscr{H}_G^0$ stated prior to the lemma above ($K_{\mathbb{R}} = K_\infty$-finite and $\mathscr{Z}_{\mathbb{R}} = \mathscr{Z}_\infty$-finite in the component $x_\infty$, $K^\infty$-finite and compactly supported in the component $x^\infty$, of the variable $x = x_\infty x^\infty$), as well as a condition of moderate growth (slowly increasing). With these constraints, we have therefore to treat $\pi$ as a representation of the "group algebra"

$$\mathscr{H} = \mathscr{H}_\infty \otimes \mathscr{H}^\infty = \mathscr{H}_{\mathbb{R}} \otimes \mathscr{H}^{\mathbb{R}}$$

of $G(\mathbb{A})$ (not to be confused with the Hilbert space $\mathscr{H}_G$ above). The factor $\mathscr{H}^\infty$ is the convolution algebra of locally constant functions of compact support on $G(\mathbb{A}^\infty)$, while $\mathscr{H}_\infty$ is the convolution algebra of distributions on the real group $G_{\mathbb{R}} = G_\infty$ that are supported on the maximal compact subgroup $K_{\mathbb{R}} = K_\infty$. They are often called Hecke algebras, and they act by right convolution on the space of automorphic forms. The algebra $\mathscr{H}_{\mathbb{R}} = \mathscr{H}_\infty$ is less well known, but it is quite elegant. It streamlines the archimedean actions as a $(\mathfrak{g}_{\mathbb{R}}, K_{\mathbb{R}})$-module ($\mathfrak{g}_{\mathbb{R}}$ is the Lie algebra of $G_{\mathbb{R}}$) into that of a single convolution algebra, which thus becomes parallel to the nonarchimedean action. It seems to have been first introduced by Flath [72], and was a basic part of the definitions in [35].

In practice, one generally wants to quantify the spaces of automorphic forms under consideration. Suppose that $\tau$ is the idempotent element in $\mathscr{H}_\infty$ attached to a finite set of irreducible characters on $K_\infty$, that $J$ is an ideal of finite codimension in $\mathscr{Z}_\infty$, that $K$ is any open compact subgroup of $G(\mathbb{A}^\infty)$ and that $N$ is a positive, slowly increasing function of $G_\infty = G(\mathbb{R})$. One can then write $\mathscr{A}(\tau, J, K, N)$ for the space of smooth (complex valued) functions $f$ on $G(\mathbb{A})$ with the following properties

(a)   $f(\gamma x) = f(x)$,        $\gamma \in G(\mathbb{Q}), x \in G(\mathbb{A})$,
(b)   $f * \tau = f$.
(c)   $zf = 0$,        $z \in J$,
(d)   $f(xk) = f(x)$,        $x \in G(\mathbb{A}), k \in K$,
(e)   $|f(x_\infty)| \leq c_f N(x_\infty)$,        $x_\infty \in G_\infty$,
   for some positive constant $c_f$.

Then $\mathscr{A}(\tau, J, K, N)$ is a space of automorphic forms on $G(\mathbb{A})$, which is stabilized by the natural subalgebra of $\mathscr{H}$ attached to $\tau$ and $K$.

   It is the nonarchimedean component of this subalgebra, the Hecke algebra

$$\mathscr{H}(G(\mathbb{A}^\infty), K) = \mathscr{H}(K \setminus G(\mathbb{A}^\infty)/K)$$

of compactly supported, $K$-biinvariant functions on $G(\mathbb{A}^\infty)$, that is particularly relevant to number theory. At the centre of its study are the unramified local factors

$$\mathscr{H}_v = \mathscr{H}(G_v, K_v) = \mathscr{H}(K_v \setminus G_v / K_v)$$

at which $K_v$ is a *hyperspecial* maximal compact subgroup of $G_v = G(\mathbb{Q}_v)$. These *spherical* Hecke algebras are abelian. What they revealed soon after Langlands' work on Eisenstein series became a critical part of his next great discovery.

   The notion of an automorphic form owes much to Harish-Chandra [89], as well as to Langlands. It is presented in [35, (1.3), (4.2)]. Borel and Jacquet followed this with the formal definition [35, (4.6)] of an automorphic representation as an irreducible constituent of a space of automorphic forms. This was supplemented with an equivalent formulation by Langlands [145, Proposition 2] as an irreducible constituent of a representation of $G(\mathbb{A})$ parabolically induced from a cuspidal automorphic representation. The passage between them was provided by Langlands' theory of Eisenstein series.

   The spectral theory of Eisenstein series was initiated by Selberg [205], [206], [207]. He established versions of the Langlands Main Theorem for various noncompact quotients

$$\Gamma \setminus X_+ \cong \Gamma \setminus \mathrm{SL}(2, \mathbb{R}) / \mathrm{SO}(2, \mathbb{R})$$

of the upper half plane

$$X_+ = \{z \in \mathbb{C} : \mathrm{Im}(z) > 0\}$$

and more generally, for the line bundles on these quotients attached to characters of the stabilizer $\mathrm{SO}(2, \mathbb{R})$ of $\sqrt{-1}$. This amounts to the theory of Eisenstein series on a noncompact quotient $\Gamma \setminus \mathrm{SL}(2, \mathbb{R})$. Selberg also included classical Hecke operators for modular forms into his constructions. This amounts in turn to the theory of Eisenstein series for the adelic quotient $\mathrm{SL}(2, \mathbb{Q}) \setminus \mathrm{SL}(2, \mathbb{A})$.

   Selberg regarded Eisenstein series as a step towards a trace formula for noncompact quotient. In adelic terms, his aim was to find a concrete formula for the trace of the operator $R_{\text{disc}}(f)$ obtained by restricting the right convolution operator

$$R(f) = \int_{\text{SL}(2,\mathbb{A})} f(x)R(x)\,dx, \quad f \in C_c^\infty(\text{SL}(2,\mathbb{A})),$$

on $L^2(\text{SL}(2,\mathbb{Q}) \setminus \text{SL}(2,\mathbb{A}))$ to the discrete spectrum. By definition

$$R_{\text{disc}}(f) = R(f) - R_{\text{cont}}(f)$$

where $R_{\text{cont}}(f)$ is the restriction of $R(f)$ to the continuous spectrum. The purpose of Eisenstein series was to provide an explicit construction for $R_{\text{cont}}(f)$. We shall return to this topic in Section 6.

# 3 *L*-functions and class field theory

The next period in Langlands' work is often identified with his 1967 letter to Weil [132]. It was subsequently expanded [138] to the series of far reaching conjectures and their consequences that became known as the Langlands program. Most immediately striking perhaps was Langlands' discovery of the long sought nonabelian class field theory. It will be a topic for the next section. In this section, we shall prepare the way. We shall review the theory of *L*-functions, and its relation to abelian class field theory. At the end, we shall then describe the hints Langlands found in Eisenstein series of what was lying ahead.

The last section might have seemed rather technical to a nonspecialist. We shall try to make amends in this section by moving at a more leisurely pace. In particular, we shall begin our discussion with the almost universally familiar notion of a Dirichlet series.

We recall that a *Dirichlet series* is an infinite series of the form

$$\sum_{n=1}^{\infty} a_n n^{-s}$$

for complex coefficients $a_n$ and a complex number $s$. If the coefficients satisfy a bound

$$|a_n| \le Cn^a, \quad n \in \mathbb{N},$$

for positive numbers $C$ and $a$, the series converges absolutely to an analytic function of $s$ in the right half plane $\mathrm{Re}(s) > a + 1$. The original model is of course the Riemann zeta function

$$\zeta(s) = \sum_{n=1}^{\infty} n^{-s}.$$

It converges to an analytic function of $s$ in the right half plane $\mathrm{Re}(s) > 1$. It also has an analytic continuation to a meromorphic function of $s \in \mathbb{C}$, whose only singularity is a simple pole at $s = 1$, and which satisfies a functional equation relating its values at $s$ and $1 - s$. The Riemann zeta function also has an Euler product. By the fundamental theorem of arithmetic, it can be represented as a product

$$\zeta(s) = \prod_{p} \left(1 - p^{-s}\right)^{-1} = \prod_{p} \left( \sum_{k=1}^{\infty} (p^k)^{-s} \right)$$

of Dirichlet series attached to prime numbers.

An *L-function* is a Dirichlet series with supplementary properties. There seems to be no universal agreement as to the definition, but let us say that an *L*-function is a Dirichlet series that converges in some right half plane, and that has an Euler product of the general form

$$L(s) = \prod_{p} \left( \sum_{k=1}^{\infty} c_{p,k} p^{-ks} \right),$$

for complex numbers $c_{p,k}$. We will not insist on analytic continuation and fundamental equation, simply because this has not been established for many of the $L$-functions that arise naturally, even though it is widely expected to hold.

In algebraic number theory, $L$-functions are used to encode arithmetic data. The coefficients $c_{p,k}$ in the Euler product turn out to provide a natural way to represent fundamental properties of the prime numbers. The essential examples are the $L$-functions of E. Artin, which were constructed from data that govern class field theory.

The goal of class field theory, as its name suggests, is to classify fields. More precisely, one would like to classify the Galois extensions of a given number field $F$. Let us review the problem in elementary terms.

We have been working in this paper with the base field $F = \mathbb{Q}$ for simplicity. Suppose then that $K$ is a finite Galois extension of $\mathbb{Q}$, which we assume to be *monogenic*. This means that we can represent $K$ as the splitting field of a monic, irreducible, integral polynomial

$$f(x) = x^n + a_{n-1}x^{n-1} + \cdots + a_1 x + a_0, \quad a_i \in \mathbb{Z}.$$

We then consider the factorization of $f(x)$ modulo a variable prime $p$. It is best to exclude the finite set of ramified primes for which $f(x)$ has repeated factors modulo $p$, a set $S$ in which we also include the archimedean place $\infty$ to be consistent with earlier notation. For each $p \notin S$, the corresponding factorization,

$$f(x) \equiv f_1(x) \cdots f_r(x) \,(\mathrm{mod}\ p)$$

of $f$ is then a product of distinct irreducible polynomials $f_i(x)$, whose degrees give us an (ordered) partition

$$\Pi_p = (n_1, \ldots, n_r), \quad n_i = \deg(f_i),$$

of $n$. We thus obtain a mapping from unramified primes $p \notin S$ to partitions $\Pi_p$ of $n$. (If $K/\mathbb{Q}$ is not monogenic, $f(x)$ need not be monic. In this case, we simply enlarge $S$ to include all prime divisors of the leading coefficient $a_n$.)

The interest in this mapping is in its implication for the Galois group

$$\Gamma_{K/\mathbb{Q}} = \mathrm{Gal}(K/\mathbb{Q})$$

of $K$ over $\mathbb{Q}$. This group is given by $f(x)$ as a conjugacy class of subgroups of the symmetric group $S_n$. Suppose for a moment that it equals the full symmetric group. The conjugacy classes of $\mathrm{Gal}(K/\mathbb{Q})$ are then parametrized by partitions of $n$. If $\Phi_p$ is the conjugacy class of $\Pi_p$ in $\mathrm{Gal}(K/\mathbb{Q})$, the factorization of $f(x)$ modulo $p$ then gives us a mapping

$$p \to \Phi_p \tag{12}$$

from primes $p \notin S$ to conjugacy classes $\Phi_p$ in $\Gamma_{K/\mathbb{Q}}$. In general, as we have said, $\mathrm{Gal}(K/\mathbb{Q})$ is determined by $f$ only as a conjugacy class of subgroups of $S_n$. However, a basic construction in algebraic number theory attaches a *canonical* conjugacy class in $\mathrm{Gal}(K/\mathbb{Q})$ to the unramified prime $p$, the Frobenius class

$$\mathrm{Frob}_p = \Phi_p.$$

Its implication for the factorization of $f$ is that among the conjugacy classes in $\Gamma_{K/\mathbb{Q}}$ attached to a partition $\Pi_p$, there is one that is canonical, the Frobenius class of $p$.

In these concrete terms, the problem for class field theory would be to characterize the prime $p$ factorization data of $f(x)$ in independent terms. More precisely, for any Galois extension $(K/\mathbb{Q})$, any partition $\Pi$ of $n$, and any conjugacy class $c$ in $\mathrm{Gal}(K/\mathbb{Q})$ that maps to the conjugacy class of $\Pi$ in $S_n$ say, can one characterize the fibre

$$P^S(c, K/\mathbb{Q}) = \{p \notin S \,:\, \Phi_p = c\}$$

of $c$ in some independent way? Of special interest is the case that $c$ is the trivial class 1 in $\mathrm{Gal}(K/\mathbb{Q})$. In this case, the fibre

$$\mathrm{Spl}(K/\mathbb{Q}) = \mathrm{Spl}^S(K/\mathbb{Q}) = P^S(1, K/\mathbb{Q}),$$

the set of primes that *split completely* in $K$, is just the set of primes for which $f(x)$ breaks into linear factors modulo $p$. Its importance is in the fact, which can be obtained as a direct consequence of the Tchebotarev density theorem for example, that the mapping

$$K/\mathbb{Q} \to \mathrm{Spl}(K/\mathbb{Q}),$$

from finite Galois extensions of $\mathbb{Q}$ to subsets of prime numbers $\{p\}$, is *injective*. In other words, the set $\mathrm{Spl}(K/\mathbb{Q})$ of primes represents a "signature" for the extension $(K/\mathbb{Q})$, in the sense that it characterizes it completely. The problem for class field theory would then be to characterize the image of this mapping in some independent fashion. This would provide a classification of Galois extensions of $\mathbb{Q}$.

These remarks are primarily for motivation, since the subject is inevitably more subtle. Nevertheless, such considerations were behind the development of (abelian) class field theory in the early part of the twentieth century. They also led E. Artin to define the $L$-functions that bear his name, as a way to encode the data provided by the conjugacy classes $\{\Phi_p\}$ in $\mathrm{Gal}(K/\mathbb{Q})$. The coefficients in $L$-functions are of course numbers. The simplest way to attach numbers to a conjugacy class is to embed the underlying Galois group into a general linear group, and then take the coefficients of the resulting characteristic polynomials. An Artin $L$-function for a Galois extension $(K/\mathbb{Q})$ therefore depends on the choice of a representation

$$r\colon \mathrm{Gal}(K/\mathbb{Q}) \to \mathrm{GL}(n, \mathbb{C}).$$

In its simplest form, it is defined as an Euler product

$$L^S(s, r) = \prod_{p \notin S} L_p(s, r), \tag{13}$$

with local factors

$$L_p(s, r) = \det(1 - r(\Phi_p)p^{-s})^{-1},$$

built out of unramified Frobenius conjugacy classes $\Phi_p$ in $\mathrm{Gal}(K/\mathbb{Q})$. It converges for $\mathrm{Re}(s) > 1$ to an analytic function of $s$.

Artin showed that the Euler product $L^S(s, r)$ actually has analytic continuation to a meromorphic function of $s \in \mathbb{C}$, with a functional equation that relates its values at $s$ and $1-s$. More precisely, he defined Euler factors $L_v(s, r)$ at the finite set of places $v \in S$, including the archimedean place $v = \infty$. He then showed that the completed product

$$L(s, r) = L_S(s, r) L^S(s, r) = \prod_v L_v(s, r) \tag{14}$$

satisfies the precise functional equation

$$L(s, r) = \varepsilon(s, r) L(1-s, r^\vee), \tag{15}$$

for the contragredient representation

$$r^\vee = {}^t r(g^{-1}), \quad g \in \mathrm{Gal}(K/\mathbb{Q}),$$

of $r$, and a certain rather mysterious monomial

$$\varepsilon(s, r) = ab^s, \quad a \in \mathbb{C}^\times, b > 0.$$

Given the analytic continuation, Artin then made the following remarkable conjecture.

**Conjecture (Artin)** *Suppose that $r$ is irreducible. Then $L(s, r)$ is an entire function of $s \in \mathbb{C}$, unless $r$ is the trivial 1-dimensional representation, in which case $L^S(s, r) = L^\infty(s, r)$ is the Riemann zeta function.*

Artin's proof of the analytic continuation and functional equation was intimately tied to abelian class field theory. This is the study of *abelian* extensions $(K/F)$ of a number field $F$, which is to say finite Galois extensions $(K/F)$ with abelian Galois group

$$\Gamma_{K/F} = \mathrm{Gal}(K/F).$$

We shall say a few words about this fundamental subject, continuing for the moment to assume that $F = \mathbb{Q}$ in order to be as concrete as possible.

Consider the special case of Artin's construction in which the Galois group $\Gamma_{K/F} = \Gamma_{K/\mathbb{Q}}$ is abelian, and the representation $r$ is irreducible. In other words, $r$ is a 1-dimensional character on $\Gamma_{K/F}$. If $p$ lies outside the finite set $S$ of ramified places, the associated Frobenius class $\Phi_p$ in $\Gamma_{K/\mathbb{Q}}$ is simply an element in this abelian group. We can use the nonzero complex numbers

$$\{r(\Phi_p) : p \notin S\}$$

to define a character $\chi^S$ on the locally compact group

$$\mathrm{GL}(1, \mathbb{A}^S) = (\mathbb{A}^*)^S = \{x^S \in \prod_{p \notin S} \mathbb{Q}_p^* : |x_p|_p = 1 \text{ for almost all } p\} = \prod_{p \notin S}^{\sim} \mathbb{Q}_p^*$$

by setting

$$\chi^S(x^S) = \prod_{p \notin S} \chi_p(x_p) = \prod_{p \notin S} r(\Phi_p)^{v_p(x_p)}, \tag{16}$$

for the valuation

$$v_p(x_p) = -\log_p(|x_p|_p), \quad x_p \in \mathbb{Q}_p^*.$$

Class field theory asserts that $\chi^S$ is the restriction to $\left(\mathbb{A}^S\right)^*$ of a uniquely determined automorphic representation of $\mathrm{GL}(1)$. In other words, there is a unique continuous, complex character $\chi$ on the quotient

$$C_F = F^* \backslash \mathbb{A}_F^* = \mathrm{GL}(1, F) \backslash \mathrm{GL}(1, \mathbb{A}_F), \quad F = \mathbb{Q},$$

whose restriction to the image of $(\mathbb{A}^*)^S$ equals $\chi^S$. Moreover, with a minor adjustment in the definition, the mapping $r \to \chi$ becomes an isomorphism of abelian groups. We can think of this property as a fundamental arithmetic law of nature. We restate it formally as follows.

Let

$$\Gamma_F^{\mathrm{ab}} = \mathrm{Gal}(F^{\mathrm{ab}}/F), \quad F = \mathbb{Q},$$

be the Galois group of the maximal abelian extension $F^{\mathrm{ab}}$ of $F = \mathbb{Q}$ (in some fixed algebraic closure $\overline{F}$ of $F$). It is an inverse limit

$$\Gamma_F^{\mathrm{ab}} = \varprojlim_K \Gamma_{K/F}$$

over finite abelian extensions, and hence a compact, totally disconnected group. A continuous 1-dimensional character $r$ in $\Gamma_F^{\mathrm{ab}}$ has cofinite kernel, and can therefore be identified with a character on the Galois group $\Gamma_{K/F}$ of a finite extension $(K/F)$. This in turn maps to a character $\chi^S$ on the group of $S$-idèles

$$I_F^S = (\mathbb{A}_F^*)^S = \{x_p \in F_p^* : p \notin S\} \subset I_F = \mathbb{A}_F^*, \quad F = \mathbb{Q}.$$

for a finite set of valuations $S$ outside of which $r$ is unramified. The following big theorem, in which we have taken $F = \mathbb{Q}$, is then the central assertion of class field theory.

**Global Reciprocity Law** *For any $r$, the character $\chi^S$ on $I_F^S$ descends to a unique character $\chi$ on the idèle class group*

$$C_F = F^* \backslash I_F = F^* \backslash \mathbb{A}_F^*, \quad F = \mathbb{Q}.$$

*The resulting mapping $r \to \chi$ becomes an isomorphism from the group of characters $r$ on $\Gamma_F^{\mathrm{ab}}$ onto the group of characters $\chi$ of finite order on $C_F$. The dual (Artin) mapping*

$$\theta_F : C_F \to \Gamma_F^{\mathrm{ab}}$$

*is therefore a continuous surjective homomorphism, whose kernel is the connected component $C_F^0$ of $1$ in $C_F$.*

This assertion is a culmination of work by mathematicians over many years, from the law of quadratic reciprocity of Gauss to the full reciprocity law of Artin, which is a more precise version of the assertion. For the Artin reciprocity law also characterizes the preimage of each finite quotient $\Gamma_{K/F}$ of $\Gamma_F^{\mathrm{ab}}$ as the cokernel of the norm mapping $N_{K/F}$ from $C_K$ to $C_F$. In general, the importance of the reciprocity law lies especially in the fact that it applies, as stated, if $F$ is any number field. We need only replace prime numbers $p$ in $\mathbb{Q}$ by prime ideals in $F$ in the definitions above.

Needless to say, the proof of the reciprocity law is deep. The original argument had a large analytic component, based on the properties of certain $L$-functions. A purely algebraic proof based on the cohomology of Galois groups came later. For the field $F = \mathbb{Q}$, the proof is actually much easier than in the general case. This is reflected in the Kronecker–Weber theorem, which characterizes $\mathbb{Q}^{\mathrm{ab}}$ explicitly as the field generated by the complex numbers

$$\{ e^{\frac{2\pi i}{n}} : n \in \mathbb{N} \}.$$

A similarly explicit result holds for any imaginary quadratic extension $F$ of $\mathbb{Q}$. For in this case, the Kronecker Jugendtraum characterizes $F^{\mathrm{ab}}$ in terms of special values of elliptic functions attached to elliptic curves over $\mathbb{Q}$ with complex multiplication in $F$. We recall that Hilbert's twelfth problem was to characterize the maximal abelian extension $F^{\mathrm{ab}}$ of any $F$ in terms of special values of natural analytic functions. Little progress has since been made on this problem, apart from the 1955 generalization of the Jugendtraum by Shimura and Taniyama to a totally complex extension of a totally real field.

How did Artin use the reciprocity law to prove his functional equation? We must first recall that Hecke has earlier attached an $L$-function $L^S(s, \chi)$ to any (quasi)character $\chi$ on the group $C_F = F^* \setminus \mathbb{A}_F^*$. We are using adelic notation here, as we did in our discussion of the reciprocity law, even though it was only later that Chevalley introduced the group of idèles $I = \mathbb{A}_F^*$. Hecke called $\chi$ a Grössencharacter, as a generalization of a Dirichlet character and of its analogue introduced by Weber for the number field $F$, which we can now take to be arbitrary. As a character on $C_F$ (rather than what Hecke would have considered a generalized ideal class character), $\chi$ has an unramified part $\chi^S$ on the locally compact group

$$\left( \mathbb{A}_F^S \right)^* = \widetilde{\prod_{v \notin S} F_v^*}.$$

It takes the form

$$\chi^S(x^S) = \prod_{v \notin S} \chi_v(x_v) = \prod_{v \notin S} c_v^{v(x_v)}$$

for complex numbers $c_v \in \mathbb{C}^*$ attached to $\chi$, and with $v$ being the normalized valuation of $F_v^*$. This of course is parallel to (16), where we had $F = \mathbb{Q}$. In particular, $S$ is a finite set of valuations of $F$ that contains both the set $V_\infty$ of archimedean places and the set of finite places $v$ at which $\chi_v$ is ramified. The unramified Hecke $L$-function is then the infinite product

$$L^S(s, \chi) = \prod_{v \notin S} L_v(s, \chi) = \prod_{v \notin S} (1 - c_v q_v^{-s})^{-1},$$

where $q_v$ is the order of the residue field $F_v^*$. After introducing this function, Hecke defined Euler factors $L_v(s, \chi)$ at the remaining places $v \in S$, and then proved that the completed product

$$L(s, \chi) = L_S(s, \chi) L^S(s, \chi) = \prod_v L_v(s, \chi)$$

satisfies the functional equation

$$L(s, \chi) = \varepsilon(s, \chi) L(1 - s, \chi^\vee) \tag{17}$$

where $\chi^\vee(x) = \chi(x^{-1})$ and

$$\varepsilon(s, \chi) = ab^s, \quad a \in \mathbb{C}^*, b > 0.$$

The central tenet of (global) class field theory can be formulated as an assertion that every abelian Artin $L$-function $L^S(s, r)$ over $F$ equals a Hecke $L$-function $L^S(s, \chi)$. This follows from the earlier definition of the mapping $r \to \chi$, the definition of the $L$-functions themselves, and of course, the Global Reciprocity Law. Therefore $L^S(s, r)$ inherits all the properties established by Hecke for $L^S(s, \chi)$. It has a completion that equals $L(s, \chi)$, and therefore has analytic continuation, and functional equation (17). This is a fundamental fact, that to this day has no direct proof. It was the main step in Artin's derivation of the functional equation (15) for an arbitrary (nonabelian) representation

$$r \colon \operatorname{Gal}(E/F) \to \operatorname{GL}(n, \mathbb{C}).$$

The other step was a formal decomposition

$$L^S(s, r) = \prod_i L^S(s, r_i)^{a_i},$$

for representations $r_i$ of cyclic Galois groups $\Gamma_i \subset \operatorname{Gal}(K/F)$ and rational numbers $a_i$, obtained by Artin by first proving a character theoretic decomposition

$$r = \sum_i a_i \operatorname{Ind}(\Gamma, \Gamma_i; r_i), \quad \Gamma = \operatorname{Gal}(K/F),$$

of $r$, and then using the compatibility of the $L$-functions with representation theoretic operations such as induction and direct sums. He was then able to use this to construct his completed $L$-function $L(s, r)$ from those of Hecke. The last step was to use the resulting formal identity of quotients

$$\frac{L(s,r)}{L(1-s,r^\vee)} = \prod_i \left( \frac{L(s,r_i)}{L(1-s,r_i^\vee)} \right)^{a_i} = \prod_i (a_i b_i^s)^{a_i}$$

to establish the identity

$$\frac{L(s,r)}{L(1-s,r^\vee)} = ab^s, \quad a \in \mathbb{C}^*, b > 0.$$

This is the functional equation (17).

We refer the reader to the article [53] of Cogdell on Artin $L$-functions for discussion of this and other questions. The earlier history of class field theory is clearly presented in the articles [252] and [93]. We also note that the argument sketched above becomes quite clear after the later proof of the Brauer Induction Theorem. This establishes that the rational numbers $a_i$ may in fact be taken to be integers. We should point out, however, that this improvement gives little information about the Artin conjecture stated above. To rule out any poles of the $L$-function $L(s,r)$, we would need to control the zeros of the abelian $L$-functions $L^S(s,r_i)$ with $a_i < 0$. This of course would be a tall order.

It was with the thesis [236] of Tate that the work of Hecke was put into the adelic form that we have followed here. This made many of Hecke's arguments more transparent. For example, the functional equation (17) is seen to be a natural consequence of the Poisson summation formula on the locally compact abelian group $\mathbb{A}_F$, with respect to the discrete subgroup $F$. It is at this point that a new ingredient enters the theory, a nontrivial additive character $\psi$ on $\mathbb{A}_F/F$, with decomposition $\psi = \widetilde{\prod_v} \psi_v$ into additive characters on the groups $F_v$. Its role is to identify the Pontryagin dual of $\mathbb{A}_F$ with $\mathbb{A}_F$ itself. The adelic formulation of [236] also leads to an important new way to see the final result. Tate showed that the global $\varepsilon$-factor in (17) has a canonical decomposition

$$\varepsilon(s,\chi) = \prod_{v \in S} \varepsilon_v(s,\chi_v,\psi_v) \tag{18}$$

into local $\varepsilon$-factors

$$\varepsilon_v(s,\chi_v,\psi_v) = \varepsilon(\chi_v,\psi_v) q_v^{-n_v(s-\frac{1}{2})}$$

that depend only in the localizations $\psi_v$ and $\chi_v$ of $\psi$ and $\chi$, for nonzero complex numbers

$$\varepsilon(\chi_v,\psi_v) = \varepsilon(\tfrac{1}{2},\chi_v,\psi_v),$$

integers $n_v = n(r_v,\psi_v)$ known as conductors, and integers $q_v$ given by the residual degree of $F_v$ if $v$ is nonarchimedean, and 1 if $v$ is archimedean. The local $\varepsilon$-factors in the functional equation turned out to have a fundamental role in the local Langlands program.

We can now turn to the work of Langlands. The first thing to mention is the volume [81] in which Godement and Jacquet generalize Tate's thesis from $\mathrm{GL}(1)$ to $\mathrm{GL}(n)$. This was a major step forward, which was not due to Langlands. However, its possible existence was clearly part of his thinking right from the beginning.

Langlands was in regular communication with Godement in the 1960s, and he mentions the extension of Tate's thesis on p. 31–32 of his fundamental paper [138] as an essential premise for his conjectures. In fact, the special case with $n = 2$ was established in the volume [103] of Jacquet and Langlands and was announced in the original article [138]. We shall return to this briefly in Section 5, but we might as well stay with the more general case here for the motivation.

The problem in [81] was to attach an $L$-function $L(s, \pi)$ to every automorphic representation

$$\pi = \widetilde{\bigotimes_v} \pi_v$$

of $G(\mathbb{A}_F)$, for the group $G = \mathrm{GL}(n)$ over a number field $F$, and to prove that this function has analytic continuation and functional equation. The initial step will by now be quite familiar. Given $\pi$, we fix a finite set of valuations $S \supset S_\infty$ of $F$ such that the local constituent $\pi_v$ of $\pi$ is unramified for any $v$ outside of $S$. This means that the restriction of $\pi_v$ to the maximal compact subgroup $G(\mathfrak{O}_v)$ of $G(F_v)$ contains the trivial one-dimensional representation. For any such $v$ there is a general bijection

$$\pi_v \to c(\pi_v),$$

from unramified representations of $G(F_v)$ onto semisimple conjugacy classes in the complex group $G(\mathbb{C})$. We shall discuss this last point again, for more general groups, in the next section. Semisimple conjugacy classes in $G(\mathbb{C}) = \mathrm{GL}(n, \mathbb{C})$ can of course be identified with complex diagonal matrices, taken only up to permutation of their entries. The unramified global $L$-function of $\pi$ is then defined as the Euler product

$$L^S(s, \pi) = \prod_{v \notin S} L_v(s, \pi),$$

where

$$L_v(s, \pi) = L(s, \pi_v) = \det(1 - c(\pi_v)q_v^{-s})^{-1}.$$

The local problem in [81] was to attach local $L$-functions

$$L_v(s, \pi) = L(s, \pi_v)$$

and the $\varepsilon$-factors

$$\varepsilon_v(s, \pi, \psi) = \varepsilon(s, \pi_v, \psi_v) = \varepsilon(\pi_v, \psi_v)q_v^{-n_v(s - \frac{1}{2})}$$

to the remaining valuations $v \in S$. This is naturally much more subtle for the nonabelian group $\mathrm{GL}(n)$, but the basic ideas resemble those for $\mathrm{GL}(1)$ in [236]. The same is true of the global problem of analytic continuation and functional equation. The essential tool was again the Poisson summation formula, this time for the additive group $\mathfrak{g}(\mathbb{A}_F) = \mathrm{M}_n(\mathbb{A}_F)$ of $(n \times n)$-adelic matrices, with respect to the discrete subgroup $\mathfrak{g}(F) = \mathrm{M}_n(F)$ of rational matrices. In the end, the authors constructed the local $L$-functions and $\varepsilon$-factors above for places $v \in S$, such that the full Euler

product

$$L(s,\pi) = L_S(s,\pi)L^S(s,\pi) = \prod_v L_v(s,\pi)$$

has analytic continuation, and satisfies the functional equation

$$L(s,\pi) = \varepsilon(s,\pi)L(1-s,\pi^\vee), \tag{19}$$

for

$$\varepsilon(s,\pi) = \prod_{v \in S} \varepsilon_v(s,\pi,\psi). \tag{20}$$

This is the general analogue for $\mathrm{GL}(n)$ of the functional equation (17) established by Tate for $\mathrm{GL}(1)$. We note that the global solution in [81] was established only in the case that $\pi$ is a *cuspidal* automorphic representation. However, the general case follows from this and the properties of Langlands' Eisenstein series. (See [33] and [102, §6].)

With the Artin $L$-functions of degree $n$ and the automorphic $L$-functions of $\mathrm{GL}(n)$, our exposition has acquired a certain symmetry. Given the Global Reciprocity Law for $\mathrm{GL}(1)$, a reader could well ask whether every Artin $L$-function is an automorphic $L$-function. If so, we would have a general extension of the reciprocity law to what we could regard as nonabelian class field theory. We would also have a proof of the Artin conjecture stated above. For it is a fundamental consequence of the harmonic analysis used by Jacquet and Godement (and Tate and Hecke for $n = 1$) that the automorphic $L$-functions $L(s,\pi)$ are *entire*, apart from certain obvious exceptions related to unramified 1-dimensional automorphic representations. Would this then be the final word on the subject?

There are three points to consider in regards to the last question. One would be the uncomfortable prospect of having to prove such a general nonabelian reciprocity law, given the historical difficulty in establishing just the abelian theory. We could expect that nonabelian class field theory, whatever form it might take, would be difficult. It would be reassuring to think that the problem at least has some further structure. A second point concerns this last possibility. Suppose that $r'$ is an irreducible Galois representation of degree $n'$, and that $\rho'$ is an irreducible $n$-dimensional representation of $\mathrm{GL}(n',\mathbb{C})$. The composition

$$\rho : \Gamma_F \xrightarrow{r'} \mathrm{GL}(n',\mathbb{C}) \xrightarrow{\rho'} \mathrm{GL}(n,\mathbb{C})$$

is then a Galois representation (often irreducible) of degree $n$. The Frobenius classes that define the $L$-functions satisfy the following relation

$$\rho(\Phi_v) = (\rho' \circ r')(\Phi_v), \quad v \notin S.$$

How would this structure be reflected in the corresponding automorphic representations? Finally, the work of Harish-Chandra has taught us that representations should be studied uniformly for all reductive groups. If some interesting phenomenon is discovered for one group, or one family of groups such as $\{\mathrm{GL}(n)\}$, it should be investigated for all groups. What are the implications of this for automorphic *L*-functions?

These considerations were undoubtedly part of the thinking of Langlands that led up to the Principle of Functoriality. However, perhaps the most decisive hints were in his theory of Eisenstein series. They came from *L*-functions he discovered in the global intertwining operators

$$M(w, \lambda) \colon \mathscr{H}_P \to \mathscr{H}_{P'}, \quad w \in W(\mathfrak{a}_P, \mathfrak{a}_{P'}), \tag{21}$$

defined by (7). (We have written $w$ here because we want to reserve $s$ for the complex variable of an *L*-function.)

Suppose for a moment that $G$ equals the group $\mathrm{SL}(2)$ over $\mathbb{Q}$, and that $P$ is the Borel subgroup of upper triangular matrices, that $\phi$ lies in the one-dimensional space of constant functions in $\mathscr{H}_p$, and that $w$ is the non-trivial Weyl element in $W(\mathfrak{a}_P, \mathfrak{a}_P)$. Then $M(w, \lambda)$ is a scalar multiple of $\phi$ given for a suitable $\lambda$ by a convergent adelic integral over

$$N(\mathbb{A}) = \left\{ \begin{pmatrix} 1 & x \\ 0 & 1 \end{pmatrix} : x \in \mathbb{A} \right\}.$$

It is not hard to evaluate. Recall that $\lambda \in \mathfrak{a}_{P,\mathbb{C}}^*$ is a complex valued linear form on $\mathfrak{a}_P$, a real 1-dimensional vector space we can in turn identify with the Cartan subalgebra

$$\left\{ \begin{pmatrix} u & 0 \\ 0 & -u \end{pmatrix} : u \in \mathbb{R} \right\}$$

of the Lie algebra of $\mathrm{SL}(2, \mathbb{R})$. The mapping

$$\lambda \to s = \lambda \begin{pmatrix} \frac{1}{2} & 0 \\ 0 & -\frac{1}{2} \end{pmatrix}$$

identifies $\lambda$ with a complex number $s \in \mathbb{C}$. It can then be shown that

$$M(w, \lambda) = \frac{L(s)}{L(s+1)},$$

for the completed Riemann zeta function

$$L(s) = L_\infty(s)\zeta(s).$$

This simple, well known formula is suggestive. If the Riemann zeta function is at the heart of the global intertwining operator for $\mathrm{SL}(2)$, something new and interesting must surely be contained in the operators for groups of larger rank. (For further background, see [78, §2].)

In his investigation of the more general intertwining operators in [141], Langlands discovered some completely new *L*-functions. Suppose that $G$ is a split, simple group, say over $\mathbb{Q}$, and that

$$P = MN \supset B$$

is a standard maximal parabolic subgroup. Suppose also that

$$\pi = \widetilde{\bigotimes_{v}} \pi_v = \pi_\infty \otimes \left( \widetilde{\bigotimes_{p}} \pi_p \right)$$

is a cuspidal automorphic representation of $M(\mathbb{A})$ that is unramified at every place $v$ of $\mathbb{Q}$. The parametrization of unramified representations $\pi_v$ of a local group $M_v = M(\mathbb{Q}_v)$ is not difficult, but it was not particularly well known at the time. One can identify them either with complex valued homomorphisms

$$\mathscr{H}(K_{M,v} \backslash M_v / K_{M,v}) \to \mathbb{C}$$

of the unramified Hecke algebra (under convolution) of functions in $C_c^\infty(M_v)$ biinvariant under a suitable maximal compact $K_{M,v}$, or with their induction parameters given by Weyl orbits of unramified characters on the Borel subgroup $B_v \cap M_v$ of $M_v$. Langlands investigated them in terms of Hecke algebras, and made a remarkable observation. For any $v$, the unramified representations $\pi_v$ of $M_v$ are parametrized by the semisimple conjugacy classes $c(\pi_v)$ in a different group, the *complex dual group* $\widehat{M}$ of $M$.

The dual group $\widehat{G}$ of $G$ is a complex simple group whose Coxeter–Dynkin diagram is dual to that of $G$, in the sense that the directions of any arrows are reversed, and whose maximal torus

$$\widehat{T} = X^*(T) \otimes \mathbb{C}^*, \quad T \subset M \subset G,$$

is dual to that of $G$. The roots $\{\alpha^\vee\}$ of $(\widehat{G}, \widehat{T})$ are the co-roots of $(G, T)$, and the corresponding fundamental dominant weights $\{\varpi^\vee\} = \{\varpi_\alpha^\vee\}$ form the dual basis of the set of co-roots $\{(\alpha^\vee)^\vee\} = \{\alpha\}$ of $(G, T)$. Finally, to the maximal parabolic subgroup $P = MN$ of $G$, there corresponds a maximal parabolic subgroup $\widehat{P} = \widehat{M}\widehat{N}$ of $\widehat{G}$, whose Levi component $\widehat{M}$ represents the dual group of $M$. We assume for simplicity that the unipotent radical $\widehat{N}$ of $\widehat{P}$ is abelian, and we write $r$ for the adjoint representation of $\widehat{M}$ in the Lie algebra $\widehat{\mathfrak{n}}$ of $\widehat{N}$. We shall say more about the dual group $\widehat{G}$, and its more sophisticated companion the $L$-group $^L G$, in the next section, but we will still not give the precise construction. For this, we simply refer the reader to the article [233] by Springer.

Given $G$, $P$, $\pi$ and $v$, we thus have a semisimple conjugacy class $c(\pi_v)$, according to what Langlands discovered in [139]. In general, $\widehat{M}$ is not a general linear group, so it does not have a characteristic polynomial. However, Langlands used $r$ to define

$$L_v(s, \pi, r) = L(s, \pi_p, r) = \det(1 - r(c(\pi_p))p^{-s})^{-1}$$

explicitly for $v = p$ nonarchimedean. He also defined

$$L_v(s, \pi, r) = L(s, \pi_\infty, r) = \prod_{\alpha^\vee} \left( \pi^{-(s+c_\alpha)/2} \Gamma\left(\tfrac{s+c_\alpha}{2}\right) \right)$$

implicitly for $v = \infty$ archimedean, for complex numbers $c_\alpha = c_\alpha(\pi_\varpi, r)$, and roots $\alpha^\vee$ of $(\widehat{G}, \widehat{T})$ that are not roots of $(\widehat{M}, \widehat{T})$. He then showed that the Euler product

$$L(s, \pi, r) = \prod_v L_v(s, \pi, r)$$

converged absolutely to an analytic function on some right half plane, which he conjectured had analytic continuation with functional equation

$$L(s, \pi, r) = L(1 - s, \pi, r^\vee), \tag{22}$$

for the contragredient representation $r^\vee(x) = {}^t r(x)^{-1}$ of $r$.

This was suggested by Langlands' theory of Eisenstein series, specifically the operators $M(w, \lambda)$ in (21). Since $P$ is maximal, there will be a unique nontrivial Weyl element $w \in W(\mathfrak{a}_P, \mathfrak{a}_{P'})$ (for a standard parabolic subgroup $P' \subset B$). Since $\pi$ is unramified, $\mathcal{H}_P$ and $\mathcal{H}_{P'}$ both have canonical, one-dimensional subspaces of constant functions, which are preserved both under translation by $w$ and by the operators $M(w, \lambda)$ themselves. We can therefore identify $M(w, \lambda)$ with a complex valued scalar function, as was the case for $G = \mathrm{SL}(2)$ above. Let $\varpi_P^\vee = \varpi_\alpha^\vee$ be the fundamental dominant weight for $(\widehat{G}, \widehat{T})$, attached to the simple root $\alpha$ of $(G, T)$ that is nontrivial on the split component $A_P = A_{M_P}$ of $P$. The mapping

$$\lambda \to s = \lambda(\varpi_P^\vee)$$

then identifies $\lambda$ with a complex number $s \in \mathbb{C}$. The main result of [139] is the formula

$$M(w, \lambda) = \frac{L(s, \pi, r^\vee)}{L(s + 1, \pi, r^\vee)}, \quad \lambda \to s, \tag{23}$$

for the scalar valued restriction of the operator (21). Langlands proved it for $G$, $P$ and $\pi$ (and, slightly modified, if $\widehat{N}$ is nonabelian) by the formula of Gindikin and Karpelevic for $G_\infty$ [80], and analogues he derived for the $p$-adic groups $G_p$.

The formula (23) not only motivated the introduction of many new automorphic $L$-functions, but also raised interesting new questions on Artin $L$-functions. Given $G$ and $P$, consider a continuous homomorphism

$$\rho \colon \Gamma_{\widehat{\mathbb{Q}}} \to \widehat{M}.$$

One could ask whether there is an automorphic representation $\pi$ of $M(\mathbb{A})$, and a corresponding $L$-function $L(s, \pi, r^\vee)$ such that

$$L(s, r^\vee \circ \rho) = L(s, \pi, r^\vee).$$

Langlands' conjectural answer to this will be taken up in the next section.

# 4 Global Functoriality and its implications

The Principle of Functoriality is the centre of the Langlands program. It postulates deep relationships among automorphic representations on different groups $G$. These in turn tie fundamental arithmetic data from number theory to equally fundamental spectral data from harmonic analysis. The global relationships also suggest local relationships that should lead to a local classification of representations. These should in turn give rise to local $L$-functions and $\varepsilon$-factors. We shall discuss Global Functoriality in this section, and Local Functoriality in the next.

The general functoriality conjecture was introduced in the revolutionary paper [138], a preprint of which was written shortly after the letter [132] to Weil. Langlands had to anticipate a number of basic properties of automorphic representations to be able even to formulate the conjecture. The name *functoriality* itself did not appear in the paper. Nor did the term *automorphic representation*, which seems to have been first introduced by A. Borel in his Bourbaki lecture [32, (5.1)].

In the original paper [138], Langlands spoke simply of an irreducible representation $\pi$ of $G(\mathbb{A})$ that "occurs in" $L^2(G(F) \setminus G(\mathbb{A}))$. This amounts to an informal definition of an automorphic representation of $G(\mathbb{A})$, with the understanding that it includes the nonunitary representations obtained by analytic continuation of their internal parameters into the complex domain. He also took for granted that any such representation could be obtained uniquely as a restricted tensor product

$$\pi = \widetilde{\bigotimes_v} \pi_v, \tag{24}$$

where for any $v$, $\pi_v$ is an irreducible representation of $G(\mathbb{Q}_v)$. This brought harmonic analysis into an area that was already broad, but that had applied mainly to complex analysis, number theory and arithmetic geometry. The formal definition was given later by Borel and Jacquet [35] in terms of automorphic forms. As we noted in §2, this was accompanied by the somewhat more direct characterization by Langlands in terms of induced cuspidal representations [145]. The tensor product decomposition was established at the same time by Flath [72].

Langlands' paper [139] on Euler products, discussed very briefly at the end of the last section, was a precursor to the functoriality paper [138] we are discussing now. Both papers required some ad hoc background to elementary properties of automorphic representations that had yet to be developed. These properties are by now well understood. We shall review a few of them here from the perspective of [138], even if this entails some repetition from our last section.

We are assuming for the moment that $G$ is a connected reductive algebraic group over the field $F = \mathbb{Q}$. A fundamental concept introduced in [138] was Langlands' notion of the *dual group* $\widehat{G}$ of $G$, and more generally, the associated *L-group*

$$^LG = \widehat{G} \rtimes \mathrm{Gal}(K/F), \tag{25}$$

where $K$ is a suitably large finite Galois extension of $F$. Once again, the names came later, from Borel [32] in the case of $L$-group, and Kottwitz [119] for the dual group. We recall from Section 3 that $\widehat{G}$ is a complex, connected, reductive algebraic Lie group, with maximal torus

$$\widehat{T} = X^*(T) \otimes \mathbb{C}^*$$

dual to the maximal torus $T$ of $G$, in the sense that

$$X^*(\widehat{T}) = \operatorname{Hom}(X^*(T), \mathbb{Z}).$$

Langlands took care to make this construction rigid by among other things fixing Borel subgroups $B \supset T$ and $\widehat{B} \supset \widehat{T}$ of $G$ and $\widehat{G}$ respectively. The supplementary structure is conveniently accommodated by letting $G$ represent a *based root datum* and $\widehat{G}$ represent the canonical *dual based root datum* [233]. With this understanding, any outer automorphism $\alpha \in \operatorname{Out}(G)$ of $G$ comes with a canonical dual outer automorphism $\widehat{\alpha} \in \operatorname{Out}(\widehat{G})$ of $\widehat{G}$. As a reductive group over $F = \mathbb{Q}$, $G$ comes with a homomorphism

$$\operatorname{Gal}(K/F) \to \operatorname{Out}(G). \tag{26}$$

The $L$-group (25) is defined as the semidirect product of $\widehat{G}$ with $\operatorname{Gal}(K/F)$ under the corresponding dual homomorphism

$$\operatorname{Gal}(K/F) \to \operatorname{Out}(\widehat{G})$$

As predicted in [138], the Langlands $L$-group has turned out to be exactly the right object to accommodate the parameters of automorphic representations.

We are taking for granted here some knowledge of the structure of reductive algebraic groups over number fields. The introductory article [233] quoted above is perhaps the most convenient reference. It is the first paper in the two volume proceedings from the 1977 AMS summer symposium in Corvallis, Oregon, which was the natural successor to the 1964 Boulder conference we discussed in Section 1. Its goal was primarily to bring the subsequent work of Langlands to a broader audience.

Langlands formulated the Principle of Functoriality (minus the name) in [138], in its full generality. However, it is easier to recognize its beauty and power if we first describe a special case, the unramified part of functoriality, for split groups $G$ over $\mathbb{Q}$. This means that the image of the outer twisting homomorphism (26) is trivial. In particular, the semidirect product (25) that gives the $L$-group is actually a direct product. It also allows us for many purposes to take $K$ equal to our base field $F = \mathbb{Q}$, and therefore to take $^L G = \widehat{G}$.

For the given split group, suppose that $\pi$ is an automorphic representation of $G(\mathbb{A})$, with decomposition (24) into local constituents. The formal definition [35] of automorphic representation includes a weak continuity condition. This implies that there is a finite set $S$ of valuations of $\mathbb{Q}$ containing the archimedean place $v = \infty$ such that for any $p \notin S$, $\pi_p$ is an *unramified* representation of $G(\mathbb{Q}_p)$. Unramified in turn means that the restriction of $\pi_p$ to the maximal compact subgroup $K_p = G(\mathbb{Z}_p)$ of $G(\mathbb{Q}_p)$ contains the trivial one-dimensional representation. (It is understood here that as a split group over $\mathbb{Q}$, $G$ comes with a suitable $\mathbb{Z}$-scheme structure.) What makes automorphic representations $\pi$ interesting is the fact, established now in some cases and conjectured by Langlands in others, is that for many $\pi$, there are deep and fundamental relationships among its unramified constituents $\{\pi_p : p \notin S\}$.

To better appreciate the phenomenon, we recall that there is a simple classification of the unramified representations $\{\pi_p\}$ of a split $p$-adic group $G(\mathbb{Q}_p)$. It takes the form of a bijective mapping

$$\pi_p \to c(\pi_p) \tag{27}$$

in which $c(\pi_p)$ ranges over the semisimple conjugacy classes in the complex dual group $\widehat{G}$ of $G$. For as Langlands pointed out in [138], any semisimple element $c_p \in \widehat{G}$ determines an unramified, one-dimensional quasi-character on a Borel subgroup $B(\mathbb{Q}_p)$ of $G(\mathbb{Q}_p)$. The corresponding induced representation $\tilde{\pi}_p$ of $G(\mathbb{Q}_p)$ then depends only on the conjugacy class $\tilde{c}_p$ of $c_p$ in $\widehat{G}$. A given unramified representation $\pi_p$ of $G(\mathbb{Q}_p)$ thus occurs as the unique irreducible constituent of $\tilde{\pi}_p$, for a unique semisimple conjugacy class $\tilde{c}_p = c(\pi_p)$ in $\widehat{G}$. We therefore have a mapping

$$\pi \to c^S(\pi) = \{c_p(\pi) = c(\pi_p) : p \notin S\} \tag{28}$$

from automorphic representations $\pi$ of $G(\mathbb{A})$ to families of semisimple conjugacy classes in $\widehat{G}$.

The semisimple conjugacy classes $c_p(\pi)$ in $\widehat{G}$ are concrete objects. They can be described in terms of complex numbers. For example, if $G$ equals $\mathrm{GL}(n)$, the semisimple classes in $\widehat{G} = \mathrm{GL}(n, \mathbb{C})$ are given by nonsingular complex diagonal matrices, taken up to permutation of their entries. These are in turn parametrized by their characteristic polynomials, or if one prefers, elements in the space $\mathbb{C}^{n-1} \times \mathbb{C}^*$ of coefficients defined by the characteristic polynomial. A general automorphic representation therefore gives us a family $c^S(\pi)$ of objects $c_p(\pi)$ with natural complex parameters. It is in the relations among these complex parameters, as $p$ varies, that the most fundamental interest of automorphic representations lies.

We can now state the global Principle of Functoriality. In the most basic case under present consideration, it applies to a pair of split groups $G'$ and $G$ over the field $F = \mathbb{Q}$, together with a homomorphism

$$\rho' \colon {}^L G' \to {}^L G$$

between the corresponding dual groups $\widehat{G}' = {}^L G'$ and $\widehat{G}$ and ${}^L G$. The homomorphism defines a mapping of semisimple conjugacy classes in the two groups, which we also denote by $\rho'$.

**Conjecture (Langlands' Principle of Global Functoriality)** *Given $G'$, $G$ and $\rho'$, suppose we also have an automorphic representation $\pi'$ of $G'(\mathbb{A})$. Then there is an automorphic representation $\pi$ of $G(\mathbb{A})$ such that*

$$c^S(\pi) = \rho'(c^S(\pi')).$$

*In other words, there is a finite set S of primes outside of which $\pi'$ and $\pi$ are both unramified, and for which*

$$c_p(\pi) = \rho'(c_p(\pi')), \quad p \notin S.$$

We can remove the finite set $S$ from the assertion by defining an equivalence relation on the set of families $c^S$ of semisimple conjugacy classes in $\widehat{G}$. We write $c^S \sim c_1^{S_1}$ for two such families $c^S$ and $c_1^{S_1}$ if

$$c_p = c_{1,p},$$

for almost all $p \notin S \cup S_1$. For any $\pi$ in the set

$$\Pi_{\mathrm{aut}}(G) = \Pi(G)$$

of automorphic representations of $G(\mathbb{A})$, we then write $c(\pi)$ for the equivalence class that contains the family $c^S(\pi)$. We thus obtain a surjective mapping

$$\pi \to c(\pi), \quad \pi \in \Pi_{\mathrm{aut}}(G),$$

from $\Pi_{\mathrm{aut}}(G)$ onto the set

$$\mathscr{C}_{\mathrm{aut}}(G) = \{c(\pi) : \pi \in \Pi_{\mathrm{aut}}(G)\}.$$

If $G = \mathrm{GL}(n)$, the restriction of the mapping to the subset $\Pi_{\mathrm{cusp}}(G)$ of cuspidal automorphic representations is injective. This is the theorem of strong multiplicity 1 for $\mathrm{GL}(n)$ [189], [108]. It gives a bijection from $\Pi_{\mathrm{cusp}}(G)$ onto the set

$$\mathscr{C}_{\mathrm{cusp}}(G) = \{c(\pi) : \pi \in \Pi_{\mathrm{cusp}}(G)\},$$

and thus allows us to identify a cuspidal automorphic representation of $\mathrm{GL}(n)$ with (the equivalence class of) a family of complex characteristic polynomials. In general, however, the mapping is not injective, for a variety of reasons that range from obvious to deep. A full understanding of the fibres of the mapping remains an important unsolved problem.

Consider the category whose objects are split reductive groups $G$ over $\mathbb{Q}$, and for which the morphisms from $G'$ to $G$ are the complex homomorphisms from $\widehat{G}'$ to $\widehat{G}$. To any object $G$, we can associate the set of (equivalence classes of) families $\mathscr{C}_{\mathrm{aut}}(G)$. For any morphism $\rho' : \widehat{G}' \to \widehat{G}$, we can associate the function

$$\mathscr{C}_{\mathrm{aut}}(\rho') : c' \to c = \rho'(c'), \quad c' \in \mathscr{C}_{\mathrm{aut}(G')},$$

postulated by the Principle of Functoriality. The correspondence $\mathscr{C}_{\mathrm{aut}}$ is then a functor from the category of groups we have just defined to the category of sets. This is essentially the origin of the term functoriality.

In his original paper [138], Langlands formulated what would become global functoriality in the setting of a general (connected) reductive algebraic group $G$ over a number field $F$. We should say something about the general analogue of the assertion of functoriality above for split groups.

A general group $G$ over $F$ can be constructed from a canonical quasi-split group $G^*$ over $F$ by an inner twist of its Galois action, a standard technique in the theory of algebraic groups. The group $G^*$ is obtained in turn from a canonical split group $G^{*,\mathrm{spl}} = G^{\mathrm{spl}}$ by an outer twist of its Galois action, attached to a homomorphism from a finite Galois group $\Gamma_{K/F} = \mathrm{Gal}(K/F)$ into the group $\mathrm{Out}(G^{\mathrm{spl}})$ of outer homomorphisms of $G^{\mathrm{spl}}$. It is the dual of this action on $\widehat{G^*}$ that gives the $L$-group

$$^L G^* = \widehat{G^*} \rtimes \Gamma_{K/F}, \quad \Gamma_{K/F} = \mathrm{Gal}(K/F)$$

described for $F = \mathbb{Q}$ earlier in this section. Since the dual $\widehat{G}$ of $G$ equals the dual $\widehat{G^*}$ of $G$ (as complex groups with actions of $\Gamma_{K/F}$), the $L$-groups $^L G = \widehat{G} \rtimes \Gamma_{K/F}$ and $^L G^* = \widehat{G^*} \rtimes \Gamma_{K/F}$ of $G$ and $G^*$ are equal. In particular, the general version of functoriality postulated by Langlands in [138] includes the case of the identity map $\rho$ from $^L G$ to $^L G^*$, and hence a correspondence from the automorphic representations of $G$ to those of the quasi-split group $G^*$. These days, this question is generally regarded as part of endoscopy, a separate theory proposed later by Langlands that seeks among other things to describe the precise nature of this correspondence. We are therefore free to consider functoriality in the more restricted context of quasi-split groups. (We shall discuss Langlands' theory of endoscopy in Section 10.).

Suppose then that $G$ is a quasi-split group over a number field $F$, which for the moment we may again take to be $\mathbb{Q}$. Our discussion above for a split group carries over with little change. In particular, an automorphic representation $\pi$ of $G(\mathbb{A})$ has a decomposition (24) such that $\pi_p$ is unramified for any $p$ outside a finite set $S$ [72]. The classification of unramified representations for split groups extends, with one minor adjustment. In the quasi-split case, the bijection $\pi_p \to c_p$ takes the unramified representations of $G(\mathbb{Q}_p)$ onto the set of $\widehat{G}$-orbits (under conjugation) in $^L G = \widehat{G} \rtimes \Gamma_{K/F}$ that map to the Frobenius class $\Phi_p$ in $\Gamma_{K/F}$. This really is a generalization of the bijection from the special case of split groups. For we recall in general that $K/\mathbb{Q}$ can be any Galois extension through which the Galois action on $\widehat{G}$ factors. In particular, if $G$ is split, we could take $^L G = \widehat{G} \times \Gamma_{K/F}$, a direct product of $\widehat{G}$ with the Galois group $\Gamma_{K/F} = \mathrm{Gal}(K/\mathbb{Q})$ of any convenient finite Galois extension $K/\mathbb{Q}$. In this case, $c_p(\pi)$ would just be the product of a semisimple conjugacy class in $\widehat{G}$ with an element in the Frobenius class $\Phi_p$ in $\Gamma_{K/F}$.

The assertion of Langlands' global functoriality conjecture then carries over to the quasi-split group $G$ as stated, with one proviso. The mapping

$$\rho' \colon {}^L G' \to {}^L G$$

between $L$-groups must be an $L$-homomorphism, by which we mean that it commutes with the projections of $^L G'$ and $^L G$ onto $\Gamma_{K/F}$. With this understanding, we can also extend the category we have defined to quasi-split groups, and the functor $\mathscr{C}_{\mathrm{aut}}$ from this category to sets. Finally, everything we have described in this section remains valid if we allow the base field $F$ to be an arbitrary number field. We took $F = \mathbb{Q}$ in order that the technical background might seem a little more concrete.

There is one further matter that requires comment. We could have referred to the Langlands conjecture above as "unramified global functoriality", since the term functoriality by itself often connotes a correspondence $\pi' \to \pi$ of automorphic representations that applies to the ramified local places $v \in S$ as well as the unramified $v \notin S$. This is certainly how Langlands introduced it in [138] (without the name). The resulting global assertion is more complicated, and presupposes Langlands' local functoriality. Langlands actually began with a question on general automorphic $L$-functions, both local and global. He then presented local and global functoriality as questions motivated by a desire to understand the $L$-functions. We shall discuss these matters in the next section.

In this section, we have presented unramified global functoriality as the primary assertion in order to give a concrete statement, and also to follow the natural progression begun with our discussion of abelian class field theory in the last section. We must now define the corresponding unramified automorphic $L$-function.

Suppose again that $G$ is a quasi-split reductive group over a number field $F$. An automorphic representation $\pi$ of $G(\mathbb{A}_F)$ gives a family

$$c^S(\pi) = \{c_v(\pi) : v \notin S\}$$

of semisimple conjugacy classes in $^LG$. To attach an automorphic $L$-function to $\pi$, we would need a family of semisimple conjugacy classes in a general linear group $\mathrm{GL}(n, \mathbb{C})$ rather than in $^LG$. We therefore fix a finite dimensional representation

$$r\colon {}^LG \to \mathrm{GL}(n, \mathbb{C})$$

in addition to the automorphic representation $\pi$ of $G(\mathbb{A})$. We then define the associated unramified automorphic $L$-function as the Euler product

$$L^S(s, \pi, r) = \prod_{v \notin S} L_v(s, \pi, r) = \prod_{v \notin S} \det(1 - r(c_v(\pi))q_v^{-s})^{-1}, \tag{29}$$

with $q_v$ being the degree of the residue class field of $F_v$. Once again, the product converges absolutely to an analytic function in some right half plane. The analogy with the Artin $L$-function (13) is clear. We can think of it as the same definition, but $\mathrm{Gal}(K/F)$ replaced by the group $^LG$, along with the extra structure provided by $\pi$.

In addition to defining automorphic $L$-functions and introducing the Principle of Functoriality, Langlands sketched the following four applications in his seminal paper [138].

**1. Analytic continuation and functional equation.** Langlands pointed out that the analytic continuation and functional equation of a general automorphic $L$-function would follow from functoriality and the special case that $G = \mathrm{GL}(n)$ and $r = St(n)$, the standard $n$-dimensional representation of $GL(n)$. This special case (at least for cuspidal $\pi$) was established soon afterwards by Godement and Jacquet [81], as we saw in the last section.

**2. Artin $L$-functions.** We have noted that quasi-split groups are the natural setting for functoriality. The Galois factor $Gal(K/\mathbb{Q})$ is then an essential part of the $L$-group $^L G$. In particular, the construction naturally includes the seemingly trivial case that $G$ is the 1-element group $\{1\}$. Its $L$-group will then be an arbitrary finite Galois group $Gal(K/\mathbb{Q})$, while $r$ becomes simply an $n$-dimensional representation of $Gal(K/\mathbb{Q})$. The associated automorphic $L$-function $L(s, \pi, r)$ (with $\pi$ being of course the trivial 1-dimensional automorphic representation of $G$) is then just the general Artin $L$-function $L^S(s, r)$. The Principle of Functoriality can thus be interpreted as an identity

$$L^S(s, r) = L^S(s, \pi, St(n)) \tag{30}$$

between a general Artin $L$-function and a standard automorphic $L$-function for $GL(n)$. This represents a general and entirely unexpected formulation of nonabelian class field theory. It identifies purely arithmetic objects, Artin $L$-functions, with objects associated with harmonic analysis, automorphic $L$-functions, thereby implying that the arithmetic $L$-functions have meromorphic continuation and functional equation, and that they are essentially entire. Functoriality would thus include a proof of the Artin conjecture stated in the last section. Abelian class field theory amounts to the special case that the dimension $n$ of $r$ equals 1. Its original aim was to establish that abelian $L$-functions are the Hecke–Tate $L$-functions attached to the automorphic representations of $GL(1)$, and thereby have analytic continuation and functional equation.

**3. Generalized Ramanujan conjecture.** The generalized Ramanujan conjecture asserts that a unitary cuspidal automorphic representation $\pi = \widetilde{\bigotimes}_v \pi_v$ of $GL(n)$ is *locally tempered*. This means that the character

$$f_v \to \mathrm{tr}(\pi(f_v)), \quad f_v \in C_c^\infty(GL(n, F_v)),$$

of each local constituent $\pi_v$ of $\pi$ is tempered, in the sense that it extends to a continuous linear form on the Schwartz space $\mathscr{C}(GL(n, F_v))$ of $GL(n, F_v)$ defined by Harish-Chandra. We recall that the classical Ramanujan conjecture applies to the case that $n = 2$, and that $\pi$ comes from the cusp form of weight 12 and level 1. It was proved by Deligne [63], who established more generally (for $n = 2$) that the conjecture holds if $\pi$ is attached to any holomorphic cusp form. (The case that $\pi$ comes from a Maass form remains an important open problem.) Langlands observed that functoriality, combined with expected properties of the correspondence $\pi' \to \pi$, would imply the generalized Ramanujan conjecture for $GL(n)$. His representation theoretic argument is strikingly similar to Deligne's geometric proof [63].

**4. Sato–Tate conjecture.** The Sato–Tate conjecture for the distribution of numbers $N_p(E)$ of solutions (mod $p$) of an elliptic curve $E$ over $\mathbb{Q}$ has a general analogue for automorphic representations. Suppose for example that $\pi$ is a (unitary) cuspidal automorphic representation of $GL(n)$. The generalized Ramanujan conjecture of 3. above asserts that the conjugacy classes, represented by diagonal $S_n$-orbits

$$c_p(\pi) = \begin{pmatrix} c_{p,1}(\pi) & & 0 \\ & \ddots & \\ 0 & & c_{p,n}(\pi) \end{pmatrix} / S_n,$$

have eigenvalues of absolute value 1. The generalized Sato–Tate conjecture describes their distribution in the maximal torus $U(1)^n$ of the maximal compact subgroup $U(n)$ of the dual group $GL(n, \mathbb{C})$. If $\pi$ is *primitive* (a notion that requires functoriality even to define, as we will describe in Section 9), the distribution of these classes should be given by the weight function in the Weyl integration formula for the unitary group $U(n)$. Langlands sketched a rough argument for establishing such a result from general functoriality. Clozel, Harris, Shepherd-Barron and Taylor followed this argument in their proof of the original Sato–Tate conjecture, but using base change for $GL(n)$ and deformation results in place of functoriality. (See [91], [238].)

# 5 Local Functoriality and early results

Within a short time of his introduction [138] of functoriality, Langlands presented some striking results that in addition to the interest they held in their own right, offered evidence for the general principle. They are contained in four[3] major works: a long, unpublished manuscript [137] related to the local properties of Artin $L$-functions, a classification of representations of algebraic tori [153], the properties of Euler products from Eisenstein series [139] mentioned in Section 3, and the monograph [103] of Jacquet and Langlands on GL(2). Each of these represents a deep and original contribution to the subject, even if like functoriality itself, each one may have been ahead of its time.

We shall describe the four contributions in sequence in the latter part of the section. A later contribution, Langlands' monograph [149] on base change for GL(2), with its dramatic applications to Artin's conjecture, will be treated separately in Section 7. In the first part of this section, we shall discuss the local aspects of the Principle of Functoriality. We should also take the opportunity to add some further remarks about the structure of Langlands' fundamental paper [138].

It is important to remember that the Langlands conjectures were so revolutionary in 1970 that they took many years to be accepted (or even noticed) by the general mathematical public. The paper [138] was dense and difficult. It consisted of seven questions (not conjectures), three local and three global. The other question, which came first, was both local and global. It concerned the possibility of defining (completed) automorphic $L$-functions. It was actually presented as the central problem, with six supplementary questions related to functoriality being a strategy for attacking it. This is the opposite of our exposition here, in which (unramified) global functoriality was presented as the centre of the Langlands program. However, the logic of [138] is compelling. As we described in the last section, it represents a (nonabelian) analogue for automorphic representations of Artin's use of abelian class field theory to establish analytic continuation and functional equation for his (abelian) $L$-functions.

Suppose that $G$ is a reductive group over a number field $F$. We then have the global $L$-group

$$^{L}G = \widehat{G} \rtimes \mathrm{Gal}(K/F),$$

where $(K/F)$ is a finite Galois extension over which $G$ splits. For any valuation $v$ of $F$, we have a local Galois extension $(K_v/F_v)$, for which $\mathrm{Gal}(K_v/F_v)$ represents a conjugacy class of subgroups of $\mathrm{Gal}(K/F)$. We also have the local $L$-group

$$^{L}G_v = \widehat{G} \rtimes \mathrm{Gal}(K_v/F_v),$$

which comes with a conjugacy class of $L$-embeddings $^{L}G_v \hookrightarrow {}^{L}G$. At this point, we are free to take motivation from the special case of automorphic $L$-functions for $\mathrm{GL}(n)$ discussed near the end of Section 3, and the unramified automorphic

---

[3] We could have added Langlands' classification of representations of real groups [151] mentioned in Section 1, as fifth work, but we are postponing its further discussion until Section 10.

$L$-functions $L^S(s,\pi,r)$ from Section 4. In particular we fix a finite dimensional representation

$$r\colon {}^L G \to \mathrm{GL}(n,\mathbb{C}),$$

a nontrivial additive character

$$\psi\colon F^* \setminus \mathbb{A}_F^* \to \mathbb{C}^*,$$

and an automorphic representation $\pi$ of $G$, objects that come with natural localizations $r_v$, $\psi_v$ and $\pi_v$. Langlands' first question was then as follows.

**Question 1 (Langlands [138])** *Is it possible to define local L-functions*

$$L_v(s,\pi,r) = L(s,\pi_v,r_v)$$

*and local $\varepsilon$-factors*

$$\varepsilon_v(s,\pi,r,\psi) = \varepsilon(s,\pi_v,r_v,\psi_v) = \varepsilon(\pi_v,r_v.\psi_v)q_v^{-(n_v - \frac{1}{2}s)},$$

*for the valuations $v \in S$, with $(G_v, r_v, \pi_v, \psi_v)$ ramified (or archimedean), such that the full Euler product*

$$L(s,\pi,r) = L_S(s,\pi,r)L^S(s,\pi,r) = \prod_v L_v(s,\pi,r)$$

*has analytic continuation, and satisfies the functional equation*

$$L(s,\pi,r) = \varepsilon(s,\pi,r)L(1-s,\pi,r^\vee), \tag{31}$$

*for*

$$\varepsilon(s,\pi,r) = \prod_{v \in S} \varepsilon_v(s,\pi,r,\psi)?$$

It is remarkable to see Langlands' comments immediately following his statement of Question 1 [138, p. 31–32]. He devotes one sentence each to: the motivation [139] he took from Eisenstein series, which we mentioned in Section 3, and which led to the Langlands–Shahidi method; his anticipation of the extension of the thesis of Tate from $\mathrm{GL}(1)$ to $\mathrm{GL}(n)$, established shortly thereafter [81]; and most striking of all, his description of his subsequent questions on functoriality as the nonabelian analogue of "the idea that led Artin to the general (abelian) reciprocity law". It is like reading fifty years of past and future history unfold in one short paragraph.

Having already stated the unramified version of global functoriality, we shall not restate all of Langlands' remaining questions . Questions 2 and 3 concern the local and global correspondence between the representations of a non-quasi-split inner twist of $G$ and those of $G$, which are now regarded as part of the conjectural theory of endoscopy. It is the next two Questions 4 and 5 that introduce functoriality. To state the first of these as local functoriality, we take quasi-split groups $G_v'$ and $G_v$ over a localization $\mathbb{Q}_v$ of $\mathbb{Q}$, together with a local $L$-homomorphism

$$\rho'_v \colon {}^L G'_v \to {}^L G_v$$

between their $L$-groups.

**Conjecture (Langlands' Principle of Local Functoriality)** *There is a natural correspondence*

$$\pi'_v \to \pi_v, \quad \pi'_v \in \Pi(G'_v),$$

*from the (equivalence classes of) irreducible representations of $G'_v$ to those of $G_v$ such that*

$$L(s, \pi_v, r_v) = L(s, \pi'_v, \rho'_v \circ r_v)$$

*and*

$$\varepsilon(s, \pi_v, r_v, \psi_v) = \varepsilon(s, \pi'_v, \rho'_v \circ r_v, \psi_v),$$

*for every finite dimensional representation $r_v$ of ${}^L G_v$ and every nontrivial additive character $\psi_v$ on $F_v$.*

Notice that unlike the statement of unramified Global Functoriality in Section 4, this conjecture is not rigid as stated, in that it does not characterize the correspondence. The property on the $L$- and $\varepsilon$- factors, which is based on a suitable answer to Question 1, represents a condition it must satisfy. Another condition is the global expectation that if $G'_v$, $G_v$ and $\rho'_v$ are localizations of global objects $G'$, $G$ and $\rho$ for each $v$, and if $\pi' = \bigotimes_v \pi'_v$ is an automorphic representation of $G'(\mathbb{A}_F)$, one can then choose a representation $\pi_v \in \Pi(G_v)$ in the image of the correspondence $\pi'_v \to \pi_v$ for each $v$ such that the product $\pi = \bigotimes_v \pi_v$ is an automorphic representation of $G$. Langlands actually imposes a stronger condition in Question 5, his general form of Global Functoriality. Namely, if $\pi_v$ is *any* representation in the image of the correspondence $\pi'_v \to \pi_v$ for each $v$, is the product an automorphic representation of $G(\mathbb{A}_F)$? He soon realized that this was asking for a little too much. Langlands' later theory of endoscopy, which has now been established in a number of cases [23] and which we will discuss in Section 10, characterizes both the correspondence $\pi'_v \to \pi_v$, and which products attached to a given $\pi'$ are actually automorphic representations of $G(\mathbb{A})$.

Questions 6 and 7 represent an important extension of functoriality to Weil groups. My understanding is that Langlands first took his conjectures to Weil in the form of Questions 1–5, as an automorphic analogue of Artin $L$-functions and the Artin reciprocity law. Weil pointed out that he had some years earlier introduced a natural generalization of Artin $L$-functions. Langlands was no doubt happy to find that his questions/conjectures extended seamlessly to Weil's $L$-functions, and that they were in fact the richer for it. Of the four major works that followed, and that we will be discussing presently, three depend intimately on the extensions of Galois groups Weil had discovered in abelian class field theory.

In the interest of making this report seem more concrete, we have sometimes avoided the most general setting. Let us now try to be a little more efficient. From this point on, we shall work over an arbitrary base field $F$ of characteristic 0, local or global. Everything discussed earlier in the case of $\mathbb{Q}$ then carries over as stated with $F$ in place of either $\mathbb{Q}$ or one of its completions $\mathbb{Q}_v$. We can also work over the

absolute Galois group

$$\Gamma_F = \varprojlim \Gamma_{K/F}, \quad \Gamma_{K/F} = \mathrm{Gal}(K/F),$$

instead to the finite Galois group $\Gamma_{K/F}$, if we agree that earlier homomorphisms $r$ and $\rho'$ are to be continuous. This is of course because the kernel of $r$ and $\rho'$ in the totally disconnected, compact group $\Gamma_F$ is a normal subgroup of $\Gamma_K$ of finite index, so that the quotient

$$\Gamma_F / \Gamma_K = \Gamma_{K/F}$$

becomes the Galois group of the finite Galois extension $K/F$. We shall freely adopt this convention in referring to past discussion, usually without further comment.

*Weil groups* are variants $W_{K/F}$ and $W_F$ of the local and global Galois groups $\Gamma_{K/F}$ and $\Gamma_F$ (with $K/F$ still being a finite Galois extension). We set $C_F$ equal to the multiplicative group $F^*$ if $F$ is local, and to the idele class group $F^* \backslash \mathbb{A}_F^*$ if $F$ is global. The relative Weil group $W_{K/F}$ is then an extension

$$1 \to C_K \to W_{K/F} \to \Gamma_{K/F} \to 1$$

in both cases. It is a locally compact group, which comes with an isomorphism $r_{K/F} \colon C_F \to W_{K/F}^{\mathrm{ab}}$ in addition to its projection onto $\Gamma_{K/F}$, while if $K = F$, $W_{K/F}$ obviously equals $C_F$.

In either the local or global setting, consider a field $E$ with $F \subset E \subset K$. Then $W_{K/E}$ is a subgroup of finite index in $W_{K/F}$, and there is a bijection

$$W_{K/F} / W_{K/E} \cong \Gamma_{K/F} / \Gamma_{K/E}$$

of finite coset spaces. If $K/E$ is Galois, the subgroup $W_{K/E}$ is then normal in $W_{K/F}$, and the quotient on the left is isomorphic to the Galois group $\Gamma_{E/F}$ that equals the quotient on the right. To obtain the Weil group $W_{E/F}$ as a quotient in this case, we need to take the commutator subgroup $W_{K/E}^c$ instead of $W_{K/E}$. We then obtain an isomorphism

$$W_{K/F} / W_{K/E}^c \cong W_{E/F}.$$

In particular, we have a continuous projection from $W_{K/F}$ onto $W_{E/F}$.

The family

$$\{W_{K/F} : K/F \text{ Galois}\}$$

is thus an inverse system. The *absolute Weil group* of $F$ is the corresponding inverse limit

$$W_F = \varprojlim_K W_{K/F}.$$

It is a locally compact group, equipped with a continuous homomorphism

$$\phi_F \colon W_F \to \Gamma_F$$

with dense image, and an isomorphism

$$r_F : C_F \to W_F^{\mathrm{ab}}.$$

If $F = \mathbb{C}$, $W_F = W_{\mathbb{C}/\mathbb{C}} = \mathbb{C}^*$. If $F = \mathbb{R}$, $W_F = W_{\mathbb{C}/\mathbb{R}}$ is the group generated by $\mathbb{C}^*$ and an element $w$, subject to the relations $w^2 = -1$, and $wzw^{-1} = \bar{z}$ for any $z \in \mathbb{C}^*$. If $F$ is a local nonarchimedean field, $W_F$ turns out to be the dense subgroup of elements $\Gamma_F$ whose image in the quotient

$$\Gamma_{F,\mathrm{un}} \cong \Gamma_F / I_F = \widehat{\mathbb{Z}}$$

of $\Gamma_F$ by its inertia subgroup $I_F$ is the dense subgroup $\mathbb{Z}$ of $\widehat{\mathbb{Z}}$. Finally, if $F$ is global, $W_F$ is given by a more complicated blend of these properties.

We refer the reader to the beginning [237, §1] of the article of Tate in the Corvallis proceedings. He takes the absolute Weil groups $\{W_F\}$ as the basic objects, equipped with mappings $\phi_F$ and $r_F$ as above that satisfy four axioms. He then defines the relative Weil groups as quotients

$$W_{K/F} = W_F / W_K^c.$$

Tate also comments very briefly [237, (1.2)] on how class field theory is implicit in the existence of Weil groups, specifically the algebraic derivation of local and global class field theory in terms of Galois cohomology.

Observe that any continuous, complex finite dimensional representation of $\Gamma_F$ pulls back to a unique continuous representation of $W_F$. Representations of Weil groups are then more general than those of Galois groups. It is for this reason that mathematicians have taken to working with the absolute Weil form

$$^L G = \widehat{G} \rtimes W_F$$

of the $L$-group of a quasi-split group $G$ over $F$, rather than the absolute Galois form $\widehat{G} \rtimes \Gamma_F$, or the original, less canonical form $\widehat{G} \rtimes \Gamma_{K/F}$. It is of course understood that the action of $W_F$ on $\widehat{G}$ is through the pullback to $\Gamma_F$ of a suitable finite quotient $\Gamma_{K/F}$. For the assertions of functoriality, one would want to take the absolute Weil form of both $L$-groups $^L G'$ and $^L G$, with $\rho'$ being an $L$-homomorphism from $^L G'$ to $^L G$, in the obvious sense that it commutes with the two projections onto $W_F$. One also has to insist that the $\rho'$-image in $\widehat{G}$ of any element in $W_F$ be semisimple, since $W_F$ is only locally compact. We shall again feel free to extend past discussions to this general context, without comment.

Langlands' final two Questions 6 and 7 apply to the Weil form of local and global functoriality, or rather the special case in which the first group $G'$ equals the trivial group $\{1\}$. In other words, $\phi = \rho'$ is an $L$-homomorphism from $W_F$ to $^L G$. (Langlands also replaces the complex dual group of $G$ by a compact real form, in anticipation perhaps of his later comments on the generalized Ramanujan conjecture.) The two questions ask for a correspondence from $L$-homomorphisms $\phi : W_F \to {}^L G$ to representations of $G(F)$ if $F$ is local and to automorphic representations of $G(\mathbb{A}_F)$ if $F$ is global. This is where we see the advantage of the role of the Weil group. For there are many more such homomorphisms than there would be for the Galois group $\Gamma_F$. The case of a local field $F$ is of particular importance. Indeed, the irreducible representations of $G(F)$ obtained in this way are believed to account for most (though not all) such representations. In fact, Question 6 has evolved into what is now known as the *local Langlands classification* (or the *local Langlands correspondence*). It has been established in a significant number of cases, for $G = \mathrm{GL}(n)$ in [92], [94] and [202], and for quasi-split classical groups in [23] and [181], even as it remains conjectural in general. (For the archimedean case, it is the original form of Question 6 that is relevant. It gives the Langlands classification for the real group, which we discussed in Section 1, and to which we will return later.)

We turn now to the four works mentioned in the beginning of the section. With many other things to discuss, we shall have to be brief, despite the importance of the results. We shall give a short description in each case, with further comments as needed later in the report.

**1. Artin $L$-functions [137].** The reader might have noticed an irregularity in our claim of symmetry in Section 3 between $n$-dimensional Artin $L$-functions and the automorphic $L$-functions for $\mathrm{GL}(n)$ of Godement and Jacquet. The $\varepsilon$-factor in Artin's functional equation (15) is a nonconstructive global function, while its automorphic counterpart in (19) has a finite product decomposition (20) into purely local functions. The problem was to find a corresponding local decomposition for Artin $\varepsilon$-factors and of course, their generalizations for Weil groups.

There were good reasons for trying to do this. Langlands was thinking of the possibility of classifying $\Pi(G)$, the set of irreducible representations of a real or $p$-adic group $G(F)$. The key to this would be local functoriality, as stated earlier in this section, but with $G' = \{1\}$ and with the local Weil group $W_F$ in place of $\Gamma_F$. In other words, $\Pi(G)$ should be closely tied to the $L$-homomorphisms

$$\phi \colon W_F \to {}^L G.$$

An essential condition would then be that for every representation $r$ of ${}^L G$, the local $L$- and $\varepsilon$-functions $L(s, \pi, r)$ and $\varepsilon(s, \pi, r, \psi)$ for a representation $\pi$ corresponding to $\phi$, which were conjectured for any $G$ in Question 1 and established for $G = \mathrm{GL}(n)$ and $r = \mathrm{St}(n)$ by Godement–Jacquet, would match independently defined functions $L(s, r \circ \phi)$ and $\varepsilon(s, r \circ \phi, \psi)$ for $\phi$.

The local $L$-functions were already part of the Artin/Weil definition, but the construction of local $\varepsilon$-factors $\varepsilon(s, r \circ \phi, \psi)$ turned out to be very difficult. It was studied with some success by B. Dwork [70], [131] for Artin $L$-functions. Langlands used Dwork's results in his investigation of Weil $L$-functions. His goal was to characterize the local $\varepsilon$-factors

$$\varepsilon(s, r, \psi), \qquad r \colon W_F \to \mathrm{GL}(n, \mathbb{C}), \qquad F \text{ local},$$

as the unique family of functions that satisfies several natural conditions as $F$, $r$ and $\psi$ vary. This question gave rise in turn to four concrete, but very complex lemmas on Gaussian sums. My understanding is that Langlands obtained a complete solution, but that it was not all contained in the manuscript [137], which itself was not published. (See [137, (3.4.1)], and Langlands' later comments on [137] in his 1969 letter [136] to Deligne.)

Langlands' proof was purely local. Some time later, Deligne [62] found a striking global argument that led to a much shorter proof. The result was the desired canonical construction of local $\varepsilon$-factors

$$\varepsilon(s, r, \psi) = \varepsilon(r, \psi) q_F^{-n_F(s - \frac{1}{2})},\tag{32}$$

for integers $n_F = n(r, \psi)$ and $q_F$, and any nontrivial additive character $\psi$ on the multiplicative group $F^*$ of the local field $F$. The global result remains the functional equation (15) for any Artin/Weil $L$-function

$$L(s, r),\qquad\qquad r\colon W_F \to \mathrm{GL}(n, \mathbb{C}),\qquad F \text{ global},\tag{33}$$

but with the global $\varepsilon$-factor in (15) now having a product decomposition

$$\varepsilon(s, r) = \prod_{v \in S} \varepsilon(s, r_v, \psi_v)$$

into purely local $\varepsilon$-factors.

**2. Automorphic representations of tori [153].** This paper gives a classification of representations (local and global) for *abelian* reductive algebraic groups, which is to say, for (algebraic) tori. Recall that a *torus* is an algebraic group that is isomorphic to a product $\mathrm{GL}(1)^k$ of multiplicative groups. A torus over $F$ (our local or global field) comes also with an outer twisting over $F$, namely a homomorphism from a finite Galois group $\mathrm{Gal}(K/F)$ into the group of automorphisms of this product. The correspondence

$$T \to X^*(T) = \mathrm{Hom}(T, \mathrm{GL}(1))$$

gives an antiisomorphism from the category of algebraic tori over $F$ that split over the finite Galois extension $K$, and the category of torsion free $\mathrm{Gal}(K/F)$-modules of finite rank. (See [233], for example.)

The paper [153] here is not nearly so complex as the last one [137]. It is still very interesting, a verification of functoriality in the simplest nontrivial case, and an elegant illustration of some of the ideas of Langlands implicit in [138]. It applies to the Weil form of functoriality in our description of Questions 6 and 7, specifically the special case that $G$ is an algebraic torus over $F$, a local of global field (of characteristic 0), and that $G'$ as before is the trivial group $\{1\}$. Functoriality then concerns the $L$-homomorphisms

$$\phi\colon W_F \to {}^L T = \widehat{T} \rtimes W_F.$$

The purpose of [153] is to establish functoriality in this case, and to establish further that the resulting correspondence is surjective.

To be more precise, we fix a local or global field $F$, and we write $\Phi(T)$ for the set of $\widehat{T}$-conjugacy classes of $L$-homomorphisms $\phi$. As we have done earlier, we also write $\Pi(T)$ for the set of equivalence classes of irreducible representations of $T(F)$ if $F$ is local, and of automorphic representations of $T(\mathbb{A}_F)$ if $F$ is global. In the global case we have the localization mappings $\phi \to \phi_v$ and $\pi \to \pi_v$, from $\Phi(T)$ to $\Phi(T_v)$ and $\Pi(T)$ to $\Pi(T_v)$ respectively. Theorem 2 is the main result of [153]. It asserts that there is a canonical mapping

$$\Phi(T) \to \Pi(T), \quad \phi \to \pi,$$

which is a bijection if $F$ is local, and a surjection if $F$ is global. In the global case, the mapping commutes with the associated localizations, and its fibres are the local equivalence classes in $\Phi(T)$, relative to the equivalence relation $\phi' \sim \phi$ in $\Phi(T)$ if $\phi'_v \sim \phi_v$ in $\Phi(T_v)$ for every localization $F_v$ of $F$.

Langlands describes his proofs as "exercises in class field theory". Observe that a mapping $\phi$ in $\Phi(T)$ is completely determined by the projection of its image onto $\widehat{T}$. This leads to an isomorphism from $\Phi(T)$ (as an abelian group) onto the Galois cohomology group $H^1(W_F, \widehat{T})$ of continuous 1-cocycles from $W_F$ to $\widehat{T}$, modulo continuous coboundaries. Such is the stuff of class field theory, in its algebraic form. It would indeed be a good exercise to study the proofs of Langlands, streamlined perhaps according to Tate–Nakayama duality [234].

**3. Euler products [139].** We discussed this paper at the end of Section 3, as historical motivation for the Principle of Functoriality. It also represents concrete evidence for functoriality, specifically the properties of $L$-functions in Question 1.

Suppose that $G$, $P$, $\pi$ and $r$ are as in the formula (23) for the global intertwining operator (21) in terms of what was then the new $L$-function $L(s, \pi, r^\vee)$. We write (23) as

$$L(s, \pi, r^\vee) = L(s+1, \pi, r^\vee) M(w, \lambda),$$

where $M(w, \lambda)$ now represents a meromorphic scalar valued function of the image $s \in \mathbb{C}$ of $\lambda$, according to the notation in Section 3 and the first assertion (a) of the Main Theorem in Section 2. We see from this that if $L(s, \pi, r^\vee)$ is meromorphic in a right half place $\mathrm{Re}(s) > b$, it continues to a meromorphic function of $\mathrm{Re}(s) > b-1$. Since we know that the Euler product for $L(s, \pi, r^\vee)$ converges in some right half plane, we conclude that $L(s, \pi, r^\vee)$ does have analytic continuation to a meromorphic function of $s$ in the complex plane.

What of the proposed functional equation

$$L(s, \pi, r) = L(1-s, \pi, r^\vee)$$

conjectured in Question 1 of [138]? The analogue of the formula (23) for the intertwining operator

$$M(w^{-1}, w\lambda)\colon \mathscr{H}_{P'} \to \mathscr{H}_P, \quad w^{-1} \in W(\mathfrak{a}_{P'}, \mathfrak{a}_P),$$

is

$$M(w^{-1}, w\lambda) = \frac{L(-s, \pi, r)}{L(-s+1, \pi, r)},$$

as one sees from the definitions [139, p. 47]. In general, the operators $M(w,\lambda)$ satisfy their own functional equation (9). In the case at hand, this becomes

$$M(w^{-1}, w\lambda)M(w,\lambda) = M(w^{-1}w,\lambda) = M(1,\lambda) = 1.$$

This gives the identity

$$\frac{L(-s,\pi,r)}{L(-s+1,\pi,r)}\frac{L(s,\pi,r^\vee)}{L(s+1,\pi,r^\vee)} = 1$$

of meromorphic functions of $s$. On the other hand, if we substitute the conjectured functional equation for the first factor in the left, the product on left becomes

$$\frac{L(1+s,\pi,r^\vee)}{L(s,\pi,r^\vee)}\frac{L(s,\pi,r^\vee)}{L(s+1,\pi,r^\vee)},$$

which is also equal to 1. In other words, the conjectural functional equation for $L(s,\pi,r)$ is compatible with the established functional equation for $M(w,\lambda)$, but is not implied by it. Nonetheless, this represents further evidence from Eisenstein series for Langlands' theory of automorphic $L$-functions and the Principle of Functoriality.

These observations were part of Langlands' article [139]. He did not assume that the dual unipotent radical $\hat{N}$ was abelian, although $M$ still represents a maximal Levi subgroup. For $G$, $P = MN$ and $\pi$ as above, but without this assumption on $\hat{N}$, the Lie algebra of $\hat{N}$ has a decomposition

$$\hat{\mathfrak{n}} = \bigoplus_{i=1}^{k}\hat{\mathfrak{n}}_i,$$

for

$$\hat{\mathfrak{n}}_i = \{U \in \hat{\mathfrak{n}} : \mathrm{ad}(\varpi_P^\vee)U = iU\},$$

and for $1 \leq k \leq 6$ (as is well known). The generalization of (23) is the formula

$$M(w,\lambda) = \prod_{i=1}^{k}\frac{L(is,\pi,r_i^\vee)}{L(is+1,\pi,r_i^\vee)}, \quad \lambda \to s,$$

where $r_i$ is the adjoint representation of $\hat{M}$ in $\hat{\mathfrak{n}}_i$. (We have been following the notation of [78, (1.2.5.3)] which is slightly simpler than that of Langlands at the end of Section 5 of [139].) If one can show that the $L$-functions $L(s,\pi,r_i^\vee)$ for $M$ have analytic continuation for $1 \leq i < k$, the argument above establishes that the same is true for $L(s,\pi,r_k^\vee)$. The extent to which this is possible is governed by the relevant Coxeter–Dynkin diagrams. At the end of the article [139], Langlands gives an extended table of such diagrams, which establish that for all but three simple groups $M$, there is at least one nontrivial representation $r$ of $\hat{M}$ for which the function $L(s,\pi,r)$ has meromorphic continuation.

Some years later, Shahidi began a sustained study that greatly expanded the theory [211], [214], [212], [213]. (See also [45], [44]). With an assumption on Whittaker models (well known for general linear groups $M$), he was able to treat general cuspidal automorphic representations $\pi \in \Pi_{\text{cusp}}(M)$ (without the condition that they be unramified everywhere). The resulting theory, known now as the Langlands–Shahidi method, has established functional equations in addition to the meromorphic continuation for many of the $L$-functions $L(s, \pi, r)$ attached to Eisenstein series. A further examination of some special cases led to a remarkable application to functoriality. Partly in collaboration with H. Kim (the case $n = 5$ below), Shahidi established functoriality for any $\pi' \in \Pi_{\text{cusp}}(G')$, where $G'$ equals GL(2), $G$ equals GL(4) or GL(5), and

$$\rho' \colon \widehat{G}' = \mathrm{GL}(2, \mathbb{C}) \to \widehat{G} = \mathrm{GL}(n, \mathbb{C}), \quad n = 4, 5,$$

is the symmetric cube or fourth power representation. This was thought to have been inaccessible, and has led to significant improvements in the bounds required by Ramanujan's conjecture for $\pi'$. (See [215], and the references therein.)

**4. Automorphic representations of** GL(2) **[103].** The 348-page monograph of Jacquet–Langlands was a partial exception to a remark from the beginning of this section. It really consists of two parts, the first twelve Sections 1–12 culminating in a striking application to Artin $L$-functions, and the last two Sections 15–16 on the comparison of representations of $G = \mathrm{GL}(2)$ with those of an inner twist, the multiplicative group $G'$ of a quaternion algebra. The two intermediate Sections 13–14 are in some sense transitional. They extend the methods from Tate's thesis from GL(1) to GL(2) in Section 13, anticipating the later volume [81] of Godement–Jacquet for GL($n$) we mentioned in our Section 3, and from GL(1) to $G'$ in Section 14, in anticipation of the comparisons in Sections 15 and 16. The first part of the monograph was widely read by mathematicians at the time, and quickly became a basic part of their thinking. The second part was slower to be taken up, but the principal result (as opposed perhaps to its proof) did soon find important applications among number theorists.

The main theme of the first part is an extension of Hecke theory to all automorphic representations of GL(2). Hecke was the first to attach $L$-functions to what amounted to a special class of cuspidal automorphic representations $\pi$ for GL(2). He showed that these $L$-functions have analytic continuation, with functional equation, to *entire* functions on $\mathbb{C}$. Hecke was of course working with classical holomorphic modular forms of weight $k$ for, let us say SL$(2, \mathbb{Z})$. It was to forms in this space that he attached his $L$-functions. Hecke also introduced the operators $\{T_n\}$ on the space that bear his name, and he proved that for any simultaneous eigenform of these operators, the associated $L$-function has an Euler product. The local components of the corresponding automorphic representation $\pi = \bigotimes_v \pi_v$ are then characterized at the $p$-adic places by the natural relation between the conjugacy classes $c(\pi_p)$ and the eigenvalues of $T_p$, and at the archimedean valuation $\infty$ by the requirement that $\pi_\infty$ be a discrete series representation corresponding to $k$ in the parametrization of Harish-Chandra.

Of particular relevance to [103] is the converse theorem Hecke then established. He showed that any Dirichlet series with certain properties, the main ones being an Euler product, the analytic continuation to an entire function of $s$, and an appropriate functional equation, gives rise to a cuspidal automorphic $L$-function $L(s, \pi)$. Unlike the original theorem, which was actually for modular forms of level $N$ (that is, for congruence subgroups $\Gamma_0(N)$, or equivalently automorphic representations $\pi$ with ramified components $\pi_p$ for $p|N$), Hecke's converse theorem really was restricted to forms for the full modular group $\Gamma_0(1) = \mathrm{SL}(2, \mathbb{Z})$. It took thirty years for it to be extended. In 1967, Weil [250] established a converse theorem for modular forms of level $N$, but with more sophisticated requirements on the given Dirichlet series. In so doing, he was able to make what became known as the Shimura–Taniyama–Weil conjecture on the modularity of elliptic curves considerably more precise.

The goal of the first part, Sections 1–12 of [103], was to extend the converse theorem of Hecke and Weil for $G = \mathrm{GL}(2)$ to any number field $F$ and to any cuspidal automorphic representation $\pi \in \Pi_{\mathrm{cusp}}(G)$. This was a bigger task than one might perhaps imagine. If $F$ is an imaginary quadratic extension of $\mathbb{Q}$, there are no archimedean discrete series representations $\pi_\infty \in \Pi_2(G_\infty)$ of $G(F_\infty)$, and the corresponding modular forms are all Maass forms. Even if $F = \mathbb{Q}$, one wants a theory that includes all Maass forms as well as holomorphic modular forms. This leaves no choice but to extend the classical results for modular forms to a full theory for automorphic representations of the adele group $G(\mathbb{A}_F)$. In particular, one must first establish a robust theory for the irreducible representations $\pi_v$ of the local components $\mathrm{GL}(F_v)$ of $G(\mathbb{A}_F)$.

Chapter I (Sections 1–8) was devoted to the local theory. For the archimedean local fields $F_v = \mathbb{R}$ or $F_v = \mathbb{C}$, the authors used basic results of Harish-Chandra to classify the irreducible representations $\pi_v$ of $G(F_v)$ (Sections 5–6). In the earlier sections (1–4), they studied the irreducible representations $\pi_v$ of the nonarchimedean groups $G(F_v)$ through Weil representations, Kirillov models and Whittaker models [249], [100], [110], [79]. They constructed various examples of representations in this case, in what amounts to be a partial classification. In all cases, they constructed local $L$-functions $L(s, \pi_v)$ and $\varepsilon$-factors $\varepsilon(s, \pi_v, \psi_v)$ for each $\pi_v$.

Chapter II (Sections 9–12) is concerned with global Hecke theory. In Theorem 11.1, the authors proved that for any $\pi \in \Pi_{\mathrm{cusp}}(G)$, the global $L$-function

$$L(s, \pi) = \prod_v L(s, \pi_v) \tag{34}$$

has analytic continuation to an *entire* function of $s$ and a functional equation

$$L(s, \pi) = \varepsilon(s, \pi) L(1 - s, \pi^\vee) \tag{35}$$

with

$$\varepsilon(s, \pi) = \prod_v \varepsilon(s, \pi_v, \psi_v). \tag{36}$$

In fact, they verified these properties for a larger family of $L$-functions $L(s, \omega \otimes \pi)$ in Corollary 11.2, where $\omega$ ranges over quasicharacters on $C_F = F^* \backslash \mathbb{A}_F^*$, following Weil's extension of Hecke's converse theorem.

The general converse theorem is Theorem 11.3 of [103]. In common with its predecessors, it asserts that the necessary conditions are sufficient. More precisely, any representation $\pi = \otimes_v \pi_v$ of $G(\mathbb{A}_F)$ that satisfies the conclusions of Theorem 11.1 and Corollary 11.2, together with two further necessary conditions, is actually a cuspidal automorphic representation of $G = \mathrm{GL}(2)$. One of the extra conditions is a bound, imposed in Theorem 11.3 to ensure that the Euler product for $L(s, \pi)$ converges in a right half plane. The other is a requirement that the local constituents $\pi_v$ of $\pi$ all be infinite dimensional. This is needed for the construction of a global Whittaker model for $\pi$ [103, Proposition 9.2], which in turn is a foundation for the desired embedding of $\pi$ into the space of cusp forms. (The last condition is known to hold for any $\pi \in \Pi_{\mathrm{cusp}}(G)$, a property for which the reader can consult the hints in the last paragraph of Section 14.)

The culmination of what we are calling the first part of [103] (Chapters I and II) is the application in Section 12 of this converse theorem to Artin $L$-functions. We are speaking now of the generalizations of these objects to continuous representations of the global Weil group $W_F$ (rather than the global Galois group $\Gamma_F$). The fundamental question is whether any irreducible two dimensional representation

$$r \colon W_F \to \widehat{G} = \mathrm{GL}(2, \mathbb{C})$$

corresponds to a cuspidal automorphic representation $\pi \in \Pi_{\mathrm{cusp}}(G)$, according to precepts of local and global functoriality.

The authors begin in Section 12 by appealing to [250]. This yields a global $L$-function $L(s, r)$ and $\varepsilon$-factor $\varepsilon(s, r)$ that satisfy analogues of the properties (34), (35) and (36), or rather their extensions from Corollary 11.2, but with the proviso that the global $L$-function need not be entire. The converse theorem actually applies to a given representation $\pi = \otimes \pi_v$ of $G(\mathbb{A}_F)$, so one first needs a way to construct the local components $\pi_v = \pi(r_v)$ of $\pi$ from the local components $r_v$ of $r$. This is exactly the form taken by the partial local classification from Chapter I. For any $r_v$, the authors prove that there is at most one $\pi_v = \pi(r_v)$ such that, in addition to some obvious requirements, the corresponding families of local $L$-functions and $\varepsilon$-factors are equal. (See p. 395 of [103] for a precise statement, which was extracted from the local results in Sections 4–6.) For the given two-dimensional representation $r$ of $W_F$, Theorem 12.2 of [103] then asserts that if both $L(s, \omega \otimes r)$ and $L(s, \omega^{-1} \otimes r^\vee)$ are entire functions for every quasicharacter $\omega$ on $C_F$, the representation $\pi_v = \pi(r_v)$ exists for every $v$, and $\pi = \otimes \pi_v$ is a cuspidal automorphic representations of $G(\mathbb{A})$.

There is a striking interpretation of this theorem. Recall the conjecture of Artin, which asserts that the $L$-function $L(s, r)$ attached to any irreducible representation $r$ of dimension greater than 1 is entire. Artin's conjecture therefore implies functoriality for $r$ when its dimension equals 2. Since this conjecture has been with us for almost a century, and is widely held to be true, we thus have strong evidence for functoriality for the case $r \to \pi$. Langlands cites Theorem 12.2, together with the partial local classification $r_v \to \pi_v$ based on $L$- and $\varepsilon$-factors, as one of the two principal contributions of [103] to the understanding of functoriality.

The other principal contribution cited by Langlands is the local and global correspondence

$$\pi' = \bigotimes_v \pi'_v \to \pi = \bigotimes_v \pi_v, \tag{37}$$

from the representations of the multiplicative group $G'$ of a quaternion algebra $M'$ to those of $G = \mathrm{GL}(2)$. Its apotheosis is the comparison in Section 16, but there was considerable preparation laid down earlier in the volume. In particular, the authors introduced the local correspondence $\pi'_v \to \pi_v = \pi(\pi'_v)$ at the same time and in the same spirit as their partial Galois/Weil correspondence $r_v \to \pi_v = \pi(r_v)$ described above. That is, they attached local $L$- and $\varepsilon$-factors $L(s, \omega_v \otimes \pi'_v)$ and $\varepsilon(s, \omega_v \otimes \pi'_v, \psi_v)$ to any $\pi'_v$, and then defined $\pi_v = \pi(\pi'_v)$ as the unique representation whose $L$- and $\varepsilon$-factors (which had already been constructed) match those of $\pi'_v$. (See p. 469 of [103] for the precise statement, taken from Sections 4 and 5.)

The local factors attached to $G'$ were of necessity defined in terms of the Lie algebra of $G'$, which is just the underlying four-dimensional quaternion algebra $M'$ over $F$. These amount to the analogue for $G'$ of the local construction from Tate's thesis. On the other hand, the local factors for $G$ were defined, according to the requirements of the converse from Hecke theory, in terms of the Mellin transform attached to a one-dimensional vector space. The analogue for $G$ of Tate's thesis would be based on a different transform attached to the Lie algebra of $G$, the four-dimensional matrix algebra $M$ over $F$. It would have in many ways been more natural. In Section 13, the authors investigate the local factors for $G$ that arise in this manner, and show in Theorem 13.1 that they do actually match the earlier local factors for $G$ defined by Hecke theory. In Section 14, they combine the global methods of Tate's thesis with the local factors already in place for $G'$. Theorem 14.1 asserts that if $\pi'$ is an automorphic representation of $G'(\mathbb{A}_F)$, the corresponding $L$-function $L(s, \pi')$ satisfies analogues of (34), (35) and (36), apart from the requirement that it be entire. But Corollary 14.3 then asserts that so long as $\pi'$ is not a one-dimensional representation attached a quasicharacter $\chi$ on $C_F$, the global $L$-functions $L(s, \omega \otimes \pi')$ and $L(s, \omega^{-1} \otimes (\pi')^\vee)$ are indeed all entire. Finally in Theorem 14.4 at the end of Section 14, the authors establish the essential property of the global correspondence (37). If we also take for granted the hints at the end of the section, as we did in the discussion of the converse theorem above, we can conclude that if $\pi'$ is not the one dimensional representation attached to a quasicharacter $\chi$ on $C_F$, its global image $\pi$ lies in $\Pi_{\mathrm{cusp}}(G)$.

In Theorem 15.2, the authors deal with the local correspondence $\pi'_v \to \pi_v$. They prove that it maps $\Pi(G'_v)$ *injectively* onto the relative discrete series $\Pi_2(G_v)$ of $G_v$, a term that carries the same meaning at any local place $v$ as that defined for $v = \mathbb{R}$ at the end of §2. Here, they use the local constructions of Tate for $\mathrm{GL}(2)$, and in particular, their own conclusion from Theorem 13.1 that the Tate local factors for $\pi_v$ match the Hecke local factors in terms of which the local correspondence was defined. (The statement of Theorem 15.2 was actually for nonarchimedean $v$, but only because the authors describe the image more explicitly as the set of representations $\pi_v$ that are either supercuspidal or "'special'", the latter being representations in the relative discrete series that are subquotients of induced representations. The result for archimedean $v$ is simpler, and is contained in the discussion in Section 5.) The authors supplement this result with an important formula for the characters of representations, regarded as locally integrable, conjugacy invariant functions $\Theta(\pi_v, \cdot)$ on the group. In Proposition 15.5, they establish that corresponding characters satisfy the identity

$$\Theta(\pi_v', b_v') = -\Theta(\pi_v, b_v), \quad \pi_v' \to \pi_v, \tag{38}$$

where $b_v' \to b_v$ is the natural bijection between the regular elliptic conjugacy classes on the two groups.

This is as far as Hecke theory goes. However, the authors still had more to say. In Section 16, they used the Selberg trace formula, something entirely different, to give a remarkable characterization of the image of the global correspondence (37). Combined with their characterization of the local correspondence $\pi_v' \to \pi_v$ we have just described, the final result is the assertion that the global correspondence $\pi' \to \pi$ from $G'$ to $G$ is a bijection, from the set of representations $\pi' \in \Pi(G')$ not attached to a quasicharacter $\chi$ on $C_F$ onto the set of representations $\pi \in \Pi_{\mathrm{cusp}}(G)$ with the property that for every $v$ such that $G_v'$ is not split, $\pi_v$ lies in the relative discrete series $\Pi_2(G_v)$ of $G_v$.

This is what Langlands cited as the second principal contribution of [103] to the understanding of functoriality. Since it characterizes automorphic representations of $G'$ as a natural subset of the cuspidal automorphic representations of its quasi-split inner form $G = \mathrm{GL}(2)$, it suggests that functoriality can be formulated purely in terms of quasi-split groups, as we have done here. As we have noted, the automorphic representations of groups that are not quasi-split are now usually treated as part of a separate theory of Langlands, the theory of endoscopy. But we could still regard the global comparison between $G'$ and $G$ as the first result in endoscopy, established years before the theory was formally proposed. We shall postpone our description of the actual comparison to the next section, as part of a general discussion of the trace formula.

Incidentally, the complementary set of automorphic representations of $G'$, those attached to quasicharacters $\chi$ of $C_F$, are just the one-dimensional automorphic representations of $G'$. They are naturally bijective with the set of one-dimensional automorphic representations of $G$, which represent the complement of $\Pi_{\mathrm{cusp}}(G)$ in the set $\Pi_2(G)$ of representations that comprise the full automorphic (relative) discrete spectrum of $G$. However, the formal correspondence between these two sets is different from (37). Its analogue for general groups was introduced [13], [18] as an attempt to describe the representations in the automorphic discrete spectrum of any group that do not satisfy the general analogue of Ramanujan's conjecture. We shall return to it in later sections.

# 6 Trace formula and first comparison

Selberg introduced his trace formula in 1956 [205]. He might actually have discovered it earlier, but since he published very little of his work, it is difficult to say. There are actually two formulas. One is the identity for compact quotient $\Gamma \setminus H$, whose elegant proof Langlands reconstructed in [133], as we discussed in Section 1. See also [79]. The other applies to many noncompact quotients $\Gamma \setminus H$ of rank 1. The rank here means the number of degrees of freedom one has to approach infinity in $\Gamma \setminus H$. This is more difficult. It requires among other things the theory of Eisenstein series, which we discussed in Section 2, and which Selberg had introduced for this express purpose.

In particular, Selberg established an explicit trace formula for $\Gamma \setminus \mathrm{SL}(2, \mathbb{R})$ (in the sense we discussed in Section 1), where $\Gamma$ equals $\mathrm{SL}(2, \mathbb{Z})$, or a congruence subgroup of $\mathrm{SL}(2, \mathbb{Z})$, or more generally, any discrete subgroup of $\mathrm{SL}(2, \mathbb{R})$ with a reasonable fundamental domain. He also established extended formulas that included the traces of the supplementary Hecke operators attached to a congruence subgroup of $\mathrm{SL}(2, \mathbb{Z})$. As we have noted, these operators are an integral part of the adelic framework that has since been adopted. Selberg used his formulas to prove striking estimates for the closed geodesics on the Riemann surface attached to $\Gamma \subset \mathrm{SL}(2, \mathbb{R})$, taken from the geometric side, and for the eigenvalues of the Laplacian on the Riemann surface, taken from the spectral side.

Langlands' interest in the trace formula was different. He saw it as an opportunity to study functoriality. Global Functoriality postulated deep reciprocity laws between automorphic representations for pairs of groups $G'$ and $G$. The trace formula for any one group, especially insofar as it existed for noncompact quotient, was clearly a complex identity. But might it be possible to compare the formulas for $G'$ and $G$, without having to evaluate the various terms in either case explicitly?

This brings us to the second part [103, §15–16] of the volume of Jacquet–Langlands. We put it aside in the previous section in order that it might serve us here as a simple introduction to the general comparison of trace formulas. On p. 516–517 of the volume, the authors stated the adelic version of Selberg's trace formula for the group $\mathrm{GL}(2)$. This was perhaps the first time the formula was stated in full, Selberg having limited his publications to specializations of the formula. In addition, it represents a two-fold extension of the formula, from a modular quotient of the upper half plane to a modular quotient of $\mathrm{GL}(2, \mathbb{R})$, and then to the adelic quotient $\mathrm{GL}(2, F) \setminus \mathrm{GL}(2, \mathbb{A})$ for a number field $F$. (In fact, the statement in [103] applies to any global field $F$, but as always for us, $F$ will be of characteristic zero.) The authors also gave a clean and concise sketch of the proof. Detailed proofs of the adelic version of Selberg's formula for $\mathrm{PSL}(2)$ [69] and groups $G$ of $F$ rank 1 [9] appeared later.

Suppose for a moment that $G$ is an arbitrary reductive group over a number field $F$. We could then take $f$ to be a function in the natural space[4] of test functions

$$C_c^\infty(G(\mathbb{A}_F)) = \varinjlim_S \left( C_c^\infty(G_\infty) \otimes C_c^\infty(G_S^\infty) \otimes \mathbf{1}^S \right)$$

on $G(\mathbb{A})$. Right convolution by $f$ in the Hilbert space $L^2(G(F) \setminus G(\mathbb{A}_F))$ obviously converges. It is easily seen to be an integral operator, with kernel

$$K(x,y) = \sum_{\gamma \in G(F)} f(x^{-1}\gamma y), \quad x,y \in G(\mathbb{A}). \tag{39}$$

If $G(F) \setminus G(\mathbb{A}_F)$ is compact, $R(f)$ is of trace class, by standard methods in functional analysis. Its trace is given by the identity (1) from Section 1, but with $G(F)$ and $G(\mathbb{A}_F)$ in place of $\Gamma$ and $G$. Indeed, the spectral extension on the right hand side of (1) equals the trace of $R(f)$, by definition, while the left hand side amounts to the geometric expression for the trace derived originally by Selberg, and by Langlands in [133].

Suppose however that $G(F) \setminus G(\mathbb{A}_F)$ is not compact. Then $R(f)$ is not of trace class. The problem is with the continuous spectrum $L_{\mathrm{cont}}^2(G(F) \setminus G(\mathbb{A}_F))$. The restriction $R_{\mathrm{cont}}(f)$ of $R(f)$ to this invariant subspace is no more of trace class than would be the convolution operator on $L^2(\mathbb{R})$ of a function in $C_c^\infty(\mathbb{R})$. One must subtract the contribution of $R_{\mathrm{cont}}(f)$ to the kernel $K(x,y)$ of $R(f)$ to obtain an operator that is better behaved. This is the role of Langlands' general theory of Eisenstein series. The point is that Eisenstein series provide a spectral formula for $K(x,y)$. As a formal consequence of Langlands' Main Theorem from Section 2, one obtains a spectral expansion

$$K(x,y) = \sum_P n_P^{-1} \int_{i\mathfrak{a}_P^*} \sum_{\phi \in \mathscr{B}_P} E(x, \mathscr{I}_P(\lambda, f)\phi, \lambda) \overline{E(y, \phi, \lambda)} \, d\lambda \tag{40}$$

in which $\mathscr{B}_P$ is an orthonormal basis of the Hilbert space $\mathscr{H}_P$ on which the induced representation $\mathscr{I}_P(\lambda)$ acts, for the kernel to accompany the simpler geometric expression (39). (See for example [22, (7.6)].) The multiple integral converges absolutely, even though there are no reasonable pointwise estimates for the integrand. The argument for this, which I learned from Langlands, and is due to Selberg, is to combine the Schwartz inequality with the fact that for a positive definite function $f = h * h^*$, the diagonal value $K(x,x)$ of the kernel bounds that of the integral in (40) over any compact subset of the domain. (See for example [22].) We shall return to these matters when we discuss the general trace formula in Section 10.

---

[4] Here $S \supset S_\infty$ is a finite set of valuations of $F$, while $C_c^\infty(G_\infty)$ is the ordinary space of test functions on the archimedean component $G_\infty = G(F_\infty)$, and $C_c^\infty(G_S^\infty)$ is the space of locally constant, complex valued functions of compact support on the "ramified" nonarchimedean component $G_{S,\infty}^\infty = G(F_S^\infty)$. We write $\mathbf{1}^S$ for the characteristic function of a suitable natural compact subgroup $K^S$ of the remaining "unramified" component $G^S = G(\mathbb{A}^S)$ of $G(\mathbb{A})$.

For the group $G = \mathrm{GL}(2)$ in §16 of [103], there are two terms indexed by $P$ in (40). One is given by the Borel subgroup $P = P_0 = B$ of upper triangular matrices in $G$. The other corresponds to $P = G$. It is the kernel of the operator $R_{\mathrm{disc}}(f)$ obtained by restricting $R(f)$ to the relative discrete spectrum.

The group $\mathrm{GL}(2)$ has a split, one-dimensional centre, the group $Z = A_G \cong \mathrm{GL}(1)$ of scalar matrices, so to have a trace at all, one must make the usual minor adjustment. However, instead of either restricting $f$ to the subgroup $G(\mathbb{A})^1$ of $G(\mathbb{A})$ or projecting it onto a function invariant under a subgroup $A_{G,\infty}^+ \cong \mathbb{R}^+$ of $A_G(F_\infty)$, according to the discussion in Section 2, the authors take $f$ to be $\eta^{-1}$-equivariant, for a character $\eta$ on the full adelic quotient $Z(F) \backslash Z(\mathbb{A}_F)$ of the centre. The kernel $K(x, y)$ is easy to adjust. For the original function $f$, the integral

$$\int_{Z(F) \backslash Z(\mathbb{A}_F)} K(zx, y) \eta(z) \, dz$$

becomes the kernel for the $\eta^{-1}$-equivariant function

$$x \to \int f(zx) \eta(z) \, dz.$$

The formulas (39) and (40) for the kernel are essentially unchanged, even if $G$ is a general group. They can be taken as stated so long as we understand that:

(i) $f$ now belongs to the space $C_c^\infty(G(\mathbb{A}_F), \eta^{-1})$ of $\eta^{-1}$-equivariant test functions.

(ii) $R(f)$ is the right convolution over $Z(\mathbb{A}_F) \backslash G(\mathbb{A}_F)$ of $f$ on the space $L^2(G(F) \backslash G(\mathbb{A}_F), \eta)$ of square integrable, $\eta$-equivariant functions on $G(F) \backslash G(\mathbb{A}_F)$.

(iii) The integrals in (40) are really taken over the kernel $i\mathfrak{a}_P^{*,G}$ in $i\mathfrak{a}_P^*$ of the canonical linear projection of $i\mathfrak{a}_P^*$ into $i\mathfrak{a}_G^*$, in which the original notation holds if we treat the dependence of $\mathscr{I}_P(\lambda, f)$ on the image of $\lambda$ in $i\mathfrak{a}_G^*$ as a Dirac distribution.

The two formulas (39) and (40) for $K(x, y)$ become $\eta^{-1}$-equivariant functions in $x$ and $y$ on $G(F) \backslash G(\mathbb{A}_F)$, making their diagonal values at $y = x$ functions on $Z(\mathbb{A}_F) G(F) \backslash G(\mathbb{A}_F)$. The trace of the restriction $R_{\mathrm{disc}}(f)$ of $R(f)$ to the $\eta$-discrete spectrum becomes the integral over this set of the expression with $P = G$ in (40).

For our group $G = \mathrm{GL}(2)$, the trace of $R_{\mathrm{disc}}(f)$ is thus the integral over $x = y$ in $Z(\mathbb{A}_F) G(F) \backslash G(\mathbb{A}_F)$ of the difference between (39) and the expression with $P = B$ in (40). Neither of these last two functions is integrable. To obtain a trace formula, one must see how their nonintegrable parts cancel, and then find an explicit formula for the integral of what remains. As we noted, Jacquet and Langlands gave a sketch of the process in §16 of [103]. The answer they obtained for the trace of $R_{\mathrm{disc}}(f)$ is given as a sum of the eight terms (i)–(viii) on p. 516–517. We shall say a few brief words about the argument, to give ourselves some general perspective.

The value of (39) at $y = x$ is the more transparent of the two expressions (and would of course be the only expression to consider in case of compact quotient). It can again be written as

$$\sum_{\{\gamma\}} \sum_{\delta \in G_\gamma(F) \backslash G(F)} f(x^{-1} \delta^{-1} \gamma \delta x),$$

where $\{\gamma\}$ is a set of representatives of $G(F)$-conjugacy classes in $A_G(F) \backslash G(F)$. If we proceed formally as if $G$ were a group with $Z(\mathbb{A}_F) G(F) \backslash G(\mathbb{A}_F)$ compact, in which all multiple integrals are absolutely convergent, we could write the integral of this function as

$$\sum_{\{\gamma\}} \int_{Z(\mathbb{A}_F) G(F) \backslash G(\mathbb{A}_F)} \sum_{\delta \in G_\gamma(F) \backslash G(F)} f(x^{-1} \delta^{-1} \gamma \delta x) \, dx, \tag{41}$$

which in turn is equal to

$$\sum_{\{\gamma\}} \int_{Z(\mathbb{A}_F) G_\gamma(F) \backslash G_\gamma(\mathbb{A}_F)} \cdot \int_{G_\gamma(\mathbb{A}_F) \backslash G(\mathbb{A}_F)} f(x^{-1} \gamma x) \, dx,$$

and hence also to

$$\sum_{\{\gamma\}} \mathrm{vol}(Z(\mathbb{A}_F) G_\gamma(F) \backslash G_\gamma(\mathbb{A}_F)) \cdot \int_{G_\gamma(\mathbb{A}_F) \backslash G(\mathbb{A}_F)} f(x^{-1} \gamma x) \, dx. \tag{42}$$

This would match the left hand side of Selberg's trace formula for compact quotient, which we quoted from [133] as (1) in Section 1. But in the case $G = \mathrm{GL}(2)$ at hand, the integrals do not all converge. There are four kinds of terms in (41), two good and two bad.

If $\gamma$ is a scalar matrix, it can be represented by 1. In this case, the corresponding integral in (41) gives a contribution

$$\mathrm{vol}(Z(\mathbb{A}_F) G(F) \backslash G(\mathbb{A}_F)) \cdot f(1). \tag{43}$$

to (42). This is the term (i) on p. 516 in [103]. If the characteristic polynomial of $\gamma$ is irreducible over $F$, its eigenvalues generate a quadratic extension $E = E_\gamma$ of $F$, and $G_\gamma(F)$ is isomorphic to $E^*$. In this case the corresponding integral in (41) converges. The contribution

$$\sum_{\{\gamma\}} \mathrm{vol}(\mathbb{A}_F^* E_\gamma^* \backslash \mathbb{A}_{E_\gamma}^*) \cdot \int_{G_\gamma(\mathbb{A}_F) \backslash G(\mathbb{A}_F)} f(x^{-1} \gamma x) \, dx \tag{44}$$

of these integrals to (42) gives the term (ii) on p. 516 of [103]. The term (iii) in [103] vanishes, being a sum over the empty set of inseparable quadratic extensions of the number field $F$. The other terms in (41) are bad. If the characteristic polynomial of $\gamma$ is a product of two *distinct* linear factors over $F$, the corresponding integral in (41) diverges. The explanation is that in the associated summand in (42), the integral converges but the volume coefficient is infinite. The remaining terms in (41) are attached to nontrivial unipotent elements, namely the complement of the scalar matrix 1 from (43) in the classes $\{\gamma\}$, with characteristic polynomial being the square of a linear factor. The explanation here is the other way around, where the associated volume is finite, and the integral is what diverges. The two kinds of bad classes $\{\gamma\}$ ultimately contribute to the trace formula as the respective terms (iv) and (v) on p. 514 of [103].

The other half of the trace formula concerns the spectral expansion of the kernel, specifically the value at $y = x$ of the summand of $P = B$ in (40). Following Selberg's basic ideas, Jacquet and Langlands multiplied the difference between the value of $y = x$ at (39) and this spectral function by the characteristic function of a large compact subset

$$\{x : H_P(x) \leq \log c_1\}$$

of $Z(\mathbb{A}_F)G(F) \setminus G(\mathbb{A}_F)$, defined in terms of the usual fundamental domain for any large positive number $c_1$. This effectively led them to a cancellation of the noninte-grable parts of each function. More precisely, they showed that for each of the two truncated functions in the difference, the integral is a sum of three linear forms in $f$: an explicit distribution that is independent of $c_1$, the product of $\log c_1$ with a simpler distribution, and a distribution that appears to be quite complicated, but which approached 0 as $c_1$ approaches infinity. The two multiples of $\log c_1$ are easily seen to cancel. The terms that approach 0 can be ignored. This leaves only the distribu-tions that are independent of $c_1$. They are equal to the sum of the terms (i)–(v) on p. 516 of [103] from the geometric kernel we have discussed above, and the terms (vi)–(viii) on p. 517 from the spectral kernel.

This completes our general remarks on the trace formula for GL(2) in [103]. It is a curious fact, observed a couple of years later, that with a minor change in the truncation process that takes into account the noncuspidal discrete spectrum (part of the term with $P = G$ in (40)), each of the two distributions that approach 0 actually vanishes if $\log c_1$ is sufficiently large. This turned out to be a general phenomenon, which carried over to the later trace formula of an arbitrary group $G$. In the general case, the truncation parameter $T_1 = \log c_1$ for GL(2) is replaced by a vector $T$ in the positive chamber $\mathfrak{a}_{P_0}^+$ that is far from the walls. It gives rise to a uniform truncation operation on the diagonal values $(x, x)$ of each of the two expansions (39) and (40) of $K(x, x)$. The integrals of the resulting two functions of $x$ in turn come with their own expansions, whose terms are polynomials in $T$ if $T$ is far from the walls of $\mathfrak{a}_{P_0}^+$ (in a sense that depends only on the support of $f$). The general trace formula comes from this. It is the identity of distributions obtained by setting the polynomial variable $T$ equal to a canonical point $T_0$ (often 0), which is determined by the choice of a suitable maximal compact subgroup $K_0 \subset G(\mathbb{A}_F)$. (See [12, §1–2] and [22, Theorem 9.1]).

There is another matter, which is a little more disturbing. Given a locally compact group $H$, one says that a continuous linear form $F$ in $h \in C_c(H)$ is *invariant* if

$$F(h^u) = F(h), \qquad h^u(x) = h(u^{-1}xu), \, u \in H,$$

for every $h$ and $u$. In the case of $H = \mathrm{GL}(2, \mathbb{A}_F)$, the trace $\mathrm{tr}(R_{\mathrm{disc}}(f))$ is an invariant distribution. One might therefore expect all of the terms (i)–(viii) in the formula for $\mathrm{tr}(R_{\mathrm{disc}}(f))$ also to be invariant distributions. They are not. The problem is that the truncation operation interferes with the symmetry under conjugation. In particular, it renders the geometric terms (iv) and (v) and the spectral terms (viii) noninvariant. However, there is a natural "renormalization" process, which converts these terms to invariant distributions, and which applies uniformly to any reductive group $G$. We refer the reader to the formula (2) of the introduction of [15] for a general statement of the final *invariant trace formula*, and to [22, §22–23] for a description of the general correction process.

These refinements are not actually needed for the comparison in §16 of [103]. Suppose again that $G$ is a general reductive group over $F$, with standard parabolic subgroup $P = MN$, and that $f_w$ is a function in $C_c^\infty(G_w)$, for a valuation $w$ of $F$. It is convenient to define a supplementary function $f_{w,M}$ on $\Pi(M_w)$ (the set of irreducible representations of $M_w = M(F_w)$) by setting

$$f_{w,M}(\sigma_w) = \mathrm{tr}(\mathscr{I}_P^G(\sigma_w, f_w)), \quad \sigma_w \in \Pi(M_w),$$

the character of the representation of $G_w$ parabolically induced from $\sigma_w$. We then say that function $f = \prod_v f_v$ in $C_c^\infty(G(\mathbb{A}_F))$ is *cuspidal* at $w$ if $f_{w,M} = 0$ for every parabolic subgroup $P \neq G$. The term really ought to be "invariantly cuspidal" so as not to conflict with the earlier notion of a cuspidal function from Section 2. However, the context should rule out any future confusion. In any case, the relevance of the property is that if $f$ is cuspidal at two places $v_1$ and $v_2$, the general invariant trace formula simplifies dramatically [15, Theorem 7.1]. It reduces to something close to what it would be if $G(F) \setminus G(\mathbb{A}_F)$ were compact. In the case $G = \mathrm{GL}(2)$ of present interest, the result for any such $f$ is that the terms (iii)–(viii) all vanish, and hence that the trace of $R_{\mathrm{disc}}(f)$ equals the sum of the terms (i) and (ii). Thus, for the group $\mathrm{GL}(2)$ (or for any group of semisimple $F$-rank 1), one does not need the invariant trace formula for this simplification. It follows easily from the basic trace formula, as Jacquet and Langlands point out for $\mathrm{GL}(2)$ in §16.

We can now begin our brief account of the comparison theorem in §16 of [103] for the global correspondence of representations (37). We should first say something on the structure of the multiplicative group $G'$ of a quaternion algebra. For a broader view of this, we extract what we need from the structure of general groups.

As we suggested in Section 2, there is a remarkable classification theory for general (connected) reductive groups $G$ over local and global fields $F$. I take the liberty of representing it as a diagram of four steps.

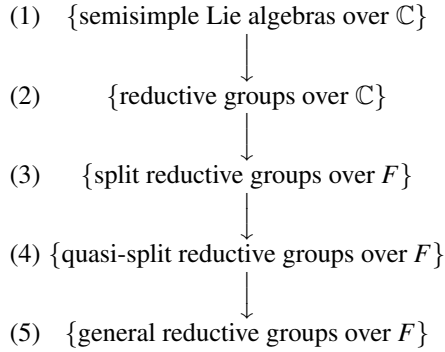

(1)   {semisimple Lie algebras over $\mathbb{C}$}

(2)      {reductive groups over $\mathbb{C}$}

(3)   {split reductive groups over $F$}

(4) {quasi-split reductive groups over $F$}

(5)   {general reductive groups over $F$}

Complex semisimple Lie algebras are bijective with Coxeter–Dynkin diagrams. These in turn are finite disjoint unions of connected diagrams, of which there are four infinite families and five exceptional diagrams. This is the starting point.

The first step $(1) \to (2)$ follows from the theory of covering groups, which is easy to manage in this setting. It is founded on our explicit knowledge of the finite abelian groups $Z(G_{\mathrm{sc}}(\mathbb{C}))$, where $G_{\mathrm{sc}}(\mathbb{C})$ ranges over the simple, simply-connected complex groups attached to the connected diagrams, and the fact that any complex reductive group is a quotient of a finite direct product of groups $G_{\mathrm{sc}}(\mathbb{C})$ and a complex torus, by a finite central subgroup. The second step $(2) \to (3)$ is a bijection, obtained from a special case of the basic construction of Chevalley groups. The third step $(3) \to (4)$ is given by *outer twistings* of the Galois action on the given split group. Its main ingredient is a homomorphism (of finite image) from the Galois group $\Gamma_F = \mathrm{Gal}(\overline{F}/F)$ to the group of automorphisms of the underlying diagram. The last step $(4) \to (5)$ is based on class field theory. It is given by *inner twistings* of the Galois action on the given quasi-split group $G^*$ over $F$, or in purely algebraic terms, elements in the Galois cohomology *set* $H^1(\Gamma_F, G^*_{\mathrm{ad}}(\overline{F}))$. These objects are a part of abelian class field theory, for they reduce to Galois cohomology *groups* $H^*(\Gamma_F, X)$, attached to abelian $\Gamma_F$-modules $X$, to which Tate–Nakayama duality applies. To view these matters in terms of Langlands dual groups, we refer the reader to [121, §1 ($F$ local) and §2 ($F$ global)].

This description of the general classification goes well beyond what is needed for quaternion groups, but it does represent an implicit foundation for various topics that will arise later. For the quaternion groups $G'$, the classification is simple. In particular, one needs only the last step $(4) \to (5)$, since the underlying quasi-split group is the split group $G = \mathrm{GL}(2)$. The conclusions are as follows. If $F$ is local, with $F \neq \mathbb{C}$, there is exactly one quaternion group up to isomorphism. (If $F = \mathbb{C}$, there are no quaternion groups.) If $F$ is global, the isomorphism classes of quaternion groups $G'$ over $F$ are parametrized by *finite, even, nonempty* sets $S$ of valuations of $F$ (with $F_v \neq \mathbb{C}$ for each $v \in S$). For any such $S$, $G'$ is characterized by the property that for *any* valuation $v$, $G'_v$ equals the unique quaternion group over $F_v$ if $v$ lies in $S$, and is equal to $\mathrm{GL}(2)_v$ otherwise.

Suppose then that $G'$ is a quaternion group over $F$ with $F$ global, and with centre $Z' = Z(G')$. The characteristic polynomial gives a natural bijection $b'_v \to b_v$ between the regular elliptic conjugacy classes in $G'_v = G'(F_v)$ and $G_v = G(F_v)$, as in the character formula (38). If $G'_v$ does not split, any regular, semisimple class is elliptic. For the group $G_v$, however, there are also regular hyperbolic conjugacy classes $a_v$ in $G_v$. We recall that according to the general definitions, a semisimple conjugacy class $\gamma = \gamma_v$ in $G_v$ (or $G'_v$) is *regular* if the centralizer $G_\gamma$ is a torus, and that the regular, semisimple classes in $G_v$ form an open dense set. The class $\gamma$ is elliptic if the quotient $Z(F_v) \setminus G_\gamma(F_v)$ is compact, and hyperbolic otherwise.

For the comparison of trace formulas, Jacquet and Langlands chose matching functions $f' \in C_c^\infty(G'(\mathbb{A}), \eta)$ and $f \in C_c^\infty(G(\mathbb{A}), \eta)$ in the relevant spaces of $\eta$-equivariant test functions. More precisely, given a function $f' = \prod f'_v$ for $G'$, they chose a function $f = \prod f_v$ for $G$ by defining the local factors $f_v$ as follows. If $G'$ is split at $v$, they simply identified $f_v$ with $f'_v$ under any of the isomorphisms from $G'(F_v)$ to $G(F_v)$ in the $G(F_v)$-conjugacy class of isomorphisms attached to $G'_v$ as an inner twist. If $G'_v$ is not split, they took $f_v$ to be any function such that

$$\mathrm{tr}(f_v(\pi_v)) = \delta(\pi_v) \cdot \mathrm{tr}(f'_v(\pi'_v)), \tag{45}$$

where $\pi_v$ is any (irreducible) tempered representation of $G(F_v)$, and $\delta(\pi_v) = 1$ if $\pi_v$ is the image of $\pi'_v$ under the local correspondence $\pi'_v \to \pi_v$ in (37), and $\delta(\pi_v) = 0$ otherwise. This is a spectral relation. One sees from the character formula (38) that it is equivalent to the geometric relation

$$\int_{G_\gamma(F_v) \setminus G(F_v)} f_v(x_v^{-1} \gamma x_v) \, dx_v = \varepsilon(\gamma) \int_{G'_{\gamma'}(F_v) \setminus G'(F_v)} f'_v(x_v^{-1} \gamma' x_v) \, dx_v \tag{46}$$

in which $\gamma$ is any regular conjugacy class in $G(F_v)$, and $\varepsilon(\gamma)$ equals $(-1)$ if $\gamma = b_v$ is the image of an (elliptic) class $\gamma' = b'_v$ in $G'(F_v)$, and equals 0 if $\gamma = a_v$ is a hyperbolic class. The function $f_v$ is not uniquely determined by either of these relations. However, if $I_v$ is any invariant distribution on $G(F_v)$, $I_v(f_v)$ is uniquely determined. We also note that if $G'_v$ is not split, the function $f_v$ is cuspidal, in the sense above.

For the group $G'$, the quotient $Z'(\mathbb{A}_F)G'(F) \setminus G'(\mathbb{A}_F)$ is compact. The trace formula for $G'$ is therefore the Selberg trace formula for compact quotient, discussed above and in Section 1. The trace of the operator $R(f') = R_{\mathrm{disc}}(f')$ is accordingly equal to the sum of the two terms

$$(Z'(\mathbb{A}_F)G'(F) \setminus G'(\mathbb{A}_F)) \cdot f'(1) \tag{47}$$

and

$$\sum_{\{\gamma\}} \mathrm{vol}(\mathbf{A}_F^* E_\gamma^* \setminus \mathbf{A}_{E_{\gamma'}}^*) \int_{G_{\gamma'}(\mathbb{A}_F) \setminus G'(\mathbb{A}_F)} f'((x')^{-1} \gamma' x') \, dx', \tag{48}$$

given by the analogues of (43) and (44), which is to say, of the terms (i) and (ii) in Chapter 16 of [103] for $G'$. As for the group $G$, the trace of $R_{\mathrm{disc}}(f)$ equals the larger sum of terms (i)–(viii) from [103]. But since the local factor $f_v$ is cuspidal at any place $v$ at which $G'_v$ does not split, and since the set $S$ of such places is nonempty and even, $f$ is cuspidal at (at least) two places $v$. Therefore, as the authors observe on p. 523, the terms (iv)–(viii) vanish for the given $f$. The term (iii) vanishes for any number field, so this leaves only the sum of (43) and (44), the terms (i) and (ii) for $G$.

The sum (48) is over regular, globally elliptic conjugacy classes $\gamma'$ in $G'(F)$, which is to say, classes whose characteristic polynomial is irreducible over $F$. The same is true of the sum over $\gamma$ in (44). As in the local case, the characteristic polynomials then give a bijection $\gamma' \to \gamma$ between the two indices of summation. Moreover, it follows from (46) (and the fact that $S$ is even) that the summands of $\gamma'$ and $\gamma$ are equal. Therefore, the sums (48) and (44) are equal, as the authors note in p. 524. All that remains of the two trace formulas are the terms (43) and (47). It follows that

$$\begin{aligned}&\text{tr}(R_{\text{disc}}(f')) - \text{tr}(R_{\text{disc}}(f))\\&= \text{vol}(Z'(\mathbb{A}_F)G'(F) \setminus G'(\mathbb{A}_F)) \cdot f'(1) - \text{vol}(Z(\mathbb{A}_F)G(F) \setminus G(\mathbb{A}_F)) \cdot f(1).\end{aligned}$$

With further simple arguments, the authors refine this formula into three identities between its basic components. For example, an approximation obtained by varying $f'$ (and its corresponding image $f$) shows that each side of the last formula vanishes. Orthogonality relations for discrete series, with the corresponding Plancherel formulas for $f'(1)$ and $f(1)$, establish that $f'(1) = f(1)$. From this it follows that

$$\text{vol}(Z'(\mathbb{A}_F)G'(F) \setminus G'(\mathbb{A}_F)) = \text{vol}(Z(\mathbb{A}_F)G(F) \setminus G(\mathbb{A}_F)) \tag{49}$$

and that

$$\text{tr}(R_{\text{disc}}(f')) = \text{tr}(R_{\text{disc}}(f)) \tag{50}$$

(See [103, p. 524–525].)

The volume identity (49) is of independent interest. The authors pointed out that while it was well known, the methods used to obtain it were not. They suggested that a similar comparison of the trace formula (unknown at the time) of an arbitrary group with that of a quasi-split inner form might be used to establish a similar identity in general. This could then be combined with Langlands' early formula [135] for the corresponding volume of a Chevalley group (or rather, its imputed generalization to quasi-split groups). The goal would be to give a general proof of Weil's conjecture on Tamagawa numbers. As we noted in Section 1, this is exactly what happened, with the subsequent work of Lai and Kottwitz (and the final step for $E_8$ by Chernousov). Weil's conjecture, incidentally, is the elegant volume formula

$$\tau(G) = \text{vol}(G(F) \setminus G(\mathbb{A})) = 1,$$

for any simply connected group $G$ over $F$, taken with respect to the canonical Tamagawa measure on $G(\mathbb{A})$. It is understood that our discussion has been for Haar measures on $G'(\mathbb{A})$ and $G(\mathbb{A})$ that are compatible, in the sense that they coincide under transfer by the inner twist, even though they do not have to be the Tamagawa measures.

It is of course the other identity (50) that is the main point. We have been discussing the proof sketched in [103, §16] of Theorem 16.1, the basic characteriza-tion of the global correspondence $\pi' \to \pi$. As we noted in Section 5, it is equivalent to the assertion that the mapping $\pi' \to \pi$ is a bijection between the subsets $\Pi_{\text{cusp}}(G') = \Pi(G')$ and $\Pi_{\text{cusp}}(G)$ described at the end of Section 5. Having already shown the mapping to be injective, the authors established the surjectivity assertion by combining (50) with some elementary functional analysis. (See [103, p. 403–503].)

This at last completes our discussion of §16 from the monograph of Jacquet and Langlands. It represents the earliest comparison of adelic trace formulas. I hope that our rather extended treatment of it will serve as an introduction to the increasingly complex comparisons that followed. We note that special cases had been established earlier by Shimizu [225], [226], [227], with comparisons of Selberg's original (non-adelic) trace formulas [205], [206], [207].

# 7 Base change

In their introduction to [103], the authors wrote of the comparison we have just reviewed,

*"... the theorem of §16 is important and its proof is such a beautiful illustration of the power and ultimate simplicity of the Selberg trace formula and the theory of harmonic analysis on semi-simple groups that we could not resist adding it. Although we are very dissatisfied with the methods of the first fifteen paragraphs we see no way to improve on those of §16. They are perhaps the methods with which to attack the question left unsettled in §12."*

For a short time afterwards, Langlands was apparently unhappy with the last sentence, possibly thinking that it was premature. It is ironic then that he found an supportive answer only a few years later.

It came from base change. Suppose again that $G = \mathrm{GL}(2)$ over the number field $F$, and that $E$ is a cyclic extension of $F$ of prime degree $\ell$. The *restriction of scalars* functor then attaches a quasi-split reductive group $G_E^0 = \mathrm{Res}_{E/F}(G)$ over $F$ to the pair $(G, E)$ such that $G_E^0(F) = G(E)$. (The reason for the superscript 0 will become clear presently.) The $L$-group of $G_E^0$ can be taken to be

$$
{}^L G_E^0 = \widehat{G}_E^0 \rtimes \Gamma_{E/F} = \underbrace{\mathrm{GL}(2,\mathbb{C}) \times \cdots \times \mathrm{GL}(2,\mathbb{C})}_{\ell} \rtimes \mathrm{Gal}(E/F),
$$

where the cyclic group $\Gamma_{E/F} = \mathrm{Gal}(E/F)$ acts by permutation on the product of groups $\mathrm{GL}(2,\mathbb{C})$. Taking ${}^L G = \widehat{G} \times \Gamma_{E/F}$ as the $L$-group of $G$, we then have an $L$-homomorphism

$$
\rho : g \times \sigma \to \underbrace{(g, \ldots, g)}_{\ell} \rtimes \sigma, \quad g \in \widehat{G} = \mathrm{GL}(2,\mathbb{C}), \sigma \in \Gamma_{E/F},
$$

from ${}^L G = \widehat{G} \times \Gamma_{E/F}$ into ${}^L G_E^0 = \widehat{G}_E^0 \rtimes \Gamma_{E/F}$. Functoriality for $\rho$ is what is known as *base change*, a problem that may be formulated in this way if $G$ is any quasi-split group over $F$, and $E/F$ is any finite extension.

Langlands had earlier [138] proposed base change for $\mathrm{GL}(2)$ as one of several natural questions chosen to illustrate the difficulty of functoriality. His unexpected solution of the problem in 1975 followed new ideas of Saito and Shintani. (The work of both Shintani and Langlands was published only later, Shintani [231] in the 1979 Corvallis proceedings, and Langlands in his 1980 monograph [149].) Motivated by the special case solved by Shintani, Langlands established a general correspondence $\pi \to \pi_E$ from the local and global representations of $G$ to those of $G_E^0$. The solution represents a new comparison of trace formulas, considerably more sophisticated than that of the local and global correspondence $\pi' \to \pi$ of (37) for quaternion groups. As such, it amounts to something beyond a proof of functoriality for this case. Like the correspondence $\pi' \to \pi$, and in common with what one might hope for in any new comparison of trace formulas, the method allows also for a characterization of the *image* of the functorial correspondence.

We draw on the introduction of [149] for a brief history of the problem. Doi and Naganuma treated cases with $F = \mathbb{Q}$, $E$ a real quadratic extension, and with the archimedean component $\pi_\infty$ of $\pi$ being in the discrete series [68]. Jacquet extended these results to more general $F$ and $\pi$ in [101], but as in [68], without characterizing the image of the correspondence. It was Saito [193] who introduced the twisted trace formula with respect to $E/F$, a new kind of trace formula that he was then able to compare with that of GL(2). He was thus able to establish base change for more general cyclic extensions $E/F$ of prime degree, and also to characterize its image. Finally, Shintani formulated these ideas in terms of representation theory and adele groups (rather than the classical framework of holomorphic modular forms). This allowed him to deal with critical problems related to the construction of a local correspondence $\pi_v \to \pi_{v,E}$ at the ramified places of $\pi$. However, Shintani was still restricted to representations $\pi$ whose local archimedean constituents $\pi_v$ belonged to Harish-Chandra's discrete series, a condition that means that the local test function $f_v$ is chosen to be cuspidal. The problem for him was in the complexity of the trace formulas that would otherwise have to be compared. We now have some idea of this difficulty, having seen the complexity of the full trace formula for GL(2), with its eight terms (i)–(viii) from [103, §16]. The restriction under which Shintani worked is essentially the condition that $f$ be cuspidal at two places, which reduces that trace formula to the simple terms (i) and (ii), as in the comparison for quaternion groups from [103].

Langlands understood how to work with the full trace formula. By comparing the twisted trace formula for $G_E^0$ with the unrestricted trace formula for $G = \mathrm{GL}(2)$, he was able to establish the general base change correspondence $\pi \to \pi_E$ for any extension $E/F$ of prime order $\ell$ and any $\pi$. It was only after having removed all the restrictions that had gone before that he was able to establish its remarkable applications to Artin $L$-functions.

Having discussed the trace formula for GL(2) and its comparisons from [103] in some detail, we shall be content with a briefer summary of the twisted trace formula. Even for GL(2), the twisted trace formula and its comparison with the ordinary trace formula are rather more technical. A reader could bypass our summary, which leads directly to the spectral comparison identity (54) (or (55), written in the notation of Langlands), and proceed to the subsequent review of the properties of the resulting base change lifting (56).

The origin of the twisted trace formula is the (algebraic) automorphism $\sigma_E$ of the quasisplit group $G_E^0$ over $F$. Its action on $G_E^0(F)$ is given by the Galois automorphism $\sigma$ of $G(E)$. It also acts on $G_E^0(\mathbb{A}_F)$, and on the quotient $G_E^0(F) \setminus G_E^0(\mathbb{A}_F)$, and indeed, on the adelic quotient $M_E(F) \setminus M_E(\mathbb{A}_F)$ of any $\sigma_E$-stable subgroup $M_E$ of $G_E^0$ over $F$. Langlands introduced a $\sigma_E$-stable character $\xi_E$ on the adelic quotient $Z_E(F) \setminus Z_E(\mathbb{A}_F)$ of the centre $Z_E = Z(G_E)$ of $G_E^0$. This is the setting of the trace formula for GL(2) discussed in the last section. In particular, for any $\xi_E^{-1}$-equivariant test function $f_E^0 \in C_c^\infty(G_E^0(\mathbb{A}_F), \xi_E^{-1})$, we have the operator $R_{\mathrm{disc}}(f_E^0)$ in the $\xi_E$-discrete spectrum

$$\mathrm{L}^2_{\mathrm{disc}}(\xi_E) = \mathrm{L}^2_{\mathrm{disc}}(G_E^0(F) \setminus G_E^0(\mathbb{A}_F), \xi_E)$$

whose trace is the object of the trace formula. But it is the twisted trace of $R_{\mathrm{disc}}(f_E^0)$ that is of interest here.

We form the semidirect product

$$G_E^+ = G_E^0 \rtimes \langle \sigma_E \rangle,$$

and write

$$G_E = G_E^0 \rtimes \sigma_E$$

for the connected component attached to a fixed generator $\sigma_E$ of $\Gamma_{E/F}$. The action

$$(\sigma_E \phi)(y) = \phi(\sigma_E^{-1}(y)), \quad y \in G_E^0(F) \setminus G_E^0(\mathbb{A}_F),$$

of $\sigma_E$ on the Hilbert space $L^2(\xi_E)$ gives a canonical extension $R_E$ of the representation $R_E^0$ of $G_E^0(\mathbb{A}_F)$ to the group generated by $G_E(\mathbb{A}_F) = \sigma_E \cdot G_E^0(\mathbb{A}_F)$. In particular, for any test function $f_E \in C_c^\infty(G_E(\mathbb{A}_F), \xi_E)$, we have a unitary operator

$$R_E(f_E) = \int_{Z_E(\mathbb{A}_F) \setminus G_E(\mathbb{A}_F)} f_E(x) R_E(x)\, dx$$

on $L^2(\xi_E)$. Our notation here is slightly different from that of Langlands, as has sometimes been the case in past discussion, but it makes no difference in the argument.

The *twisted trace formula* is an explicit formula for the trace of the operator $R_{\mathrm{disc}}(f_E)$. Its proof is very similar to that of the standard trace formula. For a start, $R_{\mathrm{disc}}(f_E)$ is an integral operator on $L^2(\xi_E)$, whose kernel is given by either the geometric expansion (39) or the spectral expression (40). Both expansions are actually valid as stated for any $G$ (taken to be a connected component of a nonconnected reductive group $G^+$ over $F$), so long as we understand that $x$ and $y$ are variables in $G^0(\mathbb{A})$, while in (40), $P$ ranges over standard parabolic *subsets* of $G$ over $F$. (See [16, §1].) Langlands sketched the proof of the (twisted) trace formula in §10 of [149], following the derivation of the standard trace formula for GL(2) from [103]. There is a minor difference here in the standard trace formula for GL(2) used in [103]. Instead of the character $\eta$ on $Z(F) \setminus Z(\mathbb{A}_F)$ in [103], Langlands here takes a character $\xi$ on the subgroup

$$Z(F) \cap N_{E/F}(Z(\mathbb{A}_E)) \setminus N_{E/F}(Z(\mathbb{A}_E))$$

of finite index, and then takes $\xi_E$ to be the pullback of $\xi$ to $Z_E(F) \setminus Z_E(\mathbb{A}_F)$ under the norm map $N_{E/F}$.

To compare the trace formulas, Langlands introduced a local correspondence

$$f_E = \prod_v f_{E,v} \to \prod_v f_v = f \tag{51}$$

of global test functions. For a given valuation $v$ of $F$, let $\gamma_v$ be a regular orbit in $G_{E,v} = G_E(F_v)$ under conjugation by the group $G_{E,v}^0 = G_E^0(F_v)$. Then its $\ell^{\mathrm{th}}$ power $\gamma_v^\ell$ is a subset of $G_{E,v}^0$, which Langlands intersected with the subgroup $G_v = G(F_v)$. He then proved that the mapping

$$\gamma_v \to \delta_v = G_v \cap \gamma_v^\ell$$

is a well defined injection from the regular $G_{E,v}^0$-orbits in $G_{E,v}$ to the regular conjugacy classes in $G_v$. (See [149, §4 and §8].) Our $G_{E,v}^0$-conjugacy in $G_{E,v}$ is the same as $\sigma_E$-conjugacy in $G_{E,v}^0$, while the power $\gamma_v^\ell$ becomes the norm $N_{E_v/F_v}$ in the group $G^0(E_v)$. (My apologies for reversing the notation for $\gamma$ and $\delta$ in [149], for the sole purpose of having the formulas (39) and (40) for the kernel extend as stated to a twisted group. It is actually quite appropriate, for it helps to unify more sophisticated notions of endoscopy.) Langlands then defined the correspondence $f_{E,v} \to f_v$ from local test functions $f_{E,v} \in C_c^\infty(G_{E,v}, \xi_{E,v})$ to functions $f_v \in C_c^\infty(G_{v,\mathrm{reg}}, \xi_v)$ in terms of orbital integrals

$$\mathrm{Orb}(\gamma_v, f_{E,v}) = \int_{G_{E,\gamma_v}^0(F_v) \backslash G_E^0(F_v)} f_{E,v}(x_v^{-1} \gamma_v x_v) \, dx_v$$

and

$$\mathrm{Orb}(\delta_v, f_v) = \int_{G_{\delta_v}(F_v) \backslash G(F_v)} f_v(x_v^{-1} \delta_v x_v) \, dx_v$$

by setting

$$\mathrm{Orb}(\delta_v, f_v) = \begin{cases} \mathrm{Orb}(\gamma_v, f_{E,v}), & \text{if } \gamma_v \to \delta_v, \\ 0, & \text{if there is no such } \gamma_v. \end{cases} \tag{52}$$

This determines the orbital integrals of $f_v$ uniquely, but of course not the function itself.

The goal of Langlands was to compare the terms in the two trace formulas, for any pair of matching global functions $f_E$ and $f$ in (51). There were three serious problems, and what we might call a curiosity, to be resolved along the way, none of which arose in the comparison [103, §16] for quaternion groups.

The first problem is obvious. One needs to know that for $v$ and $f_{E,v}$, the smooth function $f_v$ on the open subset $G_{v,\mathrm{reg}}$ of regular elements in $G_v$ can be chosen to have an extension (a priori unique) to a test function in the space $C_c^\infty(G_v, \xi_v)$. Langlands solved the problem in §6 and §8 of [149]. He called the set of regular orbital integrals of any function $f_v \in C_c^\infty(G_v, \xi_v)$ a *Harish-Chandra family*, and the set of functions of points $\delta_v$ obtained as above (with the vanishing condition) from a function $f_{E,v}$ in $C_c^\infty(G_{E,v}, \xi_{E,v})$ a *Shintani family*. In Lemma 6.2 of [149], he proved that any Shintani family is a Harish-Chandra family, and conversely, that any Harish-Chandra family with the appropriate vanishing condition is a Shintani family. His solution required a careful comparison of the singularities of the two kinds of functions at points $\delta_v$ and $\gamma_v$ near the boundary of their common domain, while keeping track of the invariant measures used to define the orbital integrals. The case that $E$ splits at $v$ was treated separately in §8 of [149]. It is much simpler. For in this case, there is a canonical choice for the function $f_v$, as a convolution of $\ell$ different functions in $C_c^\infty(G_v, \xi_v)$.

The second problem is what later became known as the fundamental lemma. It amounts to more explicit version for spherical functions of the local transfer mapping we have just discussed. In its most basic form it applies to the characteristic function $\mathbf{1}_v$ of $G(\mathcal{O}_v)$ in $G(F_v)$ and the characteristic function $\mathbf{1}_{E,v}$ of $G_E(\mathcal{O}_v)$ in $G_E(F_v)$, at any unramified valuation $v$ for $E/F$. The assertion is that $\mathbf{1}_v$ represents the image of $\mathbf{1}_{E,v}$ under the transfer mapping $f_{E,v} \to f_v$, or in other words, that the orbital integrals of $\mathbf{1}_{E,v}$ and $\mathbf{1}_v$ correspond as above. This can be regarded as the remaining step in the proof of the global correspondence (51). Namely, if $f_E$ lies in $C_c^\infty(G_E(\mathbb{A}_F), \xi_E)$, the function $f$ is itself globally smooth, in the sense that it lies in the space $C_c^\infty(G(\mathbb{A}_F), \xi)$. The fundamental lemma here was established by Langlands in §4 of [149] for general spherical functions. As a series of relations among the vertices in certain bounded subsets of the Bruhat–Tits tree for $\mathrm{SL}(2, F_v)$, it appeared at the time to be an interesting but purely combinatorial question.

The third problem concerns the more complicated "parabolic terms" in the trace formula for $G$ (and $G_E$). We recall that they vanished in the comparison of [103, §16], because the test function $f$ was cuspidal at two places. No such restriction was permitted here for what Langlands had in mind. He was able to handle the comparison of these terms by an extended analytic argument over the final three Sections 9–11 of [149]. We shall give a short description of it, if only to give a reader the chance to look up and compare the corresponding terms, (i), (ii), (iv), (v), (vi), (vii) and (viii) from [103, p. 516–517] and (10.9), (10.8), (10.12), (10.15), (10.30), (10.28) and (10.29) from [149, §10], in the two trace formulas.

We write the difference of the two trace formulas schematically as

$$
\begin{aligned}
\mathrm{tr}(R_{\mathrm{disc}}(f_E)) - \mathrm{tr}(R_{\mathrm{disc}}(f)) = {} & [(10.9) - (\mathrm{i})] + [(10.8) - (\mathrm{ii})] + [(10.12) - (\mathrm{iv})] \\
& + [(10.15) - (\mathrm{v})] + [(10.30) - (\mathrm{vi})] + [(10.28) - (\mathrm{vii})] + [(10.29) - (\mathrm{viii})].
\end{aligned}
$$

The first two square brackets contain the elliptic terms from the two trace formulas. The differences they represent are each 0. The remaining five square brackets contain the supplementary parabolic terms. Two of them, the brackets containing (vi) and (vii), still represent invariant distributions (in $f_E$ and $f$). The second of these vanishes. The first one does not, but it does represent a discrete linear combination of characters (in $f_E$ and $f$). We transfer it to the left hand side of the formula, and thereby write

$$
\begin{aligned}
\mathrm{tr}(R_{\mathrm{disc}}(f_E)) - \mathrm{tr}(R_{\mathrm{disc}}(f)) - [(10.30) - (\mathrm{vi})] \\
= [(10.12) - (\mathrm{iv})] + [(10.15) - (\mathrm{v})] + [(10.29) - (\mathrm{viii})].
\end{aligned}
\tag{53}
$$

The right hand side now consists entirely of parabolic noninvariant distributions in $f_E$ and $f$. It is where the real analytic work of Langlands began. He first wrote the sum of the noninvariant **geometric** terms, those in the brackets containing (iv) and (v), as a sum over the discrete group $F^*$ of values of a certain function in $\mathbb{A}_F^*$. He then showed that this function (apart from some manageable error terms we shall ignore here), although not infinitely differentiable at the archimedean places, was smooth enough to apply the (multiplicative) Poisson summation formula for the discrete group $F^*$ of $\mathbb{A}_F^*$. Langlands then studied the formula thus obtained from (53), roughly speaking, as an identity in the local spectral parameters of the matching test functions $f_E \to f$. The left hand side is discrete in this sense, but thanks to Poisson summation, the terms on the right hand side all become at least partly continuous. Langlands' sophisticated analytic argument allowed him to deduce at length that each side of (53) vanishes. The formula

$$\mathrm{tr}(R_{\mathrm{disc}}(f_E)) + [(\mathrm{vi}) - (10.30)] = \mathrm{tr}(R_{\mathrm{disc}}(f)) \qquad (54)$$

thus followed.

The curiosity we mentioned is the anomaly $[(\mathrm{vi}) - (10.30)]$ in the last formula. It is interesting because it comes from simple but essential terms in the two trace formulas. The contribution in each case comes from the "discrete part" of the (spectral side) of the trace formula. It represents a term that comes from Eisenstein series, and is not part of the automorphic discrete spectrum. (These terms are not be confused with the one-dimensional automorphic representations that represent the noncuspidal part of the automorphic discrete spectrum.) The question is one of interpretation. What role do they play in the final formula (54)? We shall give the answer presently.

We have completed our discussion of the base change comparison of the two trace formulas. Our notation differs somewhat from [149][5]. Langlands observes that the extra term $[(\mathrm{vi}) - (10.30)]$ in (54), which vanishes unless the prime $\ell$ equals 2, is a twisted character in $f_E$. The left hand side of (54) is therefore itself a twisted character in $f_E$, which Langlands wrote as

$$\mathrm{tr}(R(\phi)R(\sigma)) = \mathrm{tr}(R_{\mathrm{disc}}(f_E)) + [(\mathrm{vi}) - (10.30)].$$

He also wrote

$$\mathrm{tr}(r(f)) = \mathrm{tr}(R_{\mathrm{disc}}(f))$$

for the character in $f$ on the right hand side of (54). The identity (54) is therefore

$$\mathrm{tr}(R(\phi)R(\sigma)) = \mathrm{tr}(r(f)), \qquad (55)$$

in his notation. It was stated by Langlands as Theorem 11.1 of [149]. His proof, which includes the arguments from Chapters 9 and 10 of [149] that we have tried to summarize, was completed finally on p. 211 of Chapter 11 (at the end of the first paragraph there). The rest of Langlands' last chapter, Propositions 11.4–11.5 and Lemmas 11.6–11.7, was then given to deriving the properties of the base change

---

[5] Langlands wrote $\phi$ for a test function on $G(\mathbb{A}_E) = G_E^0(\mathbb{A}_F)$ rather than our test function $f_E$ on $G_E(\mathbb{A}_F) = G_E^0(\mathbb{A}_F)\sigma_E$. He then took $R(\phi)$ to be the operator $R_{\mathrm{disc}}(f_E)$ on

$$\mathrm{L}^2(G_E^0(F) \backslash G_E^+(\mathbb{A}_F)) = \mathrm{L}_{\mathrm{disc}}^2(G_E^0(F) \backslash G_E^0(\mathbb{A}_F) \rtimes \Gamma_{E/F})$$

rather than the the operator $R_{\mathrm{disc}}(f_E)$ on

$$\mathrm{L}^2(G_E^0(F) \backslash G_E^0(\mathbb{A}_F)) = \mathrm{L}_{\mathrm{disc}}^2(G_E^0(F) \backslash G_E^0(\mathbb{A}_F)\sigma_E).$$

This accounts for the integer $\ell$ he inserted in the definition of p. 199 of [149].

lifting

$$\pi = \bigotimes_v \pi_v \to \pi_E = \bigotimes_v \pi_{E,v}, \quad \pi \in \Pi(G), \tag{56}$$

of automorphic representations, stated in Chapter 2 of [149], and dual to the transfer $f_{E,v} \to f_v$ of functions in (51). The means for this were of course provided by the comparison identity (54) or (55).

Suppose that $\pi_v \in \Pi(G_v)$, and that $\pi_{E,v} \in \Pi(G^0_{E,v})$ is $\sigma_E$-stable, in the sense that the representation

$$(\sigma_E \pi_{E,v})(x_v) = \pi_{E,v}(\sigma_E^{-1}(x_v)), \quad x_v \in G^0_{E,v},$$

is equivalent to $\pi_{E,v}$. The theory of Whittaker models then gives a canonical intertwining operator $I_{E,v}$, with

$$(\sigma_E \pi_{E,v})(x_v) = I_{E,v} \circ \pi_{E,v}(x_v) \circ I_{E,v}^{-1}.$$

We can say that $\pi_{E,v}$ is a *local lifting* of $\pi_v$ if its twisted character $\Theta(\pi_{E,v} \times I_{E,v}, \cdot)$ matches the character $\Theta(\pi_v, \cdot)$ of $\pi_v$ under the norm mapping. In other words,

$$\Theta((\pi_{E,v} \times I_{E,v}), x_v \times \sigma_E) = \Theta(\pi_v, N_v x_v), \quad N_v = N_{E_v/F_v}, \tag{57}$$

whenever $N_v x_v$ is regular. (See [149, p. 11 and Definition 6.1] as well as [27, p. 51].) The automorphic representation $\pi_E$ in (56) would then be a *global lifting* of $\pi$ if for each $v$, $\pi_{E,v}$ is a local lifting of $\pi_v$. The word "lifting" here is used interchangeably with the phrases "base change lifting" or "base change transfer", or even just "base change". The word *transfer* best describes the phenomenon in general settings.

We should remind ourselves that the base change of automorphic representations is what corresponds to the restriction of Galois (or Weil group) representations. That is, if $\pi$ is the functorial image $r \to \pi$ of a two dimensional representation $r$ of $W_F$, then $\pi_E$ would be the functorial image $r_E \to \pi_E$ of the restriction $r_E$ of $r$ to the subgroup $W_E$ of index 2 in $W_F$. This follows from the functorial interpretation of base change given at the beginning of this discussion. It can also be proved directly from a comparison of global $L$-functions, the strong multiplicity one theorem for GL(2), and the local Langlands classification of representations in terms of $L$-functions and $\varepsilon$-factors.

Langlands showed that every $\pi_v$ has a unique local lifting $\pi_{E,v}$, and that every $\pi$ has a unique global lifting $\pi_E$. Therefore (56) is a well defined mapping of automorphic representations from $\Pi(G)$ to $\Pi(G^0_E)$, whose restriction is easily seen to map the subset $\Pi_{\text{temp}}(G)$ to $\Pi_{\text{temp}}(G^0_E)$. (Recall that we defined $\Pi_{\text{temp}}(G)$ somewhat informally as the subset of globally tempered automorphic representations in $\Pi(G)$. It is the set of representations that occur in the spectral decomposition of $L^2(G(F) \setminus G(\mathbb{A}))$, which by the theory of Eisenstein series, is the more concrete set of irreducible representations

$$\{\mathscr{I}^G_P(\sigma) : \sigma \in \Pi_2(M), \quad P = MN\},$$

induced parabolically from unitary, automorphic representations in the (relative) automorphic discrete spectrum of $M$.) This represents a proof of functoriality for cyclic (prime order) base change. It is the fundamental assertion from among the various local and global properties of base change that Langlands derives from the comparison identity (55), and for which the reader can consult from the two lists in [149, §2].

The most important of the remaining properties is the characterization of the image of the mapping. This is the analogue for base change of the problem solved for quaternion groups by Jacquet and Langlands by the comparison of trace formulas in [103, §16]. To describe it, we assume that $\pi_E$ belongs to the subset $\Pi_{\text{temp}}(G_E^0)$ of $\Pi(G_E^0)$. With this assumption, Langlands proves that $\pi_E$ is a base change lift if and only if it is $\sigma_E$-stable, in which case its preimage is a finite subset of $\Pi_{\text{temp}}(G)$. Moreover, if $\pi_E$ belongs to the subset

$$\Pi_1(G_E^0) = \Pi_{\text{cusp},2}(G_E^0) = \Pi_{\text{cusp}}(G_E^0) \cap \Pi_2(G_E^0)$$

of cuspidal unitary representations, its preimage is the set

$$\{\pi \otimes \omega_{E/F}^k : 1 \leq k \leq \ell\}$$

of order $\ell$, where $\pi$ lies in the associated subset

$$\Pi_1(G) = \Pi_{\text{cusp},2}(G) = \Pi_{\text{cusp}}(G) \cap \Pi_2(G)$$

of cuspidal unitary representations of $G$, and the $\omega_{E/F}$ is the class field character of order $\ell$ associated to $E/F$. Conversely, suppose that $\pi$ is a representation in the subset $\Pi_1(G)$ of the domain $\Pi_{\text{temp}}(G)$. Then its base change image $\pi_E$ lies in the subset $\Pi_1(G_E^0)$ if and only if *either* $\ell \neq 2$ *or* $\ell = 2$ but $\pi \otimes \omega_{E/F}$ is not equivalent to $\pi$. In this case, $\pi$ becomes one of the $\ell$ representations in the fibre of $\pi_E$.

The remaining case that

$$\pi \cong \pi \otimes \omega_{E/F}, \quad \pi \in \Pi_1(G), \ell = 2,$$

is of special interest. It is *dihedral*, in the sense that $\pi$ is the image $r \to \pi$ under functoriality of an irreducible induced representation $r$ of the Weil group attached to the quadratic extension $E/F$. In other words, $r$ is induced from a character $\chi_E$ on the subgroup $C_E$ of index 2 in $W_{E/F}$, with $\sigma_E \chi_E \neq \chi_E$. Then $L(s,r)$ equals the entire, abelian $L$-function $L(s, \chi_E)$. In fact it is easy to check that, conversely, any dihedral representation $r$ satisfies all the necessary conditions of Theorem 12.2 of [103], and therefore corresponds to a cuspidal automorphic representation $\pi \in \Pi_1(G)$. The character $\chi_E$ of $C_E$ can of course be interpreted as an automorphic representation of $\text{GL}(1)_E$. With this interpretation, the mapping $\chi_E \to \pi$ is sometimes called *automorphic induction*.

The dihedral representations $\pi$ are the source of the extra term on the left side of (54). On the one hand, $\pi$ contributes less than expected to the discrete spectrum in the trace formula of $G$, since the fibre of $\pi_E$ consists of $\pi$ alone (rather than a set of order 2). But on the other hand, its base change lifting

$$\pi_E = \mathrm{Ind}_{B_E^0}^{G_E^0} \begin{pmatrix} \chi_E & 0 \\ 0 & \sigma_E \chi_E \end{pmatrix}$$

is an induced representation, and does not contribute at all to the cuspidal discrete spectrum of $G_E^0$ in the twisted trace formula for $G_E^0$. The extra term in (54) measures this discrepancy. It comes entirely from the explicit Eisenstein terms (vi) and (10.30) in the "discrete parts" of the two trace formulas. In the case of $G = \mathrm{GL}(2)$, or even that of $\mathrm{GL}(n)$ [27, §3.6], one can calculate the discrepancy independently of the two trace formulas. In more complex situations, however, one must undertake a full, direct computation of "discrete parts" of the relevant trace formulas. (See [23, §4].)

To summarize, base change represents a new case of functoriality, with the pair $(G_E^0, G) = (\mathrm{Res}_{E/F} \, \mathrm{GL}(2), \mathrm{GL}(2))$ in the role of a general pair $(G, G')$ from the statements in §4 and §5. But it comes with more information than would a general case. This includes the characterization of its image and the other properties we have just described, by virtue of its origin in a comparison of trace formulas. It was these supplementary properties in particular that led Langlands to the spectacular applications to Artin's conjecture and functoriality for certain two dimensional representations $\rho$ of $W_F$. They were established in §3 in [149].

There are four classes of irreducible representations

$$r \colon W_{K/F} \to \mathrm{GL}(2, \mathbb{C}),$$

for the global Weil group

$$1 \to C_K \to W_{K/F} \to \Gamma_{K/F} \to 1$$

attached to a finite Galois extension $K/F$ of number fields. Their images are dihedral, tetrahedral, octahedral and icosahedral (in the sense of geometric symmetry described below). The ultimate goal would be to show for each $r$ that $L(s, r)$ equals $L(s, \pi)$, for a cuspidal automorphic representation $\pi \in \Pi_1(G) = \Pi_{\mathrm{cusp},2}(G)$.

If the image of $r$ is dihedral, we can arrange that $K/F$ is a quadratic extension, and that $r(C_K)$ is not central in $\widehat{G} = \mathrm{GL}(2, \mathbb{C})$. It follows that $r$ is the irreducible representation induced from a character $\chi_K$ on $C_K$, and is dihedral in the sense above. There is consequently an automorphic representation $\pi \in \Pi_1(G)$ with

$$L(s, \chi_K) = L(s, r) = L(s, \pi),$$

as desired.

In the remaining cases, the image $r(C_K)$ consists of scalar matrices, since it is easy to see that $r$ would otherwise be dihedral. The composition

$$W_{K/F} \to \mathrm{GL}(2, \mathbb{C}) \to \mathrm{PGL}(2, \mathbb{C}) \xrightarrow{\sim} \mathrm{SO}(3, \mathbb{C})$$

is then a proper orthogonal representation of the Galois group $\Gamma_{K/F}$, which by contracting $K$ if necessary, we can assume is faithful. As a finite subgroup of $\mathrm{SO}(3, \mathbb{C})$, $\Gamma_{K/F}$ becomes the group of rigid proper motions of a tetrahedron, octahedron or icosahedron, or in algebraic terms, the group $A_4$, $S_4$ or $A_5$. Nothing was known about the Artin conjecture in any of these cases before Langlands' base change. He was able to use base change to establish functoriality for any tetrahedral $\rho$. This was the first progress in Artin's conjecture in fifty years.

Langlands' argument, which is the content of Section 3 of [149], is both striking and suggestive. It is also quite compressed. We shall review it for the tetrahedral case, in which he obtains complete results. We will then say a few words about Tunnell's extension of Langlands' argument that also leads to complete results in the octahedral case. For this discussion, we shall follow the standard practice of writing $\pi = \pi(r)$ for the functorial image in $\Pi_1(G)$, if it exists, of an irreducible representation $r$ of $W_{K/F}$. We shall also write $r_E$ for the restriction of $r$ to a subgroup $W_{K/E}$ of $W_{K/F}$, as we have been doing, and $\pi_E = BC_{E/F}(\pi)$ for the base change image in $\Pi_{\text{temp}}(G_E^0)$ of a representation $\pi \in \Pi_{\text{temp}}(G)$.

Suppose that $r$ is tetrahedral. Since $r(C_K)$ is contained in the group of scalar matrices in $\mathrm{GL}(2,\mathbb{C})$, given that $r$ is not dihedral, $r$ maps the corresponding quotient $\Gamma_{K/F}$ of $W_{K/F}$ into the group $\mathrm{SO}(3,\mathbb{C}) \cong \mathrm{PGL}(2,\mathbb{C})$. Its image is equal to the tetrahedral group $A_4$, a group of order 12 with normal subgroup

$$V_4 = \{1, (12)(34), (13)(24), (14)(23)\} \cong (\mathbb{Z}/2\mathbb{Z}) \times (\mathbb{Z}/2\mathbb{Z})$$

of index 3. We can identify this subgroup in turn with the bijective image in $\mathrm{PGL}(2,\mathbb{C})$ of the set

$$\left\{ \begin{pmatrix} 1 & 0 \\ 0 & 1 \end{pmatrix}, \begin{pmatrix} 1 & 0 \\ 0 & -1 \end{pmatrix}, \begin{pmatrix} 0 & 1 \\ 1 & 0 \end{pmatrix}, \begin{pmatrix} 0 & 1 \\ -1 & 0 \end{pmatrix} \right\}.$$

Let $E$ be the Galois extension of $F$ of degree 3 in $K/F$ fixed by the subgroup $V_4$ of $A_4$. One sees easily that the restriction $r_E$ of $r$ to $W_{K/E}$ is dihedral. Indeed, the image of $r_E$ is a semidirect product

$$\left\{ \begin{pmatrix} z & 0 \\ 0 & z\varepsilon \end{pmatrix} \right\} \rtimes \left\langle \begin{pmatrix} 0 & 1 \\ 1 & 0 \end{pmatrix} \right\rangle,$$

where $\begin{pmatrix} z & 0 \\ 0 & z \end{pmatrix}$ ranges over the nontrivial image of $C_K$ and $\varepsilon$ ranges over the image of the quadratic character attached to the subgroup $\left\{ \begin{pmatrix} 1 & 0 \\ 0 & \pm 1 \end{pmatrix} \right\}$ of $\Gamma_{E/F}$. It therefore has a functorial image $\pi_E = \pi(r_E)$ in $\Pi_1(G_E^0)$. It is also easy to see that $\sigma_E \pi_E$ is isomorphic to $\pi_E$, by inspection of the action of $\sigma_E$ on the normal subgroup $\Gamma_{K/E}$ of $\Gamma_{K/F}$. Therefore $\pi_E = BC_{E/F}(\pi)$ is the base change lift of a cuspidal automorphic representation $\pi \in \Pi_1(G)$ of $\mathrm{GL}(2)$ over $F$. This last representation is uniquely determined only as an element in the set

$$\{\pi \otimes \omega_{E/F}^k : 1 \leq k \leq 3\}. \tag{58}$$

Now the determinant $\omega_r$ of $r$ and the central character $\omega_\pi$ of $\pi$ are both characters on $C_F$, which pull back under the norm mapping to the same character on $C_E$. They therefore differ by a uniquely determined power of the class field character $\omega_{E/F}$. Replacing $\pi$ by its product with this power of $\omega_{E/F}$, which is to say the unique element in (58), we can assume that $\omega_r = \omega_\pi$, and hence that $\pi$ is uniquely determined by $r$. We would expect that $\pi$ equals the functorial image $\pi(r)$ of $r$, but perhaps surprisingly at first glance, we do not yet have enough information to prove it. In pointing this out, Langlands noted that the properties of global base change we have described above establish that if $r$ *does* have a functorial image, it must necessarily be equal to $\pi$. (See [149, p. 25], where Langlands writes $\pi_{\mathrm{ps}}(r) = \pi_{\mathrm{pseudo}}(r)$ for $\pi$.)

What is missing? By strong multiplicity 1 for GL(2) [189], the representation $\pi \in \Pi_1(G)$ is uniquely determined by the family

$$\{c_v(\pi) = c(\pi_v) : v \notin S\}$$

of semisimple conjugacy classes in $\widehat{G} = \mathrm{GL}(2,\mathbb{C})$, for a finite set $S \supset S_\infty$ of valuations of $F$. We can of course choose $S$ so that $r$ as well as $\pi$ is unramified at any $v \notin S$. We then also have the family

$$\{r_v(\Phi) = r(\Phi_v) : v \notin S\},$$

where $\Phi = \Phi_v$ is the Frobenius class in $W_{K/F}$. Then $\pi$ equals the functorial image $\pi(r)$ of $r$ if and only if the two families are equal. We know that $\pi_E = \pi(r_E)$, and hence that

$$c(\pi_{E,w}) = r_E(\Phi_w)$$

for any $w$ outside the set $S_E$ of valuations of $E$ over $S$. Moreover, if $w$ lies above $v$, it follows from the definition of base change that $r_E(\Phi_w) = r(\Phi_v)^{n(w)}$ and $c(\pi_{E,w}) = c(\pi_v)^{n(w)}$, where $n(w)$ equals the degree $[E_w : F_v]$. Therefore

$$c(\pi_v)^{n(w)} = r(\Phi_v)^{n(w)}.$$

There are two possibilities for $v$. If $v$ splits completely in $E$, $n(w) = 1$, and

$$c(\pi_v) = r(\Phi_v),$$

as required. Otherwise $v$ is inert, in which case $n(w) = 3$, and we have only the relation

$$r(\Phi_v)^3 = c(\pi_v)^3.$$

Therefore, if $r(\Phi_v) = \begin{pmatrix} a_v & 0 \\ 0 & b_v \end{pmatrix}$, for numbers $a_v, b_v \in \mathbb{C}^*$, then $c(\pi_v)$ is conjugate to $\begin{pmatrix} \xi_1 a_v & 0 \\ 0 & \xi_2 b_v \end{pmatrix}$ with $\xi_1^3 = \xi_2^3 = 1$. But the central character $\omega_\pi$ of $\pi$ equals the determinant of $c(\pi_v)$, which is therefore equal to the determinant $\omega_r$ of $r$. It follows that

$$c(\pi_v) = \begin{pmatrix} \xi a_v & 0 \\ 0 & \xi^2 b_v \end{pmatrix}, \quad v \notin S, \tag{59}$$

for a complex number $\xi = \xi_v$ with $\xi^3 = 1$.

This was as far as the purely base change argument went. It still remained to be shown that $\xi = 1$. The means to do so came from two new cases of functoriality established shortly before base change, both related to the diagram

$$
\begin{array}{ccc}
PGL(2,\mathbb{C}) & = & SO(3,\mathbb{C}) \\
\uparrow & & \downarrow \\
GL(2,\mathbb{C}) & \longrightarrow & SL(3,\mathbb{C}) \\
& \phi \searrow & \downarrow \\
& & GL(3,\mathbb{C}),
\end{array}
$$

for the dual groups $\widehat{G}_1 = \mathrm{PGL}(2,\mathbb{C})$, $\widehat{H}_1 = \mathrm{SL}(3,\mathbb{C})$ and $\widehat{H} = \mathrm{GL}(3,\mathbb{C})$ of $G_1 = \mathrm{SL}(2)$, $H_1 = \mathrm{PGL}(3)$ and $H = \mathrm{GL}(3)$ respectively, and for $\phi$ the adjoint representation of $\mathrm{GL}(2,\mathbb{C})$. The first was due to Jacquet, Piatetskii-Shapiro and Shalika [104]. It implied functoriality for the 3-dimensional Galois representation $\sigma = \phi \circ r$. In other words there is cuspidal automorphic representation $\pi^1 = \pi(\sigma)$ of $\mathrm{GL}(3)$ such that

$$c(\pi_v^1) = \sigma(\Phi_v) = \phi(r(\Phi_v))$$

for almost all $v$. The second was due to Gelbart and Jacquet [77]. It implied functoriality for $\phi$, and the cuspidal automorphic representation $\pi = \pi_{\mathrm{ps}}(r)$ of $\mathrm{GL}(2)$ we described above. Namely, there is a cuspidal automorphic representation $\pi^2$ of $\mathrm{GL}(3)$ such that

$$c(\pi_v^2) = \phi(c(\pi_v))$$

for almost all $v$. Notice that these two families of conjugacy classes are obtained by composing those of $r$ and $\pi$, the ones we are trying to see are equal, by $\phi$. It was therefore expected $\pi^1$ be equivalent to $\pi^2$. But even this required something else.

Langlands established the equivalence of $\pi^1$ and $\pi^2$ by using a fundamental criterion of Jacquet, Piatetskii-Shapiro and Shalika [104], based on the Rankin–Selberg $L$-functions

$$L(s, \pi_1 \times \pi_2), \quad \pi_i \in \Pi_1(\mathrm{GL}(n_i)), i = 1, 2,$$

they had recently constructed [106]. These are the automorphic $L$-functions for the group $\mathrm{GL}(n_1) \times \mathrm{GL}(n_2)$ attached to the tensor product representation of degree $n_1 n_2$. As in the special case of Tate [236] (with $n_1 = n_2 = 1$) and Godement–Jacquet [81] (with $n_1 = n$ and $n_2 = 1$), the authors established the analytic continuation and functional equation, and what is relevant here, a criterion for $L(s, \pi_1 \times \pi_2)$ to have a pole. It is that $L(s, \pi_1 \times \pi_2)$ is entire unless $n_1 = n_2$ and $\pi_2$ equals the contragredient $\pi_1^\vee$ of $\pi_1$, in which case $L(s, \pi_1 \times \pi_2)$ has a simple pole at $s = 1$.

Langlands applied the criterion with $n_1 = n_2 = 3$, $\pi_1 = \pi^2$ and $\pi_2 = (\pi^1)^\vee$. The condition implies that $\pi^2$ is equivalent to $\pi^1$ if and only if

$$L(s, \pi_v^1 \times (\pi_v^1)^\vee) = L(s, \pi_v^2 \times (\pi_v^1)^\vee)$$

for almost all $v$, by strong multiplicity 1. That is, if and only if

$$\det(1 - |\varpi_v|^s \cdot c(\pi_v^1) \otimes {}^t c(\pi_v^1)^{-1}) = \det(1 - |\varpi_v|^s \cdot c(\pi_v^2) \otimes {}^t c(\pi_v^1)^{-1})$$

for almost all $v$, where $\varpi_v$ is a uniformizing parameter for $v$. This would of course be implied by the desired equality of the classes $c(\pi_v^1) = \phi(r(\Phi_v))$ and $c(\pi_v^2) = \phi(r(\pi_v))$, but it is in fact something that can be checked directly. Langlands did so in [149, p. 27], appealing implicitly to the fact he had noted earlier that $\sigma$ is the representation of $\Gamma_{K/F} = A_4$ induced from a character $\theta$ of order 3 on the subgroup $\Gamma_{K/F} = V_4$. This established that $\phi(c(\pi_v))$ equals $\phi(r(\Phi_v))$ for almost all $v$, and hence that $\pi^1$ equals $\pi^2$.

The final step was to combine this with (59). The result is that the conjugacy class

$$\phi(c(\pi_v)) = \mathrm{Ad} \begin{pmatrix} \xi a_v & 0 \\ 0 & \xi^2 b_v \end{pmatrix} = \begin{pmatrix} \xi^2 a_v^2 & 0 & 0 \\ 0 & \xi^3 a_v b_v & 0 \\ 0 & 0 & \xi^4 b_v^2 \end{pmatrix} = \begin{pmatrix} \xi^2 a_v^2 & 0 & 0 \\ 0 & a_v b_v & 0 \\ 0 & 0 & \xi b_v^2 \end{pmatrix}$$

is the same as

$$\phi(r(\Phi_v)) = \mathrm{Ad} \begin{pmatrix} a_v & 0 \\ 0 & b_v \end{pmatrix} = \begin{pmatrix} a_v^2 & 0 & 0 \\ 0 & a_v b_v & 0 \\ 0 & 0 & b_v^2 \end{pmatrix} \sim \begin{pmatrix} b_v^2 & 0 & 0 \\ 0 & a_v b_v & 0 \\ 0 & 0 & a_v^2 \end{pmatrix}.$$

As Langlands then argued on p. 28 of [149], this implies that either $\xi = 1$, from which (59) then becomes

$$c(\pi_v) = \begin{pmatrix} \xi a_v & 0 \\ 0 & \xi^2 b_v \end{pmatrix} = r(\Phi_v)$$

as required, or that

$$a^2 = \xi b^2.$$

Taking square roots, one sees that this in turn implies either that $a_v = \xi^2 b_v$, which again gives $\xi = 1$ as required, or that $a_v = -\xi^2 b_v$. But if the very last condition holds, we get

$$\phi(r(\Phi_v))^3 = \begin{pmatrix} a_v & 0 \\ 0 & b_v \end{pmatrix}^3 = \begin{pmatrix} -b_v^3 & 0 \\ 0 & b_v^3 \end{pmatrix} \neq 1.$$

This would imply that the class $\phi(r(\Phi_v))$ has order 6, which is impossible, since $r(\Phi_v)$ lies in the dihedral group $A_4$.

It thus follows without exception that $c(\pi_v)$ equals $r(\Phi_v)$ for almost all $v$ The theorem of strong multiplicity 1 then yields the following result.

**Theorem (Langlands [149])** *If $F$ is a number field and $r$ is a two-dimensional representation of the Weil group $W_F$ of $F$ of tetrahedral type, there is a cuspidal* *automorphic representation $\pi$ of* GL$(2, \mathbb{A}_F)$ *such that* $\pi = \pi(r)$. *In particular, the L-function*

$$L(s, r) = L(s, \pi)$$

*is entire.*

The argument we have reviewed is by any account a highly sophisticated proof. We note in passing that the two supplementary cases of functoriality [104], [77] needed to complete the argument did not use the trace formula in their proof. They were both established in the spirit of [103, Theorem 12.2], with an extension to GL(3) by Piatetskii-Shapiro of the converse theorem of Hecke and Weil. However, they would probably also be consequences of two later applications of the trace formula. The first would be automorphic induction from GL(1) to GL(3) of the character $\theta$ on $\Gamma_{K/E}$ above. This is a special case of general construction founded on base change for GL$(n)$ [27] that we will mention at the end of the section. The second would be a special case of the general endoscopic classification in [23]. It is the correspondence between automorphic representations of the group Sp$(2) = $ SL$(2)$, with dual group PGL$(2, \mathbb{C}) = $ SO$(3, \mathbb{C})$, and self-dual automorphic representations of GL(3). We will discuss the general endoscopic classification in Section 10.

Langlands' theorem for tetrahedral representations was the foundation for its extension by Tunnell to octahedral representations. Langlands himself treated some octahedral representations over $F = \mathbb{Q}$ at the end of §3 of [149], using a converse result of Deligne and Serre [67]. Shortly thereafter, Tunnell used something different, a case of nonnormal base change for GL(2) [105] to extend Langlands' result for tetrahedral representations to general octahedral representations. We shall add a few remarks on Tunnell's proof [240].

Suppose that $K/F$ is a Galois extension of number fields, and that $r$ is a faithful, two-dimensional representation of the Galois group $\Gamma_{K/F}$ of octahedral type. In other words, the image of $r$ in PGL$(2, \mathbb{C}) \cong $ SO$(3, \mathbb{C})$ is equal to the octahedral group $S_4$. There is one respect in which the octahedral case is simpler. It is that the *binary octahedral group*, the two-fold "covering group" given by its preimage in Spin$(3, \mathbb{C})$ is a direct product $(S_4 \times \mathbb{Z}/2\mathbb{Z})$. The binary tetrahedral group, on the other hand, is the nonsplit extension $S_4$ of $A_4$. It was for this reason that we took $r$ to be a tetrahedral representation of the Weil group $W_{K/F}$ earlier. In the octahedral case, we are free to let $r$ simply be a representation of the Galois group $\Gamma_{K/F}$, as we have just done.

Let $F'/F$ be the quadratic extension in $K/F$ that is fixed by $A_4$, regarded as a normal subgroup of the Galois group $\Gamma_{K/F} \cong S_4$. Then the restriction $r' = r_{F'}$ of $r$ to $\Gamma_{K/F'} \cong A_4$ is an (irreducible) subrepresentation of tetrahedral type. We also choose a 2-Sylow subgroup $Q_8$ of Gal$(K/F)$ (from the three conjugate subgroups of $S_4$ of order 8), and take $L$ to be the fixed field of this group in $K$. Finally, we let $E$ be the composite of $L \cdot F'$. We then have the following master diagram of fields, with corresponding Galois groups indicated by parentheses, for which I am indebted to W. Casselman. It can perhaps serve as a mnemonic for the complex arguments we have described. (See also p. 174 of [240].)

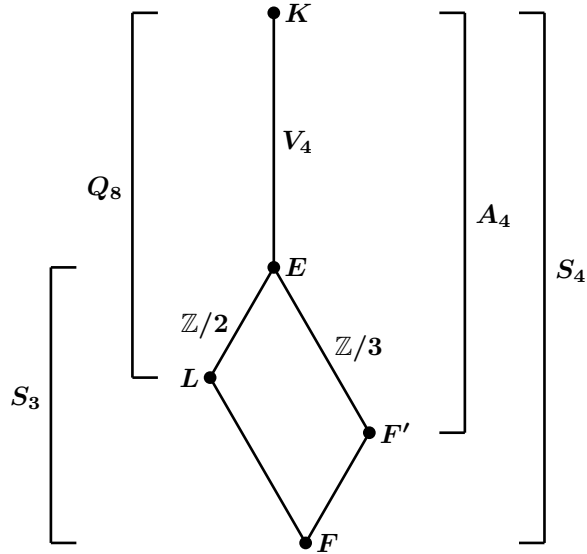

According to Langlands' tetrahedral theorem above, there is a cuspidal automorphic representation $\pi' = \pi(r')$ of $\mathrm{GL}(2, \mathbb{A}_{F'})$ that is the functorial image of $r'$. There are then exactly two cuspidal automorphic representations $\pi_1$ and $\pi_2$ of $\mathrm{GL}(2, \mathbb{A}_F)$ whose base change liftings $BC_{F'/F}(\pi_i)$ are equal to $\pi'$. They are related by $\pi_2 = \pi_1 \otimes \omega_{F'/F}$. On the other hand, the quaternion group $Q_8$ of order 8 is nilpotent, and hence monomial. It then follows from the converse theorem [103, Theorem 12.2] that there is a cuspidal automorphic representation $\pi_L = \pi(r_L)$ of $\mathrm{GL}(2, \mathbb{A}_L)$ that is the functorial image of the restriction $r_L$ of $r$ to $\Gamma_{K/L}$.

Combining these ingredients with the nonnormal base change theorem of [105], Tunnell was able to show that there is a *unique* index $i = 1, 2$ such that the base change lift $BC_{L/F}(\pi_i)$ of $\pi_i$ equals $\pi_L$ [240, Lemma, p. 174]. Following the arguments described on p. 175 of [240], he then reached the following conclusion.

**Theorem (Langlands–Tunnell [240])** *If $F$ is a number field and $r$ is a two-dimensional representation of the Galois group $\Gamma_F$ of octahedral type, there is a cuspidal automorphic representation $\pi$ of $\mathrm{GL}(2, \mathbb{A}_F)$ such that $\pi = \pi(r)$. In particular, the L-function*

$$L(s, r) = L(s, \pi)$$

*is entire.*

Motivated by Langlands' paper for $\mathrm{GL}(2)$, L. Clozel and I were later able to establish base change for the group $\mathrm{GL}(n)$ [27], again for a cyclic extension $E/F$ of prime order $\ell$. The general argument followed that of Langlands, but the comparison of trace formulas was more complicated. Recall that in our discussion of $\mathrm{GL}(2)$, there were three problems to be solved in the comparison. The first was to show that in the local transfer for functions, the image of a test function in $C_c^\infty(G_{E,v}, \xi_{E,v})$ could be chosen to lie in $C_c^\infty(G_v, \xi_v)$. Its analogue for $\mathrm{GL}(n)$ was established in [27, §1.3] by straightforward methods of local descent. The second problem was the explicit form of local transfer for the special case of nonarchimedean spherical functions (known later as the fundamental lemma). Its analogue for $\mathrm{GL}(n)$ was established in [27, §1.4]. The authors were able to rely here on a new observation [120] of Kottwitz, for the important special case of the unit functions $f_{E,v} = \mathbf{1}_{E,v}$ and $f_v = \mathbf{1}_v$. It took the form of a bijection between the summands in the two finite series that define the orbital integrals (twisted and ordinary), $\mathrm{Orb}(\gamma_v, \mathbf{1}_{E,v})$ and $\mathrm{Orb}(\delta_v, \mathbf{1}_v)$, in which the summands themselves were equal. Kottwitz was thus able to prove the required equality

$$\mathrm{Orb}(\gamma_v, \mathbf{1}_{E,v}) = \mathrm{Orb}(\delta_v, \mathbf{1}_v), \quad \gamma_v \to \delta_v$$

of the two functions without evaluating either one of them explicitly.

The relative ease with which the problems of local transfer and fundamental lemma were solved for $\mathrm{GL}(n)$ does not reflect how difficult they turned out to be in more general situations. However, the third problem, that of comparing the supplementary parabolic terms in the two trace formulas, does contain many of the difficulties encountered in general endoscopy [21], [179], [180]. Its solution, which took up much of Chapter 2 of [27], can be regarded as a blueprint for the situation in general. It is different from the global analytic argument used by Langlands for $\mathrm{GL}(2)$. The final result was a series of term by term identities (the geometric Theorem A and the spectral Theorem B, in §2.5 and §2.9 of [27]) between the constituents of the invariant trace formula for $\mathrm{GL}(n)$ (a special case of [15]) and their analogues for the twisted invariant trace formula for $\mathrm{GL}(n)_E$. Theorem B then led directly to the analogue for $\mathrm{GL}(n)$ of the spectral comparison identity (54) for $\mathrm{GL}(2)$.

The consequences of the base change comparison for $G = \mathrm{GL}(n)$ were derived in Chapter 3 of [27]. We shall state them for automorphic representations $\pi \in \Pi_{\mathrm{temp}}(G)$ and $\pi_E \in \Pi_{\mathrm{temp}}(G_E)$, even though they were derived in [27] only for the special case of "induced from cuspidal" automorphic representations, as opposed to automorphic representations induced from representations in the full discrete spectrum. The notion of a global lifting $\pi_E$ of $\pi$ is defined as for $\mathrm{GL}(2)$ above. Then

(i) Any $\pi$ has a unique global base change lifting $\pi_E$.

(ii) A given $\pi_E$ is a base change lift if and only if it is $\sigma_E$-stable, which is of course to say that the representation $\sigma_E \pi_E$ is equivalent to $\pi_E$, in which case its preimage in $\Pi_{\mathrm{temp}}(G)$ is finite.

(iii) Suppose that $\pi_E$ lies in the subset $\Pi_1(G_E) = \Pi_{\mathrm{cusp},2}(G_E)$ of cuspidal representations in $\Pi_{\mathrm{temp}}(G_E)$, and that it is $\sigma_E$-stable. Then its preimage under base change is a set of order $\ell$ of the form

$$\{\pi \otimes \omega_{E/F}^k : 1 \le k \le \ell\},$$

for some $\pi \in \Pi_1(G)$.

(iv) Suppose that $\pi_E$ belongs to $\Pi_1(G_E)$, but is *not* $\sigma_E$-stable. Then the induced representation

$$\Pi_E = \mathrm{Ind}_{P_E^0}^{G_E^0}(\pi_E \otimes \sigma_E \pi_E \otimes \cdots \otimes \sigma_E^{\ell-1}\pi_E),$$

which is an automorphic representation in $\Pi_{\mathrm{temp}}(\mathrm{GL}(n\ell)_E)$ by virtue of the theory of Eisenstein series, is $\sigma_E$-stable, and is a base change lifting of exactly one representation $\Pi \in \Pi_1(\mathrm{GL}(n\ell))$. The resulting map $\pi_E \to \Pi$ is an $n$-dimensional form of *automorphic induction*.

As an application of these results, the authors established a functorial correspondence $r \to \pi(r)$ from the irreducible $n$-dimensional representations $r$ of a *nilpotent* Galois group over $F$ to cuspidal automorphic representations $\pi = \pi(r)$ in $\Pi_1(G)$. However, since Artin's conjecture was already known for nilpotent groups, this gave no new analytic information about their $L$-functions.

Base change for $\mathrm{GL}(n)$ was used in the proof of R. Taylor, in partial collaboration with Clozel, M. Harris and N. Shepherd-Barron, in the proof of the Sato–Tate conjecture for many elliptic curves over any totally real field $F$, as we noted the end of Section 4 [238], [51], [91]. My understanding is that base change served as a substitute for functoriality in the early argument suggested by Langlands in [138].

I mention finally a recent preprint [52] of Clozel and M. S. Rajan, in which they characterize the image and fibres of solvable base change for $\mathrm{GL}(n)$. It will be interesting to see what applications come of it.

# 8 Shimura varieties

The theory of Shimura varieties has developed into a vast field, with fundamental ties to automorphic representations. I cannot do justice to the subject, and to Langlands' major contributions to it, especially in one section. I will do my best to describe some of the basic ideas and problems, with emphasis on their ties to automorphic representations. As a result, the section will perhaps be less technical, and no doubt less complete, than our earlier ones.

Shimura varieties are algebraic varieties whose complex points come from Hermitian, arithmetic, locally symmetric spaces. They are to algebraic geometry what general arithmetic locally symmetric spaces are to Riemannian geometry. The classic example is the quotient

$$\Gamma(N) \backslash X_+ = \{\gamma \in \mathrm{SL}(2,\mathbb{Z}) : \gamma \equiv 1 (\mathrm{mod}\, N)\} \backslash \{z \in \mathbb{C} : \mathrm{Im}(z) > 0\}$$

of the upper half plane, or more correctly, a certain disjoint union of such quotients. They can be compactified by adding a finite set of points, thereby becoming a (disjoint union of) compact Riemann surfaces. We shall give the formal definition of a Shimura variety $S_K$ presently, noting here only that it is defined canonically over a certain number field $E$, known as the *reflex field* of $S_K$.

Langlands' interest in Shimura varieties was in their relations with automorphic forms. As was the case for number fields, Shimura varieties also come with Galois representations, and with automorphic representations to which they are supposed to be related. The Galois representations are more complicated, in the sense that they take values in a general linear group over $\mathbb{Q}_\ell$ (or rather, some finite extension $L_\lambda$ of $\mathbb{Q}_\ell$), instead of over $\mathbb{C}$. (By convention, $\ell$ is a prime to be distinguished from $p$, another prime that might be engaged elsewhere.) However, the associated automorphic representations are expected to be more manageable, in the sense that their archimedean components $\pi_\infty$ should typically be square integrable (as in a holomorphic modular form of weight $\geq 2$), rather than an induced representation (as in a Maass form). It is for this reason that the reciprocity laws between $\ell$-adic representations and automorphic forms tend to be more concrete.

Langlands' first paper on Shimura varieties (or indeed, any aspect of algebraic geometry) is the fundamental article [140] on the basic 1-dimensional Shimura varieties $S_K$ attached to the group $\mathrm{GL}(2)$. The 1972 Antwerp conference at which Langlands gave the lectures was timely, coming soon after his long monograph [103] with Jacquet on the automorphic representations of $\mathrm{GL}(2)$. His article for these proceedings could be considered the second member of a trilogy, which also includes the monograph [149] on base change for $\mathrm{GL}(2)$ we have just discussed. Each of its three members is dedicated to a different but complementary side of the representation theory of $\mathrm{GL}(2)$.

Langlands' goal in [140] was to establish reciprocity laws between the two-dimensional $\lambda$-adic representations that arise from the étale cohomology of $S_K$ and the automorphic representations of $\mathrm{GL}(2)$. This is essentially equivalent to showing that the $L$-functions of these $\lambda$-adic representations, defined at least at the unramified places in the same way as Artin $L$-functions, are equal to automorphic $L$-functions for $\mathrm{GL}(2)$. In particular, they would have analytic continuation with functional equation to entire functions of $s \in \mathbb{C}$. This in turn should then lead ultimately to explicit formulas for the Hasse-Weil zeta functions of the varieties $S_K$ [146]. (Langlands' paper [146] is actually devoted to the compact Shimura varieties associated to various quaternion algebras, rather than the noncompact varieties $S_K$ attached to $\mathrm{GL}(2)$ over $\mathbb{Q}$.) Not surprisingly, perhaps, the power to establish such things comes from the Selberg trace formula for $\mathrm{GL}(2)$. More precisely, the results are consequences of an intricate comparison of Selberg's formula with a completely different formula, the Lefschetz trace formula, originally for a (nonsingular, projective) algebraic variety over a finite field. Let us first say something about the general theory of Shimura varieties, after which we can return to our discussion of [140] and other papers of Langlands.

In general, a Shimura variety $S_K$ is a quasiprojective variety with some auxiliary data, which is attached to a certain reductive group $G$ over $\mathbb{Q}$, and which as we have noted is naturally defined over an associated number field $E$. The subject owes its existence to the efforts of Goro Shimura over many years. Among other things, he studied the auxiliary data to be attached to various groups $G$, gave a conjectural formulation of the reflex field of definition $E$ in terms of these data, and proved the conjecture in some cases. He also studied the internal arithmetic objects of $S_K$ (as did M. Eichler at about the same time). In some cases he was able to establish reciprocity laws between these objects and automorphic forms. His most complete results were for what are now known as Shimura curves, where

$$G = \mathrm{Res}_{F/\mathbb{Q}}\, G_F$$

for the multiplicative group $G_F$ of a quaternion algebra over $F$ (as in Sections 6 and 7 here), with the requirement that $F$ be a totally real field that splits at exactly one archimedean place. In particular, $G = G_{\mathbb{Q}}$ could be the multiplicative group of a quaternion algebra that splits over $\mathbb{R}$. (See [230], and the references there.)

Deligne made a study of Shimura's work, on which he reported to Bourbaki in 1971 [61]. He reformulated Shimura's constructions in adelic terms. In this setting, a Shimura variety amounts to a family of complex varieties $S = \{S_K\}$ attached to a *Shimura datum* $(G, X)$, with $G$ being a reductive group over $\mathbb{Q}$ and $X$ a $G(\mathbb{R})$-conjugacy class of homomorphisms

$$h \colon \mathscr{R}(\mathbb{R}) = \mathbb{C}^* \to G(\mathbb{R})$$

defined over $\mathbb{R}$, for the $\mathbb{R}$-torus

$$\mathscr{R} = \mathrm{Res}_{\mathbb{C}/\mathbb{R}}(\mathbb{G}_m), \qquad \mathbb{G}_m = \mathrm{GL}(1),$$

which satisfy the natural conditions (a), (b) and (c) on p. 213–214 of [144]. In particular, if $K_h$ is the centralizer of $h(\mathbb{C})$ in $G(\mathbb{R})$ for some $h \in X$, the quotient

$$G(\mathbb{R})/K_h \cong X$$

is required to have a complex Hermitian structure. The subscripts $K$ range over open, compact, subgroups of $G(\mathbb{A}^\infty)$. For any such $K$, the associated variety has complex points parametrized by the adelic coset space

$$S_K(\mathbb{C}) = G(\mathbb{Q}) \setminus (X \times G(\mathbb{A}^\infty)/K) = G(\mathbb{Q}) \setminus G(\mathbb{A})/K_h K.$$

This space may or may not be compact. In any case, each of its (finitely many) connected components has a complex embedding into projective space, according to the Bailey–Borel compactification. Therefore, $S_K(\mathbb{C})$ is a complex quasiprojective manifold, and hence the set of complex points of a quasiprojective algebraic variety over $\mathbb{C}$. (See [61, §1.8], [144, §4], [176].)

The *Shimura variety* attached to the datum $(G,X)$ is formally taken to be the inverse limit

$$S = S(G,X) = \varprojlim_K S_K.$$

It is a complex proalgebraic variety, with complex points

$$S(\mathbb{C}) = G(\mathbb{Q}) \setminus (X \times G(\mathbb{A}^\infty)) \cong G(\mathbb{Q}) \setminus G(\mathbb{A})/K_h,$$

on which the group

$$G(\mathbb{A}^\infty) = \{x \in G(\mathbb{A}) : x_\infty = x_\mathbb{R} = 1\}$$

acts algebraically by right translation. The quasiprojective variety attached to any $K \subset G(\mathbb{A}^\infty)$ then equals the quotient

$$S_K = S_K(G,X) = S/K.$$

It is also often called a Shimura variety, as we have already done above.

The simplest examples are given by the case that $G = T$ is a torus. For $X$ then consists of a single point $h$. For an open compact subgroup $U$ of $T(\mathbb{A}^\infty)$, the set $S_U(\mathbb{C})$ is then finite, and the corresponding Shimura variety $S_U = S_U(T,h)$ is zero-dimensional. If $(G,X)$ is a general Shimura datum, a *special pair* for $(G,X)$ is a pair $(T,h)$ with $T \subset G$ and $h \in X$. Then if $U = K \cap T(\mathbb{A}^\infty)$, for an open compact subgroup $K \subset G(\mathbb{A}^\infty)$, $S_U(\mathbb{C})$ is a finite subset of $S_K(\mathbb{C})$, consisting of what are called *special points* for $S_K$. Similar notions apply to the proalgebraic complex varieties attached to the inverse limits over $U$ and $K$.

Suppose that $(G,X)$ is a Shimura datum. We first note that there are two simple homomorphisms

$$\mu, w \colon \mathbb{G}_m \to \mathscr{R},$$

which are basic to the general theory of Hodge structures [175, p. 214]. The first is defined over $\mathbb{C}$ by

$$\mu(z) = (z,1), \qquad z \in \mathbb{G}_m(\mathbb{C}) \cong \mathbb{C}^*.$$

The second is defined over $\mathbb{R}$ by

$$w(r) = r^{-1}, \qquad r \in \mathbb{G}_m(\mathbb{R}) \cong \mathbb{R}^*.$$

(In both cases, we want to keep in mind that

$$\mathscr{R}(\mathbb{R}) = \mathbb{C}^* \cong \{(z, \bar{z}) : z \in \mathbb{C}^*\} \subset \{(z_1, z_2) \in (\mathbb{C}^*)^2\} \cong \mathscr{R}(\mathbb{C}),$$

and that it is in terms of these isomorphisms that the maps are defined.) Now suppose that $h \in X$. We then also have the two homomorphisms, the *cocharacter*

$$\mu_h = h \circ \mu \colon \mathbb{G}_m \to G, \quad z \to h(z, 1),$$

which is defined over $\mathbb{C}$, and the *weight*

$$w_X = w_h = h \circ w \colon \mathbb{G}_m \to G, \quad r \to h(r)^{-1},$$

which is defined over $\mathbb{R}$. The cocharacter gives rise to a highest weight $\widehat{\mu}$ for $\widehat{G}$, which leads to a finite dimensional representation $r = r_X$ of ${}^L G$. This in turn is used to form the automorphic $L$-functions $L(s, \pi, r)$ that should ultimately be a part of the automorphic formula for the Hasse–Weil zeta function of the Shimura varieties associated to $(G, X)$. The weight of $h$ depends only on the $G(\mathbb{R})$-conjugacy class $X$ of $h$, since its image lies in the centre of $G$ ([144, condition (a), p. 213]). One can therefore write $w_X$ in place of $w_h$ as above and call it the *weight* of $(G, X)$. The homomorphisms $\mu_h$ and $w_h$ are foundations for the role of Hodge structures in the moduli of Shimura varieties ([144, §4], [175, §1]).

The cocharacter $\mu_h$ of $h \in X$ also determines the reflex field $E(G, X)$ of the datum $(G, X)$. Let $C$ be the $G(\mathbb{C})$-conjugacy class of $\mu_h$. Then $E(G, X) \subset \mathbb{C}$ is the field of definition of $C$. One can show that the intersection $C \cap G(\overline{\mathbb{Q}})$ is a $G(\overline{\mathbb{Q}})$-conjugacy class of homomorphisms $\mathbb{G}_m \to G$ over $\overline{\mathbb{Q}}$, and hence that $E(G, X) \subset \overline{\mathbb{Q}}$ is also its field of definition (see [173, Proposition 4.6(c)]). In particular, $E(G, X)$ is a number field. Finally, if $T$ is any maximal torus of $G$, the intersection $C \cap T(\overline{\mathbb{Q}})$ is a (finite) orbit in $T(\overline{\mathbb{Q}})$ under the Weyl group $W$ of $(G, T)$. It follows that $E(G, X)$ is the field of definition of this Weyl-orbit, namely the subfield of $\overline{\mathbb{Q}}$ fixed by the subgroup of elements in $\Gamma_{\mathbb{Q}} = \mathrm{Gal}(\overline{\mathbb{Q}}/\mathbb{Q})$ that stabilize the fixed finite subset $C \cap T(\overline{\mathbb{Q}})$ of $T(\overline{\mathbb{Q}})$.

We should also say something about the canonical model of $S = S(G, X)$. In general, a *model* of $S$ over a subfield $k$ of $\mathbb{C}$ is a scheme $M$ over $k$, endowed with a right action of $G(\mathbb{A}^\infty)$ over $k$, and a $G(\mathbb{A}^\infty)$-equivariant isomorphism

$$S \xrightarrow{\sim} M \times_k \mathbb{C}.$$

There could be many models of $S$ over $k$.

The *canonical model* is a model $M$ over the reflex field $E(G, X)$ with special points that satisfy a certain reciprocity law based on class field theory. If $(T, h)$ is a special pair for $(G, X)$, the reflex field $E(T, h)$ of the associated Shimura variety is contained in the reflex field $E = E(G, X)$ of $S$. For any such pair, and any element $a \in G(\mathbb{A}^\infty)$, we write $[h, a]$ for the point in

$$M(\mathbb{C}) \cong S(\mathbb{C}) = G(\mathbb{Q}) \subset X \times G(\mathbb{A}^\infty)$$

attached to the product $ha$. We also write

$$r(T,h)\colon \mathbb{A}_E^* \to T(\mathbb{A}^\infty), \quad \mathbb{A}^\infty = \mathbb{A}_{\mathbb{Q}}^\infty,$$

for the composition of elementary maps

$$\mathbb{A}_E^* \xrightarrow{Res_{E/\mathbb{Q}}(\mu_h)} T(\mathbb{A}_E) \xrightarrow{N_{E/\mathbb{Q}}} T(\mathbb{A}) \to T(\mathbb{A}^\infty),$$

obtained from the restriction of scalars functor applied to $\mu_h\colon \mathbb{G}_m \to T$, the norm map from $T(\mathbb{A}_E)$ to $T(\mathbb{A})$, and the projection $\mathbb{A} \to \mathbb{A}^\infty$ onto the finite adeles of $\mathbb{Q}$. The model $M$ is then a *canonical model* for $S$ if for every $(T,h)$, and any $a \in G(\mathbb{A}^\infty)$, the following two conditions are met.

(i) The point $[h,a]$ in $M(\mathbb{C})$ is defined over the maximal abelian extension $E^{\mathrm{ab}}$ of $E = E(T,h)$.

(ii) The special points $[h,a]$ satisfy the reciprocity law

$$\theta_E(s)[h,a] = [h,r(s)a], \qquad s \in \mathbb{A}_E^*, r = r(T,h), \tag{60}$$

where $\theta_E(s)$ is the image of $s$ in $\Gamma_E^{\mathrm{ab}} = \mathrm{Gal}(E^{\mathrm{ab}}/E)$ under the Artin map (from the Artin reciprocity law stated in Section 3).

The idea of a canonical model is remarkable. What makes it especially deep and interesting is the presence of the Artin map from abelian class field theory. The phenomenon was discovered by Shimura, who proved its existence in a number of cases. Deligne [61] established it in other cases, and as I understand it, the general case was established "somewhat independently" by Milne [173] and Borovoi [37] in the course of proving a conjecture of Langlands from [144]. (The later article [144] of Langlands will be our main topic for the next section.) As the name suggests, the canonical model $M$ of a Shimura variety $S$ is unique, up to a unique isomorphism. (See [175, Corollary 3.6].) For this reason, it is customary to identify $M$ with $S$, and then simply to regard $S = S(G,X)$ and its quotients $S_K = S_K(G,X)$ as varieties over $E(G,X)$.

The basic example is the Shimura datum $(G,X)$ in which $G = \mathrm{GL}(2)$ and $X$ is the set of $G(\mathbb{R})$-conjugates of the $\mathbb{R}$-homomorphism

$$z = a + ib \to \begin{pmatrix} a & -b \\ b & a \end{pmatrix}, \qquad z \in \mathbb{C}^*,$$

from $\mathscr{R}(\mathbb{R}) \cong \mathbb{C}^*$ to $G(\mathbb{R})$. The corresponding Shimura varieties $S_K$ were the objects Langlands studied in [140]. Before that, reciprocity laws between Frobenius classes and Hecke operators on modular curves were established by Eichler and Shimura [71], [228], [230]. It was an extension (to higher weight) of these results that Deligne [58] used to reduce the Ramanujan conjecture (for holomorphic modular forms) to the last of the Weil conjectures, which he later established in 1974 [63]. The congruence relations used by Eichler and Shimura do not generalize easily beyond modular curves, whereas the trace formulas extend in principle to arbitrary Shimura varieties. However, as in the cases of the Jacquet–Langlands correspondence and of base change, any attempt to exploit the trace formula presents an entirely new set of difficulties. Ihara [97] was the first to study the comparison in some cases, apparently following a suggestion of Sato. It was at this point that Langlands began his investigations.

Not surprisingly, the problems Langlands set out to solve were in the same spirit as those from [103, §16] and [149]. In particular, he wanted to extend the reciprocity laws to ramified places $p$. As in the Jacquet–Langlands extension of Shimizu [227] and the Shintani extension of Saito [193], this entailed reformulating the comparison in purely adelic terms. He also wanted to lay down the results for the basic noncompact varieties $S_K$ attached to $\mathrm{GL}(2)$. We recall that the noncompactness in [103, §16] was not a problem, since the test function $f$ for $\mathrm{GL}(2)$ was cuspidal at *two* places, forcing the extra terms on the geometric side of the trace formula to vanish. For the Langlands extension of Shintani [231] in base change, however, the noncompactness was very much a problem, since $f$ could not be assumed to be cuspidal at any places. Its resolution required Section 9 of [149], titled "The Primitive State of our Subject Revealed", which we discussed briefly at the end of our last section. In the case here, the problem for Langlands was to extend the Lefschetz formula to the open Shimura varieties $S_K$ for $\mathrm{GL}(2)$, and to compare the results with the Selberg formula. This was a serious task, and the main reason for the length of [140]. We do note that the difficulties for the Selberg formula were halfway between those of [103] and [149]. For this case, the test function $f$ is taken to be the cuspidal at *one* place.

The $\lambda$-adic representations of $\Gamma_{\mathbb{Q}} = \mathrm{Gal}(\overline{\mathbb{Q}}/\mathbb{Q})$ are on the étale cohomology groups $H^1_{\mathrm{et}}(S_K)$ of $S_K$. What does this have to do with automorphic representations? We have mentioned the later paper [145] of Langlands, in which he elucidated the precise relationship between automorphic representations and automorphic forms. The latter objects are well named. They are closely related to differential forms on $S_K(\mathbb{C})$, typically with values in a locally constant sheaf $\mathscr{F}$. They can therefore be used to construct de Rham cohomology groups $H^1_c(S_K(\mathbb{C}), \mathscr{F}(\mathbb{C}))$, which in turn lead to an interpretation of these groups in terms of automorphic representations.

Langlands studied the de Rham cohomology in §2 of [140]. In §3, he discussed the implications for its structure, particularly as it relates to Hecke operators. Hecke operators act on the analytic cohomology $H^1_c(S_K(\mathbb{C}), \mathscr{F}(\mathbb{C}))$ because they act on $S_K(\mathbb{C})$ as an analytic manifold. They act on the arithmetic cohomology $H^1_{\mathrm{et}}(S_K, \mathscr{F})$ because they act on $S_K$ as an algebraic variety over the reflex field $E(G, h) = \mathbb{Q}$. The intertwining operators between these two actions, combined with the theorem of strong multiplicity 1 for $\mathrm{GL}(2)$, led formally at the end of §3 to a general bijective correspondence

$$\pi \to \sigma(\pi), \tag{61}$$

between (certain) automorphic representations $\pi$ of $G$ and (certain) two-dimensional $\lambda$-adic representations $\sigma$ of $\Gamma_{\mathbb{Q}} = \mathrm{Gal}(\overline{\mathbb{Q}}/\mathbb{Q})$. The reciprocity problem was then to describe the correspondence explicitly. Langlands was able to formulate this as a precise conjecture at the end of §4. The rest of the paper [140] was devoted to the proof of two (out of three) cases of the conjecture.

We are treating [140] as the foundation of Langlands' contributions to arithmetic geometry. We shall elaborate a little further on this basic paper before turning more briefly to his subsequent articles, and to some of the new ideas they represent.

We should begin by considering the de Rham cohomology for the complex variety $S_K(\mathbb{C})$ treated in [140, §2]. Langlands was of course working with the Shimura variety attached to GL(2), but since the ideas have natural and interesting generalizations, we assume for the moment that $S$ is attached to an arbitrary Shimura datum $(G,X)$. The space $S_K(\mathbb{C})$ is typically noncompact, so one must account for the behavior of differential forms at infinity. Langlands takes the image of the cohomology of compact support $H^*_c(\cdot)$ in the full de Rham cohomology. In general, it is better to work with $L^2$ (de Rham)-cohomology $H^*_{(2)}(\cdot)$, as has been the custom since the introduction of intersection cohomology, with its role in Zucker's conjecture. We write $\mathbb{A}_{\mathrm{fin}} = \mathbb{A}^\infty$, as is convenient, and take the open compact subgroup $K \subset G(\mathbb{A}_{\mathrm{fin}})$ to be small enough so that $S_K(\mathbb{C})$ is nonsingular.

We write $(\xi, V) = (\xi, V_\xi)$ for a fixed irreducible, finite dimensional rational representation of $G$, following notation from Langlands' later article [144, §4] (rather than his notation $(\mu, L)$ from [140, §2]). Then

$$\mathscr{F} = \mathscr{F}_\xi = V_\xi \times_{G(\mathbb{Q})} (X \times G(\mathbb{A}_{\mathrm{fin}})/K),$$

the space of $G(\mathbb{Q})$-orbits in $V \times (X \times G(\mathbb{A}_{\mathrm{fin}})/K)$ under the action

$$\gamma : v \times (x,h) \to (\xi(\gamma),v) \times (\gamma x, h), \qquad \gamma \in G(\mathbb{Q}), h \in G(\mathbb{A}_{\mathrm{fin}})/K,$$

is a locally constant sheaf $\mathscr{F}_\xi(\mathbb{C})$. One is interested in the $L^2$-cohomology

$$H^*_{(2)}(S_K(\mathbb{C}), \mathscr{F}) = \bigoplus_{d=0}^{2n} H^d_{(2)}(S_K(\mathbb{C}), \mathscr{F}), \qquad n = \dim S_K(\mathbb{C}), \tag{62}$$

of $S_K(\mathbb{C})$ with coefficients in $\mathscr{F}$. We recall that this is the cohomology of the complex of $\mathscr{F}$-valued, smooth differentiable forms $\omega$ on $S_K(\mathbb{C})$ such that both $\omega$ and $d\omega$ are square integrable.

The graded, complex vector space (62) has a spectral decomposition

$$\bigoplus_\pi \Big( m_2(\pi) \cdot H^*(\mathfrak{g}_\mathbb{R}, K_\mathbb{R}; \pi_\mathbb{R} \otimes \xi) \otimes \pi^K_{\mathrm{fin}} \Big), \tag{63}$$

where $\pi = \pi_\mathbb{R} \otimes \pi_{\mathrm{fin}}$ ranges over automorphic representations of $G$ with archimedean and nonarchimedean components $\pi_\mathbb{R} = \pi_\infty$ and $\pi_{\mathrm{fin}} = \pi^\infty$ as indicated, while $m_2(\pi)$ is the multiplicity with which $\pi$ occurs in the $L^2$-discrete spectrum (with appropriate central character determined by $\xi$)[6], and $\pi^K_{\mathrm{fin}}$ is the finite dimensional space

---

[6] If $\xi$ is nontrivial, this definition requires further comment. In general, one takes functions $\phi$ on $G(\mathbb{Q}) \backslash G(\mathbb{A})$ that are $\xi(z)^{-1}$-equivariant under translation by any $z \in Z(\mathbb{R})^0$. But since $\xi$ is not generally unitary on $Z(\mathbb{R})^\circ$, one must also scale these functions by a fixed function on $G(\mathbb{Q}) \backslash G(\mathbb{A})$ whose restriction to $Z(\mathbb{R})^\circ$ equals the character $|\xi(z)|$. We shall discuss the case of GL(2) presently, following [140, p. 379]

of $K$-invariant vectors for $\pi_{\text{fin}}$. The informal proof of this decomposition is essentially a consequence [36, VII] of the definition [36, §1.5.1] of the remaining factor, the graded, finite-dimensional, complex vector space $H^*(\mathfrak{g}_{\mathbb{R}}, K_{\mathbb{R}}; \cdot)$ of $(\mathfrak{g}_{\mathbb{R}}, K_{\mathbb{R}})$-cohomology, in which $\mathfrak{g}_{\mathbb{R}}$ is the Lie algebra of $G(\mathbb{R})$, and $K_{\mathbb{R}}$ is the stabilizer of a chosen point in $X$. The formal proof for $L^2$-cohomology is in [34].

For $(\mathfrak{g}_{\mathbb{R}}, K_{\mathbb{R}})$-cohomology, we recall that a $(\mathfrak{g}_{\mathbb{R}}, K_{\mathbb{R}})$-module is a (semisimple, locally $K_{\mathbb{R}}$-finite) complex, $K_{\mathbb{R}}$-module $M$, with an action of $\mathfrak{g}_{\mathbb{R}}$ that is compatible with the adjoint action of $K_{\mathbb{R}}$ on $\mathfrak{g}_{\mathbb{R}}$. As we noted in §2, this is the same thing as a module over the real Hecke algebra $\mathscr{H}_{\mathbb{R}}$. However, we retain the Lie algebra formulation here, to emphasize its relation with differential forms. In general, the $(\mathfrak{g}_{\mathbb{R}}, K_{\mathbb{R}})$-cohomology $H^*(\mathfrak{g}_{\mathbb{R}}, K_{\mathbb{R}}; M)$ of $M$ is the $(\mathfrak{g}_{\mathbb{R}}, K_{\mathbb{R}})$-variant of the usual Lie algebra cohomology (See [36, §I.2.2].) In the case at hand, the product $\pi_{\mathbb{R}} \otimes \xi$ in (62) stands for the $(\mathfrak{g}_{\mathbb{R}}, K_{\mathbb{R}})$-module

$$M = V(\pi_{\mathbb{R}}, K_{\mathbb{R}}) \otimes V_{\xi}$$

where $M = V(\pi_{\mathbb{R}}, K_{\mathbb{R}})$ stands for the space of $K_{\mathbb{R}}$-finite vectors in the space $V(\pi_{\mathbb{R}})$ on which $\pi_{\mathbb{R}}$ acts. It is easily seen to vanish unless the infinitesimal character and central character of $\pi_{\mathbb{R}}$ equal those of $\xi$. (See [36, Theorem I.5.3].)

Consider now the case of [140], the Shimura variety attached to the group $G = \text{GL}(2)$ above. We take the representation $\xi$ of $G$ as

$$\xi = \xi_k \otimes (\det)^m,$$

where

$$\xi_k = \text{sym}^{k-1}(\text{St})$$

is the $(k-1)$-symmetric power of the standard representation of $G$ (of dimension $k$), and det is the 1-dimensional determinant representation, for integers $k \in \mathbb{N}$ and $m$, with $0 \leq m < k$. The character of $\xi$ at a diagonal matrix then equals

$$\begin{aligned}
\text{tr}\left( \xi \begin{pmatrix} \alpha & 0 \\ 0 & \beta \end{pmatrix} \right) &= (\alpha^{k-1} + \alpha^{k-2}\beta + \cdots + \alpha\beta^{k-2} + \beta^{k-1})(\alpha\beta)^m \\
&= (\alpha^{k+m-1}\beta^m + \alpha^{k+m-2}\beta^{m+1} + \cdots + \alpha^m\beta^{k+m-1}) \\
&= \frac{\alpha^n\beta^m - \alpha^m\beta^n}{\alpha - \beta}, \qquad n = k + m,
\end{aligned}$$

as on p. 389 of [140]. For each such $\xi$, the cohomology (61) is supported in degrees 0, 1 and 2, but as usual, it is the middle degree $d = 1$ that is most interesting. In this case, there is precisely one irreducible representation $\pi_{\mathbb{R}} = \pi_{\mathbb{R}}(\xi)$ that contributes to the $(\mathfrak{g}_{\mathbb{R}}, K_{\mathbb{R}})$-cohomology in (62). It is characterized by the properties

(i) $\pi_{\mathbb{R}}(z) = \xi(z^{-1}I), \qquad z \in Z(\mathbb{R})^{\circ}$,

(ii) the representation

$$x_{\mathbb{R}} \to |\xi(\det x_{\mathbb{R}})|^{\frac{1}{2}} \pi_{\mathbb{R}}(x_{\mathbb{R}}), \qquad x_{\mathbb{R}} \in G(\mathbb{R}), \tag{64}$$

is unitary[7]

These conditions imply that the vector space

$$H^1(\pi_\mathbb{R}, \xi) = H^1(\mathfrak{g}_\mathbb{R}, K_\mathbb{R}; \pi_\mathbb{R} \otimes \xi)$$

has dimension 2. (See [140, p. 388–389 and Theorem 2.10]). There is thus a decomposition

$$H^1_{(2)}(S_K(\mathbb{C}), \mathscr{F}) = \bigoplus_{\{\pi \,:\, \pi_\mathbb{R} = \pi_\mathbb{R}(\xi)\}} \left( m_2(\pi) \left( H^1(\pi_\mathbb{R}, \xi) \otimes \pi_{\text{fin}}^K \right) \right) \qquad (65)$$

of the first cohomology group. The multiplicity $m_2(\pi)$ here needs to be interpreted according to the last footnote 6. For it is understood that $\pi$ acts on the space of functions $h$ on $G(\mathbb{Q}) \setminus G(\mathbb{A})$ with

$$h(zx) = \xi(z)^{-1}h(x), \qquad z \in Z(\mathbb{R})^\circ, x \in G(\mathbb{A}),$$

such that the function

$$|\xi(\det x)|^{\frac{1}{2}}h(x)$$

is square integrable on $Z(\mathbb{R})^\circ G(\mathbb{Q}) \setminus G(\mathbb{A})$. The representation

$$(R(y)h)(x) = |\det y|^{\frac{1}{2}}h(xy), \qquad x, y \in G(\mathbb{A}),$$

of $G(\mathbb{A})$ on this space is then unitarily equivalent to the regular representation of $G(\mathbb{A})$ on $L^2(Z(\mathbb{R})^\circ G(\mathbb{Q}) \setminus G(\mathbb{A}))$. (See [140, p. 379].) This gives the interpretation of the multiplicities $m_2(\pi)$. The existence of the correspondence (61) then follows as above from the global properties of $\lambda$-adic (étale) cohomology.

Langlands' conjectural reciprocity law for $S_K$ is an explicit description of the local properties of the correspondence $\pi \to \sigma = \sigma(\pi)$. The domain consists of the automorphic representations $\pi$ that give nonzero summands in (65), while the image is the corresponding set of 2-dimensional $\lambda$-adic representations $\sigma$. For any such $\pi$, Langlands sets

$$\pi'(x) = |\det x|^{-\frac{1}{2}}\pi(x), \qquad x \in G(\mathbb{A}). \qquad (66)$$

His conjecture asserts roughly that the local components $\pi'_p$ of $\pi'$ at primes $p \neq \ell$ are images under the local Langlands correspondence[8] of Langlands parameters

---

[7] The representation in (ii) lies in the relative discrete series of $\mathrm{GL}(2, \mathbb{R})$. Even though we mentioned these objects first early in Section 1, we have so far avoided describing Harish-Chandra's general parametrization. In the case of $G = \mathrm{GL}(2, \mathbb{R})$ here, the relative discrete series with trivial central character are parametrized by the positive integers $\{k\}$. The representations $\{\pi_\mathbb{R}\}$ satisfying (i) and (ii) are therefore indexed by the same pairs $(k, m)$ that parametrize $\{\xi\}$.

[8] As we have noted before, the conjectural local Langlands correspondence has its origins in the Local Functoriality conjecture. For $G = \mathrm{GL}(2)$, it asserts a bijection from conjugacy classes of homomorphisms

$$\phi_v : L_{F_v} \to \mathrm{GL}(2, \mathbb{C}), \qquad (67)$$

known now as local Langlands parameters, and irreducible representations $\pi_v$ of $\mathrm{GL}(2, \mathbb{Q}_v)$, where attached to the local components $\sigma_p$ of $\sigma$. To set this up, he fixes an embedding of $\overline{\mathbb{Q}}$ into $\overline{\mathbb{Q}}_\ell$, and then takes on a supplementary conjectural assertion that for every element $s$ in the group

$$W_{\mathbb{Q}_p} = W_{\overline{\mathbb{Q}}_p/\mathbb{Q}_p} \subset \Gamma_{\mathbb{Q}_p} = \Gamma_{\overline{\mathbb{Q}}_p/\mathbb{Q}_p},$$

the trace of $\sigma_p(s)$ lies in the subfield $\overline{\mathbb{Q}}$ of $\overline{\mathbb{Q}}_\ell$. He then uses this in §4 (p. 403–405) to convert the $\overline{\mathbb{Q}}_\ell$-valued homomorphism

$$\sigma_p = \sigma_p(\pi)\colon \Gamma_{\mathbb{Q}_p} \to \mathrm{GL}(2, \overline{\mathbb{Q}}_\ell)$$

to a complex valued homomorphism

$$\phi_p' = \phi_p'(\pi)\colon W_{\mathbb{Q}_p} \times \mathrm{SU}(2, \mathbb{C}) \to \mathrm{GL}(2, \mathbb{C}).$$

The conjecture stated at the end of §4 then includes the supplementary assertion, with the resulting precise statement being that for any $\pi$ and $p$, $\pi_p'$ is the image under the local Langlands correspondence of the complex valued homomorphism $\phi_p'$ thus constructed.

We have largely followed the notation of Langlands from [140]. We should add a comment on the correspondence $\pi \to \pi'$ of automorphic representations in (66). It represents a transition from automorphic data to arithmetic data. Suppose for simplicity that $\xi = 1$. Consider then the archimedean parameters

$$\phi_\mathbb{R}, \phi_\mathbb{R}' \colon W_\mathbb{R} \to \mathrm{GL}(2, \mathbb{C})$$

attached to the local components $\pi_\mathbb{R}$ and $\pi_\mathbb{R}'$ of $\pi$ and $\pi'$, or rather, just the restriction of these parameters to the subgroup $\mathbb{C}^*$ of index 2 in $W_\mathbb{R}$. We can assume that the image of the restriction of each parameter lies in the group of diagonal matrices. When composed with the standard two dimensional representation of $\mathrm{GL}(2, \mathbb{C})$, the parameters $\phi_\mathbb{R}'$ gives a two dimensional representation of the group $\mathbb{C}^* = \mathbb{S}(\mathbb{R})$, which amounts to a real Hodge structure of weight 2. It follows that

$$\phi_\mathbb{R}'(z) = \begin{pmatrix} z^{-1} & 0 \\ 0 & \bar{z}^{-1} \end{pmatrix}.$$

On the other hand, the preimage $\phi_\mathbb{R}$ of $\phi_\mathbb{R}'$ should be bounded, since it is supposed to be attached to a *tempered* automorphic representation $\pi$. To see that this is so, we note that for the archimedean parameters, (66) implies the identity

---

$$L_{\mathbb{Q}_v} = \begin{cases} W_{\mathbb{Q}_v}, & \text{if } v \text{ is archimedean,} \\ W_{\mathbb{Q}_v} \times \mathrm{SU}(2), & \text{if } v \text{ is } p\text{-adic,} \end{cases}$$

such that $\pi_v$ is tempered if and only if the image of $\phi_v$ is bounded. (We shall discuss the general case in Section 10.)

$$\phi'_{\mathbb{R}}(z) = \left| \det \begin{pmatrix} z & 0 \\ 0 & \bar{z} \end{pmatrix} \right|^{-\frac{1}{2}} \phi_{\mathbb{R}}(z),$$

and hence that

$$\phi_{\mathbb{R}}(z) = |z\bar{z}|^{\frac{1}{2}} \phi'_{\mathbb{R}}(z) = \begin{pmatrix} (z/\bar{z})^{-\frac{1}{2}} & 0 \\ 0 & (\bar{z}/z)^{-\frac{1}{2}} \end{pmatrix}$$

The remaining Sections 5–7 of [140] were devoted to a new comparison of trace formulas, culminating in the proof of a significant part of the conjecture. Having spent our two last sections on the two comparisons from [103] and [149], we shall not discuss the details of the remaining comparison here, even though the Lefschetz trace formula from arithmetic geometry is quite different. Section 5 of [140] consists of some calculations on the terms in the Selberg trace formula suitable for the new comparison. Section 6 is a general description of the Selberg trace formula for GL(2), along the lines of our own discussion in the last two sections. Section 7 contains the comparison of the Selberg and Lefschetz trace formula. Langlands then uses this is prove his conjectures on the local parameters $\phi'_p$ in two cases. The first is that the two-dimensional representation $\phi'_p$ is reducible. It includes the case of good reduction, where $\phi'_p$ is a direct sum of two *unramified* quasicharacters on $W_{\mathbb{Q}_p}$, to which most if not all earlier work on the subject had been confined. The second is of "multiplicative reduction", where the homomorphism $\phi'_p$ is *special*, in the sense that it is nontrivial on the subgroup SU(2) of $L_{\mathbb{Q}_p}$. The remaining case is of "additive reduction", in which $\phi'_p$ is an irreducible two-dimensional representation of the subgroup $W_{\mathbb{Q}_p}$ of $L_{\mathbb{Q}_p}$, and the corresponding representation $\pi'_p$ of GL$(2, \mathbb{Q}_p)$ is *supercuspidal*. Langlands did not study this case. However, it was resolved soon afterwards, with an extension of the methods of Langlands by Carayol [43]. (See [152].)

We have emphasized the Langlands Antwerp paper [140] on GL(2) over some of his later contributions to Shimura varieties for a couple of reasons. One we have already mentioned is that it joins the volumes [103] and [149] discussed in the past two sections as the third member of his GL(2)-trilogy. These works were designed to illustrate Langlands' revolutionary ideas on functoriality in the simplest of cases. They are interrelated, and the influence they each have individually is amplified when they are taken together. Another reason is that [140] is essentially complete (given the later addition of Carayol for supercuspidal representations of GL(2)). With its focus on the noncompact variety $S_K(\mathbb{C})$ and on the reciprocity laws at places of bad reduction, it has served as a model for subsequent work on Shimura varieties, which continues to this day.

Langlands' Corvallis article [144] on motives will be our topic for the next section. Of his other articles, the most influential might be his remarkable conjectural description [142] of the set of points on a general Shimura variety modulo a prime of good reduction. It has been a foundation for a great deal of the work in in the subject since then. Other papers include the overview [143] of the question of the Hasse–Weil zeta function, and a more technical article [146] that answers the question for some simple Shimura varieties $S_K(G, X)$ for groups $G$ related to GL(2). There is also a short article [148] on the general question of the reduction at a prime of bad reduction, the longer article [85] (with Harder and Rapoport) on the Tate conjecture for a Shimura variety attached to a Hilbert–Blumenthal surface (with $G$ equal to $\mathrm{Res}_{F/\mathbb{Q}}\,\mathrm{GL}(2)$, for a real quadratic field $F$), and the important paper [164] with Rapoport that includes a motivic refinement of the conjecture from [142]. I have not studied these last three papers, interesting as they are, and will not have much further comment on them. A final paper [152] on Shimura varieties contains some later comments of Langlands, ostensibly on the problem from [143] of calculating the Hasse–Weil zeta function, but with observations on the other papers as well. It is very informative.

What is particularly far reaching in Langlands' later papers on Shimura varieties is the emergence of two fundamental phenomena that were also beginning to govern his work in the basic theory of automorphic forms [127], [150]. His discovery that they would have a parallel, central role in the theory of Shimura varieties seems to have been completely unexpected, perhaps because they were not critical in the Shimura varieties of small dimension that has been studied up until then. One is the question of the fundamental lemma, which we have already seen in the context of cyclic base change. It arises in the local geometric terms in the Selberg and Lefschetz trace formulas at the unramified place $p$ at which one is trying to establish the reciprocity law. The other is a broader phenomenon, which includes the fundamental lemma, and is now known as *endoscopy* rather than Langlands' original term *L-indistinguishability*.[9] This affects most of the terms in the two trace formulas one is trying to compare. For the regular elliptic orbital integrals on the geometric side of the Selberg formula, it is a reflection of the fact that two elements $\gamma_1$ and $\gamma_2$ in $G(\mathbb{Q})$ over whose $G(\mathbb{A})$-conjugacy classes one would like to integrate a test function $f$, might be conjugate over $G(\mathbb{C})$ but not over $G(\mathbb{Q})$. The theory of endoscopy will be the topic of Section 10.

Before beginning any comparison, one must first understand the explicit form taken by the Grothendieck–Lefschetz trace formula when applied to a Shimura variety $S_K$. Recall that the original Lefschetz formula [167] is an identity

$$\sum_x i(\phi, x) = \sum_{k=0}^{n} (-1)^k \,\mathrm{tr}(H^k(\phi)) \tag{68}$$

attached to any suitable mapping $\phi$, from say a compact manifold $M$ of dimension $n$ to itself. The (spectral) right hand side is an alternating sum of traces of the operators

$$H^k(\phi, \mathbb{Q})\colon H^k(M, \mathbb{Q}) \to H^k(M, \mathbb{Q}), \quad 0 \le k \le n, \tag{69}$$

on ordinary (Betti) cohomology attached to $\phi$. The (geometric) left hand side is a sum over the fixed points $x$ of $\phi$ in $M$ of certain indices $i(\phi, x)$.

---

[9] "L-indistinguishability" is a better description of what is going on. However, the fact that it has more than twice the number of syllables than does "endoscopy", together with the increasing demands being placed on mathematicians' time, may have forced the change!

Motivated by this classical formula and the Weil conjectures [248], Grothendieck introduced his version [98] for an arithmetic variety. At its simplest, it applies to a nonsingular projective variety $X$ over a finite field $\mathbb{F}$ of characteristic $p$. It is the analogue of (68) for $X$, with $\phi$ replaced by some power $\Phi$ of the Frobenius endomorphism, and $H^k(M, \mathbb{Q})$ replaced on the spectral side by $H^k(M, \mathbb{Q}_\ell)$, the $\ell$-adic (étale) cohomology of $X$ at a prime $\ell \neq p$. The geometric side becomes a sum over the finite set

$$X(\mathbb{F}'), \qquad \mathbb{F}' = (\overline{\mathbb{F}})^{\Phi}$$

of points in $X(\overline{\mathbb{F}})$ fixed by $\Phi$.

Suppose for simplicity that the Shimura reflex field $E$ of $S_K = S_K(G, X)$ equals $\mathbb{Q}$, and as usual, that $K$ is small enough that $S_K(\mathbb{C})$ is nonsingular. We fix a number field $L$ that is sufficiently large in a sense that depends in $K$, together with an embedding $L \subset \mathbb{C}$. The finite dimensional complex representation $\xi$ will then be assumed to be defined over $L$.

To proceed, it is necessary to have an explicit description of the set of points $S_K(\mathbb{F}_p)$ at an unramified prime $p$. More precisely, for any $p$ such that

$$K = K_p K^p, \qquad K^p \subset G(\mathbb{A}_{\mathrm{fin}}^p),$$

for an unramified maximal compact subgroup $K_p \subset G(\mathbb{Q}_p)$, one wants a suitable $\mathbb{Z}_p$-scheme structure on $S_K$, and a description in terms of $G$ of the set $S_K(\overline{\mathbb{F}}_p)$ of points on $S_K$ over the algebraic closure $\overline{\mathbb{F}}_p$, equipped with an action of the (geometric) Frobenius endomorphism $\mathrm{Frob}_p$. The purpose of Langlands' paper [142] was to give a conjectural such formula under very general conditions on $S_K$. He arrived at it after a study of the special case of Shimura varieties of PEL (polarization, endomorphism ring, level structure) type. PEL varieties are a rather small subset of all Shimura varieties. (See [176, §9].) However, they still represent a major generalization of Shimura curves. They had been introduced by Shimura [229], but so far as I know, without any particular interest in the Lefschetz formula.

PEL varieties parametrize abelian varieties with additional structure. They include the varieties attached to $\mathrm{GL}(2)$ [140] we discussed above, which parametrize elliptic curves $E$. They also include Siegel modular varieties, the higher dimensional generalizations attached to general symplectic groups $G = \mathrm{GSp}(2n)$.

Siegel modular varieties $S_K = S_K(G, X)$ parametrize abelian varieties $A$. More precisely, the set of complex points

$$S_K(\mathbb{C}) = G(\mathbb{Q}) \setminus (X \times G(\mathbb{A}_{\mathrm{fin}})/K), \quad G = \mathrm{GSp}(2n),$$

in $S_K$ is bijective with the set of $G(\mathbb{Q})$-orbits of pairs $(A, g)$, where $A$ is a principally polarized abelian variety over $\mathbb{C}$ up to isogeny, and $g$ is a $K$-level structure (which is to say a coset in $G(\mathbb{A}_{\mathrm{fin}})/K$). What makes the problem treated by Langlands in [142] more accessible in this case is that $S_K(G, X)$ has a canonical model, not just over the reflex field $E = \mathbb{Q}$, but also over the ring $\mathbb{Z}$, and that elevates $S_K$ to the role of a universal modular variety over $\mathrm{Spec}(\mathbb{Z})$. As I understand it, this means that there is an isomorphism class of families

$$A_K \to S_K$$

of principally polarized abelian varieties over $\mathrm{Spec}(\mathbb{Z})$, equipped with a $K$-level structure, which is *universal* in the sense that the set of all such families $A'_K \to S'_K$ over any $\mathbb{Z}$-scheme $S'_K$ is in bijection with the $\mathrm{Spec}(\mathbb{Z})$-morphisms $\phi \colon S'_K \to S_K$ under the pullback mapping

$$A'_K = \phi^* A_K.$$

I will not try to define these various terms here. But the reader can refer to [142] for the discussion of a related case, motivated by Kronecker's Jugendtraum and Hilbert's twelfth problem.

The point is that in representing a functor, $S_K$ allows one to identify the set of points

$$S_K(\overline{\mathbb{F}}_p) = \{\overline{\phi}_p \colon \mathrm{Spec}(\overline{F}_p) \to S_K\}$$

with families $\overline{A}_p = \overline{\phi}_p^* A_K$ over $\mathrm{Spec}(\overline{F}_p)$. Equipped with the action of the Frobenius endomorphism, this set amounts in turn to a classification of isogeny classes of $n$-dimensional abelian varieties over $\overline{\mathbb{F}}_p$ with $K^p$-level structure. Such objects have been well understood for some time according to Honda–Tate theory [235], [96], and can be described explicitly. These modular properties extend to the locally constant, $\lambda$-adic sheaf $\mathscr{F}_K = \mathscr{F}_{K,\xi}$ on $S_K$ attached to any finite dimensional representation $\xi$ of $G$ over $\mathbb{Q}$. They are also compatible with the Hecke correspondence $f^p$ on $S_K$ defined by right translation on $S_K(\mathbb{C})$ by any element $g^p$ in $G(\mathbb{A}_{\mathrm{fin}}^p)$. (See [124, p. 375].) The classification of abelian varieties over $\overline{\mathbb{F}}_p$ then leads to a description of the set of fixed points of the composition

$$\Phi_p \circ f^p, \qquad \Phi_p = \Phi_{p,j} = (\mathrm{Frob}_p)^j, j \in \mathbb{N}, \tag{70}$$

acting as a correspondence on the set $S_K(\overline{\mathbb{F}}_p)$. This result can be regarded as an explicit description of of the main (elliptic) part of the geometric side of the Grothendieck–Lefschetz trace formula, for the case at hand.

These are the ideas that led Langlands in [142] to conjecture a general such formula for any Shimura variety. It was a bold step, which was refined and made clearer in his later paper [164] with Rapoport. In fact, the conjectured formula was subject to further evolution, leading up to what might be its final version in §3 in the paper [123] of Kottwitz. My understanding is that this formula is still far from proved in general. The case of general PEL-Shimura varieties has itself turned out to be a challenging problem. Progress was made by Milne [172] and Zink [256], while the special case of Siegel modular varieties discussed here was established in [123] and [174]. The proof for general PEL varieties was completed in [124], [174] and [191]. (See [50] and the introduction of [124].) As an aside, I have found it difficult at times to sort out what has been established from the literature, owing no doubt to my own imperfect grasp of the technical complexities of the subject. (See the introduction to [171], which is itself quite complex.)

With a formula (either proven or conjectured) for the geometric side of the arithmetic (Grothendieck–Lefschetz) trace formula, it would then be possible to study its comparison with the geometric side of the automorphic (Arthur–Selberg) trace formula. In particular, one could consider the problems of endoscopy that had been emphasized by Langlands in his papers following [140].

Kottwitz took up the fundamental lemma in [118]. In this paper, he was able to reduce it to a more familiar problem, the identity for twisted spherical functions required for cyclic base change. This was the problem solved by Langlands in the special case of GL(2). In a subsequent paper [120], Kottwitz reduced the problem further. For the special case of the unit function in the Hecke algebra of spherical functions, he reduced the twisted fundamental lemma to its original version for orbital integrals. The special case of GL($n$) of this last reduction was, incidentally, an essential ingredient of the proof of base change for GL($n$) in [27].

In each paper, Kottwitz was able to treat the problem as a natural combinatorial identity. In fact, in each case he observed that two finite series were equal simply because there was a term by term matching of their summands. However, the original fundamental lemma has turned out to be much more subtle. Even its statement draws upon the deeper notions from endoscopy, such as stability, endoscopic groups and transfer factors, that will be part of our discussion in Section 10. In particular, it is not just a combinatorial problem. Its ultimate proof by Ngô Bao Châu, drawing on the work of Waldspurger, came considerably later [186]. It required among other things, his remarkable observation that the Hitchin fibration over the field of meromorphic functions on a compact Riemann surface matches the geometric side of the trace formula (in characteristic $p$) over the global field of functions on a smooth projective curve in characteristic $p$.

The fundamental lemma is an important ingredient for our understanding of the automorphic properties of the problem. However, the full comparison requires something more, what can be called the stabilization of the Lefschetz trace formula. As such, it becomes one of a number of such constructions, beginning with the stabilization [21] of the basic automorphic trace formula (for a quasisplit group $G^*$), its analogue [21] for an inner form $G$ of $G^*$, the stabilization of the general twisted trace formula [179], [180], the expected stabilization of a metaplectic trace formula, possible stabilizations of various relative trace formulas, and who know what else. In all such constructions, the basic ingredient is the stable trace formula for a quasisplit group $G^*$, which is defined by an inductive process from the basic automorphic (Arthur–Selberg) trace formula for $G^*$. In the other cases, the comparison is not just with this stable trace formula for $G^*$, but rather a linear combination of stable trace formulas for a collection of quasisplit groups $G'$, known as endoscopic groups (attached to the problem at hand).

In [123, §4, §7], Kottwitz stabilized the elliptic part of the Lefschetz trace formula, which is to say, his proposed formula from §3 of [123]. In so doing, he appealed to his earlier results in [121] and [119], as well as Langlands' original monograph [150] on stabilization. We shall describe his construction in a way that can be compared with our general discussion of endoscopy in Section 10.

To simplify the notation, we shall allow $G$ to represent the given Shimura datum $S_K(G, X)$, as well as the reductive group over $\mathbb{Q}$ that is its essential ingredient. This is similar to a convention from the theory of endoscopy, in which $G'$ is often used to represent a full endoscopic datum $(G', s', \mathscr{G}', \xi')$ attached to any $G$ under consideration, as well as the quasisplit group $G'$ that is its main component. While we are at it, let us also write $f$ here in place of the triplet $(\xi, \Phi_p, f^p)$ in keeping track of the dependence of the Lefschetz trace formula on these quantities. We cannot quite think of $f$ as a test function on $G(\mathbb{A})$, which we would have if we were dealing with the automorphic trace formula. However, for any endoscopic datum $G'$ attached to $G$, there is a transfer mapping

$$f \rightarrow f' = f'_\infty \cdot f'_p \cdot (f^p)'$$

from triplets $f$ to test functions $f'$ on $G'(\mathbb{A})$, defined by Kottwitz in [123, §7]. We shall write $\Lambda^p_{\mathrm{ell}}(f)$ for the conjectural elliptic part of the Lefschetz trace formula to exhibit its dependence on $p$ as well as $f$. It is by definition given by the formula [123, (3.1)].

The stabilization in [123] is then an expansion

$$\Lambda^p_{\mathrm{ell}}(f) = \sum_{G'} \iota(G, G') \widehat{S}_{\mathrm{ell}}'(f') \tag{71}$$

of this linear form in terms of corresponding *stable*, linear forms[10] in various *stable* (automorphic) trace formulas[10]. The groups $G'$ on the right hand side represent equivalence classes of elliptic endoscopic data[10] for $G$, according to the convention above. For any such $G'$, $S'_{\mathrm{ell}}$ is the elliptic part of the geometric side of the stable trace formula

$$S'_{\mathrm{geom}}(\cdot) = S'_{\mathrm{spec}}(\cdot)$$

for $G'$. The corresponding test function $f'$ for $G'$ is the transfer of $f$ mentioned above, defined by the general transfer factors Langlands and Shelstad. It is determined only up to the values it takes when paired with a stable distribution on $G'(\mathbb{A})$. With this proviso in mind, we have written $\widehat{S}_{\mathrm{ell}}'(f')$ for the *uniquely determined pairing* of the stable distribution $S'_{\mathrm{ell}}(\cdot)$ with the function $f'$. Finally, the coefficients $\iota(G, G')$ are constructed in an elementary manner from $G$ to $G'$. (See [123, Theorem 7.2].)

We should perhaps pause to remind ourselves of the ultimate goal. It is founded on the two interpretations of cohomology and the fundamental data they support. On the one hand, we have the $L^2$-de Rham cohomology $H^*_{(2)}(S_K(\mathbb{C}), \mathscr{F})$, and its decomposition (63) in terms of automorphic representations. On the other hand, we have the intersection cohomology $IH^*(\overline{S}_K(\mathbb{C}), \mathscr{F})$ of the Baily–Borel compactification $\overline{S}_K(\mathbb{C})$, equipped with the correspondences defined by Hecke operators. Zucker's conjecture [257], established around 1988 [168], [198], asserts that these two complex graded vector spaces are isomorphic. We are assuming that the sheaf $\mathscr{F}$, assigned as it is to the representation $\xi$ of $G$, is defined over the chosen number field $L$. Like ordinary Betti cohomology, one can vary the coefficient field of $IH^*(\overline{S}_K(\mathbb{C}), \mathscr{F})$, taking it here to be $L$. If $\lambda$ is any finite place of $L$ not lying over $p$,

---

[10] We are asking a reader unfamiliar with these other terms to wait until (or look ahead to) Section 10, where they will be discussed in somewhat greater detail.

the tensor product

$$IH_\lambda^* = IH^*(\overline{S}_K(\mathbb{C}), \mathscr{F}_\lambda) = IH^*(\overline{S}_K(\mathbb{C}), \mathscr{F}) \otimes_L L_\lambda$$

represents a change of coefficients from $L$ to the $\lambda$-adic field $L_\lambda$. What makes it all work is that the last space is isomorphic to the $\lambda$-adic (étale) cohomology of the variety $\overline{S}_K$ at its reduction modulo $p$. In particular, this $\lambda$-adic vector space comes with a representation of $\mathrm{Gal}(\overline{\mathbb{Q}}/\mathbb{Q}) \times \mathscr{H}_K$, where $\mathscr{H}_K$ is the Hecke algebra of compactly supported $L$-valued functions in $K \backslash G(\mathbb{A}_{\mathrm{fin}})/K$. In this setting, the spectral side of the Lefschetz trace formula is the Euler number

$$\Lambda_{\mathrm{spec}}^p(f) = \sum_{k=0}^m (-1)^k \mathrm{tr}(IH_\lambda^k(\Phi_p \times f^p)), \tag{72}$$

for the operator on $IH_\lambda^k$ attached to the composition (70). (This discussion follows the beginnings of §1 of [123].)

The actual Lefschetz formula is an identity

$$\Lambda_{\mathrm{geom}}^p(f) = \Lambda_{\mathrm{spec}}^p(f). \tag{73}$$

The geometric side $\Lambda_{\mathrm{geom}}^p(f)$ includes the elliptic part $\Lambda_{\mathrm{ell}}^p(f)$ studied by Kottwitz [123]. In general, however, there are complementary terms attached to fixed points "at infinity", or more precisely, in strata of $\overline{S}_K(\mathbb{C})$ attached to proper parabolic subgroups $P$ of $G$. (The largest stratum $S_K(\mathbb{C}) \subset \overline{S}_K(\mathbb{C})$ is attached to the group $G$ itself.) These are more difficult. However, they do not occur if $S_K(\mathbb{C})$ is already compact, and are not significant in many noncompact Shimura varieties. For example, if $G = \mathrm{Res}_{F/\mathbb{Q}}(H)$, for a split reductive group $H$ over a totally real number field $F$ with more than one archimedean valuation, the nonelliptic terms on the geometric side of the automorphic trace formula vanish [15, Theorem 7.1(b)], a property that should be reflected in the corresponding terms in $\Lambda_{\mathrm{geom}}^p(f)$. However, in the basic case $G = \mathrm{GSp}(2n)$ of Siegel modular varieties we have discussed, there are complementary terms in $\Lambda_{\mathrm{geom}}^p(f)$. Their analysis is part of the work of S. Morel.

At any rate, Kottwitz' stabilization (71) of $\Lambda_{\mathrm{ell}}^p(f)$ leads us to expect a similar stabilization

$$\Lambda_{\mathrm{geom}}^p(f) = \sum_{G'} \iota(G, G') \widehat{S}_{\mathrm{geom}}^p(f') \tag{74}$$

of the full geometric side of the Lefschetz trace formula. We should point out that the linear forms on each side of (74) are defined, since the left hand side equals that of (73), while the right hand side is given by the stable trace formula of each of the groups $G'$. It is just that we cannot say that the two sides are equal unless $\Lambda_{\mathrm{ell}}^p(f)$ equals $\Lambda_{\mathrm{geom}}^p(f)$. (Even in this last case, the inequality still rests on Kottwitz' basic conjectural formula [123, (3.1)].) Similar comments apply to the stabilization

$$\Lambda_{\mathrm{spec}}^p(f) = \sum_{G'} \iota(G, G') \widehat{S}_{\mathrm{spec}}^p(f') \tag{75}$$

of the spectral side, which would follow formally from (74), when combined with the identities $\Lambda_{\mathrm{spec}}^p(f) = \Lambda_{\mathrm{geom}}^p(f)$ (the Lefschetz formula (73)) and

$$\widehat{S}_{\mathrm{spec}}(f') = \widehat{S}_{\mathrm{geom}}(f')$$

(the stable formula for $G'$). In particular, (75) would follow from Part I of [123] (Sections 1–7) if $\Lambda_{\mathrm{ell}}^p(f)$ equals $\Lambda_{\mathrm{geom}}^p(f)$ *and* the conjectured formula [123, (3.1)] is valid.

Part II of [123] (Section 8–10) is what Kottwitz calls the destabilization of (75). The ultimate goal of the theory would be to deduce reciprocity laws for $S_K$, of the kind established by Langlands for $G = \mathrm{GL}(2)$ in [140], from the identity (75). But there are not shortcuts. The only way we can expect to prove (75) is to derive it from the geometric stabilization (74), and this would rest on the imposing task of proving the explicit fixed point formula [123, (3.1)] for $S_K$, and its extension to the compactification for $S_K$. After this, there would be a new set of problems on the interpretation of the right hand side of (75). These concern a set of conjectures on the classification of the automorphic discrete spectrum that I described in [18], and applied to the cohomology of Shimura varieties in §9 if that paper.

Kottwitz assumed these conjectures in [123]. He then compared the explicit expression they yield for the right hand side of (75) with the right hand side of (72). The former concerns the characters of automorphic representations, the latter the virtual, $\lambda$-adic characters from intersection cohomology. By manipulating the terms in (75), he obtained a concrete formula [123, (10.5)] for the alternating sum of traces in (72) in terms of automorphic characters. This in turn suggested an explicit decomposition of the virtual, $\lambda$-adic representation of $\mathrm{Gal}(\overline{\mathbb{Q}}/\mathbb{Q}) \times \mathscr{H}_{\widehat{L}}$ on the étale cohomology of $\overline{S}_K$ in terms of the automorphic discrete spectrum of $G$. It is the direct sum displayed in the last paragraph of §10 of [123], but unfortunately not labeled.

This last decomposition from [123] could be regarded as the ultimate goal. It represents something beyond the character formula given by (10.5) of [123], even with all the conjectures that were taken for granted. For one thing, the integer $j$ in (10.5) has for technical reasons to be taken to be sufficiently large. With the formula established, one could presumably extend it to all $j$ by extrapolating the behavior of each side for $j$ in the restricted range. This would give the reciprocity law for the unramified places $p$. The remaining finite set of places $p$ with $\overline{S}_K$ ramified are more serious. They would require a solution for general $G$ of the problems for $\mathrm{GL}(2)$ solved by Langlands and Carayol. Their solution would presumably give reciprocity laws for $\overline{S}_K = \overline{S}_K(G, X)$ at all places. In particular, it should ultimately express the associated motivic $L$-function

$$L(s, IH^k(\overline{S}_K, \mathscr{F})), \qquad 0 \leq k \leq 2d = 2\dim(S_K),$$

of weight $k$, and the Hasse–Weil zeta function $\zeta(s, \overline{S}_K(G, X))$ of $\overline{S}_K$, explicitly in terms of automorphic $L$-functions.

In the final Part III of [123], Kottwitz established the conjectural fixed point formula [123, (3.1)] for the Siegel modular varieties $S_K = S_K(G, X)$ attached to the groups $G = \mathrm{GSp}(2n)$. This represents a grounding of sorts for the conjectured identity whose general version had itself been the foundation for Parts I and II of the article, and for our discussion above. Siegel modular varieties have been studied more recently by Morel, as we have noted. She established explicit formulas for their complementary terms in $\Lambda^p_{\mathrm{geom}}(f)$ [182], those given by boundary components in $\overline{S}_K$ attached to proper parabolic subgroups $P$ of $G$. The stabilization (74) of $\Lambda^p_{\mathrm{geom}}(f)$ was then given an explicit form in her paper [183]. This was based on a formula [17] for the associated nonelliptic geometric terms in the automorphic trace formula, in terms of the characters of Harish-Chandra's discrete series, and the stabilization of this formula treated in [82]. The only thing now standing in the way of a complete solution for Siegel modular varieties, or at least the immediate precursor [123, (10.5)] of a complete solution, is the "destabilization" of the expansion (74) of $\Lambda^p_{\mathrm{spec}}(f)$. This requires the local and global conjectures in [18] for the group $\mathrm{GSp}(2n)$. The conjectures have now been established for the group $\mathrm{Sp}(2n)$ [23]. It is obviously important to try to extend them to $\mathrm{GSp}(2n)$, even though there are indications that new difficulties arise.

Siegel modular varieties attached to the groups $G = \mathrm{GSp}(2n)$ represent a fundamental class of examples. We can think of their role in the general theory of Shimura varieties as something akin to that of the groups $\mathrm{GL}(N)$ in the general theory of automorphic forms. The groups $\mathrm{GL}(N)$ themselves are disqualified from this larger role since they are only attached to Shimura varieties when $N = 2$ (in which case $\mathrm{GL}(2) = \mathrm{GSp}(2n)$).

I have written more on Shimura varieties than I had originally intended. I wanted simply to give some sense of the foundations laid down by Langlands in this very rich field, following Shimura himself [230] and Deligne [61]. As we have seen, these include the comprehensive beginnings with $\mathrm{GL}(2)$, the conjectural formula for the points of a Shimura variety modulo $p$, and the recognition of the central role played by the fundamental lemma and endoscopy. Nonetheless, we shall return to the subject in the next section. In the second half of the section, we shall discuss the conclusions of Kottwitz from another point of view, which includes a conjectural derivation in terms of motives. This is turn will serve as a springboard for the introduction of further questions and problems.

# 9 Motives and Reciprocity

One of Grothendieck's great insights was the idea of a motive. In general terms, motives are supposed to have two simultaneous roles. On the one hand, they are to be regarded as fundamental building blocks of algebraic varieties. On the other, they represent also a universal cohomology functor for algebraic varieties, of which Betti cohomology, (algebraic) de Rham cohomology, $\ell$-adic étale cohomology and crystalline cohomology become concrete realizations. Their existence was predicated by Grothendieck on a number of conjectures [112], [113] which are still largely unproved today. Perhaps for this reason, they were not widely discussed in formal mathematical circles before 1977. However, the unwritten understanding, if it existed, came to an end with the articles of Tate [237] and Langlands [144] in 1979 Corvallis proceedings.

We are speaking of pure motives, which are the semisimple objects in the larger category of mixed motives. They would be the objects obtained from the category $(\mathrm{SProj})_F$ of smooth, projective algebraic varieties over a field $F$, as opposed to all varieties. Our object is to describe their role in [144], the article in which Langlands proposed what is now known as the *Reciprocity Conjecture*. Roughly speaking, Reciprocity is the analogue for general (smooth, projective) varieties of the Shimura–Taniyama–Weil conjecture for elliptic curves over $\mathbb{Q}$. Together with Functoriality, it is often regarded as one of the two fundamental pillars of the Langlands program.

Consider a pair of fields $(F, Q)$ of characteristic 0, equipped with a complex embedding $F \subset \mathbb{C}$ for $F$. The category $\mathrm{Mot}_{F,Q}$ of (pure) motives over $F$ with coefficients in $Q$ (or just $Q$-motives over $F$) should be a semisimple, $Q$-linear category (which is to say, an abelian category in which short exact sequences split, and in which the abelian groups $\mathrm{Hom}(M, N)$ have been enriched to $Q$-vector spaces), equipped with a functor

$$(\mathrm{SProj})_F \to \mathrm{Mot}_{F,Q} \tag{76}$$

from the category of smooth projective varieties over $F$. It would in fact be a tensor category, equipped with a tensor product structure $\otimes$ over $Q$, satisfying several natural axioms. With the complex embedding $F \hookrightarrow \mathbb{C}$ of $F$, ordinary Betti cohomology of complex varieties with $Q$-coefficients gives a functor

$$H_B = H_B^* \colon (\mathrm{SProj})_F \to (\mathrm{Vect})_Q$$

with values in the category $(\mathrm{Vect})_Q$ of graded $Q$-vector spaces. This should factor through $\mathrm{Mot}_{F,Q}$ to a $\otimes$-functor

$$H_B = H_B^* \colon \mathrm{Mot}_{F,Q} \to (\mathrm{Vect})_Q. \tag{77}$$

It is interesting to note that the two arrows in the composition

$$(\mathrm{SProj})_F \to \mathrm{Mot}_{F,Q} \to (\mathrm{Vect})_Q \tag{78}$$

of (76) and (77), paired with the corresponding fields $F$ and $Q$, illustrate the dual roles played by motives. The first arrow describes motives as building blocks of smooth projective varieties (given the grading of $\mathrm{Mot}_{F,Q}$ by weights implicit in the second one). The second arrow for its part suggests a role for $\mathrm{Mot}_{F,Q}$ as a cohomology theory for smooth projective varieties. This would be universal, in the sense that there should be similar "realizations" for all the various arithmetic cohomology theories for $(\mathrm{SProj})_F$. (See [65].) In the case of Shimura varieties, $F$ equals the reflex field $E$.

Suppose now that $F$ and $Q$ are number fields. In fact, we might as well assume that $Q = \mathbb{Q}$, since Betti cohomology comes with a restriction of scalars functor. The functor

$$H_B \colon \mathrm{Mot}_{F,\mathbb{Q}} \to (\mathrm{Vect})_{\mathbb{Q}}$$

in (77) is called a *fibre functor* for the tensor category $\mathrm{Mot}_{F,\mathbb{Q}}$ over $\mathbb{Q}$. It gives $\mathrm{Mot}_{F,\mathbb{Q}}$ the structure of a natural (neutral) *Tannakian category*. A fundamental observation of Grothendieck was that a Tannakian category is equivalent to the category of (finite dimensional) representations over the ground field ($\mathbb{Q}$ in this case) of an affine proalgebraic group. There would then be an (anti)equivalence from $\mathrm{Mot}_{F,\mathbb{Q}}$ to the category of representations of a proalgebraic group

$$\mathscr{G}_F = \mathscr{G}_{F,\mathbb{Q}} = \mathrm{Aut}^{\otimes}(H_B)$$

over $\mathbb{Q}$. This group would be an extension

$$\mathscr{G}_{F,\mathbb{Q}} \mapsto \Gamma_F,$$

of the absolute Galois group of $F$ by a connected, reductive, proalgebraic group over $\mathbb{Q}$, whose finite dimensional representations over $\mathbb{Q}$ parametrize (up to equivalence) the objects $M$ in $\mathrm{Mot}_F = \mathrm{Mot}_{F,\mathbb{Q}}$. (See [210].)

With the Shimura–Taniyama–Weil conjecture in mind, Langlands was interested in the relations between motives and automorphic representations. Treating general automorphic representations and motives as objects over the complex numbers, he considered the complexification

$$\mathscr{G}_F = \mathscr{G}_{F,\mathbb{C}} = \mathscr{G}_{F,\mathbb{Q}} \times_{\mathbb{Q}} \mathrm{Spec}(\mathbb{C})$$

of the *motivic Galois group* above. (We are identifying the $\mathbb{Q}$-group $\mathscr{G}_F = \mathscr{G}_{F,\mathbb{Q}}$ above with its associated group $\mathscr{G}_F = \mathscr{G}_{F,\mathbb{Q}} \cong \mathscr{G}_F(\mathbb{C})$ of complex points here.) It is an extension

$$\mathscr{G}_F \to \Gamma_F$$

of $\Gamma_F$ by a *complex*, reductive proalgebraic group whose finite dimensional complex representations parametrize *complex $F$-motives*. Langlands then suggested that the best (and perhaps only) way to express the relations would be through a parallel *automorphic Galois group*.

It was an audacious proposal. As far as I know, nothing of the sort had ever been imagined before. Langlands had studied the $L$-homomorphisms

$$\phi \colon W_F \to {}^L G \qquad (79)$$

of the global Weil group $W_F$ ten years earlier as a means to parametrize some automorphic representations of $G$. But by 1977, it was well understood that these objects could not account for most automorphic representations (except in the case of a torus $G = T$ [153]). There seemed to be a general feeling that there would be nothing more in this particular direction to say.

Langlands formulated the construction of an automorphic Galois group, again as an extension

$$G_{\Pi(F)} \to \Gamma_F$$

of $\Gamma_F$ by a connected, complex reductive group. Its representations of degree $n$ would parametrize the set of automorphic representations $\Pi(F)$ of $\mathrm{GL}(n)$ over $F$ that he called *isobaric*. These are the representations of $\mathrm{GL}(n, \mathbb{A}_F)$ denoted symbolically by

$$\pi = \pi_1 \boxplus \cdots \boxplus \pi_r, \qquad \pi_i \in \Pi_{\mathrm{cusp}}(\mathrm{GL}(n_i)), \qquad (80)$$

on p. 207 of [144]. The ranks $n_i$ correspond to a partition $(n_1, \ldots, n_r)$ of $n$, and $\pi$ stands for a canonical irreducible constituent of the associated induced representation from the standard parabolic subgroup of $\mathrm{GL}(n, \mathbb{A}_F)$ attached to the partition.

Langlands' idea was to try to attach a Tannakian category to the representations (80). For a start, it was necessary to have a classification of automorphic representations of $\mathrm{GL}(n)$ in terms of the isobaric representations (80), which he formulated at the bottom of p. 207 of [144]. The classification was established soon afterwards by Jacquet and Shalika [107], using their theory with Piatetskii-Shapiro of Rankin–Selberg $L$-functions. This gave $\Pi(F)$ the structure of an abelian category. However, to obtain a tensor category from $\Pi(F)$ would require something much stronger, functoriality attached to the tensor product representations

$$\mathrm{GL}(n_i, \mathbb{C}) \times \mathrm{GL}(n_j, \mathbb{C}) \to \mathrm{GL}(n_i n_j, \mathbb{C})$$

of dual groups. This is one of the deepest cases of functoriality, and is still far from being resolved. The final ingredient would be a suitable fibre functor

$$\Pi(F) \to (\mathrm{Vect})_{\mathbb{C}}.$$

It is hard to imagine any construction of this last ingredient, but one could hope that it would be a part of the eventual proof of the cases of functoriality above. At any rate, this would essentially make the tensor category $\Pi(F)$ into a Tannakian category, of which the complex automorphic Galois group $G_{\Pi(F)}$ would then be a consequence.

The group $G_{\Pi(F)}$ would be a replacement for the Weil group, in its role (79) for parameters of automorphic representations. As we have noted, the $n$-dimensional representations of $G_{\Pi(F)}$ would parametrize all isobaric representations of $G = \mathrm{GL}(n)_F$. This would include all of the globally tempered representations $\Pi_{\mathrm{temp}}(G)$, which we recall are the irreducible representations of $G(\mathbb{A}_F)$ that occur in the decomposition of $L^2(G(F) \backslash G(\mathbb{A}_F))$. More generally, if $G$ is any connected, quasisplit group over $F$, the algebraic $L$-homomorphisms

$$\phi \colon G_{\Pi(F)} \to {}^L G$$

over $\Gamma_F$ would parametrize disjoint global packets of automorphic representations of $G$ ($L$-packets) whose union contains $\Pi_{\mathrm{temp}}(G)$. The representation theory of a general reductive group $G$ will thus be more complicated than that of $\mathrm{GL}(n)$, and will probably be best understood, through the theory of endoscopy, in terms of the theory for its quasisplit inner twist.

Langlands assumed the existence of $G_{\Pi(F)}$, and turned to the problem of relating motives to automorphic representations. His proposed solution of the problem, after the first bold step of postulating the existence of $G_{\Pi(F)}$, was elegant and simple. It was to conjecture the existence of a surjective canonical mapping

$$G_{\Pi(F)} \to \mathscr{G}_F \tag{81}$$

of complex, proalgebraic groups over $\Gamma_F$. Among other things, the mapping would be compatible with local data attached to each of the two groups, which we will discuss presently in slightly different guise. The surjectivity of the mapping (81) was not stated explicitly in [141], but it was clearly a part of Langlands' thinking, in the Tate conjecture [152] and his implicit aim of formulating a general analogue of the Shimura–Taniyama–Weil conjecture. A complex $F$-motive $M$, being identified with a finite dimensional complex representation of $\mathscr{G}_F$ pulls back to a complex finite dimensional representation $r_M$ of $G_{\Pi(F)}$. The local data in $G_{\Pi(F)}$ and $\mathscr{G}_F$ would then yield an identity

$$L(s, r_M) = L(s, M) \tag{82}$$

of $L$-functions. In particular, the motivic $L$-functions on the right would inherit the analytic continuation and functional equation from the standard automorphic $L$-functions on the left hand side. It is the conjectured mapping (81) that is known today as Langlands' *Reciprocity Conjecture*.

A number of cases of the Reciprocity Conjecture have already been established, if we are prepared to state them directly as relations between motives and automorphic representations (that is, without the universal groups in (81)). For example, any complex $n$-dimensional representation of $\Gamma_F$ is a motive, known for obvious reasons as an *Artin motive* [185, §2]. Reciprocity in this case amounts to functoriality for Artin $L$-functions. As we discussed in Section 7, it was established by Langlands for complex two-dimensional Galois representations of dihedral, tetrahedral and octahedral representations, with Galois groups isomorphic respectively to $D_{2n}$, $A_4$ and $S_4$ (the latter with Tunnell).

Shimura varieties will be another source of examples. For the Shimura varieties attached to $\mathrm{GL}(2)$, there are many complex two-dimensional motives, which should correspond to many cuspidal automorphic representations of $\mathrm{GL}(2)$. The local computations required to state Reciprocity in this case were the content of the conjecture at the end of §4 of [140]. As we discussed in Section 8, this was established later in the same article and the subsequent article of Carayol. The same phenomena occur for general Shimura varieties, although they are more subtle. For a general Shimura datum $(G, X)$, the co-character $\mu_h$ defined in Section 8 is dual to a minuscule weight $\widehat{\mu}$ of $\widehat{G}$, which in turn gives rise to the finite dimensional representation $r = r_X$ of the $L$-group $^L G$. This assigns motives to the various constituents of the $L^2$-cohomology (63) of $S_K(\mathbb{C})$ obtained from $(\mathfrak{g}_{\mathbb{R}}, K_{\mathbb{R}})$-cohomology. On the other hand, automorphic representations of $G$ could be attached to the global parameters $\phi$ in the generalization of (79) if we had the group $G_{\Pi(F)}$. The reciprocity correspondence between the two kinds of objects would then be within reach if one could establish a local reciprocity law at any place $v$ of $\mathbb{Q}$. This would be the general analogue for $(G, v)$ of the conjecture [140, §4] for $(G, v) = (\mathrm{GL}(2), p)$ (See [144, Lemma, p. 240] and [18, Proposition 9.1] for the archimedean case $v = \infty$, topics we will take up later in this section, and the discussion in [123, §8–10] reviewed in the last section for the case $v = p$.)

The most famous example of Reciprocity is of course the Shimura–Taniyama–Weil (STW) conjecture. It applies to the motive of weight 1 (corresponding to the first cohomology group $H^1$) of an arbitrary elliptic curve $E$ over $\mathbb{Q}$. The problem was to show that it corresponds to a cusp form $f$ of weight 2 over $\mathbb{Q}$ such that

$$L(s, \pi) = L(s, E),$$

where $\pi$ is the automorphic representation of $\mathrm{GL}(2)$ attached to $f$ in the $L$-function on the left, and $E$ is identified with its motive of weight 1 on the right. In other words, the local data of the two objects match, in the sense of [140, §4], discussed in the last section and above. This includes the requirement that the conductors of $f$ and $E$, independently defined nonnegative integers, also match. We recall that this last condition was the quantitative improvement [250] Weil added to the conjecture in 1967. The problem was of course much deeper than the similar question above for a two dimensional motive in a Shimura variety $S_K$ attached to $\mathrm{GL}(2)$. For in the latter case, there was all the geometric structure of $S_K$ that drove the Lefschetz trace formula, and its comparison with the Selberg trace formula.

The STW conjecture was established for semistable elliptic curves $E$ by Andrew Wiles in 1995 [253], in partial collaboration with Richard Taylor [239]. *Semistable* means that for each prime $p$, $E$ either has good reduction or multiplicative reduction. For example, it is included in the two conditions from [140, §4] under which Langlands established the local reciprocity law for $S_K$ in §7 of that paper. The authors were content to work with this restriction on $E$, since according to the results of Ribet [192], it was sufficient to establish Fermat's last theorem. Their proof of the STW conjecture relied on techniques [170] that were quite different from what we have been discussing in this article, as well as new methods developed expressly for the purpose. However there was one fundamental theorem from Langlands' work that was critical. It was the Langlands–Tunnell theorem, the reciprocity law for two dimensional representations of Galois groups isomorphic to $S_4$. With the exceptional isomorphism

$$S_4 \cong \mathrm{PGL}_2(\mathbb{F}_3),$$

the theorem was the starting point for Wiles' extended study of deformation rings and the congruence properties for modular forms. (See [184], [83], [56] for general introductions to Wiles' proof.)

Six years after the two papers of Wiles and Taylor, C. Breuil, B. Conrad, F. Diamond and R. Taylor published a proof of the STW conjecture [40] for general $E$. They built on the work of various authors to remove the ramification constraints step by step, the most difficult being various calculations associated with the prime $p = 3$. The proof of the general STW conjecture has led to the proofs of unsolvability of other Fermat-like diophantine equations. Taken on its own, it represents a milestone in number theory, the proof of a longstanding fundamental case of what is now the general Reciprocity Conjecture. (See [55].)

We return now to Langlands' proposed universal automorphic Galois group $G_{\Pi(F)}$. Shortly after the publication of [144], Kottwitz pointed out that the Langlands group would be simpler if it were taken to be in the category of locally compact groups, rather than complex proalgebraic groups [119, §12]. In this formulation, the universal group would be an extension $L_F$ of the absolute Weil group $W_F$ by a connected compact group. It would thus take its place in a sequence

$$L_F \to W_F \to \Gamma_F$$

of three locally compact groups, all having ties to the arithmetic of the global field $F$.

This represented a less severe conceptual change from what had been in place before 1979. The earlier set of global Langlands parameters, as $L$-homomorphisms of $W_F$ into the complex group ${}^L G$, would simply be enriched to the larger set of $L$-homomorphisms of $L_F$ into ${}^L G$. Local Langlands parameters would remain the same, namely $L$-homomorphisms from the locally compact group

$$L_{F_v} = \begin{cases} W_{F_v}, & \text{if } v \text{ is archimedean,} \\ W_{F_v} \times \mathrm{SU}(2), & \text{if } v \text{ is nonarchimedean,} \end{cases} \tag{83}$$

into the complex group ${}^L G_v \subset {}^L G$. (Langlands had introduced the product $W_{F_v} \times \mathrm{SL}(2,\mathbb{C})$ in [144, §4] as an equivalent version of the Weil–Deligne group [237, §4]. Kottwitz chose $\mathrm{SU}(2)$ in place of $\mathrm{SL}(2,\mathbb{C})$ in order that the bounded local parameters (with respect to their images in $\widehat{G}$) would continue to be ones that correspond to locally *tempered* representations of $G(F_v)$.) The global locally compact group $L_F$ should then be equipped with a commutative diagram

$$\begin{array}{ccccc} L_{F_v} & \longmapsto & W_{F_v} & \longrightarrow & \Gamma_{F_v} \\ \downarrow & & \downarrow & & \downarrow \\ L_F & \longmapsto & W_F & \longrightarrow & \Gamma_F \end{array} \tag{84}$$

of continuous homomorphisms for each valuation $v$ of $F$. As usual, the vertical embedding on the left in (84) would be determined only up to conjugacy, and would extend the embeddings of local Weil and Galois groups.

Kottwitz proposed a set of axioms for $L_F$, but he seems to have been thinking of a group directly related to Langlands' construction of $G_{\Pi(F)}$. Roughly speaking, $G_{\Pi(F)}$ would be regarded as the "algebraic hull" of $L_F$, a proalgebraic group over $F$ whose algebraic representations were in bijection with the continuous representations of $L_F$. In other words, if the proposed Tannakian category exists, thereby giving rise to a group $G_{\Pi(F)}$, the existence of a group $L_F$ should follow formally. The algebraic hulls of the local groups $L_{F_v}$ introduced in [144], could then revert back to the original groups (84). Finally Langlands' Reciprocity Conjecture (81) would become the existence of a continuous $L$-homomorphism

$$L_F \to \mathscr{G}_F \qquad\qquad (85)$$

Motivated by Langlands' paper, and the supplementary remarks by Kottwitz, I introduced a constructive version [20] of the group $L_F$ in 2002. It is more concrete, and it leads to a correspondingly concrete description of the motivic Galois group $\mathscr{G}_F$. It also has the advantage of not requiring a Tannakian category for its existence. The construction is still conjectural. It relies on the Principle of Functoriality, as formulated in its unramified form in Section 4. Moreover, some further conditions related to functoriality would be needed for the resulting group to have the desired properties. These were discussed somewhat tentatively in §5 of [20]. They will probably be resolved one way or another by Beyond Endoscopy, the strategy proposed by Langlands around 2000 for attacking functoriality.

In Langlands' conjectural definition of $G_{\Pi(F)}$, the basic building blocks are cuspidal automorphic representations of general linear groups. This follows the implicit definition of the motivic Galois group in terms of irreducible (complex) motives. But there are automorphic representations that could be regarded as more fundamental. These would be the cuspidal, tempered automorphic representations of quasisplit groups $G$ that are not functorial images from any smaller group. Let us be more precise.

We first recall the set

$$\mathscr{C}(G) = \mathscr{C}_{\mathrm{aut}}(G) = \{c(\pi) : \pi \in \Pi(G)\}$$

introduced in Section 4. It consists of equivalence classes of families

$$c^S(\pi) = \{c_v(\pi) = c(\pi_v) : \pi \in \Pi(G), v \notin S\}$$

of semisimple conjugacy classes in $^LG$. It is best that we then consider the subset $\mathscr{C}_{\mathrm{bdd}}(G)$ of classes $c \in \mathscr{C}(G)$ that are *bounded*, in the sense that for almost all $v$, the projection of $c_v$ onto $\widehat{G}$ meets a maximal compact subgroup of $\widehat{G}$. These would be the classes $c = c(\pi)$ attached to the globally tempered representations $\pi \in \Pi_{\mathrm{temp}}(G)$ (which we recall means that they should occur in the spectral decomposition of $L^2(G(F) \setminus G(\mathbb{A}_F))$, but with the further condition that they be locally tempered, in the sense that they satisfy the bounds from the general analogue of Ramanujan's conjecture. We are assuming functoriality. This implies that the unramified $L$-functions

$$L^S(s,c,r) = \prod_{v \notin S} L_v(s, c_v, r), \qquad c \in \mathscr{C}(G), \tag{86}$$

attached to representations $r$ of $^L G$ have analytic continuation and functional equation. The $L$-functions attached to classes $c \in \mathscr{C}_{bdd}(G)$ are the ones for which the Euler product on the right converges absolutely for $\mathrm{Re}(s) > 1$.

Suppose that $G$ is simple and *simply connected*, as well as quasisplit. In this case, we say that a class $c \in \mathscr{C}_{bdd}(G)$ is *primitive* if for any $r$,

$$\mathrm{ord}_{s=1} L^S(s, c, r) = [r \colon 1_{^L G}].$$

This amounts to asking that $c$ not be a proper functorial image $\rho'(c')$, for some $c' \in \mathscr{C}_{bdd}(G')$, and some $L$-homomorphism

$$\rho' \colon {}^L G' \to {}^L G \tag{87}$$

whose image in $^L G$ is proper. We write $\mathscr{C}_{prim}(G)$ for the set of primitive classes in $\mathscr{C}_{bdd}(G)$, for the simply connected group $G$. It is these objects, or if one prefers, corresponding automorphic representations $\pi \in \Pi_{prim}(G)$, that represent the fundamental building blocks of $L_F$. They are the smallest family in the embedded sequence

$$\mathscr{C}_{prim}(G) \subset \mathscr{C}_{cusp}(G) \cap \mathscr{C}_{bdd}(G) \subset \mathscr{C}_{bdd}(G)$$

where $\mathscr{C}_{cusp}(G)$ is the subset of classes $c = c(\pi)$ in $\mathscr{C}(G)$ for which $\pi$ is cuspidal.

The main ingredient in the construction of $L_F$ is an indexing set $\mathscr{C}_F$. It consists of isomorphism classes of pairs

$$(G, c), \qquad c \in \mathscr{C}_{prim}(G),$$

with $G$ simple and simply connected, in which $(G, c)$ is isomorphic to $(G_1, c_1)$ if there is an isomorphism of groups $G \to G_1$ over $F$, and a dual isomorphism $^L G_1 \to {}^L G$ that takes $c_1$ to $c$. Suppose that $c$ belongs to $\mathscr{C}_F$ (in the sense that it represents an isomorphism class of pairs $(G, c)$). Since the group $G$ is simply connected, the complex dual group $\widehat{G}$ is of adjoint type. We write $K_c$ for a compact real form of its simply connected cover $\widehat{G}_{sc}$. The Weil group $W_F$ acquires an action on $K_c$ from the semidirect product

$$^L G = \widehat{G} \rtimes W_F.$$

For any such $c$, there is then a natural extension

$$1 \to K_c \to L_c \to W_F \to 1 \tag{88}$$

of $W_F$ by $K_c$, which does not generally split. There are two separate constructions of this extension, for which a reader can consult [20, §4]. The second of these includes a description of localizations

$$
\begin{array}{ccc}
L_{F_v} & \longmapsto & W_{F_v} \\
\downarrow & & \downarrow \\
L_c & \longmapsto & W_F
\end{array}
\tag{89}
$$

for the group $L_c$, which is based on the local Langlands correspondence.

The extensions (88) and localizations (89) attached to elements $c \in \mathscr{C}_F$ are what is needed to construct the locally compact Langlands group $L_F$. It is defined as the fibre product

$$
L_F = \prod_{c \in \mathscr{C}_F} (L_c \to W_F)
\tag{90}
$$

over $W_F$. As such, it is an extension

$$
1 \to K_F \to L_F \to W_F \to 1
\tag{91}
$$

of $W_F$ by the compact simply connected group

$$
K_F = \prod_{c \in \mathscr{C}_F} K_c,
$$

and is hence locally compact. The required localizations (84) follow from their analogues (89) for each $L_c$.

We have appealed to functoriality in the definition of the sets $\mathscr{C}_{\mathrm{prim}}(G)$, and therefore in the indexing set $\mathscr{C}_F$ used to define $L_F$. The expectation is that $L_F$ will lead to a classification of automorphic representations. The best outcome would be that for any quasisplit $G$, the set $\Pi_{\mathrm{bdd}}(G)$ of locally tempered representations that occur in the spectral decomposition of $L^2(G(F) \setminus G(\mathbb{A}_F))$ is a disjoint union of global $L$-packets, parametrized by $\widehat{G}$-conjugacy classes of $L$-homomorphisms

$$
\phi : L_F \to {}^L G
$$

whose image in $\widehat{G}$ is bounded. This reflects what might be expected for the subset of such representations attached to parameters $\phi$ for the Weil group $W_F$. As we have said, the matter would likely be resolved with a proof of functoriality by the methods of Beyond Endoscopy. In particular, if the proposed classification above needs only minor adjustments, these ought to become clear from Beyond Endoscopy.

We should emphasize that $L_F$ represents a "thickening" of the Weil group $W_F$. The two groups should satisfy similar qualitative properties, such as those outlined for the Weil group by Tate in §1 of [237]. For example, if $E$ is a finite extension of $F$, $L_E$ would be the preimage of $\Gamma_E \subset \Gamma_F$ in $L_F$ under the composition

$$
L_F \to W_F \to \Gamma_F.
$$

In particular, this could serve as a definition of $L_F$ in terms of the group $L_{\mathbb{Q}}$.

In §6 of [20], there is also a tentative construction of the complex motivic Galois group $\mathscr{G}_F$ (which could conceivably serve at some point as an actual definition). It is modeled on the construction of $L_F$ and the version (85) of Langlands' Reciprocity homomorphism. For a quasisplit group $G$ over $F$, we could define a complex $G$-*motive* to be an $L$-homomorphism from $\mathscr{G}_F$ to $^LG$ (with the Galois form $^LG = \widehat{G} \rtimes \Gamma_F$ of the $L$-group, since $\mathscr{G}_F$ is to be regarded as a complex proalgebraic group over $\Gamma_F$). In the case that $G$ is simple and simply connected, we would also speak of *primitive* $G$-motives. They would correspond to elements $c \in \mathscr{C}_{\mathrm{prim}}(G)$ that are *algebraic* (or *motivic*). This means that if $c = c(\pi)$ for an automorphic representation $\pi$, and if

$$\phi_v \colon W_{F_v} \to {}^LG_v$$

is a Langlands parameter for $\pi_v$ at an archimedean place $v$, the composition of $\phi_v$ with any finite dimensional representation $r$ of $^LG$ whose kernel contains a subgroup of finite index in $W_F$ is of *Hodge type*. In other words, the restriction of $r \circ \phi_v$ to the subgroup $\mathbb{C}^*$ of $W_{F_v}$ is a direct sum of (quasi)characters of the form

$$z \to z^{-p} \overline{z}^{-q}, \qquad z \in \mathbb{C}^*, p, q \in \mathbb{Z}. \tag{92}$$

Our restriction on $r$ would in fact imply that each of these is of weight 0, in the sense that $p + q = 0$, and hence a character of the form

$$z \to (z/\overline{z})^q.$$

The construction of $\mathscr{G}_F$ amounts to a fibre product analogous to (90), but with two changes.

(i) The indexing set $\mathscr{C}_F$ in (90) is replaced by the subset $\mathscr{C}_{F,\mathrm{alg}}$ of algebraic indices $c \in \mathscr{C}_F$.

(ii) The diagram (88) of locally compact groups is replaced by a diagram

$$1 \to \mathscr{D}_c \to \mathscr{G}_c \to \mathscr{T}_F \to 1$$

of complex proalgebraic groups, in which $\mathscr{D}_c$ equals the simply connected complex dual $\widehat{G}_{\mathrm{sc}}$ of the group $G$ attached to $c$, and $\mathscr{T}_F$ is the Taniyama group over $\Gamma_F$ introduced by Langlands in Section 5 of [144].

The automorphic Galois group would then be given as a fibre product

$$\mathscr{G}_F = \prod_{c \in \mathscr{C}_{F,\mathrm{alg}}} (\mathscr{G}_c \to \mathscr{T}_F)$$

over $\mathscr{T}_F$. The construction would thus express $\mathscr{G}_F$ explicitly in terms of the families $\mathscr{C}_{F,\mathrm{alg}}$. Otherwise said, it builds $\mathscr{G}_F$ out of primitive $G$-motives rather than irreducible $\mathrm{GL}(n)$-motives. With this proposed construction of $\mathscr{G}_F$, the homomorphism (85) is then defined concretely on the last page 481 of [20]. It is to be regarded as an $L$-homomorphism of groups over $\Gamma_F$.

Given $\mathscr{G}_F$, and any $G$ over $F$, we write $\Phi_{\mathrm{alg}}(G)$ for the set of $G$-motives. This is a set (possibly empty) of $\widehat{G}$-conjugacy classes of proalgebraic $L$-homomorphisms

$$\Phi \colon \mathscr{G}_F \to {}^L G,$$

with respect to the projections of $\mathscr{G}_F$ and ${}^L G$ onto $\Gamma_F$. We can then identify these mappings with their restrictions to $L_F$ under Langlands' proposed reciprocity mapping $\phi \colon L_F \to \mathscr{G}_F$, since the Reciprocity Conjecture should imply that the two sets are indeed in bijection. In other words, we can also regard a parameter $\Phi \in \Phi_{\mathrm{alg}}(G)$ as $\widehat{G}$-conjugacy class of $L$-homomorphisms

$$\Phi \colon L_F \to {}^L G$$

(with respect again to the two projections onto $\Gamma_F$) that is of Hodge type. Namely, if $(r,V)$ is any finite dimensional representation of ${}^L G$, and $v$ is any archimedean place of $F$, the restriction of $\Phi_v$ to the subgroup $\mathbb{C}^*$ of $W_{F_v}$ is a direct sum of characters of the form (92). For any such $v$, we write $\Phi_{\mathrm{alg}}(G_v)$ for the obvious local analogue of this set.

The algebraic Langlands parameters $\Phi \in \Phi_{\mathrm{alg}}(G)$ generally have nonzero weights, which means that they are nontempered. However, they can be projected naturally onto tempered parameters. To see this, we first introduce the weight homomorphism

$$w \colon C_F \to L_F$$

for $L_F$. It is the mapping from $C_F = F^* \backslash \mathbb{A}_F^*$ into the preimage

$$L_F^0 = K_F W_F^0$$

in $L_F$ of the identity component $W_F^0$ of $W_F$, defined by mapping the norm $\|c\|$ for any $c \in C_F$ to the multiplicative subgroup $\mathbb{R}^+$ of $W_F^0$ described in [237, (1.4.4)]. We also have the *Tate homomorphism*

$$t \colon L_F \to C_F,$$

which can be expressed in terms of the composition of the projection $L_F \to W_F$ and $W_F \to C_F$. One sees that the image of $w$ lies in the centre of $L_F$, and that the composition $(t \circ w)$ maps any $c \in C_F$ to $\|c\|^{-2}$. For any algebraic parameter $\Phi \in \Phi_{\mathrm{alg}}(G)$, the modified Langlands parameter

$$\phi(x) = \Phi\big(x \cdot (w \circ t)(x)^{-\frac{1}{2}}\big), \qquad x \in L_F, \tag{93}$$

then lies in $\Phi_{\mathrm{temp}}(G)$. It has the property that if $v$ is an archimedean valuation of $F$ and $r_v$ is a finite dimensional representation of ${}^L G_v$, and if $(r_v \circ \Phi_v)(z)$ is a sum of algebraic characters (92), then $(r_v \circ \phi_v)(z)$ is the corresponding sum of continuous characters

$$(z/|z|^{\frac{1}{2}})^{-p}(\bar{z}/|\bar{z}|^{\frac{1}{2}})^{-q} = (z/\bar{z})^{-\frac{1}{2}(p-q)}. \tag{94}$$

We can write $\Phi_{\text{temp,alg}}(G)$ for the set of parameters $\phi \in \Phi_{\text{temp}}(G)$ obtained in this way, and

$$\Phi_{2,\text{alg}}(G) = \Phi_2(G) \cap \Phi_{\text{temp,alg}}(G)$$

for the subset of such parameters in $\Phi_2(G)$. We note that the correspondence $\phi_v \longleftrightarrow \Phi_v$ is a generalization of Langlands' dual correspondences $\phi_v \longleftrightarrow \phi_v'$ and $\pi_v \longleftrightarrow \pi_v'$ for GL(2) described in the last section.

We observe in passing that there is a parallel structure for the motivic Galois group $\mathscr{G}_F$. It was pointed out by Serre [210], who noted that $\mathscr{G}_F$ comes with a weight homomorphism

$$w \colon \mathbb{G}_m \to \mathscr{G}_F,$$

whose image lies in the centre of $\mathscr{G}_F$. It has the property that for a representation $r$ of $\mathscr{G}_F$, the *weights* of the corresponding motive are given by the decomposition of $r \circ w$ into characters of $\mathbb{G}_m$. The motivic Galois group also comes with the Tate motive

$$t \colon \mathscr{G}_F \to \mathbb{G}_m$$

of weight $(-2)$, which is to say that $t(w(x)) = x^{-2}$. These objects for $\mathscr{G}_F$ would follow immediately from their analogues for $L_F$ and the basic properties of the Taniyama group $\mathscr{T}_F$.

Langlands' statement of the Reciprocity Conjecture actually came at the beginning of the article [144]. According to his introduction in §1, the original intent was to study two specific problems in the theory of Shimura varieties. These were taken up in the final two Sections 6 and 7 of the article. The earlier sections arose from his afterthoughts on the problems, but were for reasons of exposition presented in the reverse order. Section 2 contains the proposed mapping (81), on which our discussion to this point has been based.

Section 3 of [144] contains a brief discussion of another question, automorphic representations Langlands called anomalous. Unlike what happens for GL($n$), these can include cuspidal automorphic representations that are not locally tempered. The first examples had been introduced for the group PSp(4) by Kurokawa shortly before the Corvallis conference. A second family of examples for PSp(4), discovered by Howe and Piatetskii-Shapiro, was discussed informally at the conference. These both turned out to be special cases of the representations for general groups introduced in [18], which we discussed in the last section as part of the conjectural stabilization (74) of spectral side of the Lefschetz trace formula. The relevant conjectures have now been established for quasisplit classical groups [23], [181].

In §4 of [144], Langlands returned to his original topic, the theory of Shimura varieties. This section again concerns motives, and as such, provides a bridge between the specific problems in the later sections of [144] and the formulation of the general Reciprocity Conjecture in Section 2. However, it applies to the motives as they are thought to appear on the geometric side of the Lefschetz trace formula, rather than the motives on the spectral side that would govern Reciprocity. This dual role would serve as a grand generalization of what happens for elliptic curves, which on the one hand represent moduli on the geometric side of the Lefschetz trace for-mula for $GL(2)$, and on the other, spectral objects classified by the (now resolved) STW-conjecture. An equally fundamental analogy is the dual way of representing extensions of a number field $F$ that is at the core of the automorphic trace formula for $GL(n)$. On the one hand, they come from the irreducible polynomials of degree $n$ over $F$ in the elliptic part of the geometric side, as governed by the base of the Steinberg–Hitchin fibration that we will recall in the final section, and on the other, the irreducible complex representations of $\mathrm{Gal}(\overline{F}/F)$ of degree $n$ that are conjectured to be basic components of the discrete part of the spectral side.

The conjectural formula in Section 3 of Kottwitz' paper [123], as it evolved from [142] and [164], is what allows us to understand the spectral properties of Shimura varieties $S_K$. However, its proof was restricted to cases in which $S_K(\mathbb{C})$ could be realized as a moduli space for geometric objects related to abelian varieties. Many, if not most, Shimura varieties are not PEL varieties, the basic moduli spaces of this sort. Abelian varieties are of course motives. The idea of §4 of [144], which Langlands learned from Deligne, was to treat an arbitrary Shimura variety $S_K(G, X)$ (apart from those that behave badly on the cocentre, that is with $G \neq G_0$ in the notation of [144, p. 217]), as a moduli space of *motives*. Deligne regarded Shimura varieties as parameter spaces for certain Hodge structures. His construction was predicated on a conjecture that these objects in turn parametrize uniquely determined motives. We shall review it as presented in [144] to get a further sense of the remarkable internal structure of a general Shimura variety.

As we have done with motives, we speak here of pure Hodge structures. Following the beginning of §4 of [144], we recall that a *real Hodge structure* $V$ is a finite dimensional real vector space $V_{\mathbb{R}}$ whose complexification has a decomposition

$$V_{\mathbb{C}} = \bigoplus_{p,q \in \mathbb{Z}} V^{p,q},$$

for complex subspaces $V^{p,q}$ such that $V^{q,p} = \overline{V^{p,q}}$. It is equivalent to a finite dimensional representation $\sigma$ of the group $\mathscr{R} = \mathrm{Res}_{\mathbb{C}/\mathbb{R}}(\mathbb{G}_m)$ over $\mathbb{R}$, with

$$V^{p,q} = \{v \in V_{\mathbb{C}} = V_{\mathbb{R}} \otimes \mathbb{C} : \sigma(z_1, z_2)v = z_1^{-p} z_2^{-q} v\}.$$

It is also essentially the same as an archimedean parameter

$$\phi_{\mathbb{R}} \colon W_{\mathbb{R}} \to GL(n, \mathbb{C})$$

of Hodge type, of the kind that classifies representations of $GL(n, \mathbb{R})$ of motivic significance, as we can recall from the quasicharacters (92) and the related parameters $\phi'_{\mathbb{R}}$ for $GL(2)$ in Section 8. With either interpretation, the set of real Hodge structures is a Tannakian category, with associated group $\mathscr{R}$. A real Hodge structure $V$ is of *weight $n$* if $V^{p,q} = 0$ unless $p + q = n$.

Recall also that a *rational Hodge structure* $V$ is a finite dimensional vector space $V_{\mathbb{Q}}$ over $\mathbb{Q}$ with a direct sum decomposition

$$V_{\mathbb{Q}} = \bigoplus_n V_{\mathbb{Q}}^n,$$

together with real Hodge structures of weight $n$ on the real vector spaces $V_{\mathbb{R}}^n = V_{\mathbb{Q}}^n \otimes_{\mathbb{Q}} \mathbb{R}$. The basic example is the Tate rational Hodge structure $\mathbb{Q}(1)$ of weight $-2$, in which $\mathbb{Q}(1)_{\mathbb{Q}} = 2\pi i \mathbb{Q} \subset \mathbb{C}$ and $\mathbb{Q}(1)_{\mathbb{C}} = \mathbb{Q}(1)^{-1,-1}$. Rational Hodge structures $V$ contain much more information than real Hodge structures. In particular, if real Hodge structures give representations of general linear groups $\mathrm{GL}(n, \mathbb{R})$ in terms of archimedean Langlands parameters, supplementary $\mathbb{Q}$-structures should in many cases lead to motives, and therefore the enriched structure of automorphic representations of the groups $\mathrm{GL}(n)$. A necessary condition for this, however, is that $V$ be *polarizable*, in the sense that it can be endowed with a bilinear form of the sort described on p. 215 of [144]. The situation is a generalization of the theory of abelian varieties, whereby a complex torus represents an abelian variety if and only if it has a Riemann bilinear form. The category

$$\mathrm{Hod} = \mathrm{Hod}_{\mathbb{C},\mathbb{Q}} = \{V = (V_{\mathbb{C}}, V_{\mathbb{Q}})\}$$

of polarizable rational Hodge structures is Tannakian, with fibre functor $V \to V_{\mathbb{Q}}$, and a corresponding Hodge group $\mathscr{H} = \mathscr{H}_{\mathbb{C}} = \mathscr{H}_{\mathbb{C},\mathbb{Q}}$ over $\mathbb{Q}$.

Recall finally that $\mathrm{Mot} = \mathrm{Mot}_{\mathbb{C},\mathbb{Q}}$ is the Tannakian category of $\mathbb{Q}$-motives over $\mathbb{C}$, with motivic Galois group $\mathscr{G} = \mathscr{G}_{\mathbb{C}} = \mathscr{G}_{\mathbb{C},\mathbb{Q}}$. Every object in this category comes with a polarizable rational Hodge structure, according to the properties of $\mathbb{Q}$-Betti cohomology of nonsingular complex projective varieties. There is consequently a $\otimes$-functor

$$h_{\mathrm{BH}} \colon \mathrm{Mot} \to \mathrm{Hod}$$

of categories, and a corresponding group homomorphism

$$h_{\mathrm{BH}}^* \colon \mathscr{H} \to \mathscr{G}$$

over $\mathbb{Q}$. The Hodge conjecture implies that the functor is fully faithful, which means that no data is lost from a motive (over $\mathbb{C}$) in passing to its Hodge structure.

The construction described in §4 of [144] applies to any Shimura variety $S_K = S_K(G, X)$ with $G = G_0$. It also depends on a rational representation $(\xi, V)$ of $G$. The idea is to associate to every point $x = (h, g)$ in the space

$$\mathscr{X}_K = X \times (G(\mathbb{A}_{\mathrm{fin}})/K)$$

a pair of objects $(V^x, \phi_{\mathrm{fin}}^x)$ as follows. The first component $V^x$ is the rational Hodge structure on the vector space $V_{\mathbb{Q}}^x = V_{\mathbb{Q}} = V$ given by the representation $\xi \circ h$ of $\mathscr{R}(\mathbb{R}) = \mathbb{C}^*$ on $V_{\mathbb{C}}$. The second component $\phi_{\mathrm{fin}}^x$ is the isomorphism $V_{\mathbb{A}_{\mathrm{fin}}}^x \to V_{\mathbb{A}_{\mathrm{fin}}}$ given by $v \to \xi(g)^{-1} v$, with the understanding that it be defined only up to composition by an element of $\xi(K)$. Each such pair is also implicitly fitted with an underlying family of polarizations $\mathscr{P}^x$ of $V^x$, described briefly on p. 216 of [144], with the property that if $x' = \gamma x = (\gamma h, \gamma g)$ for some $\gamma \in G(\mathbb{Q})$, there is a natural map $\gamma \colon \mathscr{P}^x \to \mathscr{P}^{x'}$.

With this machinery in place, Langlands varies the representation

$$\xi = (\xi, V) = (\xi, V(\xi)).$$

The construction then attaches to any $x$ a functor

$$\eta^x \colon (\xi, V(\xi)) \to V^x(\xi)_{\mathbb{Q}}, \qquad V(\xi) = V,$$

from the category $\mathrm{Rep}(G)$ of rational representations of $G$ to the category Hod of polarizable rational Hodge structures. It is a $\otimes$-functor (commuting with tensor products), with a matching

$$\omega_{\mathrm{Hod}} \circ \eta^x = \omega_{\mathrm{Rep}(G)}$$

of underlying fibre functors. But Hod is equivalent to the category $\mathrm{Rep}(\mathscr{H})$. The properties of Tannakian categories then provide a homomorphism $\phi^x \colon \mathscr{H} \to G$ over $\mathbb{Q}$ for which $\eta^x$ can be identified with the functor

$$\eta^x \colon (\xi, V(\xi)) \to (\xi \circ \phi^x, V(\xi)).$$

A comparison of $\phi^x$ with the given mapping $\phi^x_{\mathrm{fin}}$ then leads in [144] back to the original element $g$ in the pair $x = (h, g)$.

The conclusion reached on p. 214 on [144] is that $\mathscr{X}_K$ parametrizes pairs $(\phi, g)$, where $\phi$ is a homomorphism from $\mathscr{H}$ to $G$ over $\mathbb{Q}$, and $g$ is an element in $G(\mathbb{A}_{\mathrm{fin}})$ taken only up to right multiplication by an element in $K$. Moreover, $\phi$ is subject to the constraint that its composition $\phi \circ h$ with the canonical homomorphism $h$ from $\mathscr{R}$ to $\mathscr{H}$ lies in the set $X$ of homomorphisms from $\mathscr{R}$ to $G$ in the original Shimura datum $(G, X)$.

Langlands notes finally that this construction of Deligne comes with the hope/expectation that any homomorphism $\phi' \colon \mathscr{H} \to G$ over $\mathbb{Q}$ such that $\phi' \circ h$ lies in $X$ is a composition $\phi' = \phi \circ h^*_{\mathrm{BH}}$, for a (uniquely determined) homomorphism $\phi \colon \mathscr{G} \to G$ over $\mathbb{Q}$, which is to say, a $G$-motive over $\mathbb{C}$ with coefficients in $\mathbb{Q}$. The complex variety $S_K(\mathbb{C})$ would then parametrize equivalence classes of pairs

$$\{(\phi, g) : \phi \colon \mathscr{G} \to G, \, \phi \circ h \in X, \, g \in G(\mathbb{A}_{\mathrm{fin}})/K\},$$

where $(\phi', g')$ is equivalent to $(\phi, g)$ if

$$(\phi', g') = (\mathrm{ad}(\gamma)\phi, \gamma g)$$

for some $\gamma \in G(\mathbb{Q})$.

This completes our discussion of §4 of [144]. Langlands observed that it did not yet yield a "moduli problem in the usual sense". He also asserted that "nonetheless, there is a good deal to be learned" from the discussion we have just sketched. Indeed there is. I have yet to learn much of the subsequent history of the problem. Among other things, I am puzzled by what seems to be a paucity of later references to a construction that seems so compelling (brief as our distilled review here is). I have not yet studied the later paper [164] of Langlands and Rapoport, or the fundamental volume [66] on absolute Hodge cycles, or later papers of Milne such as [173] and [175]. I presume that this construction has had to be reformulated in terms of gerbes in order to accommodate relations among motives in different characteristic. The reduction of a moduli space modulo $p$ was of course essential to the proof of any case of the conjectural formula [142], [164], [123] for the terms on the geometric side of the Lefschetz trace formula. A recent paper of [111] of Kisin establishes the formula for Shimura varieties of abelian type. (As I understand it [176, §9], a motive is of *abelian* type if it lies in the category generated by abelian varieties; a Shimura variety $S_K = S_K(G,X)$ with rational weight $w_X$ is of *abelian type* if it is a moduli space in the sense of the construction above of abelian motives).

It was in §5 of [144] that Langlands introduced the Taniyama group $\mathscr{T}_F$. We recall that it is the replacement of the Weil group $W_F$ in the diagrams (90) and (91) for the construction we have proposed for the motivic Galois group $\mathscr{G}_F$. It is an extension

$$1 \to \mathscr{S}_F \to \mathscr{T}_F \to \Gamma_F \to 1$$

of the Galois group $\Gamma_F = \mathrm{Gal}(\overline{F}/F)$ (with $F$ embedded in $\mathbb{C}$) by the Serre group. The Serre group $\mathscr{S}$ is in turn a complex proalgebraic torus, with a continuous action of $\Gamma_{\mathbb{Q}} = \mathrm{Gal}(\overline{\mathbb{Q}}/\mathbb{Q})$. (We are writing $\mathscr{S}_F$ for the same group, but with its Galois action restricted from $\Gamma_{\mathbb{Q}}$ to the subgroup $\Gamma_F$.) It was actually defined by Langlands in §4, following Serre's construction in [208]. The Serre group should be the commutator quotient of both $\mathscr{H}$ and $\mathscr{G}$, and thereby fit into a diagram

$$\mathscr{H} \mapsto \mathscr{G} \to \mathscr{S}$$

of complex, connected groups. The (polarizable, rational) Hodge structures and (complex) motives defined by representations of $\mathscr{S}$ are then said to be of *CM-type*.

Langlands defined $\mathscr{T}_F$ by an explicit 2-cycle from $\Gamma_F$ with values in $\mathscr{S}_F$. He then used it in §6 to formulate the conjectural solution to a moduli problem he had posed at the end of §4. The problem was just one of a number of questions that would need to be resolved in order to be able to treat a general Shimura variety $S_K = S_K(G,X)$ as a moduli space, as in the special case of a PEL-variety. It concerned how the proposed parametrization of $S_K(\mathbb{C})$ by pairs $(\phi, g)$ changes under an automorphism $\tau$ of $\mathbb{C}$. Langlands observed in general terms that $(\phi, g)$ would be replaced by a pair $(\phi', g')$ attached to another Shimura datum $(G', X') = (G^{\tau, \phi}, X^{\tau, \phi})$. The problem was to describe the group $G^{\tau, \phi}$ explicitly. We will not review Langlands' conjectural resolution, obtained as part of his construction of the Taniyama group, but we recall that the conjecture was proved soon afterwards, independently by Borovoi [37] and Milne [173]. The Taniyama group itself has remained an important part of the theory. Besides its fundamental role in the motivic Galois group $\mathscr{G}_F$, and in the volume [66] on absolute Hodge cycles and elsewhere, it has also been used in a rather different context. The paper [6] introduced an extension of the Taniyama group, which in turn has become a part of a very interesting generalization [75] of the theory of motives.

The last Section 7 of [144] contains remarks of Langlands on the cohomology of general Shimura varieties $S_K = S_K(G,X)$, with a view towards the Hasse–Weil zeta function of $S_K$. These were taken up, at least implicitly, in the later paper [123] of Kottwitz discussed in the last section. I will not review them here, and in particular, their bearing on the inner twist $G^{\tau,\phi}$ of $G$ introduced by Langlands in §4, and studied further in §5 and §6 of [144].

However, I would like to draw attention to the lemma of Langlands on p. 240 of [144]. We shall use it as the starting point for a discussion related to the final formula obtained by Kottwitz in §10 of [123], which we reviewed (but did not state) at the end of the last section. We shall describe a reciprocity identity for any Shimura variety, based on the general conjectures for motives. We shall then discuss some further motivic questions suggested by the identity.

Langlands considered an archimedean parameter in the set

$$\{\phi \in \Phi(G_{\mathbb{R}}) : S_{\phi}^0 \subset Z(\widehat{G}_{\mathbb{R}})\}$$

such that the graded vector space

$$H^*(\phi,\xi) = \bigoplus_{\pi \in \Pi_{\phi}} H^*(\mathfrak{g}_{\mathbb{R}}, K_{\mathbb{R}}; \pi \otimes \xi)$$

in (63) is nonzero. He then attached two representations of the real Weil group $W_{\mathbb{R}}$ to $\phi$ with two different roles in mind, one motivic and one automorphic. The lemma asserts that the two representations are equivalent.

The "motivic" representation for $\phi$ comes from a real Hodge structure on the space $H^*(\phi,\xi)$, introduced by Langlands at the lower half of p. 239 of [144]. We note that it is to be regarded as a spectral object, in contrast to the Hodge structures from §4 of [144] we have just reviewed. As a representation of $\mathbb{C}^*$ on $H^*(\phi,\xi)$, it is defined on a subgroup of index 2 in $W_{\mathbb{R}}$. To extend it, he first considered the archimedean Weil group $W_{\mathbb{C}/\overline{E}} = W_E$, where $E$ is the Shimura field over which $S_K$ is defined, and $\overline{E}$ is its completion defined by the embedding $E \subset \mathbb{C}$ with which it is equipped. If $\overline{E} = \mathbb{R}$, Langlands observed that the action of $\mathbb{C}^*$ on $H^*(\phi,\xi)$ extends naturally to the Weil group $W_{\mathbb{R}} = W_{\mathbb{C}/\mathbb{R}}$. With this in hand, his representation of $W_{\mathbb{R}}$ can then be taken in general to be

$$H^*(\phi,\xi)^+ = \mathrm{Ind}(W_{\mathbb{R}}, W_{\overline{E}}, H^*(\phi,\xi)),$$

the representation on a space $H^*(\phi,\xi)^+$ obtained by inducing the subrepresentation of $W_{\overline{E}}$ on $H^*(\phi,\xi)$.

Langlands' other representation of $W_{\mathbb{R}}$, the one whose role would be automorphic, is given by the parameter

$$\phi_{\mathbb{R}} : W_{\mathbb{R}} \to {}^L G_{\mathbb{R}}$$

itself. We recall from the last section that $S_K = S_K(G,X)$ comes with a cocharacter $\mu = \mu_h$ for $G$, and a corresponding character $\widehat{\mu}$ that serves as the highest weight for an irreducible representation $(r, V(r))$ of $\widehat{G}$. It follows from the definition of the reflex field $E$ of $S_K$ that the Galois group $\Gamma_E$ stabilizes the $\widehat{G}$-orbit of the mi-nuscule weight $\widehat{\mu}$, and hence that $r$ extends to a representation $r_E$ of the $L$-group $^L G_E = \widehat{G} \rtimes \varGamma_E$ of $E$ such that $\varGamma_E$ acts trivially on the weight space of $\widehat{\mu}$. Langlands then introduced the representation

$$r^+ = \mathrm{Ind}(^L G, {}^L G_E, r_E) \tag{95}$$

of the $L$-group $^L G = \widehat{G} \rtimes \varGamma_{\mathbb{Q}}$ obtained by induction from the representation $r_E$ of $^L G_E$ to $^L G$. His second representation of $W_{\mathbb{R}}$ can then be defined as

$$r^+ \circ \varPhi,$$

where

$$\varPhi(w) = \phi(w)|w|^{-\frac{d}{2}}, \qquad d = \dim S_K, w \in W_{\mathbb{R}}.$$

It is not hard to see from the definitions that $\varPhi$ is algebraic, in the sense that it is the local component of a parameter in the global set $\varPhi_{\mathrm{alg}}(G)$, or equivalently, that its constituents are of the form (92). The lemma on p. 240 of [144] can be taken as the assertion that the two representations of $W_{\mathbb{R}}$ are equivalent.

There are several deeper phenomena suggested by this lemma, simple as it may be. With one representation of $W_{\mathbb{R}}$ acting on an archimedean cohomology group, and the other given explicitly in terms of an archimedean Langlands parameter, it is suggestive of the local archimedean component of the Global Reciprocity correspondence between motives and automorphic representations. We would be dealing with a specific motive here. It is represented by the de Rham cohomology (63) of $S_K(\mathbb{C})$, and can therefore simply be regarded as the motive of the Shimura variety $S_K$ over $E$. We are in fact speaking of what is known as a *realization* of the motive, specifically the Hodge realization. We can think of the lemma as a property of the local archimedean part of the mapping $\varPhi \to \phi$ from $\varPhi_{\mathrm{alg}}(G)$ to $\varPhi_{\mathrm{temp}}(G)$ in (93).

There is something else in Langlands' lemma. It suggests a broader global perspective, one that goes beyond the reflex field $E$ and the complex embedding $E \subset \mathbb{C}$. The local archimedean reflection of this phenomenon is given by the extensions $H^*(\phi)^+$ and $V(r)^+$ of the complex vector spaces $H^*(\phi)$ and $V(r)$ on which Langlands' two representations act. The global implication is that we would need to consider an extension of the motive of $S_K$, with components attached to Galois conjugates $E' \subset \mathbb{C}$ of $E \subset \mathbb{C}$. What would these motives be? Are they attached to Shimura varieties $S'_{K'}$? The conjecture of Langlands of §6 of [144], established not long afterwards [37], [173] asserted that they are.

Before I try to expand on these global implications, I should first describe the local generalization [18, Proposition 9.1] that was motivated by Langlands' lemma. This in turn relies on the local conjectures of [18, §8], which together with their global counterparts were a part of our discussion of Kottwitz' conjectural spectral (de)stabilization of the Lefschetz trace formula from the last section.

The local conjectures apply to any completion $F_v$ of a number field $F$. They concern enriched parameters

$$\psi_v \colon L_{F_v} \times \mathrm{SL}(2, \mathbb{C}) \to {}^L G_v, \qquad \psi_v \in \varPsi(G_v),$$

taken as usual up to conjugacy by $\widehat{G}$, but with the property that the image of $L_{F_v}$ projects onto a bounded subset of $\widehat{G}$. In other words, the restriction $\phi_v$ of $\psi_v$ to $L_{F_v}$ is a tempered Langlands parameter. Given $\psi_v$, one forms the centralizer

$$S_{\psi_v} = \operatorname{Cent}(\psi(L_{F_v} \times \operatorname{SL}(2, \mathbb{C})), \widehat{G})$$

in $\widehat{G}$ of its image, and the finite group

$$\mathbf{S}_{\psi_v} = S_{\psi_v}/S^0_{\psi_v} Z(\widehat{G})^{\Gamma_v},$$

often abelian, of connected components in $S_{\psi_v}$ modulo the Galois invariants in the centre of $\widehat{G}$. For each $\psi_v$, the conjectures assert the existence of a finite set $\Pi_{\psi_v}$ of representations of $G(F_v)$. This set would parametrize (in a noncanonical way) the irreducible characters

$$s_v \to \langle s_v, \pi_v \rangle, \qquad s_v \in \mathbf{S}_{\psi_v}, \, \pi_v \in \Pi_{\psi_v},$$

on the group $\mathbf{S}_{\psi_v}$. However, in contrast to the (bounded) Langlands parameters $\phi_v$, the representations $\pi_v \in \Pi_{\psi_v}$ need not be either tempered or irreducible. On the other hand they are conjectured to be unitary, and to be finite sums of irreducible representations. The local parameters $\psi_v$ and packets $\Pi_{\psi_v}$, together with their global counterparts $\psi$ and $\Pi_{\psi}$, are really a part of the theory of endoscopy. We have mentioned this term regularly in earlier sections, but we have not yet said what it is. We shall do so in the next section.

For the proposition from [18], we take $v$ to be the real valuation of $\mathbb{Q}$, and start with a parameter $\psi_{\mathbb{R}}$ in the set

$$\Psi_2(G_{\mathbb{R}}) = \{\psi_{\mathbb{R}} \in \Psi(G_{\mathbb{R}}) : S^0_{\psi_{\mathbb{R}}} \subset Z(\widehat{G}_{\mathbb{R}})\}.$$

The groups $S_{\psi_{\mathbb{R}}}$ and $\mathbf{S}_{\psi_{\mathbb{R}}}$ are then abelian. However, the situation here is slightly different from that of Langlands' lemma, given the requirement that the $\psi_{\mathbb{R}}$-image of $L_{\mathbb{R}} = W_{\mathbb{R}}$ be bounded in $\widehat{G}$. We are regarding the irreducible representation $\xi$ of $G$ as algebraic, which means that its restriction to $G_{\mathbb{R}}$ generally has nonunitary central character. This in turn forces the cohomology $H^*(\mathfrak{g}_{\mathbb{R}}, K_{\mathbb{R}}, \pi_{\mathbb{R}} \otimes \xi)$ to vanish for representations $\pi_{\mathbb{R}}$ in the packet $\Pi_{\psi_{\mathbb{R}}}$. To rectify the problem, we write

$$\xi_{\mathbb{R}}(x_{\mathbb{R}}) = \xi(x_{\mathbb{R}})|\det(x_{\mathbb{R}})|^{\alpha_{\xi}}, \qquad x_{\mathbb{R}} \in G_{\mathbb{R}} = G(\mathbb{R}),$$

where $\alpha_{\xi} \in \mathbb{R}^+$ is the number that makes the central character of $\xi_{\mathbb{R}}$ unitary. (This is related to the earlier footnote 6 and the ensuing discussion for $\operatorname{GL}(2)$. It is also implicit in the discussion on p. 61 of [20].) With this understanding, we define $\Psi_2(G_{\mathbb{R}}, \xi_{\mathbb{R}})$ to be the subset of archimedean parameters in $\Psi_2(G_{\mathbb{R}})$ such that the graded vector space

$$H^*(\psi_{\mathbb{R}}, \xi_{\mathbb{R}}) = \bigoplus_{\pi_{\mathbb{R}} \in \Pi_{\psi_{\mathbb{R}}}} H^*(\mathfrak{g}_{\mathbb{R}}, K_{\mathbb{R}}; \pi_{\mathbb{R}} \otimes \xi_{\mathbb{R}})$$

is nonzero. The representations $\pi_{\mathbb{R}} \in \Pi_{\psi_{\mathbb{R}}}$ are interesting examples of the unitary representations with cohomology classified by Vogan and Zuckerman [244].

Consider then a parameter $\psi_{\mathbb{R}}$ in the set $\Psi_2(G_{\mathbb{R}}, \xi_{\mathbb{R}})$. The representations $\pi_{\mathbb{R}} \in \Pi_{\psi}$ then give rise to a real Hodge structure. This relies on the analysis of such parameters by Adams and Johnson [2], [18, §5]. It in turn gives a representation of $\mathbb{C}^*$ on $H^*(\psi_{\mathbb{R}}, \xi_{\mathbb{R}})$, but with components of the form (94) rather than (92), for our having replaced $\xi$ by $\xi_{\mathbb{R}}$. The group $S_{\psi_{\mathbb{R}}}$ also acts on $H^*(\psi_{\mathbb{R}}, \xi_{\mathbb{R}})$, and commutes with the action of $\mathbb{C}^*$. The extra ingredient in the parameter is the group $\mathrm{SL}(2, \mathbb{C})$. It also acts on the space $H^*(\psi_{\mathbb{R}}, \xi_{\mathbb{R}})$, where it governs the grading in the manner familiar from the hard Lefschetz theorem. Since it commutes with the action of the product of $S_{\psi_{\mathbb{R}}} \times \mathbb{C}^*$, we obtain a representation of the product $S_{\psi_{\mathbb{R}}} \times \mathbb{C}^* \times \mathrm{SL}(2, \mathbb{C})$ on the space $H^*(\psi_{\mathbb{R}}, \xi_{\mathbb{R}})$. Following Langlands, we can then construct an induced representation

$$\rho_{\psi_{\mathbb{R}}}^+(s, (w, u)), \qquad (s, w, u) \in S_{\psi_{\mathbb{R}}} \times (W_{\mathbb{R}} \times \mathrm{SL}(2, \mathbb{C})),$$

of the product $S_{\psi_{\mathbb{R}}} \times (W_{\phi_{\mathbb{R}}} \times \mathrm{SL}(2, \mathbb{C}))$ on a graded vector space $H^*(\psi_{\mathbb{R}}, \xi_{\mathbb{R}})^+$ that contains $H^*(\psi_{\mathbb{R}}, \xi_{\mathbb{R}})$ and that is the analogue of the space $H^*(\phi, \xi)^+$ introduced above. This is the "motivic" representation for $\psi_{\mathbb{R}}$.

The "automorphic" representation for $\psi_{\mathbb{R}}$ is constructed as above, but in terms of the new parameter $\psi_{\mathbb{R}}$. It equals

$$\sigma_{\psi_{\mathbb{R}}}^+(s, (w, u)) = r_{\mathbb{R}}^+(\psi_{\mathbb{R}}(w, u)s), \quad (s, w, u) \in S_{\psi_{\mathbb{R}}} \times (W_{\mathbb{R}} \times \mathrm{SL}(2, \mathbb{C})), \qquad (96)$$

where $r_{\mathbb{R}}^+$ is the representation of $^L G_{\mathbb{R}}$ attached as above to the Shimura datum $(G, X)$. Proposition 9.1 of [18] amounts to the assertion

$$\rho_{\psi_{\mathbb{R}}}^+ \cong \sigma_{\psi_{\mathbb{R}}}^+ \qquad (97)$$

that the two representations of $S_{\psi_{\mathbb{R}}} \times (W_{\mathbb{R}} \times \mathrm{SL}(2, \mathbb{C}))$ are equivalent. It reduces to the lemma of Langlands in the case that the parameter $\psi_{\mathbb{R}} = \phi_{\mathbb{R}}$ lies in the subset $\Phi_2(G_{\mathbb{R}})$ of $\Psi_2(G_{\mathbb{R}})$, which is to say that its restriction to $\mathrm{SL}(2, \mathbb{C})$ is trivial. (The proposition was actually formulated and proved for the smaller representations with $\mathbb{C}^*$ in place of $W_{\mathbb{R}}$, but its extension follows from the various definitions.)

We can now turn to the global implications of these local results. We shall introduce two global representations, one motivic and one automorphic, which for the moment serve simply to help us focus our thoughts. They depend on the global conjectures from [18, §8], about we can first say a few words.

If $G$ is any reductive group over a number field $F$, we write $\Psi(G)$ for the set of $L$-homomorphisms

$$\psi \colon L_F \times \mathrm{SL}(2, \mathbb{C}) \to {}^L G$$

such that $L_F$ has bounded image in $\widehat{G}$, taken up to $\widehat{G}$-conjugacy. The domain here now includes the hypothetical global Langlands group $L_F$, in which the local groups $L_{F_v}$ are imbedded. Any $\psi$ would therefore give local parameters $\psi_v \in \Psi(G_v)$, local packets $\Pi_{\psi_v}$ of representations, and a global packet $\Pi_{\psi}$ of representations

$$\pi = \widetilde{\bigotimes_v} \pi_v, \qquad \pi_v \in \Pi_{\psi_v},$$

of $G(\mathbb{A})$, in which $\pi_v$ is required to be unramified in a certain sense for almost all $v$. The natural mappings $s_v \to s$ from the local groups $\mathbf{S}_{\psi_v}$ to the corresponding global group $\mathbf{S}_\psi = S_\psi / S_\psi^0 Z(\widehat{G})^\Gamma$ will then attach a global pairing

$$s \to \langle s, \pi \rangle = \prod_v \langle s_v, \pi_v \rangle, \qquad s \in \mathbf{S}_\psi, \pi \in \Pi_\psi,$$

on $\mathbf{S}_\psi$ to any representation in the global packet. The factors on the right will equal 1 for almost all $v$, while the product of the noncanonical local pairings $\langle s_v, \pi_v \rangle$ will become canonical. The result would be a canonical finite dimensional character $\langle \cdot, \pi \rangle$ on $\mathbf{S}_\psi$ for every representation $\pi$ in the global packet $\Pi_\psi$. Finally, suppose that $\psi$ lies in the subset $\Psi_2(G)$ of global parameters such that the connected centralizer $S_\psi^0$ is contained in $Z(\widehat{G})$. In general, there is a natural one-dimensional sign character $\varepsilon_\psi$ on $\mathbf{S}_\psi$, constructed in a simple way from symplectic root numbers attached to $\psi$. The main conjecture is then that the automorphic discrete spectrum of $G$ (taken modulo the centre) is a direct sum over $\psi \in \Psi_2(G)$ of representations $\pi \in \Pi_\psi$, taken with multiplicities equal to the multiplicities

$$m_2(\pi) = |\mathbf{S}_\psi|^{-1} \sum_{s \in \mathbf{S}_\psi} \varepsilon_\psi(s) \langle s, \pi \rangle, \qquad \pi \in \Pi_\psi, \tag{98}$$

of $\varepsilon_\psi$ in $\langle \cdot, \pi \rangle$. (See [18, §8].)

We return now to the Shimura varieties $S_K = S_K(G, X)$, with $F$ equal to either $\mathbb{Q}$ or the reflex field $E$. Our main global representation will be the one that is "automorphic". It is built in a natural way from the archimedean representation (96) introduced above and its nonarchimedean complement from the corresponding expression (63) for the $L^2$-cohomology of $S_K(\mathbb{C})$. It is the representation

$$\bigoplus_{\psi \in \Psi_2(G, \xi)} \bigoplus_{\pi \in \Pi_\psi} (\sigma_\psi^+ \otimes \pi_{\mathrm{fin}}^K)_{\varepsilon_\psi} \tag{99}$$

of $(L_\mathbb{Q} \times \mathrm{SL}(2, \mathbb{C})) \times \mathscr{H}_K(G)$, whose terms we describe as follows. The representation itself is a direct sum over global parameters

$$\psi \colon L_\mathbb{Q} \times \mathrm{SL}(2, \mathbb{C}) \to {}^L G$$

in the subset $\Psi_2(G, \xi)$ of parameters in $\Psi_2(G)$ that restrict to archimedean parameters in the subset $\Psi_2(G_\mathbb{R}, \xi_\mathbb{R})$ of $\Psi_2(G_\mathbb{R})$, and representations $\pi \in G(\mathbb{A})$ in the corresponding global packet $\Pi_\psi$. For any $\psi$ and $\pi$, $(\sigma_\psi^+ \otimes \pi_{\mathrm{fin}}^K)_{\varepsilon_\psi}$ is then the representation of $(L_\mathbb{Q} \times \mathrm{SL}(2, \mathbb{C})) \times \mathscr{H}_K$ given by the multiplicity of the sign character $\varepsilon_\psi$ on $\mathbf{S}_\psi$ in the representation

$$(\sigma_\psi^+ \otimes \pi_{\mathrm{fin}}^K)((w, u), s, f) = r^+(\psi(w, u)s) \otimes \langle s, \pi_{\mathrm{fin}}^K \rangle(\pi_{\mathrm{fin}}^K(f))$$

of $(L_{\mathbb{Q}} \times \mathrm{SL}(2, \mathbb{C})) \times \mathbf{S}_{\psi} \times \mathscr{H}_K$. In other words,

$$(\sigma_{\psi}^+ \otimes \pi_{\mathrm{fin}}^K)((w,u),f) = |\mathbf{S}_{\psi}|^{-1} \sum_{s \in \mathbf{S}_{\psi}} \varepsilon_{\psi}(s)\big(r^+(\psi(w,u)s)\big) \otimes \langle s, \pi_{\mathrm{fin}}^K\rangle\big(\pi_{\mathrm{fin}}^K(f)\big). \tag{100}$$

The representation (99) is the centre of our discussion. The essential point is that it should be equivalent to the natural representation of $(L_{\mathbb{Q}} \times \mathrm{SL}(2, \mathbb{C})) \times \mathscr{H}_K$ on our expression

$$H_{(2)}^*(S_K(\mathbb{C}), \mathscr{F})^+ = \bigoplus_{\pi} m_2(\pi)\big(H^*(\mathfrak{g}_{\mathbb{R}}, K_{\mathbb{R}}, \pi_{\mathbb{R}} \times \xi_{\mathbb{R}})^+ \otimes \pi_{\mathrm{fin}}^K\big) \tag{101}$$

for the extended $L^2$-cohomology. Indeed (101) is defined in the same way as (99), but with the archimedean representation $\sigma_{\psi_{\mathbb{R}}}^+$ of (96) replaced by the representation $\rho_{\psi_{\mathbb{R}}}^+$ on $(\mathfrak{g}_{\mathbb{R}}, K_{\mathbb{R}})$-cohomology. The equivalence of the two representations (99) and (101) would then follow from the equivalence (97) of $\rho_{\psi_{\mathbb{R}}}^+$ and $\sigma_{\psi_{\mathbb{R}}}^+$, the form (98) for the multiplicity $m_2(\pi)$ in terms of the parameter $\psi$ and the sign character, and the various definitions.

The analogue of (99) for the actual $L^2$-cohomology $H_{(2)}^*(S_K(\mathbb{C}), \mathscr{F})$ is a similar representation, but with the group $L_E$ in place of $L_{\mathbb{Q}}$. It amounts to the formula stated as [18, (9.3)] (with the local representation $\rho_{\psi_{\mathbb{R}}}$ in place of $\sigma_{\psi_{\mathbb{R}}}$), where it follows from the definition of $\rho_{\psi_{\mathbb{R}}}$ and the actual assertion of Proposition 9.1 of [18]. We can think of (99) as the representation of $L_{\mathbb{Q}}$ induced from this representation of $L_E$. This relation is in keeping with the bijections

$$L_{\mathbb{Q}}/L_E \cong W_{\mathbb{Q}}/W_E \cong \Gamma_{\mathbb{Q}}/\Gamma_E \cong \mathrm{Hom}(E, \mathbb{C})$$

that we expect of the Langlands group, and that are stated for the Weil and Galois groups on the first page of Tate's article [237].

**Remarks.** 1. The representation (99) is an interesting expression on several counts. It provides some insight into automorphic representations (acting on spaces of automorphic forms) rather than the characters in terms of which they were classified in [23]. By displaying the global parameters and their corresponding families of $\mathscr{H}_K$-modules $\{\pi_{\mathrm{fin}}^K\}$ together as tensor products, it bears a philosophical resemblance to the theta correspondence. The automorphic modules $\pi^K$ are lacking the archimedean components $\pi_{\mathbb{R}}$, but these are hidden in the representations $\sigma_{\psi}^+$, or rather the $(\mathfrak{g}_{\mathbb{R}}, K_{\mathbb{R}})$-cohomology within the equivalent representations $\rho_{\psi}^+$.

2. The finite group $\mathbf{S}_{\psi}$ is abelian, since it is a subgroup of the archimedean group $\mathbf{S}_{\psi_{\mathbb{R}}}$, which in the case of the algebraic parameters $\psi$ is itself abelian [2]. However, the other groups $\mathbf{S}_{\psi_v}$ need not be abelian, so $\langle s, \pi_{\mathrm{fin}}^K\rangle$ can be a higher dimensional character, which means that $\sigma_{\psi}^+ \otimes \pi_{\mathrm{fin}}^K$ is not strictly a representation of the component group $\{s\} = \mathscr{S}_{\psi}$. It perhaps ought to be replaced in the notation by a higher dimensional $\mathbf{S}_{\psi}$-module. But the formula (99) makes sense as stated, and in any case, the groups $\mathbf{S}_{\psi_v}$ are typically abelian, making $\langle s, \pi_{\mathrm{fin}}^K\rangle$ a one-dimensional character.

3. The sign character $\varepsilon_\psi$ is an interesting arithmetic object in its own right. That it occurs in the basic automorphic expression (99) for Shimura varieties is perhaps surprising.

4. It is the automorphic expression (99) that is closely related to the formula of Kottwitz displayed in the last paragraph of [123, §10], and that was our object of discussion at the end of the last section. His derivation of the formula by the comparison of trace formulas, even though it rests on the conjectural fixed point formula [123, (3.1)] and the conjectures of [18, §8], of course brings us closer to an actual proof of the formula than our derivation of (99) by Langlands' Functoriality and Reciprocity. Our goal for (99) has been conceptual.

5. The isomorphism of $H^*_{(2)}(S_K(\mathbb{C}), \mathscr{F})^+$ with the representation (99) could be regarded as a global counterpart of a conjectural formula of Kottwitz [190] for the representations of local groups in the cohomology of local Shimura varieties.

The automorphic representation (99) can be regarded as the primary object of this discussion, from which the others follow, and to which subsequent questions ultimately return. The representation of the product $(L_\mathbb{Q} \times \mathrm{SL}(2, \mathbb{C})) \times \mathscr{H}_K$ on the extended $L^2$-cohomology (101), while being equivalent to (99), is really secondary. For as we have noted, it is essentially a realization of a motive. The true motivic companion of (99) would be a more direct analogue. We define it formally as the algebraic representation

$$\bigoplus_\Psi \bigoplus_{\Pi \in \Pi_\Psi} (\Sigma^+_\Psi \otimes \Pi^K_{\mathrm{fin}})_{\varepsilon_\psi} \tag{102}$$

of $(\mathscr{G}_\mathbb{Q} \times \mathrm{SL}(2, \mathbb{C})) \times \mathscr{H}_K$ whose pullback to $(L_\mathbb{Q} \times \mathrm{SL}(2, \mathbb{C})) \times \mathscr{H}_K$ under Reciprocity equals (99). In particular, $\Psi$ is the preimage of the parameter $\psi \in \Psi_2(G, \xi)$ from (99), under the analogue of the mapping (93). The automorphic representations in (102), which we have written as $\Pi \in \Pi_\Psi$ in place of $\pi \in \Pi_\psi$ are an interesting part of the construction. They are representations of $\mathscr{H}_K$, as an algebra of cycles in the extended cohomology space that commutes with the motivic representation of $\mathscr{G}_\mathbb{Q} \times \mathrm{SL}(2, \mathbb{C})$, and which in turn acts as a kind of diagonalization of this algebra.

The motivic representation (102) does not display much of its internal structure. This is because we have been treating it as a representation of the complex group $\mathscr{G}_\mathbb{Q} = \mathscr{G}_\mathbb{Q}(\mathbb{C})$, the group of complex points of the underlying group $\mathscr{G}_{\mathbb{Q},\mathbb{Q}}$ over $\mathbb{Q}$. (We adopted the overlapping notation earlier to emphasize the parallel conjectural structure of the universal groups $L_F$ and $\mathscr{G}_F$. Applied to $\mathscr{G}_F$ alone, the ambiguity is harmless; we can regard $\mathscr{G}_F$ either as the group of complex points of a scheme over $\mathbb{Q}$, or more traditionally, as a complex group with underlying structure as a group defined over $\mathbb{Q}$. In the case here of $F = \mathbb{Q}$, we shall often write $\mathscr{G} = \mathscr{G}_\mathbb{Q}$ simply for the reductive group over $\mathbb{Q}$.) It is this $\mathbb{Q}$-structure that governs the arithmetic properties of (101).

Among other things, the $\mathbb{Q}$-structure is needed to complete the Hodge realization of the Shimura variety $S_K$ on the extended cohomology space $H^*_{(2)}(S_K(\mathbb{C}), \mathscr{F})^+$. We saw in the construction of the motivic representation (101) (using (97) and Langlands' earlier lemma) how to define the real Hodge structure on the space. To extend it to a rational Hodge structure, we would need to use the fact that the representation (101) of $\mathscr{G}_{\mathbb{Q}} \times \mathrm{SL}(2,\mathbb{C})$ as a fibre functor, can be defined over $\mathbb{Q}$. The $\mathbb{Q}$-Hodge structure is required in turn to be polarizable [144, p. 215]. We would want to be able to attach an explicit polarization to the motivic parameters. This should bear a simple relation to the Lefschetz structure given by the $\mathrm{SL}(2,\mathbb{C})$-component of the parameters.

The most widely studied realization functor for motives is defined by their étale cohomology and its corresponding compatible families of $\ell$-adic Galois representations. Known as the $\mathbb{A}_{\mathrm{fin}}$-realization [65], it is of obvious arithmetic importance. For the Shimura variety $S_K$, it was the main topic of our last section. As we recall, the Galois representations act on the $\ell$-adic (étale) version of the intersection cohomology $IH^*(\bar{S}_K, \mathscr{F}_\lambda)^+$ for the Baily–Borel compactification of $\bar{S}_K$. Like the Hodge realization, it depends very much on the $\mathbb{Q}$-structure of the group $\mathscr{G}_{\mathbb{Q}}$. As a matter of fact, we cannot really speak of the $\mathbb{A}_{\mathrm{fin}}$-realization of $S_K$, and the individual $\ell$-adic representations in particular, without this $\mathbb{Q}$-structure. For it is considerably more subtle than the $\mathbb{Q}$-structure on the Shimura group $G$. Without it, one has to work with compatible families of $\lambda$-adic representations, where $\lambda$ ranges over the non-archimedean completions of a finite extension $L$ of $\mathbb{Q}$ that depends on the group $K$. This was the point of view of Langlands in [140] and Kottwitz in [123].

We display the representation spaces we have discussed in the diagram

$$
\begin{array}{ccc}
H^*_{(2)}(S_K, \mathscr{F})^+ & \longleftarrow & \bigoplus_{\psi} \bigoplus_{\pi} (\sigma^+_\psi \otimes \pi^K_{\mathrm{fin}})_{\varepsilon_\psi} \qquad (99) \\
\big\updownarrow{\scriptstyle\sim} & & \big\downarrow \\
IH^*(\bar{S}_K, \mathscr{F})^+ & \longleftarrow & \bigoplus_{\psi} \bigoplus_{\Pi} (\Sigma^+_\psi \otimes \Pi^K_{\mathrm{fin}})_{\varepsilon_\psi} \qquad (102)
\end{array}
$$

The diagram is somewhat impressionistic, but it is helpful for us in keeping track of what is highly conjectural, and what is better understood and more concrete. The vertical arrow on the right is in the former category. Indeed, its domain in the upper right hand corner is given in terms of the hypothetical automorphic Galois group $L_F$, whose existence is closely related to Langlands' Principle of Functoriality. The ultimate proof of this would be the goal of Beyond Endoscopy, the recent long term program of Langlands we will discuss in §11. The codomain in the lower right hand corner depends on the hypothetical motivic Galois group $\mathscr{G}_F$, while the arrow itself is given by Langlands' Reciprocity Conjecture. There is no concrete program for the proof of this, but it will surely demand everything we can prove about Shimura varieties. The vertical arrow on the left is the isomorphism of Zucker's conjecture, which we recall has been known for twenty years. The lower horizontal arrow is given by the $\mathbb{A}_{\mathrm{fin}}$-realization of $S_K^+$, while its composition with the isomorphism on the left is the Hodge realization of $S_K^+$. The upper horizontal arrow is built out of Langlands' lemma and the proposition from [18] with which we began this discussion. It seems clear from these remarks that the automorphic representation in the upper right hand corner can indeed be regarded as the foundation of the other spaces and arrows in the diagram.

We shall conclude this discussion with a list of five problems. These represent refinements of conjectures that would enhance our understanding, as opposed to ideas that might be applied to their eventual proofs. Some of them are accessible, requiring perhaps only a little careful thought. In fact, some of these may in fact already be known. But taken together, they present a broader picture that can only serve to help us.

**Problems: 1. Realizations of $S_K$.** We are thinking again of the Hodge and $\mathbb{A}_{\text{fin}}$-realizations of the (extended) Shimura motive (102). The problem would include a more explicit description of the fibre functor, as a $\mathbb{Q}$-refinement of the complex representation (102) of $\mathscr{G}_{\mathbb{Q}} \times \mathrm{SL}(2, \mathbb{C})$ on either $H^*_{(2)}(S_K, \mathscr{F})^+$ or $IH^*(\overline{S}_K, \mathscr{F})^+$. This can be considered as a special case of the corresponding problem for the full motivic Galois group $\mathscr{G}_{\mathbb{Q}}$, which we will state as Problem 3 below. However, there are some supplementary (and simpler) questions we could ask about the case of a Shimura variety $S_K$ here.

For example, (102) is the direct sum over a finite set of parameters $\Psi \in \Psi_{\text{alg}}(G)$ (the analogue of $\Phi_{\text{alg}}(G)$ for $\Phi(G)$) of complex representations of $\mathbb{G}_{\mathbb{Q}} \times \mathrm{SL}(2, \mathbb{C})$. What is the partition of this set that gives the decomposition of the sum, as a representation over $\mathbb{Q}$? Can one answer this question without a full understanding of the $\mathbb{Q}$-structure on $\mathscr{G}_{\mathbb{Q}}$? What supplementary information might be required on the complex Hecke algebra

$$\mathscr{H}_K = \mathscr{H}(K \backslash G(\mathbb{A}_{\text{fin}})/K)$$

as a rational algebra of correspondences of cycles on $S_K$, and on its representations on the relevant space of automorphic forms $\mathscr{A}_\xi(G(\mathbb{Q}) \backslash G(\mathbb{A})/K)$. Furthermore, for each $\Psi$, the representation $\Sigma_\Psi^+$ in (102) is defined as in (99) in terms as an $L$-homomorphism of $\mathscr{G}_{\mathbb{Q}} \times \mathrm{SL}(2, \mathbb{C})$ into the complex $L$-group $^L G$. Do we need to impose a $\mathbb{Q}$-structure on $^L G$ to be able to ask these questions on the $\mathbb{Q}$-fibre functor? I have not thought about the problem, even to the point of being confident in posing the questions.

**2. Hasse–Weil zeta function of $S_K$.** This is of course a famous longstanding problem. It was posed in this context by Langlands [146, 159, 152], but I am not quite sure of its present status. The question here is of a conjectural formula, since the expressions (99) and (102) on which it could be based are hypothetical. The same is true of the similar formula at the end of [123, §10], even though its conjectural foundations are much less severe. The problem appears to be quite accessible, amounting no doubt to a careful collection of the relevant terms in either of the formulas above, but it would be well worth any time taken to fully understand it. For the special case of Picard modular surfaces, a family of Shimura surfaces attached to various forms of the unitary group in three variables, the answer is known, and has been fully proved. We refer the reader to the volume on the subject edited by Langlands and Ramakrishnan, and their summary [163] from the volume of its main result.

We note that there are really two zeta functions. One is the zeta function of $S_K$ as a variety over $E$. The other is attached to the disconnected variety over $\mathbb{Q}$ represented by the extended cohomology spaces on the left hand side of the diagram above. It would be a product of zeta functions taken over its components. Each factor would be the zeta function of a separate Shimura variety, obtained from $S_K$ by the inner twist proposed by Langlands in Sections 4–6 of [144].

I will try to return to the representations in (99) and (102) elsewhere, with supplementary details and possible extensions. The goal would be to give a precise (conjectural) formula for these zeta functions.

**3. $\mathbb{Q}$-structure on $\mathscr{G}_F$.** We have given an explicit conjectural description of the complex motivic Galois group $\mathscr{G}_F$ over a number field $F$. It is a fibre product

$$\prod_c (\mathscr{G}_c \to \mathscr{T}_F),$$

over a set $\mathscr{C}_{F,\mathrm{alg}}$ of equivalence classes of pairs $(G, c)$, of extensions

$$\mathscr{G}_c \to \mathscr{T}_F \to \varGamma_F$$

of complex simply connected groups $\mathscr{G}_c$. The problem would be to give an explicit conjectural description of the $\mathbb{Q}$-structure on $\mathscr{G}_F$ attached to the complex embedding $F \subset \mathbb{C}$ and the corresponding Betti fibre functor. (We should not forget that the two fields $F$ and $\mathbb{Q}$ here are the two fields $F$ and $Q$ in (78), and have quite different sources.) Langlands' Taniyama group $\mathscr{T}_F$ is the set of complex points of a proalgebraic group that is already defined over $\mathbb{Q}$. The problem is to extend this $\mathbb{Q}$-structure to the fibre product of groups $\mathscr{G}_c$. Might the solution be given in terms of the equivalence classes of concrete families $c = \{c_v\}$ of semisimple conjugacy classes in $^L G$ that make up the indexing set $\mathscr{C}_{F,\mathrm{alg}}$?

There would be two steps. We have been writing $\mathscr{G}_F$ for the group over $\mathbb{Q}$ whose structure we seek. Let us write $\mathscr{G}_F^{\mathrm{spl}}$ for the same group of complex points, but with the structure of a (disconnected) split group over $\mathbb{Q}$. Intermediate between $\mathscr{G}_F$ and $\mathscr{G}_F^{\mathrm{spl}}$ would be a quasisplit group $\mathscr{G}_F^* = \mathscr{G}_F^{\mathrm{qs}}$. The first step would be to describe $\mathscr{G}_F^*$ explicitly by an outer twist of the Galois action on $\mathscr{G}_F^{\mathrm{spl}}$. The second step would be to describe $\mathscr{G}_F$ as an inner twist of $\mathscr{G}_F^*$.

How would we approach the first step? The co-ordinates of the conjugacy classes $c = \{c_v\}$ ought to be algebraic numbers. This is known in many cases, where it can be established from the "finite form" of the trace formula, with the test function $f$ being cuspidal at an Archimedean place. It would follow in general from functoriality, which we have already taken as a prerequisite for this section. The Galois group $\varGamma_{\mathbb{Q}}$ would then act by permutation of the families $c$, and hence on the indices $\mathscr{C}_{F,\mathrm{alg}}$ in the fibre product that defines $\mathscr{G}_F$. It is tempting to think of using this to construct a quasisplit outer form of $\mathscr{G}_F^{\mathrm{spl}}$ over $\mathbb{Q}$.

However, important as that phenomenon may be, it is not what we want! It should lead to an outer from attached to the étale/$\ell$-adic/$\mathbb{A}_{\text{fin}}$ realization for the motivic Galois group!

Our concern here is the Betti realization. It would be attached to the fibre functor that assigns to a motive defined over $F$ its Betti cohomology with $\mathbb{Q}$-coefficients. This should give a quasisplit outer from over $F$ of the original group $G$. The problem is that it has also to depend on the embedding of $F$ into $\mathbb{C}$ in order to become a group over $\mathbb{Q}$. This leads to the second step, the description of the inner twist $\mathscr{G}_F$ of $\mathscr{G}_F^*$. The problem is to describe the associated nonabelian cohomology class in $H^1(\mathbb{Q}, G_{\text{ad}}^*) \cong H^1(\mathbb{Q}, \widehat{G})$ explicitly.

I regret not having had the time to think about the question (as well as many other things!), because I suspect that the answer is both simple and interesting. To exploit it, we would note that motivic Galois groups should behave like Weil groups, in the sense that as a complex group, $\mathscr{G}_F$ would be the preimage of the subgroup $\Gamma_F \subset \Gamma_Q$ in the projection $\mathscr{G}_Q \to \Gamma_Q$ for any field $Q \subset F$, and hence that

$$\mathscr{G}_Q/\mathscr{G}_F \cong \Gamma_Q/\Gamma_F \cong \text{Hom}_Q(F, \mathbb{C}).$$

Taking $Q$ to be the rational field $\mathbb{Q}$, we could then write $\mathscr{G}_{\mathbb{Q}}^*$ as a disjoint union of groups $\mathscr{G}_F^*$, taken over the embeddings that parameterize the different Betti fibre functors. They define an inner twist of $\mathscr{G}_{\mathbb{Q}}^*$ that depends only on $F$. The inner form $\mathscr{G}_{\mathbb{Q}}$ of $\mathscr{G}_{\mathbb{Q}}^*$ that we seek would then presumably be the direct limit over increasing fields $F$ of the inner twists defined in this way for each $F$. Note that as a reductive proalgebraic group with $\mathbb{Q}$-structure (over the group $\Gamma_{\mathbb{Q}}$), $\mathscr{G}_{\mathbb{Q}}$ is completely canonical. To finish the second step, we could simply take $\mathscr{G}_F$ to be the preimage of $\Gamma_F \subset \Gamma_{\mathbb{Q}}$ in $\mathscr{G}_{\mathbb{Q}} \to \Gamma_{\mathbb{Q}}$ attached to an embedding $F \subset \mathbb{C}$.

What we have just described is related to the conjecture of Langlands stated in §6 of [144] (and proved in [37, 173]). The conjecture applies to a Shimura variety over the reflex field $E \subset \mathbb{C}$. It attaches different Shimura varieties to different complex embeddings of $E$, each obtained from the original one by an *explicit* inner twist. What we called the extended cohomology $H_{(2)}^*(S_K, \mathscr{F})^+$ (with locally constant sheaf $\mathscr{F}$) for the resulting motive over $\mathbb{Q}$ then becomes a disjoint union of the motives of the different Shimura varieties, or rather the Hodge realizations of these motives. It would be very interesting to compare this conjecture, and its solution, with the second step above. I hope to return to some of these questions in a future paper.

**4. Realizations for $\mathscr{G}_F$.** The Hodge and $\mathbb{A}_{\text{fin}}$ realizations of a Shimura variety are a fundamental part of its theory. A full solution of Problem 3 would give us a different way to view the realizations of any motive, each based on some further structure on the group $\mathscr{G}_F$.

The $\mathbb{A}_{\text{fin}}$-realization of a motive is a compatible family of $\ell$-adic representations

$$\bigotimes_{\ell \neq p} r_\ell, \qquad \ell \notin S,$$

of $\Gamma_\mathbb{Q}$. The prime $p$ represents a $\mathbb{Q}$-rational conjugacy class $c_p$, which embed diagonally in the $\ell$-adic vector spaces. This formulation presupposes that as a representation of $\mathscr{G}_F$, the motive is defined over $\mathbb{Q}$. But without the $\mathbb{Q}$-structure in hand, one has to work implicitly with the groups $\mathscr{G}_F^{\mathrm{spl}}$. As we have noted, this is what necessitates taking $\lambda$-adic representations for the completions of a finite extension $L$ of $\mathbb{Q}$. But we are now supposing that we have the $\mathbb{Q}$-structure on $\mathscr{G}_F$ attached to a Betti fibre functor. The $\mathbb{A}_{\mathrm{fin}}$-realization then becomes more fundamental. For it would amount to a compatible family of $\ell$-adic homomorphisms from $\Gamma_\mathbb{Q}$ to $\mathscr{G}_F$ over $\mathbb{Q}$.

Similar comments would also apply to the Hodge realization. I have not thought precisely about how best to express them, but it would clearly be interesting to formulate the Hodge realization as further structure on the group $\mathscr{G}_F$. It is closely related to the *period realization*, which we will discuss in a moment. We should add that for any index $(G, c)$ in $\mathscr{C}_{F,\mathrm{alg}}$, the ramified local complements

$$\{\Phi_v \in \Phi(G_v) : v \in S\}$$

of the family $c = \{c_v : v \notin S\}$ ought to be uniquely determined by $c$ itself. In particular, $c$ would give us the archimedean parameter $\Phi_\infty$ on which the Hodge structure depends. This presumably follows from the theorem of strong multiplicity one for $\mathrm{GL}_N$, and the fact that $c$ is primitive.

There are other realizations for motives. One would like to understand them all in terms of $\mathscr{G}_F$. We shall say a few more words about one of them, the *period realization* (which I believe is the same as what is often called the *De Rham–Betti realization*). It is yet another extraordinary side of Grothendieck's vision for motives. It suggests a systematic approach to classical transcendental number theory. Even more remarkable is that it represents an extension of algebraic number theory to many of the classical transcendental numbers that have been with us as definite integrals or infinite series since the advent of calculus. The role of the classical Galois group $\Gamma_F$ is then played by its extension given by the totally disconnected group $\mathscr{G}_F(\mathbb{Q})$. The full theory has also to include mixed motives, which we will discuss very briefly as Problem 5, but the ideas are perhaps easier to sketch in terms of pure motives. (See [116] [7]).

The basic idea comes from the familiar De Rham theorem, which asserts that for a manifold $X$, the pairing

$$H_{\mathrm{DR}}^k(X, \mathbb{C}) \times H_k(X, \mathbb{C}) \to \mathbb{C}, \qquad (\phi, c) \to \int_c \phi,$$

between De Rham cohomology and complex Betti homology is nonsingular, and therefore gives an isomorphism from $H_{\mathrm{DR}}^k(X, \mathbb{C})$ to complex Betti cohomology $H_{\mathrm{B}}^k(X, \mathbb{C}) = H_k(X, \mathbb{C})^*$. Suppose now that $X$ is a nonsingular, projective algebraic variety over $\mathbb{Q}$. The Betti cohomology of $X(\mathbb{C})$ can of course have $\mathbb{Q}$-coefficients, and so becomes a graded vector space $H_{\mathrm{B}}(X)$ over $\mathbb{Q}$. A deep theorem of Grothendieck asserts that the same is true of De Rham cohomology. Namely, there is a rational graded vector space $H_{\mathrm{DR}}(X)$ over $\mathbb{Q}$ whose complexification equals $H_{\mathrm{DR}}(X, \mathbb{C})$, together with a canonical isomorphism

$$\varpi_X \colon H_{\mathrm{DR}}(X) \otimes \mathbb{C} \xrightarrow{\sim} H_{\mathrm{B}}(X) \otimes \mathbb{C}.$$

The isomorphism is that of the original De Rham theorem. It assigns a complex number

$$\langle \varpi_X(\phi), c \rangle = \int_c \phi$$

to every rational differential form $\phi$ and every rational singular cycle $c$ of a given degree on $X$. These numbers are called *periods* of $X$.

The main point is that this construction extends to the Tannakian category $\mathrm{Mot}_{\mathbb{Q}}$ of motives over $\mathbb{Q}$. The *period realization* of $\mathrm{Mot}_{\mathbb{Q}}$ is the $\otimes$-functor

$$M \to (H_{\mathrm{DR}}(M), H_{\mathrm{B}}(M), \varpi_M)$$

from $\mathrm{Mot}_{\mathbb{Q}}$ to the Tannakian category of triples

$$(V, W, \varpi), \qquad V, W \in \mathrm{Vect}_{\mathbb{Q}},$$

where $\varpi$ is an isomorphism between the complex vector spaces $V_{\mathbb{C}}$ and $W_{\mathbb{C}}$. *Grothendieck's period conjecture* represents an analogue for the period realization of the Hodge conjecture for the Hodge realization, or the Tate conjecture for the $\mathbb{A}_{\mathrm{fin}}$-realization, fundamental foundations we have not been able to discuss. It implies that the period realization is fully faithful. The actual conjecture stated in [7, 4.1.1] applies to the period torsor

$$P_{\mathrm{mot}}(M) = \mathrm{Isom}^{\otimes}(H_{\mathrm{DR}|\langle M \rangle}, H_{\mathrm{B}|\langle M \rangle}),$$

where $\langle M \rangle$ is the Tannakian subcategory of $\mathrm{Mot}_{\mathbb{Q}}$ generated by a motive $M$, and $H_{\bullet|\langle M \rangle}$ stands for the restriction of the realization $H_{\bullet}$ to $\langle M \rangle$. The period torsor is the (noncanonical) subvariety over $\mathbb{Q}$ given by a finite dimensional affine general linear group over $\mathbb{Q}$. It is a torsor under the motivic Galois group of $M$, the finite dimensional quotient $\mathscr{G}(M) = \mathscr{G}_{\mathbb{Q}}(M)$ of $\mathscr{G} = \mathscr{G}_{\mathbb{Q}}$ attached to the subcategory $\langle M \rangle$ of $\mathrm{Mot}_{\mathbb{Q}}$. Grothendieck's conjecture asserts that the canonical complex point

$$\varpi_M \in P_{\mathrm{mot}}(M, \mathbb{C})$$

is the generic point in $P_{\mathrm{mot}}(M)$. It amounts to the assertion that the smallest algebraic subvariety of $P_{\mathrm{mot}}(\mathbb{C})$ defined over $\mathbb{Q}$ and containing $\varpi_M$ is $P_{\mathrm{mot}}$ itself. It is in this form that the conjecture suggests applications to transcendental number theory (see [7, §4]).

Taking the natural extension to $M$ of the definition for $M = X$, we define the *periods* of $M$ to be the entries in the matrix of $\varpi_M$ with respect to bases of the $\mathbb{Q}$-vector spaces $H_{\mathrm{DR}}(M)$ and $H_{\mathrm{B}}(M)$. (To be canonical, we allow the bases to vary, or equivalently, we take the periods to be the $\mathbb{Q}$-vector space $\mathscr{P}(M)$ generated by the periods with respect to any fixed pair of bases.) Grothendieck's period conjecture implies that any polynomial relations with rational coefficients among the periods of $M$ are among the relations that define the variety $P_{\mathrm{mot}}(M)$. It follows easily that the algebra over $\mathbb{Q}$ generated by the periods coincides with the algebra $\mathbb{Q}[P_{\mathrm{mot}}(M)]$ over $\mathbb{Q}$. This implies in turn that $\mathbb{Q}[P_{\mathrm{mot}}(M)]$ coincides with the $\mathbb{Q}$-algebra $\mathscr{P}(\langle M \rangle)$ obtained by taking the periods of all motives in the category $\langle M \rangle$. On the other hand, the group $\mathscr{G}(M, \mathbb{Q})$ of rational points in $\mathscr{G}(M)$ acts simply transitively on the rational points $P_{\mathrm{mot}}(M, \mathbb{Q})$ in $P_{\mathrm{mot}}(M)$, and hence on the $\mathbb{Q}$-algebra

$$\mathscr{P}(\langle M \rangle) = \mathbb{Q}[P_{\mathrm{mot}}(M)].$$

Taking limits over $M$, we see finally that the group

$$\mathscr{G}(\mathbb{Q}) = \varprojlim_{M} \big( \mathscr{G}(M, \mathbb{Q}) \big)$$

of $\mathbb{Q}$-rational points in the motivic Galois group $\mathscr{G}$ acts canonically on the $\mathbb{Q}$-algebra

$$\mathscr{P} = \varinjlim_{M} \big( \mathscr{P}(\langle M \rangle) \big)$$

of motivic periods over $\mathbb{Q}$. This is clearly a generalization of Galois theory for $\overline{\mathbb{Q}}/\mathbb{Q}$, with $\mathscr{G}(\mathbb{Q})$ being an extension of the Galois group $\varGamma_{\mathbb{Q}}$ and $\mathscr{P}$ a $\mathbb{Q}$-algebra that contains the algebraic closure $\overline{\mathbb{Q}}$ of $\mathbb{Q}$ in $\mathbb{C}$. (See [7, §5.1] for further analogies with classical Galois theory.)

I have followed the short introduction [7] in saying a few words on the Galois theory of periods. I have not stated a specific problem. Let us simply ask the same question about the period realization that we posed above for the $\mathbb{A}_{\mathrm{fin}}$ and Hodge realizations. Namely, can we formulate the theory above strictly in terms of supplementary internal structure on the motivic Galois group $\mathscr{G}$ over $\mathbb{Q}$? Even if this makes sense, it would not seem to have any immediate application. But in adding to the underlying structure of $\mathscr{G}$, it would clearly give us a broader understanding.

**5. Mixed motives.** Our discussion to this point has applied only to pure motives. Grothendieck's original vision was for a broader theory of mixed motives. (See [209, p. 345].) They would be attached to varieties over $F$ that need not be either projective or nonsingular. (The case of open Shimura varieties is actually an anomaly, since $L^2$-cohomology and intersection cohomology take it back into the domain of pure motives.) The theory of mixed motives was subsequently developed by Deligne, initially through his extensive theory of mixed Hodge structures [59], [60], [64], and more recently, through other means such as those in [66] and [65]. It remains a major area of activity, encompassing many deep and fundamental concepts.

One of Grothendieck's basic tenets was the existence of a broader group, the mixed motivic Galois group. Over the number field $F$, it would be a semidirect product

$$\mathscr{G}_F^+ = \mathscr{N}_F \rtimes \mathscr{G}_F$$

of the (pure) motivic Galois group $\mathscr{G}_F$ with a proalgebraic unipotent radical $\mathscr{N}_F$. Its existence was again predicated on the theory of Tannakian categories. In particular, with a suitable fibre functor, $\mathscr{G}_F^+$ would again become a proalgebraic group over $\mathbb{Q}$. However, Grothendieck's axioms for mixed motives are deep and difficult. They generalize his standard conjectures for pure motives, which are still far from proved. My impression is that much current work in the area is to find other means to characterize mixed motives and the group $\mathscr{G}_F^+$, which are more concrete and perhaps less difficult to establish.

The problem we pose here would be to find a concrete conjectural description of $\mathscr{G}_F^+$, comparable to what we considered for pure motives in Problem 3. This would be harder than the other problems, and might seem unrealistic to some. But if we were to go ahead, there would be two possible ways to proceed. One would be to try to extend Langlands' Reciprocity Conjecture. This is the approach of Harder [84], who has studied automorphic analogues of mixed motives in terms of the (nonunitary) values of Eisenstein series. The other would be to combine a solution of Problem 3 for the group $\mathscr{G}_F$ with a description of the unipotent radical $\mathscr{N}_F$ in elementary terms. There is an explicit solution of this problem for the category of mixed Tate motives, which yields the simplest interesting mixed motivic Galois group, and is attached to the (pure) Tate motive $\mathbb{Q}(1)$ [144, p. 214]. The solution was remarkably simple, if also quite difficult to prove [65], [41]. It is perhaps a good omen.

Everything we have discussed for pure motives should extend to the theory of mixed motives. In particular, the conjectural category of mixed motives over $\mathbb{Q}$, say, would have a period resolution that adds greatly to the set of periods, a list that would then include algebraic numbers, the periods of elliptic curves over $\mathbb{Q}$ (these sets both being pure motives), the number $\pi$, values of the logarithm at rational numbers $q \notin \{-1, 0, 1\}$, special values of the gamma function, special values of the hypergeometric function, and perhaps most striking of all, the unknown values

$$\{\zeta(2n+1) : n \in \mathbb{N}\}$$

of the Riemann zeta function that have been a preoccupation of mathematicians since the time of Euler. (See [7, §5.2–5.7] and also [116] for more examples.)

There is one number that is conspicuously absent from the list. The exponential base $e$ is in fact not a period. But it is an *exponential period*, a larger (countable) class of transcendental numbers attached to what are known as *exponential motives* [75].

**Coda: Particle physics.** There is a third conjectural universal group, in addition to the automorphic and (mixed) motivic Galois groups. This was proposed by P. Cartier, who called it the *cosmic Galois group*. It would be a quotient $\mathscr{C}^+ = \mathscr{C}_{\mathbb{Q}}^+$ of the mixed motivic Galois group $\mathscr{G}^+ = \mathscr{G}_{\mathbb{Q}}^+$. The corresponding group $\mathscr{C}_{\mathbb{Q}}^+(\mathbb{Q})$ of rational points would act like a Galois group on the $\mathbb{Q}$-vector space of periods of Feynman integrals, sums of which form the amplitudes attached to Feynman graphs [42], [54]. It is apparently unknown what this quotient should be, even as $\mathscr{C}^+$ might well turn out to be the full (mixed) motivic Galois group $\mathscr{G}^+$. This group suggests a fundamental relationship between the arithmetic Langlands program and basic particle physics, of the kind perhaps that is sometimes dreamt of. (See [188, p. 503], for example.) I am hardly a disinterested observer, and my knowledge of physics is fragmentary at the very best, but I would argue as follows.

Feynman integrals have long been a foundation for the theory of fundamental particles. In principle, they ought to give the quantum probabilities for the output data, measured from collision experiments with given input data. However, the calculations have traditionally been purely numerical, and of great difficulty. The infinite sums that go into a Feynman amplitude were originally thought to provide a convergent series. However, according to my very limited understanding, they were shown by Dyson around 1950 not to converge, but rather to give only an asymptotic formula, except in the idealized case of free particles, with input Lagrangian having only kinetic (quadratic) terms. As an approximation of this asymptotic formula, the first few terms of the infinite series, taken at points close to the origin, still give astonishing good results in the case of QED (quantum electrodynamics). However, they fail in more complex experiments. It is a fundamental problem in physics to discover a more sophisticated theory for describing quantum amplitudes in general, but which would still reduce to Feynman amplitudes in simple situations.[11]

It was shortly before the year 2000 that the physicist D. Kreimer discovered the number $\zeta(3)$ among the more complex calculations of QED. He was later joined by A. Connes, and as I understand it, they soon found that many other such calculations also gave periods of mixed motives. Moreover, the Galois action of $\mathscr{G}^+(\mathbb{Q})$ on periods, or rather its restriction to the unipotent radical $\mathscr{N}(\mathbb{Q})$, seemed to be closely related to the conceptually difficult (at least for mathematicians) physical process of renormalization.

The mixed motivic Galois group $\mathscr{G}^+$ is at the heart of much of modern arithmetic geometry. However, it has generally been regarded as inaccessible. Langlands' Reciprocity Conjecture makes it much more concrete. Combined with suitable conjectural solutions for Problems 3, 4 and 5, it would impose a rich internal automorphic structure on both $\mathscr{G}^+$ and its associated Galois group $\mathscr{G}^+(\mathbb{Q})$ for periods. Put simply, Functoriality and Reciprocity would give us the automorphic Galois group $L_F$, together with its close ties to the motivic Galois group $\mathscr{G}_F$. They are the centre of the Langlands program. On the other hand, Feynman diagrams have long been central to theoretical particle physics. It is expected that there will be something more fundamental that could eventually take their place. Whatever this might turn out to be, it is also reasonable to believe that the Langlands program would be a part of it.

This completes the second of our two sections on arithmetic geometry. Some of it is clearly speculative. However, I hope that the mathematical side of it at least will hold in principle, and that any inaccuracies will require only minor adjustments. In general, I will be happy if my attempts to describe some of the broader ideas behind Langlands' work and their subsequent development are some compensation for any misstatements that might also be present.

---

[11] I thank Marco Gualtieri for illuminating conversations. Any misinterpretations are entirely my doing.

# 10  The theory of endoscopy

There were a number of natural questions arising from his ideas that Langlands thought deeply about in the decade of the 1970s. For example, the conjectural correspondence

$$\pi' = \bigotimes_v \pi'_v \to \bigotimes_v \pi_v = \pi$$

of functoriality (Questions 4 and 5 of [138]) was just that, a correspondence. Could it be reformulated as a well defined mapping? Compared to the explicit results for GL(2) in [103], the representation theory of the group SL(2) has more structure. What was the explanation? Also, with his more recently gained experience in the $\lambda$-adic representations of Shimura varieties, Langlands found some unexpected anomalies in the associated Hasse–Weil zeta functions [143]. Again, what was the explanation? And finally, in Harish-Chandra's classification of the discrete series representations for a real group $G(\mathbb{R})$, a monumental achievement that was ahead of its time, there were some unusual aspects of his formulas for their characters. Could they be related to Local Functoriality? Langlands confronted this last problem in the work that led to his classification [151], and in his later work [165] and [166] with Shelstad.

The questions all turned out to be related. The underlying phenomena eventually became part of Langlands' conjectural theory of endoscopy. We have mentioned endoscopy a number of times already, most notably in Kottwitz' conjectural stabilization of the Lefschetz trace formula in Section 8. In this section, we shall try to give a more systematic description of the theory, and of some of the progress that has come in its development.

Given his success with the trace formula for GL(2) (as described in the three applications from our Sections 6, 7 and 8), Langlands would of course have considered how these methods might be applied to other groups, and to more general cases of functoriality. There was no clear strategy on how to proceed. But he appears to have acquired a strong sense that the trace formula would ultimately lead to a solution, informed no doubt by his general theory of Eisenstein series, and perhaps also by a skepticism as to whether other possible approaches would have the power to treat the general case.

One might try to think about comparing trace formulas for two groups $G'$ and $G$ related by the $L$-homomorphism $\rho' \colon {}^L G' \to {}^L G$ of functoriality. The immediate question would be to relate the basic elliptic terms on the geometric side of each trace formula. The conjugacy classes $\gamma \in \Gamma_{\mathrm{ell,reg}}(G)$ that index these terms for $G$ have coordinates defined by the algebra of $G$-invariant polynomials on $G$. One could think of using $\rho'$ to relate these coordinates for $G$ and $G'$. However, a serious problem arises immediately. The coordinates parametrize only geometric conjugacy classes, while for most groups $G$ other than GL($n$), there can be distinct (elliptic, regular) conjugacy classes in $G(F)$ over a ground field $F \subset \mathbb{C}$ that represent the same conjugacy class in $G(\mathbb{C})$. The theory of endoscopy begins with this problem.

It brings to bear on it some sophisticated new techniques that originate with (abelian) class field theory.

Consider the example of SL(2), with $F = \mathbb{R}$. The regular elliptic elements

$$\left\{ \gamma = \begin{pmatrix} \cos(\theta) & -\sin(\theta) \\ \sin(\theta) & \cos(\theta) \end{pmatrix}, \ \delta = \begin{pmatrix} \cos(\theta) & \sin(\theta) \\ -\sin(\theta) & \cos(\theta) \end{pmatrix} \right\}, \quad \theta \in (0, \pi),$$

are not conjugate in SL$(2, \mathbb{R})$. However, the matrix $g = \begin{pmatrix} i & 0 \\ 0 & -i \end{pmatrix}$ in SL$(2, \mathbb{C})$ has the property that

$$g \gamma g^{-1} = \begin{pmatrix} i & 0 \\ 0 & i^{-1} \end{pmatrix} \begin{pmatrix} \cos(\theta) & -\sin(\theta) \\ \sin(\theta) & \cos(\theta) \end{pmatrix} \begin{pmatrix} i^{-1} & 0 \\ 0 & i \end{pmatrix} = \begin{pmatrix} \cos(\theta) & \sin(\theta) \\ -\sin(\theta) & \cos(\theta) \end{pmatrix} = \delta,$$

so the two elements are conjugate in SL$(2, \mathbb{C})$. The regular elliptic conjugacy classes $\gamma \in \Gamma_{\mathrm{reg,ell}}(\mathrm{SL}(2))$ for SL$(2, \mathbb{R})$ thus occur naturally in pairs $(\gamma, \delta)$ that map to the same conjugacy class in SL$(2, \mathbb{C})$. The dual spectral property concerns the representations $\pi \in \Pi_2(\mathrm{SL}(2))$ in the discrete series for SL$(2, \mathbb{R})$. They too occur naturally in pairs $(\pi_n^+, \pi_n^-)$, which are parametrized by the positive integers $n$. Harish-Chandra's theory of infinite dimensional characters, which we will discuss presently, shows that the two phenomena are indeed dual in a precise sense. The characters of any pair $(\pi_n^+, \pi_n^-)$, as locally integrable functions on regular conjugacy classes, differ only on the pairs $(\gamma, \delta)$, and for these, only in a simple manner.

The group $G = \mathrm{SL}(2)$ is quite special. For this case, the dual properties have formulations in terms of the real group $G'(\mathbb{R})$ with $G' = \mathrm{GL}(2)$, as well as for the complex group $G(\mathbb{C})$. Since the element $g = \begin{pmatrix} i & 0 \\ 0 & -i \end{pmatrix}$ in SL$(2, \mathbb{C})$ is the product of the central element $\begin{pmatrix} i & 0 \\ 0 & i \end{pmatrix}$ in GL$(2, \mathbb{C})$ with the matrix $g_1 = \begin{pmatrix} 1 & 0 \\ 0 & -1 \end{pmatrix}$ in GL$(2, \mathbb{R})$, $g_1 \gamma g_1^{-1}$ equals $\delta$, and the elements $\gamma$ and $\delta$ are also conjugate in GL$(2, \mathbb{R})$. With this interpretation, the dual spectral property can be regarded as a very special case of Local Functoriality. It applies to $G' = \mathrm{GL}(2)$ and $G = \mathrm{SL}(2)$, with the homomorphism

$$\rho' \colon \widehat{G}' = \mathrm{GL}(2, \mathbb{C}) \to \widehat{G} = \mathrm{PGL}(2, \mathbb{C})$$

being the natural projection. In the local classification for GL$(2, \mathbb{R})$, the representations $\pi_n' \in \Pi_2(G')$ in the relative discrete series for GL$(2, \mathbb{R})$ (with respect to a fixed central character) are parametrized by irreducible 2-dimensional representations $\phi'$ of $W_{\mathbb{R}}$. These in turn are bijective with positive integers $\{n\}$. The corresponding pairs of representations $\{\pi_n^{\pm}\} \subset \Pi_2(G)$ for SL$(2, \mathbb{R})$ are attached to composite homomorphisms

$$\phi = \rho' \circ \phi' \colon W_{\mathbb{R}} \to \mathrm{PGL}(2, \mathbb{C}), \qquad \phi' \in \Phi_2(G').$$

They consist simply of the irreducible constituents of the restriction of $\pi_n'$ to the subgroup SL$(2, \mathbb{R})$ of GL$(2, \mathbb{R})$. One can obviously think of the pair $\{\pi_n^+, \pi_n^-\}$ attached in this way to $\phi$ as a torsor under the group

$$\mathbf{S}_\phi = \mathrm{Cent}(\phi(W_\mathbb{R}), \widehat{G}) = \mathbb{Z}/2\mathbb{Z}.$$

The sets $\pi_\phi = \{\pi_n^+, \pi_n^-\}$ and $\pi_{\phi'} = \{\pi_n'\}$ are called local $L$-packets for $G$ and $G'$.

The description of $\mathbf{S}_\phi$ as a torsor is an obvious tautology in the case of $G = \mathrm{SL}(2)$, but its generalization to arbitrary groups became part of the local Langlands classification. There were still interesting questions for $\mathrm{SL}(2)$, and particularly, for groups related to $\mathrm{SL}(2)$ and its inner twists over local and global fields $F$. The paper [127] of Langlands with Labesse contains a comprehensive study of them.

To see how Langlands' ideas progress in general, suppose that $G$ is a (connected) reductive group over a local or global field of characteristic 0. Quasisplit groups again play a special role in the theory, which for questions of transfer entail an underlying, fixed inner twist

$$\psi\colon G \to G^*,$$

of $G$ to a quasisplit group $G^*$ over $F$. One also has to work with classes $\gamma \in \Gamma_{\mathrm{reg}}(G)$ that are *strongly* regular, in the sense that the centralizer $G_\gamma$ in $G$ of (any representative of) $\gamma$ is a maximal torus $T$. (Regular elements satisfy the weaker property that the identity component $G_\gamma^0$ is a maximal torus.) We may as well simplify our notation slightly by agreeing to have the subscript *reg* mean strongly regular rather than regular. The problem then is to understand the set $\Gamma_{\mathrm{reg}}(G)$ of strongly regular conjugacy classes of $G(F)$ in a given stable conjugacy class. (A strongly regular *stable* conjugacy class is by definition the intersection of $G(F)$ with a strongly regular conjugacy class in the group $G(\overline{F})$ of points over an algebraic closure of $F$.)

Suppose that $\delta \in G(F)$ is strongly regular, with centralizer the maximal torus $T \subset G$ over $F$, and that $\gamma \in G(F)$ is another element in the stable class of $\delta$. Then $\gamma$ equals $g^{-1}\delta g$, for some element $g \in G(\overline{F})$. If $\sigma$ lies in the Galois group $\Gamma_F = \mathrm{Gal}(\overline{F}/F)$, we see that

$$\delta = \sigma(\delta) = \sigma(g\gamma g^{-1}) = \sigma(g)\gamma\sigma(g)^{-1} = t(\sigma)^{-1}\delta t(\sigma),$$

where $t(\sigma)$ is the 1-cocycle $g\sigma(g)^{-1}$ from $\Gamma_F$ to $T(\overline{F})$. It is easy to check that a second element $\gamma_1 \in G(F)$ in the stable class of $\delta$ is $G(F)$-conjugate to $\gamma$ if and only if the corresponding 1-cocycle $t_1(\sigma)$ has the same image as $t(\sigma)$ in the Galois cohomology group

$$H^1(F, T) = H^1(\Gamma_F, T),$$

which is to say that $t_1(\sigma)t(\sigma)^{-1}$ is of the form $t'\sigma(t')^{-1}$, for some element $t' \in T(\overline{F})$. Conversely, a class in $H^1(F, T)$ comes from an element $\gamma_1$ if and only if it is represented by a 1-cocycle of the form $g\sigma(g)^{-1}$. The mapping $\gamma \to \{t(\sigma)\}$ therefore defines a bijection from the set of $G(F)$-conjugacy classes in the stable conjugacy class of $\delta$ to the kernel

$$\mathscr{D}(T) = \mathscr{D}(T/F) = \ker(H^1(F, T) \to H^1(F, G)).$$

The codomain $H^1(F,G)$ is only a set with distinguished element 1, since $G$ is generally not abelian. The preimage $\mathscr{D}(T)$ of this element in $H^1(F,T)$ therefore need not be a subgroup. However, $\mathscr{D}(T)$ is contained in the subgroup

$$\mathscr{E}(T) = \mathscr{E}(T/F) = \operatorname{Im}(H^1(F,T_{\mathrm{sc}}) \to H^1(F,T))$$

of $H^1(F,T)$, where $T_{\mathrm{sc}}$ is the preimage of $T$ in the simply connected cover $G_{\mathrm{sc}}$ of the derived group of $G$. This is because the canonical map $\mathscr{D}(T_{\mathrm{sc}}) \to \mathscr{D}(T)$ is surjective. If $H^1(F,G_{\mathrm{sc}}) = \{1\}$, which is the case whenever $F$ is a nonarchimedean local field [232, §3.2], $\mathscr{D}(T)$ actually equals the subgroup $\mathscr{E}(T)$. This is why Langlands worked with the groups $\mathscr{E}(T)$ in place of $\mathscr{D}(T)$, and why the simply connected group $G_{\mathrm{sc}}$ plays an important role in the theory. Langlands introduced these ideas in the initial pages of his foundational article [147], although he had discussed them widely in the years preceding it. We are following some of the discussion in [22, §27].

We have referred to the trace formula regularly in this report, especially in the last three sections. However, to understand the refinements that originate with Langlands' observations above, we need to say something more formal. Continuing with the reductive group $G$, we now take its field of definition $F$ to be global. The (invariant) trace formula for $G$ is a general identity

$$I_{\mathrm{geom}}(f) = I_{\mathrm{spec}}(f), \qquad f \in C_c^\infty(G(\mathbb{A})), \tag{103}$$

obtained by integrating the geometric and spectral expansions (39) and (40) of the kernel $K(x,y)$ over $x = y$ in $G(F) \backslash G(\mathbb{A})^1$. As we have noted, this cannot be taken literally, since neither integral converges in general. Making sense of it is a long process, but roughly speaking, one truncates the two expansions of $K(x,x)$ in a consistent way so that the integrals converge. One then observes that as functions of the of variable of truncation $T$, a vector in some translate of a positive cone $\mathfrak{a}_0^+$, the integrals are polynomials in $T$. One can then set $T$ equal to the polynomial variable at $T_0 \in \mathfrak{a}_0$, a canonical point that depends on a maximal compact subgroup $K_0 \subset G(\mathbb{A})$ and a minimal parabolic subgroup $P_0 \subset G$, both of which are part of the truncation process. The result is a natural identity

$$J_{\mathrm{geom}}(f) = J_{\mathrm{spec}}(f), \qquad f \in C_c^\infty(G(\mathbb{A})), \tag{104}$$

which is independent of the choice of $P_0$.

However, (104) is only an intermediate step. We recall that a distribution $J$ on $G(\mathbb{A})$ is said to be *invariant* if

$$J(f^y) = J(f), \qquad f \in C_c^\infty(G(\mathbb{A})), y \in G(\mathbb{A}),$$

where

$$f^y(x) = f(yxy^{-1}), \qquad x \in G(\mathbb{A}).$$

The point here is that the linear forms on each side of (104) are noninvariant. One has then to "renormalize" the identity. There is no need to describe this explicitly, although I have been told that it is in the same spirit as a similar (but more complex) operation in quantum field theory that restores the symmetry that was lost in the truncation of divergent intervals. (I may however have misunderstood this. Renormalization seems actually to be a physics analogue of the original truncation process.) In any case, it leads to the identity (103), in which each side is now an invariant distribution. Moreover, the choice of the point $T_0$ makes each side of (103) independent of $K_0$ as well as $P_0$.

We should emphasize that neither (103) nor (104) is just an abstract formula. As in the case of GL(2), each side of (103) represents a rather complex expansion into explicit invariant linear forms, one geometric and one spectral, which can all be decomposed explicitly into their local constituents ([14], [15]). For example, on the geometric side, we have the strongly regular elliptic part

$$I_{\text{ell,reg}}(f) = \sum_{\gamma \in \Gamma_{\text{ell,reg}}(G)} \text{vol}(\gamma)\text{Orb}(\gamma, f), \tag{105}$$

where

$$\text{Orb}(\gamma, f) = \int_{G_\gamma(\mathbb{A}) \backslash G(\mathbb{A})} f(x^{-1}\gamma x)\, dx$$

and

$$\text{vol}(\gamma) = \text{vol}(Z_+ G_\gamma(F) \backslash G_\gamma(\mathbb{A})).$$

Its analogue on the spectral side would be the trace

$$I_2(f) = \sum_{\pi \in \Pi_2(G)} \text{mult}(\pi)\Theta(\pi, f) \tag{106}$$

where

$$\Theta(\pi, f) = \text{tr}(\pi(f)) = \text{tr}\left(\int_{G(\mathbb{A})} f(x)\pi(x)\, dx\right),$$

and $\text{mult}(\pi)$ is the multiplicity with which $\pi$ occurs discretely in $L^2(Z_+ G(F) \backslash G(\mathbb{A}))$. We are writing $Z_+$ here for a fixed connected, central subgroup of $G(\mathbb{A})$ that is a complement to the subgroup $G(\mathbb{A})^1$, defined as in §2 in case $F = \mathbb{Q}$.

These terms have been familiar[12] since Selberg introduced his original trace formula for compact quotient. The complementary terms on each side are in some sense just as explicit, but often considerably more complex. We shall not discuss them here. In fact, we will allow ourselves to write

$$I_{\text{ell,reg}}(f) \sim I_2(f) \tag{107}$$

as a heuristic approximation of the invariant trace formula. The right hand side is what one wants to understand, and left hand side represents the means by which one hopes to investigate it. It was with this strategy that Langlands created the conjectural theory of endoscopy in the 1970s.

---

[12] In particular, the analogues $J_{\text{ell,reg}}(f)$ and $J_2(f)$ in the original (noninvariant) trace formula were already invariant, as we can recall from the discussion of the special case of GL(2), and are therefore the same as $I_{\text{ell,reg}}(f)$ and $I_2(f)$.

In the larger scheme of things, the invariant trace formula (103) is itself an intermediate step. The final goal was the stable trace formula, which came later [21]. We shall describe it in general terms for further perspective.

First of all, it is useful to keep in mind that there is a simple description of the space of invariant distributions on $G_v = G(F_v)$, for any localization $F_v$ of $F$. For it is known that this space is the closed linear span (with respect to the weak topology) of *either* the set

$$\mathrm{Orb}(\gamma_v, f_v), \qquad \gamma_v \in \Gamma_{\mathrm{reg}}(G_v), f_v \in C_c^\infty(G_v),$$

of strongly regular orbital integrals, *or* the set

$$\Theta(\pi_v, f_v), \qquad \pi_v \in \Pi_{\mathrm{temp}}(G_v), f_v \in C_c^\infty(G_v),$$

of irreducible tempered characters. We can also use the first description here to define the notion of a stable distribution. Let $\Delta_{\mathrm{reg}}(G_v)$ be the set of strongly regular stable conjugacy classes in $G_v$. For any $\delta_v$ in this set, we define the *stable* orbital integral as the associated sum

$$\mathrm{SOrb}(\delta_v, f_v) = \sum_{\gamma_v \to \delta_v} \mathrm{Orb}(\gamma_v, f_v), \qquad f \in C_c^\infty(G_v),$$

of orbital integrals over the finite set of conjugacy classes $\gamma_v$ in $\delta_v$. We then define the subspace of *stable distributions* on $G_v$ to be the closed linear span of the space of stable orbital integrals. The spectral analogue of this description should also be true, but it requires us to know what a stable (tempered) character is. An explicit description of this notion is available in many cases, but not in general. Its general formulation is one of the main goals of the local theory of endoscopy.

The stable trace formula for a quasisplit[13] group $G$ over $F$ is a refinement

$$S^G_{\mathrm{geom}}(f) = S^G_{\mathrm{spec}}(f), \qquad f \in C_c^\infty(G(\mathbb{A})), \tag{108}$$

of the invariant trace formula (103) in which each side is stable.[14] Its construction, and its role in the broader operation of *stabilization*$13^0$, is not difficult to describe in general terms. It is in fact quite similar to the stabilization of the Lefschetz trace formula we discussed briefly in Section 8.

Suppose for a moment that $F$ is a local or global field, and that $G$ is any reductive group over $F$. One of the central notions of endoscopy is the assignment to $G$ of a family of *endoscopic data* $(G', \mathscr{G}', s', \xi')$, where $G'$ is a quasisplit reductive group over $F$, $\mathscr{G}'$ is a split extension of $W_F$ by $\widehat{G}'$, $s'$ is a semisimple element in $\widehat{G}$, and $\xi'\colon \mathscr{G}' \to {}^L G$ is an $L$-homomorphism, subject to various conditions.[15] *Equivalence* of endoscopic data is also defined, by a relation that is closely related to conjugation in ${}^L G$ by elements $g \in \widehat{G}$. (See [165, (1.2)].)

---

[13] The stable trace formula is best regarded as a phenomenon for quasisplit groups. On the other hand, *stabilization*, which we will describe presently applies to arbitrary groups.

[14] A stable distribution on $G(\mathbb{A})$ would of course be a continuous linear form that is stable on each of the factors $G_v$ of $G(\mathbb{A})$.

This is an admittedly technical part of the theory, but the basic idea is simple enough. Its origins in Langlands' sets $\mathscr{D}(T)$ and $\mathscr{E}(T)$, which we will review presently, are really quite remarkable. Basically one wants to attach smaller quasisplit groups $G'$ to $G$ by taking the dual group $\widehat{G}'$ to be the connected centralizer in $\widehat{G}$ of a semisimple element $s'$, and by constructing the $L$-group ${}^L G'$, which then determines $G'$ as a quasisplit group, in terms of the centralizer of $s'$ in the larger group ${}^L G$. If the derived group of $G$ is simply connected, one can identify the subgroup $\xi'(\mathscr{G}')$ of ${}^L G$ with ${}^L G'$ [147]. The general case, however, is a little more subtle, and one has to attach some auxiliary data to $G'$ that serve the same purpose. (See [19] for example.) As in Section 8, we write $G'$ to represent the full endoscopic datum $(G', \mathscr{G}', s', \xi')$. One says that $G'$ is elliptic if the image $\xi'(\mathscr{G}')$ in ${}^L G$ is contained in no proper parabolic subgroup of ${}^L G$, or equivalently, if $(Z(\widehat{G}')^{\Gamma_F})^0$ is mapped by $\xi'$ onto $(Z(\widehat{G})^{\Gamma_F})^0$. Finally, if $F$ is local, there are only finite many equivalence classes of endoscopic data, while if $F$ is global, there are finitely many classes that are unramified outside a given finite set of places.

At the centre of the theory is the endoscopic transfer $f \to f'$ of functions from $G$ to $G'$, a topic we can revisit after our brief discussions from Sections 6–8. It was defined formally by Langlands and Shelstad [165], following Shelstad's treatment of the case of real groups [221]. If $F$ is local, it is a mapping from functions $f \in C_c^\infty(G(F))$ to smooth functions $f'$ on $\Delta_{\mathrm{reg}}(G')$. The Langlands–Shelstad transfer conjecture was then formulated in [165] as the hypothesis that

$$f'(\delta') = \mathrm{SOrb}(\delta', h'), \qquad \delta' \in \Delta_{\mathrm{reg}}(G'),$$

for some function $h' \in C_c^\infty(G'(F))$. It had already been established for archimedean $F$ by Shelstad in [221], using the foundations of harmonic analysis on real groups laid out by Harish-Chandra, as we shall discuss later in this section. In the case of nonarchimedian $F$, the conjecture was reduced to the fundamental lemma by Waldspurger [245], which we have already noted was finally established by Ngô [186], with further contributions [246] from Waldspurger. If $F$ is global, and $f = \prod_v f_v$ lies in $C_c^\infty(G(\mathbb{A}))$, the global mapping is defined by setting $f' = \prod_v f_v'$. One sees that it satisfies the global version of the Langlands–Shelstad conjecture, namely that $f'$ is the image of a function $h' \in C_c^\infty(G'(\mathbb{A}_F))$, by applying the local conjecture to the local functions $f_v$, and the fundamental lemma at places $v \notin S$ for which $f_v$ is the unramified unit function.

With all of this background, we can now describe, again in quite general terms, the stabilization of the invariant trace formula (103). We are assuming that $F$ is a global field, that $G$ is any reductive group over $F$, and that $f$ is a function in $C_c^\infty(G(\mathbb{A}_F))$. The stabilization is then represented by decompositions

---

[15] The main conditions are (iv), (a) and (b), on p. 234 of [165]. They identify $\widehat{G}'$ under the restriction of $\xi'$ with the connected centralizer of $s'$ in $\widehat{G}$, and relate $\mathscr{G}'$ under $\xi'$ with the full centralizer $s'$ in ${}^L G$.

$$I_{\text{geom}}(f) = \sum_{G'} \iota(G, G') \widehat{S}'_{\text{geom}}(f') \tag{109}$$

and

$$I_{\text{spec}}(f) = \sum_{G'} \iota(G, G') \widehat{S}'_{\text{spec}}(f') \tag{110}$$

of the two sides of (103). These are entirely analogous to the decompositions (74) and (75) of the two sides of the Lefschetz trace formula. The summands are indexed by the equivalence classes of elliptic endoscopic data $G'$ for $G$, while $f'$ is the Langlands–Shelstad transfer of $f$ to $G'(\mathbb{A}_F)$. The linear forms $S'_{\text{geom}} = S^{G'}_{\text{geom}}$ and $S'_{\text{spec}} = S^{G'}_{\text{spec}}$ are the analogues for the quasisplit groups $G'$ of the linear forms on each side of (108). In particular, they are stable, and therefore have uniquely determined pairings $\widehat{S}'_{\text{geom}}(f')$ and $\widehat{S}'_{\text{spec}}(f')$ with the function $f'$, in the notation of (74) and (75). The coefficients $\iota(G, G')$ attached to $G$ and $G'$ were defined in [150], and given a particularly simple formula in [119]. The expansions (109) and (110) for arbitrary $G$, and the stable trace formula (108) for quasisplit $G$, are then proven together.

The basic strategy is quite simple. We emphasize that the linear forms $\widehat{S}'_{\text{geom}}(f')$ and $\widehat{S}'_{\text{spec}}(f')$ in the expansions depend only[16] on $G'$, as a quasisplit group, even though the coefficients $\iota(G, G')$, the function $f$, the other components $\mathscr{G}'$, $s'$ and $\xi'$ of $G'$ as an endoscopic datum, and the family $\{G'\}$ itself, all depend on $G$ as well. The summands with $G' \neq G$ in (109) and (110) can therefore be assumed inductively to have been defined.

Suppose first that $G$ is quasisplit. Then $G' = G$ is among the indices of summation in the expansions (109) and (110). We can therefore rewrite them as

$$S_{\text{geom}}(f) = I_{\text{geom}}(f) - \sum_{G' \neq G} \iota(G, G') \widehat{S}'_{\text{geom}}(f')$$

and

$$S_{\text{spec}}(f) = I_{\text{spec}}(f) - \sum_{G' \neq G} \iota(G, G') \widehat{S}'_{\text{spec}}(f').$$

This extends the inductive definition to $G$, and establishes the formula (108) from its analogues for $G'$ and the formula (103). However, there is still something serious to prove in this case. One must show that the right hand side of each of these two expansions is a stable linear form in $f$. Suppose next that $G$ is not quasisplit. Then the sums in (109) and (110) include the term with $G'$ equal to $G^*$, the quasisplit inner form of $G$. Assuming that we have already dealt with this case, we may suppose that all of the terms in (109) and (110) are defined. The remaining problem, then, is simply to establish the two identities in this case. Its proof is deep, but turns out to be quite similar to the proof of stability in the quasisplit case.

This discussion is a useful overview, but it is slightly misleading. For once again, the stable trace formula (108) is not to be regarded as just an abstract identity. Like its predecessors (104) and (103), each side of (108) represents a complex expansion into explicit linear forms, one geometric and one spectral, which are now all stable. The stable trace formula (108) is then to be understood as the identity between these two complex expansions. This is how it is proved, and how it is to be used in the applications of endoscopy. The stabilizations (109) and (110), can be regarded heuristically as an identification of the invariant trace formula (103) with a stable trace formula (represented by the summands with $G = G^*$ in (109) and (110)), modulo an unstable error term (represented by the sum over $G \neq G^*$).

---

[16] They are also the same as in the summands in the stabilization (73) and (74) of the Lefschetz trace formula.

These remarks may not be specific enough to be of much help to a general reader. Part of the reason for rehearsing them is for their application to the concrete terms in the heuristic approximation (107) of the invariant trace formula. They are essentially how Langlands constructed a conjectural but explicit stabilization of each side of (107). More precisely, he constructed a stabilization

$$I_{\mathrm{ell,reg}}(f) = \sum_{G'} \iota(G, G') \widehat{S}'_{\mathrm{ell},G\text{-reg}}(f') \tag{111}$$

of the left hand side of (107) that adheres to the general principles above, but which was obtained directly from its definition (105) in terms of orbital integrals.[17] Using this for guidance, enhanced by the results established for special cases in [127], Langlands then conjectured a partial stabilization of the right hand side of (107). Some of his ideas are contained in the expository sections of the volume [150].

To review these things, we can return to the earlier discussion of the sets $\mathscr{D}(T)$ and $\mathscr{E}(T)$ introduced by Langlands to analyze stable conjugacy classes. We are taking $G$ to be a reductive group with quasisplit inner twist $G^*$, over the local or global field $F$, with a maximal torus $T \subset G$ over $F$. The sets $\mathscr{E}(T)$ are to be regarded as geometric objects, for they are clearly founded on the terms on the left hand side of the (approximate) trace formula (107). It is in their spectral counterparts that class field theory appears, specifically in the Tate–Nakayama duality theory [234].

If $F$ is local, the theory provides a canonical isomorphism

$$H^1(F,T) = H^1(\Gamma, T) \xrightarrow{\sim} \pi_0(\widehat{T}^\Gamma)^*, \quad \Gamma = \Gamma_F = \mathrm{Gal}(\overline{F}/F),$$

of $H^1(F,T)$ with the group of characters on the component group of the $\mathrm{Gal}(\overline{F}/F)$-invariants in the dual torus $\widehat{T}$. In the case that $F$ is global, it provides a canonical isomorphism

$$H^1(F, T(\overline{\mathbb{A}}_F)/T(\overline{F})) = H^1(\Gamma, T(\mathbb{A}_{\overline{F}})/T(\overline{F})) \xrightarrow{\sim} \pi_0(\widehat{T}^\Gamma)^*.$$

Using standard techniques, specifically an application to the short exact sequence

$$1 \to T(\overline{F}) \to T(\mathbb{A}_{\overline{F}}) \to T(\mathbb{A}_{\overline{F}})/T(\overline{F}) \to 1$$

---

[17] The subscript $G$-reg on the right side of (111) denotes the subset of classes in $\Gamma_{\mathrm{ell,reg}}(G')$ that are images of strongly regular classes in $G(F)$. Its dependence on $G$ is an anomaly that would disappear if we had started on the left with the larger set $\Gamma_{\mathrm{ell}}(G)$ of all (elliptic) semisimple classes in $G(F)$.

of $\Gamma_F$-modules to the isomorphism

$$H^1(F, T(\mathbb{A}_{\overline{F}})) \xrightarrow{\sim} \bigoplus_v H^1(F_v, T)$$

provided by Shapiro's lemma, one then obtains a characterization of the diagonal image of $H^1(F, T)$ in the direct sum over $v$ of the groups $H^1(F_v, T)$. It is given by a canonical isomorphism from the cokernel

$$\mathrm{coker}^1(F, T) = \mathrm{coker}(H^1(F, T) \to \bigoplus_v H^1(F_v, T))$$

onto the image

$$\mathrm{im}\big(\bigoplus_v \pi_0(\widehat{T}^{\Gamma_v})^* \to \pi_0(\widehat{T}^{\Gamma})^*\big)$$

(See [30], [234] and [121, §1–2].)

If these results are combined with their analogues for $T_{\mathrm{sc}}$, they provide similar assertions for the subgroups $\mathscr{E}(T/F)$ of $H^1(F, T)$. In the local case, one has only to replace $\pi_0(\widehat{T}^{\Gamma_v})$ by the group $\mathscr{K}(T/F_v)$ of elements in $\pi_0(\widehat{T}/Z(\widehat{G})^{\Gamma})$ whose image in $H^1(F_v, Z(\widehat{G}))$ is trivial. In the global case, one replaces $\pi_0(\widehat{T}^{\Gamma})$ by the group $\mathscr{K}(T/F)$ of elements in $\pi_0(\widehat{T}/Z(\widehat{G})^{\Gamma})$ whose image in $H^1(F, Z(\widehat{G}))$ is locally trivial, in the sense that their image in $H^1(\Gamma_v, Z(\widehat{G}))$ is trivial for each $v$. (See [150] and [121].)

Langlands introduced these ideas to be able to construct the stabilization (111) of the strongly regular part of the trace formula. The first step was to write the left hand side as

$$I_{\mathrm{ell,reg}}(f) = \sum_{\gamma \in \Gamma_{\mathrm{ell,reg}}(G)} \mathrm{vol}(\gamma)\,\mathrm{Orb}(\gamma, f)$$

$$= \sum_{\delta \in \Delta_{\mathrm{ell,reg}}(G)} \mathrm{vol}(\delta)\,\Big(\sum_{\gamma \to \delta} \mathrm{Orb}(\gamma, f)\Big),$$

where $\gamma$ is summed in the brackets over the preimage of $\delta$ in $\Gamma_{\mathrm{ell,reg}}(G)$, and $\mathrm{vol}(\delta) = \mathrm{vol}(\gamma)$ depends only on $\delta$. We are of course assuming that $F$ is global here. The last sum over $\gamma$ looks as if it might be stable in $f$. But stability is a local concept, and there are not enough rational conjugacy classes $\gamma$ to make this sum stable in each component $f_v$ of $f$. The problem is the failure of every $G(\mathbb{A}_F)$-conjugacy class in the $G(\mathbb{A}_F)$-stable class of $\delta \in \Delta_{\mathrm{ell,reg}}(G)$ to have a representative $\gamma$ in $G(F)$. If $T = G_\delta$, the cokernel we denoted by $\mathrm{coker}^1(F, T)$ gives a measure of this failure. Langlands' construction treats the sum $\big(\sum_{\gamma \to \delta} \mathrm{Orb}(\gamma, f)\big)$ as the value at 1 of a function on the finite abelian group $\mathrm{coker}^1(F, T)$. The critical step is to expand this function according to Fourier inversion on $\mathrm{coker}^1(F, T)$. It leads naturally to the definition of endoscopic data $\{G'\}$, and finally the desired stabilization (111).

To simplify the construction, we might as well assume for the present that $G = G_{\mathrm{sc}}$. Then $T = T_{\mathrm{sc}}$, while $\mathscr{E}(T/F) = H^1(F, T)$ and $\mathscr{K}(T/F) = \pi_0(\widehat{T}^{\Gamma})$ if $F$ is either local or global. In particular, $\mathscr{K}(T/F) = \widehat{T}^{\Gamma}$ if $T$ is elliptic in $G$ over $F$.

With this condition on $G$, we apply Fourier inversion on $\operatorname{coker}^1(F, T)$. One has to keep track here of the redundancy from the sum $\gamma \to \delta$, given by the set of $G(F)$-conjugacy classes in the $G(\mathbb{A}_F)$-conjugacy class of $\delta$ (regarded as a representative in $G(F)$ of the class in $\Delta_{\mathrm{ell,reg}}(G)$). It can be seen (with an appeal to the Hasse principle for $G = G_{\mathrm{sc}}$) that this is bijective with the finite abelian group

$$\ker^1(F, T) = \ker(H^1(F, T) \to \bigoplus_v H^1(F_v, T)).$$

It then follows that

$$I_{\mathrm{ell,reg}}(f) = \sum_{\delta \in \Delta_{\mathrm{ell,reg}}(G)} \iota(T) \operatorname{vol}(\delta) \sum_{\kappa \in \widehat{T}^\Gamma} \operatorname{Orb}^\kappa(\delta, f), \tag{112}$$

where $T = G_\delta$ is the centralizer of (some fixed representative of) $\delta$, $\iota(T)$ equals the product of $(\widehat{T}^\Gamma)^{-1}$ with $|\ker^1(F, T)|$, and

$$\operatorname{Orb}^\kappa(\gamma, f) = \sum_{\{\gamma_\mathbb{A} \in \Gamma(G(\mathbb{A}_F)) : \gamma_\mathbb{A} \sim \delta\}} \operatorname{Orb}(\gamma_\mathbb{A}, f) \kappa(\gamma_\mathbb{A}).$$

The last sum is over the $G(\mathbb{A}_F)$-conjugacy classes $\gamma_\mathbb{A}$ in the stable class of $\delta$ in $G(\mathbb{A}_F)$. For any such $\gamma_\mathbb{A}$, its local component $\gamma_v$ is $G(F_v)$-conjugate to $\delta_v$ for almost all $v$, from which it follows that $\gamma_\mathbb{A}$ maps to an element $t_\mathbb{A}$ in the direct sum of the groups $H^1(F_v, T)$. This in turn maps to a point in the cokernel (111), and hence to a character in $(\widehat{T}^\Gamma)^*$. The coefficient $\kappa(\gamma_\mathbb{A})$ is the value of this character at $\kappa$.

The expression (112) is a step closer to the desired stabilization of $I_{\mathrm{ell,reg}}(f)$. In particular, it contains the origins of the endoscopic data $G'$ in (111). For suppose that $T$ and $\kappa$ are as in (112). One chooses an admissible[18] embedding $\widehat{T} \subset \widehat{G}$ of its dual group, taking then $s' \in \widehat{G}$ to be the resulting image of $\kappa \in \widehat{T}$, and $\widehat{G}' = \widehat{G}_{s'}$ the connected centralizer of $s'$. It is known that there is also an $L$-embedding

$${}^L T = \widehat{T} \rtimes W_F \to {}^L G = \widehat{G} \rtimes W_F,$$

of the $L$-group of $T$ into that of $G$, which restricts to the given embedding of $\widehat{T}$ into $\widehat{G}$. (This is a little more subtle and entails some choices to which the embedding is sensitive. See [165, (2.6)].) In any case, for a fixed such embedding, the product

$$\mathscr{G}' = {}^L T \cdot \widehat{G}'$$

is an $L$-subgroup of ${}^L G$, which commutes with $s'$, and provides a split extension

$$1 \to \widehat{G}' \to \mathscr{G}' \to W_F \to 1$$

---

[18] This means that it is the mapping assigned to a choice of some pair $(\widehat{B}, \widehat{T})$ in $\widehat{G}$, and some Borel subgroup $B$ of $G$ containing $T$.

of $W_F$ by $\widehat{G}'$. In particular, it determines an action of $W_F$ on $\widehat{G}'$ by outer automorphisms, which factors through a finite quotient of $\Gamma_F$. We take $G'$ to be a quasisplit group over $F$ for which $\widehat{G}'$, with the given action of $\Gamma_F$, is a dual group. Finally, if we let $\xi'$ be the identity $L$-embedding of $\mathscr{G}'$ into ${}^L G$, the 4-tuple $(G', s', \mathscr{G}', \xi')$ becomes an endoscopic datum for $G$. We have thus obtained a correspondence

$$(T, \kappa) \to (G', s', \mathscr{G}', \xi'),$$

from the pairs $(T, \kappa)$ taken from (112), to the endoscopic data derived from them as above.

There is another ingredient to the last correspondence. Given the pair $(T, \kappa)$, we can choose a maximal torus $T' \subset G'$ over $F$, together with an isomorphism from $T'$ to $T$ over $F$ that is admissible, in the sense that the associated isomorphism $\widehat{T}' \to \widehat{T}$ of dual groups is the composition of an admissible embedding $\widehat{T}' \subset \widehat{G}$ (as in the footnote 18) with an inner automorphism of $\widehat{G}$ that takes $\widehat{T}'$ to $\widetilde{\widehat{T}}$. Let $\delta' \in T'(F)$ be the associated preimage of the original point $\delta$. The tori $T$ and $T'$ are the centralizers in $G$ and $G'$ of $\delta$ and $\delta'$, so we can regard $\delta$ and $\delta'$ as the primary objects. They become part of a larger correspondence

$$(\delta, \kappa) \to (G', \delta') = ((G', \mathscr{G}', s', \xi'), \delta'). \tag{113}$$

Elements $\delta'$ obtained in this way are called *images from $G$* [165, (1.3)].

Now suppose that $G$ is arbitrary. The general form of the expansion (112) is derived the same way, and takes an almost identical form

$$I_{\mathrm{reg,ell}}(f) = \sum_{\delta \in \Delta_{\mathrm{ell,reg}}(G)} \iota(T, G) \operatorname{vol}(\delta) \sum_{\kappa \in \mathscr{K}(T/K)} \operatorname{Orb}^\kappa(\delta, f),$$

where

$$\iota(T, G) = |\ker(\mathscr{E}(T/F) \to \bigoplus_v \mathscr{E}(T, F_v))| |\kappa(T/F)|^{-1}$$

and $\operatorname{Orb}^\kappa(\delta, f)$ is defined as in (112). The correspondence (113) remains in place, and is easily seen to have an inverse, which in general extends to a bijection

$$\{(G', \delta')\} \xrightarrow{\sim} \{(\delta, \kappa)\}. \tag{114}$$

The domain is the set of equivalence classes of pairs $(G', \delta')$, where $G'$ is an elliptic endoscopic datum for $G$, $\delta'$ is a strongly $G$-regular, elliptic element in $G'(F)$ that is an image from $G$, and equivalence is defined by isomorphisms of endoscopic data. The range is the set of equivalence classes of pairs $(\delta, \kappa)$, where $\delta$ belongs to $\Delta_{\mathrm{ell,reg}}(G)$, $\kappa$ lies in $\mathscr{K}(G_\delta/F)$, and equivalence is defined by conjugation by elements in $G(\overline{F})$. (See [150] and [121, Lemma 9.7].) Given $(G', \delta')$, we set

$$f'(\delta') = f_G^\kappa(\delta), \tag{115}$$

where to emphasize the bijection, and to keep us mindful of its essential simplicity, we have written $f_G^\kappa(\delta)$ in place of $\mathrm{Orb}^\kappa(\delta, f)$. In other words

$$f_G^\kappa(\gamma_\mathbb{A}) = \sum_{\{\gamma_\mathbb{A} \in \Gamma(G(\mathbb{A})): \gamma_\mathbb{A} \sim \delta\}} f_G(\gamma_\mathbb{A}) \kappa(\gamma_\mathbb{A}), \qquad (116)$$

with

$$f_G(\gamma_\mathbb{A}) = \mathrm{Orb}(\gamma_\mathbb{A}, f).$$

We can then write

$$I_{\mathrm{ell,reg}}(f) = \sum_{G' \in \mathscr{E}_{\mathrm{ell}}(G)} |\mathrm{Out}_G(G')|^{-1} \sum_{\delta' \in \Delta_{\mathrm{ell}, G\text{-reg}}} \mathrm{vol}(\delta') \iota(G_\delta, G) f'(\delta')$$

for the finite group

$$\mathrm{Out}_G(G') = \mathrm{Aut}_G(G')/\mathrm{Int}(G')$$

of outer automorphisms of $G'$ as an endoscopic datum, and with the understanding that $f'(\delta') = 0$ if $\delta'$ is not an image from $G$. Langlands showed that for any pair $(G', \delta')$, the number

$$\iota(G, G') = (\iota(G_\delta, G) \iota(G'_{\delta'}, G')^{-1}) |\mathrm{Out}_G(G')|^{-1}$$

is independent of $\delta'$ and $\delta$. (Kottwitz later expressed the product in the brackets on the right as a quotient $\tau(G)\tau(G')^{-1}$ of Tamagawa numbers [119, Theorem 8.3.1].) Set

$$\widehat{S}'_{\mathrm{ell}, G\text{-reg}}(f') = \sum_{\delta' \in \Delta_{\mathrm{ell}, G\text{-reg}}} \mathrm{vol}(\delta') \iota(G'_{\delta'}, G') f'(\delta'). \qquad (117)$$

It then follows that

$$I_{\mathrm{ell,reg}}(f) = \sum_{G' \in \mathscr{E}_{\mathrm{ell}}(G)} \iota(G, G') \widehat{S}'_{\mathrm{ell}, G\text{-reg}}(f'), \qquad (118)$$

since

$$\mathrm{vol}(\delta) = \mathrm{vol}(G_\delta(F) \backslash G_\delta(\mathbb{A}_F)^1) = \mathrm{vol}(G'_{\delta'}(F) \backslash G'_{\delta'}(\mathbb{A}_F)^1) = \mathrm{vol}(\delta').$$

We have sketched how Langlands derived a version (118) of the desired formula (111). However, the term $f'(\delta')$ in (117) is defined in (115) only as a function on the rational classes $\Delta_{\mathrm{ell}, G-\mathrm{reg}}(G')$. As we have said earlier, one wants it to be the adelic stable orbital integral of a function in $C_c^\infty(G'(\mathbb{A}_F))$. But to this point, we do not yet have a well defined candidate for its *local* orbital integrals. The sum in (116) is over adelic products $\gamma_\mathbb{A} = \prod_v \gamma_v$, in which $\gamma_v$ is a conjugacy class in $G(F_v)$ that lies in the stable class of the image $\delta_v$ of $\delta$ in $G(F_v)$. It follows that if $f = \prod_v f_v$, then

$$f'(\delta') = f_G^\kappa(\delta) = \prod_v f_v^\kappa(\delta_v)$$

where

$$f_v^{\kappa}(\delta_v) = \prod_{\gamma_v \sim \delta_v} f_{v,G}(\gamma_v) \kappa(\gamma_v). \tag{119}$$

But this is not quite the local definition we are looking for. The problem is that we have been treating $\delta$ as both a stable class in $\Delta_{\text{ell,reg}}(G)$ and a representative in $G(F)$ of that class. The distinction has not mattered up until now, since $f'(\delta') = f_G^{\kappa}(\delta)$ depends only on the class of $\delta$. However, the coefficients $\kappa(\gamma_v)$ in the local functions (119) are defined in terms of the "relative position" of $\gamma_v$ and $\delta_v$, a notion that comes from the original pairing between $H^1(F_v, T)$ and $\pi_0(\widehat{T^{\Gamma_v}})$, and is sensitive to how $\gamma_v$ and $\delta_v$ are situated within their local conjugacy classes.

The solution for Langlands and Shelstad was to replace $\kappa(\gamma_v)$ with a function $\Delta_G(\delta_v', \gamma_v)$ that they called a *transfer factor*. This is the deepest part of the theory, and it is the content of the papers [165] and [166]. The function is defined as a product of $\kappa(\gamma_v)$ with some subtle factors that depend on $\delta_v'$ and $\delta_v$, but not on $\gamma_v$. The product $\Delta_G(\delta_v', \gamma_v)$ then turns out to be independent of the choice of $\delta_v$, and depends therefore only on the local stable class of $\delta_v'$ in $G'(F_v)$ and the local conjugacy class of $\gamma_v$ in $G(F_v)$. Moreover, if $\delta_v'$ is the local image of $\delta' \in \Delta_{\text{ell},G\text{-reg}}(G')$ for each $v$, the product over $v$ of the corresponding local transfer factors is equal to the coefficient $\kappa(\gamma_{\mathbb{A}})$ in (116).

Transfer factors play the role of a kernel in the local transfer of functions. After first introducing them, Langlands and Shelstad defined the transfer to $G_v'$ of a function $f_v \in C_c^{\infty}(G_v)$ on $G_v$ as an "integral transform"

$$f_v'(\delta_v') = \sum_{\gamma_v \in \Gamma_{\text{reg}}(G_v)} \Delta_G(\delta_v', \gamma_v) f_{v,G}(\gamma_v), \quad \delta_v' \in \Delta_{G\text{-reg}}(G_v'), \tag{120}$$

where

$$f_{v,G}(\gamma_v) = |D(\gamma_v)|^{\frac{1}{2}} \text{Orb}(\gamma_v, f_v), \qquad \gamma_v \in \Gamma_{G\text{-reg}}(G_v), \tag{121}$$

is now a *normalized*[19] orbital integral, and $f_v'(\delta_v') = f_v^{G'}(\delta_v')$ is the analogue for $G_v'$ of the normalized stable orbital integral

$$f_v^G(\delta_v) = |D(\delta_v)|^{\frac{1}{2}} \text{SOrb}(\delta_v, f_v), \qquad \delta_v \in \Delta_{G\text{-reg}}(G_v), \tag{122}$$

for $G_v$. The normalizing factor is the absolute value of the Weyl discriminant

$$D(\gamma_v) = D_G(\gamma_v) = \det(1 - \text{Ad}(\gamma_v))_{\mathfrak{g}_v / \mathfrak{t}_v},$$

---

[19] This is not in conflict with the global notation in (115), thanks to the product formula

$$|D(\gamma)| = \prod_v |D(\gamma_v)|_v = 1, \qquad \gamma \in G(F).$$

Langlands and Shelstad put the quotient $|D_G(\gamma_v)|^{\frac{1}{2}} |D_{G'}(\delta_v')|^{-\frac{1}{2}}$ into their transfer factor as term $\Delta_{\text{IV}}(\delta_v', \gamma_v)$ in [165, §3.6]. However, it is instructive to use it to normalize the orbital orbitals, as we will observe later in the section, even as we continue to use the notation of [165] for the transfer factor in (120) for the Lie algebras $\mathfrak{g}_v$ and $\mathfrak{t}_v$ of $G_v$ and $T_v = G_{\gamma_v}$. It was only then that they could pose their local transfer conjecture. It became part of the global conjecture (together with the fundamental lemma for the nonarchimedean unit function), which we recall was established later.

Thus, the function $f'(\delta')$ in (117) really is the stable orbital integral at $\delta'$ of a function $h'$ in $C_c^\infty(G'(\mathbb{A}))$. It is at this point that one can treat the left hand side of (117) as the pairing of a stable distribution $S'_{\mathrm{ell},G\text{-}\mathrm{reg}}$, defined by the inductive procedure[20] from (118) outlined earlier, with the function in $C_c^\infty(G'(\mathbb{A}))$, rather than just the sum over rational points $\delta'$ on the right hand side of (117).

There are still a couple of technical points that we should at least mention. The Langlands–Shelstad transfer factor depends on a choice of $L$-embedding of $^L G'$ into $^L G$. If $G_{\mathrm{der}}$ equals $G_{\mathrm{sc}}$, such embeddings exist, and to fix one, it suffices to choose an $L$-isomorphism from $\mathscr{G}'$ to $^L G'$. If not, $\mathscr{G}'$ might not be $L$-isomorphic to $^L G'$. (It is a question whether a certain 2-cocycle with values in $Z(\widehat{G})$ splits.) In this case, minor adjustments have to be made, which entail choosing a central extension

$$1 \to \widetilde{C}' \to \widetilde{G}' \to G' \to 1,$$

and taking $f'$ to be a function on $\widetilde{G}'(\mathbb{A})$ with a certain central character on $\widetilde{C}'(F) \backslash \widetilde{G}'(\mathbb{A})$. (See [165, (4.4)], which involves also taking a central extension $\widetilde{G}$ of $G$, or [22, p. 202], which is based on the adjustment made in [125] for *twisted* transfer factors.)

Another point is the subscript "$G$-reg" in the summands on the right hand side of (118). This is a minor logical violation of our general inductive definitions in the case of the stable linear form $S^G_{\mathrm{ell},\mathrm{reg}}$ on $G(\mathbb{A})$, since the stable distributions $S'$ in (118) are supposed to depend only on $G'$ (and not $G$). We have already remarked in the footnote 20 that we do not need the general inductive definitions in this concrete case. Still, to be consistent, why don't we just replace the subscripts "$G$-reg" by "reg", and replace the equality sign in (118) by the "heuristic approximation" symbol $\sim$ we have already used in (107). This in any case is philosophically sound since the complement of $\Delta_{G\text{-}\mathrm{reg}}(G')$ in $\Delta_{\mathrm{reg}}(G')$ is sparse.

The stabilization of the strongly regular, elliptic part of the trace formula becomes

$$I_{\mathrm{ell},\mathrm{reg}}(f) \sim \sum_{G' \in \mathscr{E}_{\mathrm{ell}}(G)} \iota(G,G') \widehat{S}_{\mathrm{ell},\mathrm{reg}}(f'). \tag{123}$$

We could combine this with (107), and the elementary induction arguments that precede (109) and (110). We would then obtain a stabilization

$$I_2(f) \sim \sum_{G' \in \mathscr{E}_{\mathrm{ell}}(G)} \iota(G,G') \widehat{S}_2(f') \tag{124}$$

---

[20] One does not actually need the general inductive definition in the simple case here of the $G$-regular elliptic terms. One obtains the stable distribution $S^G_{\mathrm{ell},G\text{-}\mathrm{reg}}$ ($G$ being quasisplit) directly as a line combination of stable, adelic orbital integrals from the construction of Langlands we have just described.

of the $L^2$-discrete part of the trace formula. These arguments all apply to "approximate" identities, which means that (124) is something we expect to be true. Langlands reviewed such arguments, and would then have used (124) to guess at some of the spectral implications of the theory of endoscopy [150]. These include most notably versions of his conjectural classification of representations into local and global $L$-packets.

We have concluded our discussion of the explicit stabilization (123) of the regular elliptic part (104) of the invariant trace formula. It must seem rather murky to a nonspecialist. It is helpful to think of the bijection (114) as the centre of the process. We can illustrate the transition schematically as follows.

$$I_{\mathrm{ell,reg}}(f) \qquad (104)$$

Galois cohomology $\Big\downarrow$

$$\{(G', \delta')\} \qquad (114)$$

transfer factors $\Big\downarrow$

$$\sum_{G'} \iota(G, G')\widehat{S}'_{\mathrm{ell,reg}}(f') \qquad (123)$$

Before going on, we recall that there are a number of further topics we promised to take up in this section. We may not be able to give them all the attention they deserve, but there is one critical paper that would in any case be the next step in this discussion. It is Langlands' classification of the representations of real groups [151], which together with subsequent work of Shelstad, represents an essential link[21] between the earlier work of Harish-Chandra on representation theory and the emerging theory of endoscopy. We shall take this opportunity for a short digression on the work of Harish-Chandra, in which we assume that $G$ is a reductive algebraic group over $\mathbb{R}$.

We have alluded to the construction by Harish-Chandra of the discrete series [86], [88], those representations of a real group that occur discretely[22] in $L^2(G(\mathbb{R}))$, but we have not described it. It was a climax in his long and comprehensive study of the harmonic analysis on a general (semisimple) real group $G(\mathbb{R})$.

Harish-Chandra's harmonic analysis also represents an interplay between the geometric objects and the spectral objects on $G(\mathbb{R})$. These are the orbital integrals and the irreducible characters, whose global versions became the heart of the trace formula. They were both introduced by Harish-Chandra in the early stages of his career. Both are fundamental and deep. It was of course the theory of characters that became a foundation for the discrete series.

---

[21] As it could also be argued that Langlands' manuscript on Eisenstein series represents a link between Harish-Chandra's investigations into the Plancherel formula and a future trace formula, even though Langlands' monumental volume stands on its own, and in fact also influenced the subsequent course of Harish-Chandra's work.

[22] We have been using the term *relative discrete series* (or *square integrable representations*) to describe the representations that occur discretely modulo the centre $Z(\mathbb{R})$ of $G(\mathbb{R})$. These are slightly more general. They are the representations of Levi subgroups used in the parabolic induction process that yields all the tempered representations, the ones that occur in the full spectral decomposition of $L^2(G(\mathbb{R}))$.

Recall that an irreducible unitary representation $\pi$ of $G(\mathbb{R})$ is infinite dimensional (unless it is 1-dimensional or attached to a representation of compact factor on $G(\mathbb{R})$). It was not initially clear how it could have a character, since the trace of an infinite dimensional unitary matrix $\pi(x)$ is not defined. Harish-Chandra's idea was to make systematic use of the general theory of distributions that had just been introduced by Laurent Schwartz [203], [204]. For any irreducible $\pi$, Harish-Chandra proved that the operator

$$\pi(f) = \int_{G(\mathbb{R})} f(x)\pi(x)\,dx$$

attached to a function $f \in C_c^\infty(G(\mathbb{R}))$ was of trace class, and that the linear form

$$f \to \Theta(\pi, f) = \mathrm{tr}(\pi(f))$$

was a distribution (which is to say, continuous for the usual topology on $C_c^\infty(G(\mathbb{R}))$). This is what he called the *character* of $\pi$. The proof was not particularly difficult as these things go. Much deeper was a second theorem on characters, his so called *regularity theorem*. It asserts that any (irreducible) character $\Theta(\pi)$ is a locally integrable function $x \to \Theta(\pi, x)$ on $G(\mathbb{R})$, which is to say that

$$\Theta(\pi, f) = \int_{G(\mathbb{R})} \Theta(\pi, x)f(x)\,dx, \qquad f \in C_c^\infty(G(\mathbb{R})). \tag{125}$$

The regularity theorem is really about the differential equations

$$z\Theta = \chi(\Theta, z)\Theta, \qquad z \in \mathscr{Z}_G, \tag{126}$$

satisfied by any invariant eigendistribution $\Theta$ of the centre $\mathscr{Z}_G$ of the universal enveloping algebra $\mathscr{U}_G$ of the complex Lie algebra of $G(\mathbb{R})$. This is a property that holds for any character $\Theta = \Theta(\pi)$, by an infinite dimensional version of Schur's lemma previously established by Harish-Chandra. In this case, the homomorphism

$$\chi(\Theta)\colon z \to \chi(\Theta, z), \qquad z \in \mathscr{Z}_G,$$

from $\mathscr{Z}_G$ to $\mathbb{C}^*$ is called the *infinitesimal character* of $\pi$. Like many results in this area, Harish-Chandra's argument yields not only the existence of the function $\Theta(\pi, x)$, but frequently also an interesting, explicit formula that it satisfies. The first (and easier) half of the proof uses the elliptic regularity theorem for differential equations to prove that the restriction of $\Theta(\pi, x)$ to the open, dense subset $G_{\mathrm{reg}}(\mathbb{R})$ of (strongly) regular elements in $G(\mathbb{R})$ is a (real) analytic function of $x$. The second half classifies the singularities of the *normalized* character

$$\Phi(\pi, x) = |D(x)|^{\frac{1}{2}}\Theta(\pi, x), \qquad x \in G_{\mathrm{reg}}(\mathbb{R}), \tag{127}$$

by the Weyl discriminant $D(x)$, at the hypersurfaces in the complement of $G_{reg}(\mathbb{R})$ in $G(\mathbb{R})$. It establishes that any left invariant derivative of $\Phi(\pi, x)$ remains bounded as $x$ approaches a singular hypersurface, and for many $\pi$, also gives an explicit formula for the "jump" of the function as it crosses the hypersurface. This yields an interesting boundary value problem satisfied by $\Phi(\pi, x)$ on the closures of the open connected components of $G_{reg}(\mathbb{R})$. Its solution is what provides the explicit formula for $\Phi(\pi, x)$.

Since characters are invariant distributions, their functions $\Theta(\pi, x)$ are conjugacy invariant in $x$. This can be combined with the Weyl integration formula

$$\int_{G(\mathbb{R})} h(x)\, dx = \sum_{\{T\}} |W(G(\mathbb{R}), T(\mathbb{R}))|^{-1} \int_{T_{reg}(\mathbb{R})} \left( |D(t)| \int_{T(\mathbb{R}) \backslash G(\mathbb{R})} h(x^{-1}tx)\, dx \right) dt$$

for the change of variables used to express the integral of a function $h \in C_c(G(\mathbb{R}))$ as an integral of its averages over conjugacy classes. Here $\{T\}$ is the set of $G(\mathbb{R})$-conjugacy classes of maximal tori in $G(\mathbb{R})$, $W(G(\mathbb{R}), T(\mathbb{R}))$ is the normalizer of $T(\mathbb{R})$ in $G(\mathbb{R})$ modulo its centralizer $T(\mathbb{R})$, $T_{reg} = T \cap G_{reg}$, and

$$D(t) = \det((1 - \mathrm{Ad}(t))_{\mathfrak{g}/\mathfrak{t}})$$

is again the Weyl discriminant. It then follows from (125) that

$$\Theta(\pi, f) = \sum_{\{T\}} |W(G(\mathbb{R}), T(\mathbb{R}))|^{-1} \int_{T_{reg}} \Phi(\pi, t) f_G(t)\, dt, \tag{128}$$

for the normalized character $\Phi(\pi, t)$ from (127) and the normalized orbital integrals $f_G(t)$ from (121). Harish-Chandra used this formula repeatedly in his development of the discrete series.

Here in general terms is what he proved. First of all, a (connected) reductive group $G$ over $\mathbb{R}$ has a discrete series of representations $\pi$ if and only if it has a maximal torus $T$ over $\mathbb{R}$ that is *anisotropic*, which means that $T(\mathbb{R})$ is compact. We should recall that for *any* maximal torus $T \subset G$ in $G$ over $\mathbb{R}$, we have a chain of three Weyl groups

$$W(G(\mathbb{R}), T(\mathbb{R})) \subset W_{\mathbb{R}}(G, T) \subset W(G, T),$$

in which $W(G, T)$ is the full (complex) Weyl group and $W_{\mathbb{R}}(G, T)$ is the subgroup of elements that stabilize $T(\mathbb{R})$, while $W(G(\mathbb{R}), T(\mathbb{R}))$ is as in (128), the subgroup of elements in $W_{\mathbb{R}}(G, T)$ induced from $G(\mathbb{R})$. In the case here that $T(\mathbb{R})$ is compact, $W(G, T)$ equals $W_{\mathbb{R}}(G, T)$, but $W(G(\mathbb{R}), T(\mathbb{R}))$ is generally a proper subgroup of $W_{\mathbb{R}}(G, T)$. This last circumstance is responsible for some of the complexity of discrete series representations. The second main property he established is that $\pi$ is completely determined by the restriction $\Theta(\pi, t)$ of its character to the anisotropic torus $T(\mathbb{R})$. Harish-Chandra in fact showed that this restriction $\Theta(\pi, t)$ satisfies an explicit formula that is more complicated than, but nevertheless reminiscent of, the Weyl character formula.

We recall that the Weyl character formula applies to the special case that $G$ is anisotropic, which means that $G(\mathbb{R})$ itself is compact. The discrete series accounts for all of the irreducible representations in this case, and they are all finite dimensional. According to Weyl's classification, they are parametrized by orbits $\{\chi\}$ of characters $\chi$ on $T(\mathbb{R})$ under the Weyl group $W(G,T)$. To state the Weyl character formula for any such representation, we have to choose an order on the roots $\{\alpha\}$ of $(G,T)$. This gives a corresponding set $\{\alpha > 0\}$ of positive roots as well as the associated linear form

$$\rho = \tfrac{1}{2} \sum_{\alpha > 0} \alpha \tag{129}$$

on the Lie algebra $\mathfrak{t}(\mathbb{R})$ of $T(\mathbb{R})$, and from each $W(G,\mathbb{R})$-orbit $\{\chi\}$, a unique character $\chi$ whose differential $d\chi$ lies in the closure of the associated positive chamber in the dual space $\mathfrak{t}^*(\mathbb{R})$. The Weyl character formula for the representation $\pi_\chi$ attached to the $W(G,T)$-orbit of $\chi$ is then

$$\Theta(\pi_\chi, t) = \sum_{w \in W(G(\mathbb{R}), T(\mathbb{R}))} \left( \frac{\varepsilon(w) \chi(w \cdot \exp H) e^{\rho(wH - H)}}{\Delta(\exp H)} \right),$$

for any point

$$t = \exp H$$

in $T_{\mathrm{reg}}(\mathbb{R})$, and for

$$\Delta(\exp H) = \prod_{\alpha > 0} \left( 1 - e^{-\alpha(H)} \right).$$

It is a simple matter to check that the right hand side of the formula remains the same if either $H$ or $\chi$ is replaced by a Weyl translate $wH$ or $w\chi$, for any $w \in W(G,T)$. The former property is needed for the function $\Theta(\pi_\chi, \cdot)$ on $G(\mathbb{R})$ to be invariant under conjugation, the latter for it to depend only of the $W(G,T)$-orbit of $\chi$.

Harish-Chandra's discrete series had of course to generalize this. The obvious structural difference is that it is the subgroup $W(G(\mathbb{R}), T(\mathbb{R}))$ of

$$W(G,T) = W_{\mathbb{R}}(G,T),$$

acting on the anisotropic torus $T$, that ought to reflect the $G(\mathbb{R})$-invariance of characters. The similarity is that they would still turn out to be determined by their values on $T(\mathbb{R})$. The natural guess, in retrospect at least, would be that the discrete series are parametrized by suitably defined $W(G(\mathbb{R}), T(\mathbb{R}))$-orbits. This is precisely what Harish-Chandra established, but only after years of concentrated study of the underlying harmonic analysis. All that we need to specify his classification is his formula [86, Theorem 3] and [88, Theorem 18] for their characters on $T_{\mathrm{reg}}(\mathbb{R})$. We shall state the version of it formulated by Langlands on p. 134 of [151], which is closer to the Weyl character formula, and is also compatible with Langlands' ideas on the broader classification of representations.

Given the group $G$ with anisotropic maximal torus $T$, we can choose characters $\chi$ on $T(\mathbb{R})$ and an order on the roots $\{\alpha\}$, as in the special case of anisotropic $G$ above. Langlands represents the order by the corresponding linear form $\rho$ in (129), and therefore considers pairs $(\chi, \rho)$ in which $\chi$ lies in the closure of the positive chamber in $\mathfrak{t}^*$ attached to $\rho$. Harish-Chandra's classification of discrete series is then given by a bijection

$$(\chi, \rho) \to \pi_{\chi, \rho}$$

from the $W(G(\mathbb{R}), T(\mathbb{R}))$-orbits of such pairs onto the equivalence classes of discrete series representations such that

$$\Theta(\pi_{\chi, \rho}, t) = (-1)^{q_G} \sum_{w \in W(G(\mathbb{R}), T(\mathbb{R}))} \left( \frac{\varepsilon(w) \chi(w \cdot \exp H) e^{\rho(wH - H)}}{\Delta(\exp H)} \right), \qquad (130)$$

for any point $t = \exp H$ in $T_{\mathrm{reg}}(\mathbb{R})$, where $q_G = \frac{1}{2} \dim(G(\mathbb{R})/K_{\mathbb{R}})$ is one-half of the dimension of the symmetric space attached to $G(\mathbb{R})$.

We note in passing that the values of $\Theta(\pi_{\chi, \rho}, t_1)$ on a general maximal torus $T_1 \subset G$ can be reduced according to the theory developed by Harish-Chandra to its values (130) on $T(\mathbb{R})$. This is because the singularities of the normalized character $\Phi(\pi_{\chi, \rho}, t_1)$, expressed in Harish-Chandra's jump conditions at a singular hypersurface $T_{01}$, are given in terms of the limits at points $t_{01} \in T_{01}(\mathbb{R})$ of its values $\Phi(\pi_{\chi, \rho}, t_0)$ on a maximal torus[23] $T_0$ that shares the hypersurface $T_{01}$ with $T_1$, but whose anisotropic part is of one dimension greater than that of $T_1$. Increasing the anisotropic dimension $d_a(T_1)$ of a maximal torus $T_1$ makes the corresponding character values simpler. For this reason, the solution of the boundary value problem for $\Phi(\pi_{\chi, \rho}, t_1)$ follows by decreasing induction on $d_a(T_1)$, using the differential equations (126), the jump conditions, the basic explicit formula (130) on $T(\mathbb{R})$, and the fact that $\Phi(\pi_{\chi, \rho}, t_1)$ is bounded on $T_1(\mathbb{R})$, a consequence in turn of Harish-Chandra's proof that the characters $\Theta(\pi_{\chi}, \rho)$ of discrete series are tempered distributions. (See the formulas of Harish-Chandra [86] and their simplifications in [95].) We recall from Section 8 that such formulas arose later in the invariant trace formula [15], [17], and were then used by Morel [183] for the geometric boundary terms in the Lefschetz trace formula for the Shimura varieties attached to $\mathrm{GSp}(2n)$.

Langlands arrived at his version (130) of Harish-Chandra's formula in the process of relating the irreducible representations of $G(\mathbb{R})$ to his ideas for parameterizing the local components of automorphic representations. The obvious lesson to be taken from Harish-Chandra's discrete series is that these representations occur naturally in finite sets. Each such set corresponds to an irreducible finite-dimensional representation of the complex group $G(\mathbb{C})$, or equivalently, the representation $\pi_\chi$ of the compact real form of $G(\mathbb{C})$ with highest weight $d\chi$ in the Weyl classification. It consists of representations with the same infinitesimal character, and can be parameterized by the cosets in $W(G(\mathbb{R}), T(\mathbb{R})) \backslash W(G, T)$. These finite sets became known as $L$-packets, a term we have used regularly throughout the report without actually defining it.

To see its meaning, we note that the first two sections of Langlands' paper [151] were devoted to something quite different. For any group $G$ over $\mathbb{R}$, Langlands con-

---

[23] One says that $T_1$ is a *Cayley transform* of $T_0$.

sidered local $L$-homomorphisms

$$\phi\colon W_{\mathbb{R}}\to {}^LG$$

from the Weil group to the $L$-group, with the property[24] that if the image of $\phi$ is contained in some parabolic subgroup ${}^LP$ of ${}^LG$, the corresponding parabolic subgroup $P$ of $G$ is defined over $\mathbb{R}$. He wrote $\Phi(G)$ for the set of $\widehat{G}$-conjugacy classes of such parameters $\phi$, as we have noted earlier, and $\Phi_2(G)$ for the subset of such classes such that the image of $\phi$ lies in no proper parabolic subgroup ${}^LP$ of ${}^LG$. For the local field $\mathbb{R}$, these are reasonably elementary objects. Langlands calculated them directly in terms of simple data within ${}^LG$. He then observed that for the group $G$ with anisotropic torus $T$, the set $\Phi_2(G)$ was naturally bijective with the $L$-packets of discrete series. On the other hand, for any representation

$$r\colon {}^LG\to \mathrm{GL}(n,\mathbb{C}),$$

one can attach $L$-functions $L(s,r\circ\phi)$ and $\varepsilon$-factors $\varepsilon(s,r\circ\phi,\psi)$ to parameters $\phi\in\Phi_2(G)$, according to the prescription described in [237, (3.1) and (3.3.1)]. The representations $\pi$ in the corresponding $L$-packets $\Pi_\phi$ then have the property that their $L$-functions $L(s,\pi,r)$ and $\varepsilon$-factors $\varepsilon(s,\pi,r,\psi)$ for any $r$ match those of $\phi$. This has been proven in cases where the representation theoretic functions have an independent meaning. In cases where they do not, it can be taken simply as their definition.

The Langlands classification was of course for all irreducible representations of a group $G(\mathbb{R})$, not just the discrete series. Langlands first extended his parametrization of $L$-packets of discrete series by parameters $\phi\in\Phi_2(G)$ to $L$-packets of square integrable representations (relative discrete series), the case that $G$ has a maximal torus that is anisotropic modulo the centre of $G$, at the bottom of p. 134 of [151]. This condition holds by definition for the Levi component of any *cuspidal* parabolic subgroup $P=N_PM$ of $G$. He then observed that a general parameter $\phi$ could be represented as the image in $\Phi(G)$ of a parameter $\phi_M\in\Phi_2(M)$ under the embedding ${}^LM\subset{}^LG$ attached to some cuspidal parabolic subgroup $P=N_PM$ of $G$. One can refine this embedding to a two stage embedding

$$\phi_M\to\phi_L\to\phi,\qquad \phi_M\in\Phi_2(M),$$

for parabolic subgroups $P=N_PM\subset Q=N_QL$ defined as follows. One first writes $\phi_M$ uniquely as a twist

$$\phi_M=\phi_{M,\mathrm{temp},\lambda},$$

---

[24] This condition is often called *relevance*. The reason we have not encountered it before is that we have usually been working with quasisplit groups, where the condition is automatic. We are also using the more streamlined notation of [119, §1] rather than the original formulation in §1–2 of [151]. As a matter of fact, these days one often formulates matters in terms of Vogan's *pure inner forms* [243], in which different inner forms are treated as components of the same object, and where the condition of relevance is only implicit.

where $\phi_{M,\text{temp}} \in \Phi_2(M)$ has bounded image in $\widehat{M}$, $\lambda \in \mathfrak{a}_M^*$ is a uniquely determined, real valued linear form on the real vector space $\mathfrak{a}_M$, and $\phi_{M,\text{temp},\lambda}$ is the parameter whose $L$-packet is the set of representations

$$\pi_{M,\lambda}(m) = \pi_M(m)e^{\lambda(H_M(m))}, \qquad m \in M(\mathbb{R}),$$

such that $\pi_M$ lies in the $L$-packet of $\phi_{M,\text{temp}}$. One then chooses $P$ such that $\lambda$ lies in the *closure* of the corresponding chamber $(\mathfrak{a}_P^*)^+$, and $Q \supset P$ so that $\lambda$ lies in the *open* chamber $(\mathfrak{a}_Q^*)^+$, regarded as a convex cone in the closure of $(\mathfrak{a}_Q^*)^+$. The set $\Phi(G)$ can then be identified with the set of $\widehat{G}$-orbits of triplets $(\phi_M, P, Q)$ of this form.

Langlands' goal in [151] was to define an explicit partition of $\Pi(G)$ into a disjoint union over $\phi \in \Phi(G)$ of finite $L$-packets $\Pi_\phi$. In particular, the $L$-packets for $G$ would also be indexed by $\widehat{G}$-orbits of triplets. His answer, which is not hard to describe, is an elegant reformulation of some of Harish-Chandra's fundamental results.

Given a triplet $(\phi_M, P, Q)$, we can define the $L$-packet $\Pi_{\phi_M} \subset \Pi_2(M)$ by Langlands' parametrization of Harish-Chandra's relative discrete series. For the next step, we form the parabolic subgroup $R = L \cap P$ of $L$ with Levi component $M$. For each $\pi_M \in \Pi_{\phi_M}$, we then take $\mathscr{I}_R^L(\pi_M)$, a representation of $L(\mathbb{R})$ parabolically induced from the representation $\pi_M$ in the relative discrete series of $M(\mathbb{R})$ that is unitary modulo the centre of $L(\mathbb{R})$. Let $\Pi_{\phi_L,\pi_M}$ be its set of irreducible constituents, a finite set of representations of $L(\mathbb{R})$ that are tempered modulo the centre. (The induced representation is generically irreducible, but when is it not, its irreducible constituents are of considerable interest.) Langlands defined the $L$-packet of $\phi_L$ in $\Pi(L)$ to be the set

$$\Pi_{\phi_L} = \bigcup_{\pi_M} \Pi_{\phi_L,\pi_M}, \qquad \pi_M \in \Pi_{\phi_M},$$

a union that was known to be disjoint. Finally, for any representation $\pi_L$ in $\Pi_{\phi_L}$, one can take the induced representation $\mathscr{I}_Q^G(\pi_L)$. This is a nontempered induced representation of $G(\mathbb{R})$, which in general is reducible. However, Langlands proved that it has a unique irreducible quotient $\pi(\pi_L)$. He then defined the $L$-packet of $\phi$ to be the set

$$\Pi_\phi = \{\pi = \pi(\pi_L) : \pi_L \in \Pi_{\phi_L}\}$$

of all these *Langlands quotients*, establishing at the same time that these representations were all disjoint.

This is the Langlands classification for real groups. When he introduced it in 1973, the irreducible constituents of the (essentially) tempered induced representations $\mathscr{I}_R^L(\pi_R)$ were not completely understood. However, they were classified soon afterwards by Knapp and Zuckerman [114], [115], thus providing in particular an explicit classification

$$\Pi_{\text{temp}}(G) = \coprod_{\phi \in \Phi_{\text{temp}}(G)} \Pi_\phi$$

of the irreducible *tempered* representations of $G(\mathbb{R})$ in terms of $L$-packets parameterized by the bounded Langlands parameters $\Phi_{\mathrm{temp}}(G)$. This is the special case of the general classification in which the groups $Q$ in the triplets $(\phi_M, P, Q)$ are all equal to $G$.

It was observed later in the 1970s that the Langlands classification would apply in principle also to $p$-adic groups. In this case, there is still no explicit classification of the relative discrete series $\Pi_2(M)$, or of the irreducible constituents of the induced tempered representations $\mathscr{I}_R^L(\pi_M)$ (although much is known, especially about this second question). However, the general classification, including the properties Langlands established for the quotients that bear his name, remains in force.

Langlands' classification established his conjecture of Local Functoriality for $F = \mathbb{R}$, as stated in Question 4 or 6 of [138], or as *Local Functoriality* stated in Section 5 here. In fact, it gives affirmative answers to all the questions in [138], insofar as they apply to the local field $F = \mathbb{R}$. As we noted in the special case of discrete series above, it assigns local $L$-functions and $\varepsilon$-factors to representations $\pi \in \Pi(G)$ by setting

$$L(s, \pi, r) = L(s, r \circ \phi), \qquad \pi \in \Pi_\phi,$$

and

$$\varepsilon(s, \pi, r, \psi) = \varepsilon(s, r \circ \phi, \psi), \qquad \pi \in \Pi_\phi,$$

for any $\phi \in \Phi(G)$ and $r \colon {}^L G \to \mathrm{GL}(n, \mathbb{C})$, and again for the functions on the right defined as in [237]. Moreover, the functorial correspondence $\pi' \to \pi$ of local representations is defined explicitly in terms of their $L$-packets by the transfer of parameters

$$\pi' \colon \Pi_{\phi'} \to \pi \in \Pi_\phi, \qquad \phi = \rho' \circ \phi', \phi' \in \Phi(G'),$$

for $\rho' \colon {}^L G' \to {}^L G$ as in the statement of Local Functoriality. This answers the local version of the original questions we posed at the beginning of this section, in the case that $F_v = \mathbb{R}$.

The Langlands classification for real groups suggests a spectral analogy with the theory of local endoscopy we have described. The Langlands parameters $\phi \in \Phi_{\mathrm{temp}}(G)$ ought to be analogues of (strongly) regular stable conjugacy classes $\Delta_{\mathrm{reg}}(G)$ over $\mathbb{R}$. Moreover the packets themselves ought to be analogues of the set of conjugacy classes $\gamma \in \Gamma_{\mathrm{reg}}(G)$ in a stable class $\delta$. However, there is more structure than this on the geometric side. The elements in a "geometric packet" are bijective with the explicit set $\mathscr{D}(T)$, where $T$ is the centralizer of a chosen base point $\delta$ in the stable class. We do have their dual analysis in terms of Tate–Nakayama duality, with its ties to endoscopic groups. But there is also the more refined structure given by the Langlands–Shelstad transfer factors, together with the associated transfer conjecture and its ultimate proof. Were there spectral analogues of any of these things?

The question was answered for real groups by Shelstad. We must not forget that she had first to introduce the archimedean transfer factors that became the inspiration for [165], and establish the associated transfer of functions. These remain basic links to the general theory, but they also have foundations in the work of Harish-

Chandra. We should discuss them briefly before we describe their spectral conse-
quences. In so doing, we need to step back in history, say to the year 1975. What
was available then were the basic ideas of Langlands on stable conjugacy and endo-
scopic groups, his new preprint on the classification for real groups, and of course,
the work of Harish-Chandra.

Shelstad's transfer factors are closely related to some curious factors in a refined
normalization of orbital integrals that had been forced on Harish-Chandra. He de-
fined a real group $G$ to be *acceptable* if the linear form (129) on the Lie algebra of
any maximal torus $T$ over $\mathbb{R}$ lifts to a quasicharacter $\xi_\rho$ on $T(\mathbb{C})$. This condition is
independent of the underlying system of positive roots on $T$, and holds whenever
$G_{\mathrm{der}}$ is simply connected. Harish-Chandra often worked with this assumption, with
the understanding that adjustments for the general case were easy to add separately.
(See for example [90, §8].) Under this condition, he normalized the orbital integrals
on $T_{\mathrm{reg}}(\mathbb{R})$ by setting

$$F_f(t) = \varepsilon_{\mathbb{R}}(t)\xi_\rho(t)\Delta(t)\mathrm{Orb}(t,f), \tag{131}$$

for $\Delta(t)$ as in (130) and $\varepsilon_{\mathbb{R}}(t) \in \{\pm 1\}$ the locally constant sign function

$$\mathrm{sign}\Big(\prod_{\alpha \in P_{\mathbb{R}}}(1 - \xi_\alpha(t^{-1}))\Big)$$

on $T_{\mathrm{reg}}(\mathbb{R})$ (with $P_{\mathbb{R}}$ being the set of positive real roots on $T$) [87, §22], [200, p. 5].
Since

$$|\varepsilon_{\mathbb{R}}(t)\xi_\rho(t)\Delta(t)| = |D(t)|^{\frac{1}{2}},$$

this does represent a refinement of our normalization $f_G(t)$ from (121). It was cho-
sen by Harish-Chandra to have the property that if $T(\mathbb{R})$ is compact, and $f$ is a ma-
trix coefficient of a discrete series representation, a function in his Schwartz space
$\mathscr{C}(G(\mathbb{R}))$ on $G(\mathbb{R})$ defined in [88], then $F_f(t)$ extends from $T_{\mathrm{reg}}(\mathbb{R})$ to a smooth
function on $T(\mathbb{R})$.

Before commenting on Shelstad's transfer factors, we should first include a cou-
ple of remarks that further illustrate the dual nature of orbital integrals and irre-
ducible characters. Orbital integrals satisfy differential equations

$$F_{zf}(t) = \gamma(z)F_f(t), \qquad t \in T_{\mathrm{reg}}(\mathbb{R}), z \in \mathscr{Z}_G, \tag{132}$$

where $z \to \gamma(z)$ is the Harish-Chandra homomorphism from $\mathscr{Z}_G$ to the algebra of
invariant differential operators on $T(\mathbb{R})$. Also, any left-invariant derivative of $F_f(t)$
remains bounded on $T_{\mathrm{reg}}(\mathbb{R})$ as $t$ approaches a singular hypersurface [87, Theo-
rem 3], and has an explicit formula for the jump as this function crosses the hy-
persurface [90, Theorem 1]. These are dual to the properties we have described for
irreducible characters (which are actually simpler when stated with the normaliza-
tion (127) replaced by the analogue of Harish-Chandra's refined normalization). In
particular, they lead to a boundary value problem for each function $F_f(t)$ on the clo-
sure of a connected component of $T_{\mathrm{reg}}(\mathbb{R})$. The only difference with what happens for irreducible characters is that the torus that shares the singular hypersurface with $T$ has anisotropic dimension one *less* than that of $T$. In other words, it is a Cayley transform of $T$, rather than other way around. It is consequently a decrease in the anisotropic dimension $d_a(T)$ that makes the associated orbital integrals simpler.

Shelstad's transfer factors for the real group $G$ serve as a bridge to the general Langlands–Shelstad transfer factors in [165]. These were defined in §3 of [165] as products

$$\Delta_G(\delta', \gamma) = \Delta_I(\delta', \gamma)\Delta_{II}(\delta', \gamma)\Delta_{III}(\delta', \gamma)\Delta_{IV}(\delta', \gamma), \qquad (133)$$

in which the third term comes with a further decomposition

$$\Delta_{III}(\delta', \gamma) = \Delta_{III_1}(\delta', \gamma)\Delta_{III_2}(\delta', \gamma) = \Delta_1(\delta', \gamma)\Delta_2(\delta', \gamma).$$

The Langlands–Shelstad transfer factors are complex and subtle, but they can be illuminated in their specialization to real groups, and the relations the latter bear to the quotients of Harish-Chandra's normalizing factors for $G$ and $G'$. In commenting briefly on this, we might as well assume that $G_{\mathrm{der}}$ is simply connected.

The term $\Delta_{IV}$ in (133) is the quotient of our original normalizing factor $|D(\gamma)|^{\frac{1}{2}} = |D_G(\gamma)|^{\frac{1}{2}}$ by $|D_{G'}(\delta')|^{\frac{1}{2}}$, as we agreed in the footnote 19. The term $\Delta_1 = \Delta_{III_1}$ is essentially the local form of the character $\kappa(\gamma)$ with which we began the original construction. The term $\Delta_2 = \Delta_{III_2}$ deals with the contribution $\bar{\xi}_\rho(\gamma)\xi_{\rho'}(\delta')^{-1}$ of the function $\xi_\rho(\gamma)$ in (131). It normalizes the choice of the $L$-isomorphism from ${}^L G'$ to $\mathscr{G}'$ that makes $\xi'$ an $L$-embedding of ${}^L G'$ into ${}^L G$. The term $\Delta_{II}$ in (133) addresses the contribution of the function $\varepsilon_{\mathbb{R}}(\gamma)\Delta(\gamma)$ in (131), or rather the contribution to $\Delta(\delta', \gamma)$ of its quotient by the factor $|D(\gamma)|^{\frac{1}{2}}$ that was put into $\Delta_{IV}$. Finally the term $\Delta_I$ is a sign, which is independent of $\gamma$, and compensates for various other noncanonical choices that had to go into the previous terms.

This is a necessarily superficial description of transfer factors. Shelstad actually worked with Harish-Chandra's later normalization $'F_f(\gamma)$ [90, §17] of orbital integrals, which makes sense for any $G$, and gives a more complete motivation for her work, but which is also a little more complicated to describe. From now on, let us just assume without comment that every endoscopic datum $G'$ attached to a given $G$ (over a local or global field $F$) is such that the group $\mathscr{G}'$ is $L$-isomorphic to the $L$-group ${}^L G'$. As we have noted above, this condition always holds if $G_{\mathrm{der}}$ is simply connected [147, Proposition 1].

Having introduced the transfer factors $\Delta(\delta', \gamma) = \Delta_G(\delta', \gamma)$ for the real group $G$ with endoscopic datum $(G', \mathscr{G}', s', \xi')$, Shelstad defined the transform

$$f'(\delta') = \sum_{\gamma \in \Gamma_{\mathrm{reg}}(G)} \Delta_G(\delta', \gamma) f_G(\gamma), \quad f \in C_c^\infty(G(\mathbb{R})), \, \delta' \in \Delta_{\mathrm{reg}}(G'),$$

that became the archimedean precursor of the general local transfer mapping (120). She then applied the harmonic analysis of Harish-Chandra systematically to its study, making use of the adjoint relation (128) between characters and orbital integrals, the differential equations (126) satisfied by characters and their analogues (132) for orbital integrals, and the accompanying boundary conditions in each case.

She used these techniques first to establish the real form of what became the general Langlands–Shelstad transfer conjecture. She then applied them to the spectral question raised above. Her results, which of course depend also on the Langlands classification for real groups, are as follows

We assume first that $G$ is quasisplit over $\mathbb{R}$. For every tempered Langlands parameter $\phi \in \Phi_{\text{temp}}(G)$, set

$$f^G(\phi) = \sum_{\pi \in \Pi_\phi} f_G(\pi), \qquad f \in C_c^\infty(G(\mathbb{R})), \tag{134}$$

for

$$f_G(\pi) = \Theta(\pi, f) = \text{tr}(\pi(f)).$$

The first (spectral) result of Shelstad is that the linear form $f \to f^G(\phi)$ is a stable distribution. It is called the *stable character* of $\pi$, and is clearly the spectral analogue of a stable orbital integral. Given the validity of the transfer conjecture for $G$ and $G'$ (for any $G$), this proves that the pairing $f'(\phi')$ is well defined for any tempered parameter $\phi' \in \Phi_{\text{temp}}(G')$ for the quasisplit group $G'$.

The second spectral result applies to any group $G$ over $\mathbb{R}$. It is an expansion

$$f'(\phi') = \sum_{\pi \in \Pi_\phi} \Delta_G(\phi', \pi) f_G(\pi), \tag{135}$$

for complex coefficients $\Delta_G(\phi', \pi)$ supported on the subset of $\pi \in \Pi_\phi$ in $\Pi_{\text{temp}}(G)$, in which $\phi = \xi' \circ \phi'$ is the image of $\phi'$ in $\Pi_{\text{temp}}(G)$. These coefficients are spectral analogues of Shelstad's (geometric) transfer factors, and can be called *spectral transfer factors*.

Shelstad's third spectral result makes these coefficients explicit, in a sense that depends on a chosen base point $\pi_1$ in the packet $\Pi_\phi$. Following Labesse and Langlands [127], [218], and similar definitions we have seen in our earlier sections, she defined $S_\phi = \text{Cent}(\phi(W_F), \widehat{G})$, the centralizer in $\widehat{G}$ of the image of a parameter $\phi \in \Pi_{\text{temp}}(G)$, and $\mathbf{S}_\phi = S_\phi / S_\phi^0 Z(\widehat{G})^\Gamma$, the group of connected components in $S_\phi$ modulo the Galois invariants in the centre of $\widehat{G}$. To state the result, we use the spectral analogue of the bijection (114), or rather the local spectral analogue of (114) for $F = \mathbb{R}$. Its inverse is a bijection

$$(G', \phi') \xrightarrow{\sim} (\phi, s), \tag{136}$$

from the isomorphism classes of pairs $(G', \phi')$, in which $G'$ is an endoscopic datum for $G$ and $\phi'$ lies in $\Pi_{\text{temp}}(G)$, onto isomorphism classes of pairs $(\phi, s)$, where $\phi$ is a parameter in $\Pi_{\text{temp}}(G)$ and $s$ is a semisimple element in $S_\phi$. Shelstad's third spectral result asserts that for any $(\phi, s)$, and for $\pi_1 \in \Pi_\phi$ fixed and $\pi \in \Pi_\phi$ arbitrary, the quotient of $\Delta(\phi', \pi)$ by $\Delta(\phi', \pi_1)$ depends only on the image $x$ of $s$ in $\mathbf{S}_\phi$, and that the resulting mapping

$$x \to \langle x, \pi | \pi_1 \rangle = \Delta(\phi', \pi) \Delta(\phi', \pi_1)^{-1} \tag{137}$$

is an *injection* from $\Pi_\phi$ to the group of characters on the abelian 2-group $\mathbf{S}_\phi$. This is parallel to what happens for the geometric transfer factors, where $(x, \pi, \pi_1)$ would be replaced by $(\kappa, \gamma, \gamma_1)$.

We have completed our brief review of Shelstad's work. Her results first appeared in the papers [219], [217], [220], [221], but they were expanded into a somewhat more expository treatment in the papers [222], [224], [223].

Suppose now that $G$ is a group over any local field $F$ of characteristic 0, which we take to be quasisplit. We have alluded to the conjectural local Langlands correspondence for $\Pi_{\text{temp}}(G)$ in past sections. It seems to have evolved with Langlands' ideas in the early 1970s, based on his experience with $\mathrm{GL}(2)$, and then the group $G = \mathrm{SL}(2)$. Its conjectural premises are close to the results of Shelstad for $F = \mathbb{R}$, but there is one significant difference. As we have noted earlier, the Weil group $W_F$ has to be replaced by the local Langlands group $L_F$. It thus remains equal to $W_F$ if $F$ is archimedean, but is taken to be the product $W_F \times \mathrm{SU}(2)$ if $F$ is nonarchimedean. This is to account for the Steinberg representation of $G(F)$, and more generally, the representations in $\Pi_2(G)$ whose matrix coefficients do not have compact support modulo $Z(G)$. Moreover, despite the fact that the conjectural assertions are otherwise similar to those of Shelstad for $F = \mathbb{R}$, any general proof seems unlikely to be the same. What is missing is Harish-Chandra's explicit classification of the discrete series. Without it, we do not know in general how to attach $L$-packets to parameters $\phi \in \Phi(G)$ and hence how to construct candidates $f^G(\phi)$ for the basic stable distributions.

Although perhaps already clear to the reader, it does not hurt to emphasize the two-fold nature of the final classification for real groups. There is the Langlands classification, based on Harish-Chandra's discrete series, and then there is Shelstad's endoscopic extension based on this and other work of Harish-Chandra. Given the apparent difficulty of an independent classification of general supercuspidal representations, we would hope to establish the local Langlands correspondence without having an explicit construction of the representations in the packets $\Pi_\phi$.

If $G = \mathrm{GL}(n)$, the local Langlands correspondence was proved by Harris and Taylor [92], Henniart [94] and Scholze [202]. It was established by global means, taken from the theory of Shimura varieties. Since stable conjugacy is the same as conjugacy in this case, there is no endoscopy. The local correspondence becomes a canonical bijection $\phi \to \pi_\phi$ from $\Phi_{\text{temp}}(G)$ to $\Pi_{\text{temp}}(G)$ (or between the larger sets $\Phi(G)$ and $\Pi(G)$). It is characterized by the requirement that the two kinds of local Rankin–Selberg $L$-functions and $\varepsilon$-factors attached to the representations

$$r: \mathrm{GL}(n_1, \mathbb{C}) \times \mathrm{GL}(n_2, \mathbb{C}) \to \mathrm{GL}(n_1 n_2, \mathbb{C})$$

coincide for parameters $(\phi_1, \phi_2)$ and representations $(\pi_1, \pi_2)$.

The local Langlands correspondence was established for a quasisplit symplectic or special orthogonal group over the nonarchimedean local field $F$ in Chapter 6 of [23]. The methods are again global, but here they come from the stabilization of the trace formula. For there is considerable endoscopy to contend with in this case. However, there is also a natural way to construct the basic stable distributions $f^G(\phi)$ attached to parameters $\phi \in \Pi_{\text{temp}}(G)$. They are twisted transfers from corresponding twisted invariant distributions on a general linear group $\text{GL}(N)$, relative to the standard outer automorphism $x \rightarrow {}^t x^{-1}$. This is because $G$ is a *twisted endoscopic group* for $\text{GL}(N)$, where $N = 2n+1$ if $G = \text{Sp}(2n)$ and $2n$ if $G$ equals either $\text{SO}(2n+1)$ or $\text{SO}(2n)$. The assertions are similar to those of Shelstad for real groups, but there is no need to state them at this point, since we will soon describe their generalizations that accompany the global classification. We note that similar methods were used by Mok [181] to establish the local Langlands correspondence for quasisplit unitary groups $G$.

The global endoscopic classification is deeper. There were hints of what form it should take in Langlands' paper [127] with Labesse. However, it goes beyond the global Question 7 of [138], whose local version Question 6 was a foundation for the local Langlands correspondence. This is because the global Weil group was known even for $\text{GL}(2)$ to provide only a sparse set of cuspidal automorphic representations. (The same might be said of the local Weil group $W_F$ in Question 6 for a nonarchimedian field $F$. But its extension was easily accommodated, either as the Weil–Deligne group [237], or equivalently, as the local Langlands group $L_F = W_F \times \text{SU}(2)$ we have taken here.) For a conjectural global classification, one would need to replace the global Weil group by the hypothetical global Langlands group $L_F$ discussed in the last section.

There is another ingredient that has also to be included in any global classification. It consists of the conjectures introduced in [18], [13], and discussed in the treatment of Shimura varieties in Section 8, for describing those automorphic representations in the discrete spectrum $L^2_{\text{disc}}$ that are not locally tempered. These represent the counterexamples to what would be the natural extension of Ramanujan's conjecture. For a quasisplit group $G$ over a global field $F$, the global parameters are the $L$-homomorphisms

$$\psi \colon L_F \times \text{SL}(2, \mathbb{C}) \rightarrow {}^L G \tag{138}$$

whose restriction to $L_F$ has bounded image in $\widehat{G}$, taken as usual up to $\widehat{G}$-conjugacy. As the essential global objects, these parameters force us to account also for their local analogues

$$\psi_v \colon L_{F_v} \times \text{SL}(2, \mathbb{C}) \rightarrow {}^L G_v.$$

For if one is to obtain a global classification of automorphic representations for $G$, as they occur in the discrete spectrum $L^2_{\text{disc}}$, one would at the same time have to establish a generalization of the local Langlands correspondence for these local parameters.

We assume now that $G$ is a quasisplit symplectic or special orthogonal group, this time over a global field $F$. This is the group for which the endoscopic classification of automorphic representations was established in [23]. The results, which take up the entire monograph, rest on the stabilization of the trace formula for $G$, as well as the twisted trace formula for $\text{GL}(N)$, and at some points even the twisted trace formula for $\text{SO}(2N)$. The stabilization of the ordinary (invariant) trace formula completed in [21], is of course a special case of the stabilization of the general twisted trace formula in the monographs [179], [180]. These in turn depend on other things, including the general Langlands–Shelstad transfer conjecture, its twisted analogue by Kottwitz and Shelstad [125], the fundamental lemma and its twisted analogue, and finally, a weighted fundamental lemma and its twisted analogue [247], [46], [47] required for terms in the complement of $I_{\text{ell,reg}}(f)$ in $I_{\text{geom}}(f)$. These results have now all been established, except possibly for the twisted, weighted fundamental lemma, which presumably would follow from the methods of [45] and [46].

In [23], an ad hoc substitute $\mathscr{L}_\psi$ for the hypothetical Langlands group $L_F$ was introduced in §1.4, as well as an ad hoc set of $L$-homomorphisms

$$\psi\colon \mathscr{L}_\psi \times \mathrm{SL}(2,\mathbb{C}) \to {}^L G,$$

in order to be able to formulate the global classification unconditionally. The global results were then stated in §1.5. The local results were stated in §2.3 and established in Chapter 7. Their proof depends on a special case, the actual local Langlands correspondence, established in Chapter 6 and described briefly above. It also relies on some properties of the intertwining operators between induced representations formulated in §2.3–2.4. These were established also in §7, apart from two references [A25] and [A26] from [23] that have still to be written, but which I expect will be completed soon.

As endoscopic identities for the localizations $\psi_v$ of global parameters, the local results are similar to the local Langlands correspondence for parameters $\phi_v$. In particular, they resemble our summary of Shelstad's results for real groups. However, there are also a couple of differences.

One is that spectral transfer factors here can be normalized in terms of the Whittaker models inherited from $\mathrm{GL}(N)$. This is actually a simplification. It allows us to take the analogue of the factor $\Delta(\phi', \pi_1)$ in (137) to be 1. The analogue of (137) becomes an injection

$$x_v \to \langle x_v, \pi_v \rangle = \Delta(\psi_v, \pi_v), \qquad \pi_v \in \Pi_{\psi_v}, \tag{139}$$

from a packet of representations $\Pi_{\psi_v}$ to the group of characters on the abelian 2-group $\mathbf{S}_{\psi_v} = S_{\psi_v}/S^0_{\psi_v} Z(\widehat{G})^{\Gamma_v}$. Another difference is a complication. This occurs in the construction of the packets $\Pi_{\psi_v}$. For a general local parameter $\psi_v \in \Psi(G_v)$, the packet becomes a finite set of *reducible* representations

$$\pi_v = \pi_v^1 \oplus \cdots \oplus \pi_v^k, \qquad \pi_v^i \in \Pi(G_v),$$

of $G(F_v)$. Moreover, while some of the irreducible constituents $\pi_v^i$ of these finite sums are tempered, others are not. However, they are all unitary. A second complication, minor but nonetheless interesting, concerns the coefficients $\langle x_v, \pi_v \rangle = \Delta(\psi_v', \pi_v)$ that ought to occur in the sum (135). What actually occurs are coefficients $\langle s_{\psi_v} x_v, \pi \rangle$, in which $s_{\psi_v}$ is the image in $\mathbf{S}_{\psi_v}$ of the element $s_\psi = \psi \begin{pmatrix} -1 & 0 \\ 0 & -1 \end{pmatrix}$. The analogues for $\psi_v$ of (134) and (135) then become

$$f_v^{G_v}(\psi_v) = \sum_{\pi_v \in \Pi_{\psi_v}} \langle s_{\psi_v}, \pi_v \rangle f_{G_v}(\pi_v) \tag{140}$$

and

$$f_v'(\psi_v') = \sum_{\pi_v \in \Pi_{\psi_v}} \langle s_{\psi_v} x_v, \pi_v \rangle f_{G_v}(\pi_v). \tag{141}$$

The local parameters $\psi_v \in \Psi(G_v)$ with trivial restrictions to the factor $\mathrm{SL}(2, \mathbb{C})$ are the usual Langlands parameters $\phi_v \in \Psi_{\mathrm{temp}}(G_v)$. The corresponding analogues of (140), (141) and (139) describe the local (endoscopic) Langlands correspondence for the group $G_v$. In this case, the representations $\pi_v$ in a packet $\Pi_{\phi_v}$ are irreducible and expected to be tempered, while the element $s_{\psi_v}$ in (140) and (141) equals 1. This is clearly close to Shelstad's endoscopic classification for real groups. However, as we have noted, a more complicated global proof is required because of the lack of any $p$-adic analogue of the Harish-Chandra classification of discrete series. We refer the reader also to the earlier, clearly written volume [1] by Adams, Barbasch and Vogan for *archimedean* parameters $\psi_v$, in which the conjectures were established for *any* real group, but without our defining property by twisted transfer to $\mathrm{GL}(N)$ that was needed for the global (and $p$-adic) theory.

Having stated the analogues of the Langlands correspondence for the localizations $\psi_v$ of the global parameters[25] $\psi \in \Psi(G)$, we can describe the global endoscopic classification for $G$. We write $\Psi_2(G)$ for the subset of global parameters $\psi \in \Psi(G)$ such that $S_\psi^0 = \{1\}$, which is to say that the image of $\psi$ does not lie in any proper parabolic subgroup of $^L G$. For any such $\psi$, we can then form the global packet

$$\Pi_\psi = \{\pi = \bigotimes_v \pi_v : \pi_v \in \Pi_{\psi_v}, \langle \cdot, \pi_v \rangle = 1 \text{ for } v \notin S\}$$

of representations of $G(\mathbb{A})$ that are unramified at almost every place $v$ of $F$. (The local construction is such that if the function $\langle x_v, \pi_v \rangle$ equals 1, the representation $\pi_v$ is irreducible and unramified.) For any $\pi \in \Pi_\psi$, the function

$$\langle x, \pi \rangle = \prod_v \langle x_v, \pi_v \rangle$$

is then defined. The main global result is Theorem 1.5.2 of [23]. It asserts that

$$L^2_{\mathrm{disc}}(G(F) \backslash G(\mathbb{A})) \cong \bigoplus_{\psi \in \Psi_2(G)} \left( \bigoplus_{\pi \in \Pi_\psi(\varepsilon_\psi)} m_\psi \pi \right). \tag{142}$$

Here $m_\psi$ equals 1 or 2, and

$$\varepsilon_\psi : \mathbf{S}_\psi \to \{\pm 1\}$$

---

[25] Bear in mind that $\psi$ represents one of our ad hoc parameters on a product $\mathscr{L}_\psi \times \mathrm{SL}(2, \mathbb{C})$. But the group $\mathscr{L}_\psi$ was constructed in [23, 1.4] to contain the local Langlands group $L_{F_v}$, so its restriction to $\psi_v$ is a homomorphism from $L_{F_v} \times \mathrm{SL}(2, \mathbb{C})$ to $^L G_v$.

is a linear character defined explicitly in terms of global symplectic $\varepsilon$-factors, while $\Pi_\psi(\varepsilon_\psi)$ is the subset of representations in the global packet $\Pi_\psi$ such that the character $\langle \cdot, \pi \rangle$ on $\mathbf{S}_\psi$ equals $\varepsilon_\psi$.

This completes our very brief summary of the endoscopic classification of representations of the quasisplit symplectic or special orthogonal group $G$. For more information, a reader could begin with the introduction in [23], and then go to Chapter 1 and perhaps the first few sections of Chapter 2 and Chapter 4. Similar results have been established for quasisplit unitary groups by Mok [181]. We recall also that some of the geometric implications of these matters were discussed in the last two sections.

I should say that the description I have given here is not quite correct as stated. There are technical adjustments required for the case that $G$ equals the group $\mathrm{SO}(2n)$. Suppose for example that $G$ is split over a completion $F_v$, and that

$$\phi_v^1, \phi_v^2 \colon W_{F_v} \to \mathrm{SO}(2n, \mathbb{C}) = \widehat{G}$$

is a pair of distinct, irreducible, special orthogonal representations of $W_{F_v}$ that are conjugate under the action of $\mathrm{O}(2n, \mathbb{C})$. The local results above imply that $\{\phi_v^1, \phi_v^2\}$ corresponds to a pair $\{\pi_1, \pi_2\}$ of irreducible representations in $\Pi_2(G_v)$. What they do not give is the actual bijection (from the two possible choices) between these two sets of order 2 implicit in the local Langlands correspondence. Such ambiguity has to be built into the assertions of [23]. (Their global manifestations are closely related to the multiplicity $m_\psi \in \{1, 2\}$ in (142).) I have omitted them deliberately in my summary above, in the hopes of better conveying the essence of what is going on. A reader can easily restore them from the definitions of Chapter 1 of [23]. In any case, this section has gone on long enough!

# 11 Beyond Endoscopy

The theory of endoscopy we have just discussed has the potential to establish interesting cases of functoriality. They arise from the elliptic endoscopic data $(G', s', \mathscr{G}', \xi')$ attached to a quasisplit group $G$, and the $L$-embedding

$$\rho' \colon {}^L G' \to {}^L G$$

obtained from $\xi'$ and the choice of an $L$-isomorphism from ${}^L G'$ to the group $\mathscr{G}'$.

For example, suppose that $G$ is a quasisplit classical group, and that $G'$ is a product of two quasisplit classical groups $G_1 \times G_2$ for which the canonical direct product ${}^L G_1 \times {}^L G_2$ is a maximal $L$-subgroup of ${}^L G$. Functoriality for this case follows from the results of [23] and [181]. The proof depends on the stabilization of the trace formula of $G$ [21], which among other things yields the stable trace formula for $G$. It also depends on the twisted stabilization of $\mathrm{GL}(N)$, with respect to the standard outer automorphism. We have not discussed twisted endoscopy very much, but the formal definitions are similar to those of ordinary endoscopy. The twisted stabilization of the ordinary (twisted) trace formula was established in complete generality in the two volumes [179], [180], apart from the provisio mentioned earlier on the twisted, weighted fundamental lemma. The classical groups $G$ are themselves twisted endoscopic groups for general linear groups. A consequence of this, functoriality for the natural embedding of ${}^L G$ into $\mathrm{GL}(N, \mathbb{C})$, was also a part of the results in [23] and [181]. We have discussed these matters already, at the end of the last section.

There are certainly other interesting cases of functoriality that come from endoscopy, but most of these are presently out of reach. And at any rate, the examples of functoriality attached in one way or another to endoscopy are pretty sparse compared to the general case.

Beyond Endoscopy is a strategy proposed by Langlands around 2000 for attacking the general Principle of Functoriality. The ideas represent a departure from anything that has gone before. They do involve a comparison of trace formulas, stable trace formulas in fact. However, they entail something else as well, the automorphic $L$-functions

$$L(s, \pi, r), \qquad \pi \in \Pi_2(G), r \colon {}^L G \to \mathrm{GL}(n, \mathbb{C}),$$

attached to $G$. Langlands' proposal was to refine the stable trace formula for $G$ by inserting a supplementary factor into the stable multiplicity of $\pi$ on the spectral side, namely the order of the pole of $L(s, \pi, r)$ at $s = 1$. For fixed $r$, this would vary with $\pi$, or rather the global packet of $\pi$, and according to what is expected about functoriality, would give information on the "functorial lineage" of $\pi$. We refer the reader to Section 1 of [155] for a basic introduction to the ideas, together with a number of critical examples.

To my view, this strategy of Langlands is fundamental, and of the greatest significance. It is also deep and difficult, much more so even than the theory of endoscopy.

Despite the fact that Langlands' proposal is now twenty years old, its study is still in the very early stages.

We shall generally assume for the rest of the section that $G$ is the general linear group $\mathrm{GL}(n+1)$ over the field $F = \mathbb{Q}$. The stable trace formula for $G$ is then the same as the invariant trace formula (103). The fundamental problem is to understand how the representations $\pi \in \Pi_2(G)$ in the discrete spectrum are related to functoriality, or more precisely, how they might arise as functorial images of triplets $(G', \pi', \rho')$, for $L$-homomorphisms

$$\rho' \colon {}^L G' \to {}^L G = \mathrm{GL}(n+1, \mathbb{C}).$$

It is best to discard the nontempered representations $\pi$, which one might in any case expect to be able to treat by induction. We therefore restrict our consideration to the subset $\Pi_1(G) = \Pi_{\mathrm{cusp},2}(G)$ of cuspidal automorphic representations in $\Pi_2(G)$. Langlands noted that the functorial preimages $(G', \pi', \rho')$ of $\pi \in \Pi_1(G)$ should be closely related to the poles at $s = 1$ of the $L$-functions $L(s, \pi, r)$ attached to $r$. This of course presupposes the meromorphic continuation of $L(s, \pi, r)$ to a half space $\mathrm{Re}(s) > 1 - \varepsilon, \varepsilon > 0$, something that is not known in general.

Let us assume for a moment that functoriality holds, say for all general linear groups $G = \mathrm{GL}(n+1)$ and $G' = \mathrm{GL}(m+1)$. This implies the meromorphic continuation of the $L$-functions $L(s, \pi, r)$, and allows us to set

$$m_\pi(r) = -\mathrm{ord}_{s=1} L(s, \pi, r) = \mathrm{res}_{s=1} \left( -\frac{d}{ds} \log L(s, \pi, r) \right),$$

the order of the pole of $L(s, \pi, r)$ at $s = 1$. For any $r$, we are then free to define

$$I^r_{\mathrm{cusp}}(f) = \sum_{\pi \in \Pi_1(G)} m_\pi(r) \cdot \mathrm{mult}(\pi) \cdot \Theta(\pi, f), \tag{143}$$

the contribution of the representations $\pi \in \Pi_1(G)$ to the primary spectral part of (107), but weighted by these integers. For example, if $r$ equals the trivial representation $1_G$ of ${}^L G$, $L(s, \pi, r)$ is just the completed Riemann zeta function

$$L(s, 1) = \pi^{-s/2} \Gamma(s/2) \zeta(s), \qquad \pi = 3.1416...,$$

for each representation $\pi$, which has a simple pole at $s = 1$. In this case $I^r_{\mathrm{cusp}}(f)$ is itself just the cuspidal part

$$I_{\mathrm{cusp}}(f) = \sum_{\pi \in \Pi_1(G)} \mathrm{mult}(\pi) \Theta(\pi, f) \tag{144}$$

of (107). Langlands suggested the possibility of finding a geometric expansion for $I^r_{\mathrm{cusp}}(f)$ for any $r$, thereby giving a refinement of the invariant trace formula (103).

The idea was thus to try to find a trace formula whose spectral side is the modified cuspidal expansion (143). To see how this might be possible, we note that $m_\pi(r)$ is the residue at $s = 1$ of

$$-\frac{d}{ds}(\log(L^S(s, \pi, r)))$$

$$=-\frac{d}{ds}(\log(\prod_{p \notin S} \det(1 - r(c(\pi_p))p^{-s})^{-1}))$$

$$=\sum_{p \notin S} \frac{d}{ds}(\log(\det(1 - r(c(\pi_p)))p^{-s}))$$

$$=\sum_{p \notin S} \sum_{k=1}^{\infty} \log(p) \operatorname{tr}(r(c(\pi_p))^k) p^{-ks}.$$

We can discard the terms with $k \geq 2$ in this last sum, since the function they define would be holomorphic at $s = 1$. It follows that

$$m_\pi(r) = \operatorname{res}_{s=1} \left( \sum_{p \notin S} \log(p) \operatorname{tr}(r(c(\pi_p))) p^{-s} \right).$$

A familiar application of the Wiener–Ikehara theorem[26] would then give a formula

$$m_\pi(r) = \lim_{N \to \infty} \left( |S_N|^{-1} \sum_{p \notin S_N} \log(p) \operatorname{tr}(r(c(\pi_p))) \right),$$

where

$$S_N = \{p \notin S : p \leq N\}.$$

(See [208, p. I-29].)

To exploit this last formula, one can write the test function $f \in C_c^\infty(G(\mathbb{A}))$ as a product $f_S \cdot \mathbf{1}^S$, for $f_S \in C_c^\infty(G(F_S))$ and $\mathbf{1}^S$ the characteristic function of the compact open subgroup $G(\widehat{\mathbb{Z}}^S) = \prod_{p \notin S} G(\mathbb{Z}_p)$ in $G(\mathbb{A}^S)$. One can then enrich this function by adding a factor at any $p \notin S$. We set

$$f_p^r(x) = f(x) \cdot h_p^r(x_p), \qquad x \in G(\mathbb{A}),$$

where $x_p$ is the component of $x$ in $G(\mathbb{Q}_p)$, and $h_p^r$ is the function in the unramified Hecke algebra $\mathscr{H}(G(\mathbb{Z}_p) \backslash G(\mathbb{Q}_p)/G(\mathbb{Z}_p))$ on $G(\mathbb{Q}_p)$ whose Satake transform equals

$$\widehat{h_p^r}(c_p) = \operatorname{tr}(r(c_p)),$$

for any semisimple conjugacy class $c_p$ in $\widehat{G} = \operatorname{GL}(n+1, \mathbb{C})$, which is to say that

$$\operatorname{tr}(\pi_p(h_p^r)) = \operatorname{tr}(r(c(\pi_p))),$$

---

[26] We need to assume here that $\pi$ is locally tempered, and hence that the generalized Ramanujan conjecture holds for $G$. This implies that the $L$-series $L(s, \pi, r)$ converges *absolutely* for $\operatorname{Re}(s) > 1$, and therefore satisfies conditions (a) and (b) of the theorem stated in [208]. This is really part of functoriality, on which we are basing the present motivational argument. (The assumption should really have been explicit in the discussion of these matters in [24, §2]).

for any unramified representation $\pi_p$ of $G(\mathbb{Q}_p)$. Then

$$\Theta(\pi, f_p^r) = \operatorname{tr}(\pi(f))\operatorname{tr}(\pi_p(h_p^r)) = \Theta(\pi, f)\operatorname{tr}(r(c(\pi_p))),$$

for any $\pi \in \Pi_1(G)$ such that $\pi_p$ is unramified. Combining this with the last formula for $m_\pi(r)$ and the definition (143) of $I_{\text{cusp}}^r(f)$, one sees that $I_{\text{cusp}}^r(f)$ would equal

$$\sum_{\pi \in \Pi_1(G)} \Big( \lim_{N \to \infty} |S_N|^{-1} \sum_{p \notin S_N} \log(p) \operatorname{tr}(\pi_p(h_p^r)) \Big) \operatorname{mult}(\pi) \Theta(\pi, f)$$
$$= \lim_{N \to \infty} (|S_N|^{-1} \sum_{p \notin S_N} \log(p) \Big( \sum_{\pi \in \Pi_1(G)} \operatorname{mult}(\pi) \Theta(\pi, f_p^n) \Big).$$

In other words,

$$I_{\text{cusp}}^r(f) = \lim_{N \to \infty} \Big( |S_N|^{-1} \sum_{p \notin S_N} \log(p) I_{\text{cusp}}(f_p^r) \Big). \tag{145}$$

On the other hand, we can formally rewrite the invariant trace formula (103) as

$$I_{\text{geom,temp}}(f) = I_{\text{cusp}}(f), \tag{146}$$

where

$$I_{\text{geom,temp}}(f) = I_{\text{geom}}(f) - (I_{\text{spec}}(f) - I_{\text{cusp}}(f)). \tag{147}$$

The subscript *temp* indicates that the distribution should be locally tempered. This is because the representations in $\Pi_1(G)$ would satisfy the analogue of the Ramanujan conjecture, according to our assumption and the argument of Langlands sketched at the end of [138], as we observed in footnote 26. The Dirichlet series for $L(s, \pi, r)$ would then converge *absolutely* for Re$s > 1$, an implicit condition for the Wiener–Ikehara theorem we applied above. The result would then be an $r$-trace formula

$$I_{\text{geom,temp}}^r(f) = I_{\text{cusp}}^r(f), \tag{148}$$

for any $r$, where

$$I_{\text{geom,temp}}^r(f) = \lim_{N \to \infty} \Big( |S_N|^{-1} \sum_{p \notin S_N} \log(p) I_{\text{geom,temp}}(f_p^r) \Big). \tag{149}$$

This would provide a large family of refined trace formulas from which one might try to deduce functoriality.

However, this was all under the assumption of meromorphic continuation for the $L$-function $L(s, \pi, r)$. The way to establish this, according to Langlands' conjectures, is to apply functoriality to the homomorphisms

$$\rho' = r : \widehat{G} = \operatorname{GL}(n+1, \mathbb{C}) \to \operatorname{GL}(N, \mathbb{C}).$$

We cannot very well assume functoriality when it is what we ultimately want to prove. Langlands' idea was to use the spectral expression (143) initially only to motivate the proposed limit (149). For it does tell us that the limit (149) ought to exist. Langlands' hope is that it will eventually be possible to prove independently that the limit does exist, and to express it in terms of a reasonably explicit geometric expansion. One could then work on trying to establish a spectral expansion for the limit akin to (143).

It is clearly an enormous problem. The first major step would be the formidable task of finding a geometric expansion for the difference $I_{\text{geom,temp}}(f)$ in (147). To see what this might entail, we shall consider the formal approximation (107) of the full trace formula (103) given by

$$I_{\text{ell,reg}}(f) = \sum_{\gamma \in \Gamma_{\text{ell,reg}}(G)} \text{vol}(\gamma) \text{Orb}(\gamma, f)$$

$$\sim \quad I_2(f) = \sum_{\pi \in \Pi_2(G)} \text{mult}(\pi) \Theta(\pi, f).$$

Its analogue for the cuspidal trace formula (146) is the approximation

$$I_{\text{ell,reg,temp}}(f) \sim I_{\text{cusp}}(f), \tag{150}$$

where

$$I_{\text{ell,reg,temp}}(f) = I_{\text{ell,reg}}(f) - \sum_{\pi \notin \Pi_1(G)} \text{mult}(\pi) \Theta(\pi, f), \tag{151}$$

the sum being over the complement of $\Pi_1(G)$ in $\Pi_2(G)$. We are retaining the subscript *temp* in the distribution $I_{\text{ell,reg,temp}}(f)$ to emphasize that it is supposed to be an approximation of $I_{\text{geom,temp}}(f)$. The difference between the two distributions is the sum of the supplementary geometric terms in the full trace formula (103) minus the sum of the supplementary spectral terms. We have generally avoided discussing these auxiliary terms (except for GL(2) in Sections 6 and 7), but if $G$ is not equal to GL(2), there are some locally nontempered distributions $I_M(\pi, f)$ among the supplementary spectral terms. The distribution $I_{\text{ell,reg,temp}}(f)$ therefore cannot be locally tempered. It is to be regarded as we have said, simply as a formal approximation of the distribution $I_{\text{geom,temp}}(f)$ we do expect to be tempered.

The point to be made is that the regular elliptic part $I_{\text{ell,reg}}(f)$ of the trace formulas conceals some secrets. It is very familiar, having been known (if not always in adelic form) since Selberg first introduced his trace formula for compact quotient. But despite the fact that the individual orbital integrals in its summands are locally tempered (a fact proved by Harish-Chandra for real groups in 1966 [88], and for $p$-adic groups somewhat later), their sum is not. The main obstruction is simple enough, the sum over $\pi$ in (151), but this is its spectral form. We are thus asking for an explicit, and we hope reasonably simple, geometric expansion for the difference (151). As far as I know, this natural question was never considered (except perhaps for $G = \text{GL}(2)$) before 2000. It is the essential part of what we called the first major step above. It will demand much effort, supported no doubt by experience gained from experiments in special cases.

Langlands discussed these and other ideas, with various examples, in Part I of his foundational article [155]. In Part II he examined various terms in the trace formula for GL(2). Part III of [155] is devoted to actual experiments, using computer calculations to estimate some of the quantities from Part II. Part IV contains among other things a few remarks on general groups. A significant part of the article is devoted to a topic that will have to be understood after the initial questions have been answered. It is the supplementary geometric part of the trace formula for GL(2), represented by the noninvariant terms (iv) and (v) on p. 516–517 of [103]. The contributions of these terms will be locally tempered for GL(2), as will the noninvariant contributions of the spectral terms (vi), (vii) and (viii) from [103]. It is a simple enough example for one to be able to ask what influence all of these terms might have on the limit (149) in the case of GL(2). Langlands studied them in some detail [155, §2.4, §4.3 and Appendix C], and found some interesting cancellations among their contributions to the limit.

Langlands' initial article on Beyond Endoscopy was actually the unpublished precursor [154] of [155]. This paper contains some of his first ideas on the new program, with some comparisons to the developing theory of endoscopy. It represents an informal introduction to the main article [155]. The successor to [155] was the report [156], in which Langlands reviewed some of the constructions and calculations from [155]. He then described how they could be formulated in the case of function fields, that is, global fields of positive characteristic.

The next paper was Langlands' 2010 article [74] with Frenkel and Ngô. It contains three critical suggestions for the analysis of $I_{\mathrm{ell,reg,temp}}(f)$. The paper encompasses both number fields and function fields, opening the possibility of extending to number fields techniques that had been exploited by Ngô in his recently completed proof of the fundamental lemma for Lie algebras of positive characteristic. One of these became the first of the three suggestions. It was to parameterize the classes $\gamma \in \Gamma_{\mathrm{ell,reg}}(G)$ that index the summands in $I_{\mathrm{ell,reg}}(f)$ by points in what the authors called the base of the Steinberg–Hitchin fibration. For the case of $G = \mathrm{GL}(n+1)$ here, this amounts to a parameterization of the semisimple classes $\gamma$ by their characteristic polynomials. It represents a significant change of perspective despite its simplicity.

The base of the Steinberg–Hitchin fibration for $G = \mathrm{GL}(n+1)$ is the product

$$\mathscr{A}(n) = \mathscr{B}(n) \times \mathbb{G}_m,$$

of affine $n$-space $\mathscr{B}(n)$ with the multiplicative group $\mathbb{G}_m = \mathrm{GL}(1)$. The characteristic polynomial

$$\det(\lambda I - \gamma) = \lambda^{n+1} - a_1 \lambda^n + \cdots + (-1)^n a_n \lambda + (-1)^{n+1} a_{n+1} = p_a(\lambda)$$

of a class $\gamma \in \Gamma_{\mathrm{ell,reg}}(G)$ has rational coefficients. It gives a bijection $\gamma \to a$ from $\Gamma_{\mathrm{ell,reg}}(G)$ onto the subset $\mathscr{A}_{\mathrm{irred}}(n, \mathbb{Q})$ of elements

$$a = (a_1, \ldots, a_n, a_{n+1})$$

in $\mathscr{A}(n,\mathbb{Q})$ such that $p_a(\lambda)$ is irreducible over $\mathbb{Q}$. The suggestion was thus to rewrite the sum over $\gamma \in \Gamma_{\text{ell,reg}}(G)$ in $I_{\text{reg,ell}}(G)$ as a sum over elements

$$a = (b, a_{n+1}), \qquad b \in \mathbb{Q}^n, a_{n+1} \in \mathbb{Q}^*,$$

in $\mathscr{A}_{\text{irred}}(n,\mathbb{Q})$.

The second suggestion from [74] was to try to apply the Poisson summation formula to the sum over $b \in \mathbb{Q}^n$. This is certainly not immediately possible, for a variety of reasons. The problem would be to modify the resulting expression for $I_{\text{ell,reg}}(f)$ in such a way that Poisson summation could be applied to the rearranged sum over $b$ in $\mathbb{Q}^n$ (regarded of course as a lattice in $\mathbb{A}^n$). The third suggestion was to try then to account for the nontempered representations $\pi$ in (151) directly in terms of the summands attached to the dual variables $\xi \in \mathbb{Q}^n$. The authors conjectured in particular that the contribution of the trivial one-dimensional representation $\pi_1$ of $G(\mathbb{A}^n)$ was contained in the dual summand with $\xi = 0$.

With these ideas, the hidden structure in $I_{\text{ell,reg}}(f)$ becomes more compelling. The summands on the right hand side of (151) are parameterized by among other things, unipotent conjugacy classes in $\hat{G}$. The suggestions in [74] are that their contributions to $I_{\text{ell,reg}}(f)$ might have an unexpectedly explicit form (See [26] for a conjectural description in the case of general linear groups.). The proposed phenomena become all the more intriguing for groups $G$ other than $\text{GL}(n+1)$. The authors discussed their suggestions in general terms in the first three sections of [74]. In the remaining sections, they offered some evidence.

They devoted some time to describing how best to normalize the invariant measures on the various spaces in play, the most important being the additive adelic space $\mathscr{B}(n,\mathbb{A}) \cong \mathbb{A}^n$ that would be the domain of the Fourier transform from Poisson summation. The next step was essentially to construct a function on this space from the local orbital integrals of $f$, for which the global orbital integrals in $I_{\text{ell,reg}}(f)$ represent the values on the lattice $\mathscr{B}(n,\mathbb{Q}) = \mathbb{Q}^n$. The original function itself was quite unsuitable. However, the authors used an idea of J. Getz in §4 of [74] to truncate it in a certain way so as to make it amenable to Poisson summation. They then showed that in the resulting sum of Fourier transforms over $\xi$ in the dual lattice, the contribution of the 1-dimensional representations, the most highly nontempered representations from the right hand side of (151), is indeed contained in the summand with $\xi = 0$. The sum over these Fourier transforms with $\xi \neq 0$ therefore removes these representations from the difference (151). The authors in fact showed that it is asymptotically smaller in a natural sense than the sum of the 1-dimensional representations. Estimates of this sort are exactly what are being sought, even if in this case the results are too weak to apply. However, they constitute striking evidence for the general proposals in [74], especially since they apply to any (quasisplit) group $G$.

The heart of Beyond Endoscopy would be a general comparison of trace formulas. However, this is really something for the future. For in general it could come only after the three suggestions from [74] have been successfully carried out, and indeed, only after the stable generalizations

$$S^r_{\text{geom,temp}}(f) = \lim_{N \to \infty} \left( |S_N|^{-1} \sum_{p \notin S_N} \log(p) S_{\text{geom,temp}}(f^r_p) \right)$$

and

$$S^r_{\text{geom,temp}}(f) = S^r_{\text{cusp}}(f) \tag{152}$$

of the limit (149) and the $r$-trace formula (148) have been established for any $G$. The right hand side of (152) would have to be defined as the general analogue

$$S^r_{\text{cusp}}(f) = \lim_{N \to \infty} \left( |S_N|^{-1} \sum_{p \notin S_N} \log(p) S_{\text{cusp}}(f^r_p) \right)$$

of (145), whose existence would be a consequence of the existence of the limit $S^r_{\text{geom,temp}}(f)$. (The stable analogue $S_{\text{cusp}}(f)$ of $I_{\text{cusp}}(f)$ in this limit would be a sum as in (144), but over global $L$-packets *all* of whose constituents are cuspidal, or equivalently, all of whose constituents are expected to be locally tempered.) The goal of the comparison would be to provide information about the stable distributions $S^r_{\text{cusp}}(f)$ akin to the stable analogue of the right hand side of (143). We reiterate that a priori, we would know nothing about the right hand side of (143). The role of this formula was only to motivate the limit (145) with which the heart of the argument would begin. Once one has versions of the formulas (143) for various $r$, and comparisons of them with analogues for other groups $G'$, one could finally begin what would presumably be the last stage of the argument, the search for confirmation of functoriality.

Our experience with endoscopy can inform what might be the new comparison. The "Beyond Endoscopic" comparison suggested by Langlands in [155] will have some features in common with endoscopic comparison (better known as stabilization), and some features that are quite different. In general, we could imagine a "Beyond Endoscopic datum" attached to a quasisplit group $G$ over a number field $F$ simply as a pair $(G', \rho')$, for a reductive group $G'$ over $F$, and an $L$-homomorphism

$$\rho' \colon {}^L G' \to {}^L G$$

whose restriction to $\widehat{G}'$ is an embedding. As stated, this is much broader than an endoscopic datum, to the extent that $\widehat{G}'$ is not required to be the connected centralizer in $\widehat{G}$ of a semisimple element $s' \in \widehat{G}$. It is also incomplete, in the sense that $(G', \rho')$ should at least be replaced by a triplet $(G', \mathscr{G}', \xi')$ that satisfies the same conditions as an endoscopic datum $(G', s', \mathscr{G}', \xi')$ (without the point $s'$), while $\rho'$ could represent the choice of an $L$-isomorphism from $\mathscr{G}'$ to $G'$. Langlands proposed a transfer

$$f = \prod_v f_v \to \prod_v f'_v = f', \qquad f \in C^\infty_c(G(\mathbb{A})),$$

of functions from $C^\infty_c(G(\mathbb{A}))$ to functions in $C^\infty_c(G'(\mathbb{A}))$, to be known as stable transfer since it should depend only on the stable orbital integrals of *both* $f$ and $f'$. This would be quite different from endoscopic transfer.

If we assume the local Langlands correspondence for both $G$ and $G'$ at the places $v$ of $F$, something that could well be available before this stage of development of Beyond Endoscopy, stable transfer would be easy to define. For any function $f_v \in C_c^{\infty}(G_v)$, it would be enough to specify the value $f_v'(\phi_v') = f_v^{G'}(\phi_v')$ of any tempered local Langlands parameter $\phi_v' \in \Phi(G_v')$. We would do so by setting

$$f_v'(\phi_v') = f_v^G(\rho_v' \circ \phi_v'),\tag{153}$$

for the localization $\rho_v'$ of $\rho'$. The global transfer $f'$ at a function $\prod_v f_v$ in $C_c^{\infty}(G(\mathbb{A}_F))$ would then simply be defined as the product $\prod_v f_v'$ of local transfers. With this definition, the real problem would then be to determine the orbital integrals of $f_v'$ in terms of those of $f_v$. In other words, one would like to determine the value $f'(a_v')$ at an $F_v$-valued point $a_v' \in \mathscr{A}(G_v', F_v)$ of the Steinberg–Hitchin base for $G_v'$ in terms of the values $f_v^G(a_v)$ of $f_v$ at $F_v$-valued points $a_v \in \mathscr{A}(G_v, F_v)$ of the Steinberg–Hitchin base for $G_v$. According to some version of the Schwartz kernel theorem, we would expect to be able to write

$$f_v'(a_v') = \int_{\mathscr{A}(G, F_v)} \Delta(a_v', a_v) f_v^G(a_v)\, da_v,$$

for some integral kernel $\Delta(a_v', a_v)$. The integral is to be understood as the pairing of the "Schwartz function" $f_v^G(\cdot)$ on $\mathscr{A}(G_v, F_v)$ with the "tempered distribution" $\Delta(a_v', \cdot)$, notions that would have to be suitably interpreted. The qualitative difference between this and the Langlands–Shelstad transfer factor $\Delta(\delta_v', \gamma_v)$ is manifest.

With a theory of Beyond Endoscopic transfer, one could then finally consider comparing stable trace formulas. There are many possibilities, and it is hard to know in advance which of these might be best. For example, one might try to compare the $r$-trace formula of $G$ with the stable trace formula of $G'$, or more likely, a linear combination over Beyond Endoscopic data $G'$ of stable trace formulas of $G'$. One proposed model for such a comparison, founded on the structure of Langlands' automorphic Galois group $L_F$ proposed in Section 9, was described in [24, §2] and [25, §4], where it was called the *primitization* of the $r$-trace formula. However, what is needed now is more understanding of the many fundamental questions that would have to be answered first, if not in general then at least for natural examples.

This last discussion is intended as background for [158], the article on singularities and transfer Langlands wrote as a continuation of his fundamental paper with Frenkel and Ngô. Our remarks on stable transfer and comparison of trace formulas are taken largely from the early pages of [158]. However, the ostensible purpose of the paper was to explore them in detail for the group $G = \mathrm{SL}(2)$, and the dihedral groups $G'$ of elements of norm 1 in quadratic extensions $E$ of $F$, together with the degenerate case of $G' = \mathrm{GL}(1)$. Langlands proved the existence of stable transfer $f_v \to f_v'$ in each of these cases, for localizations $F_v$ of $F$ with $q_v \neq 2$. To do so, he relied on general properties of stable characters for $\mathrm{SL}(2)$ in [127], and explicit formulas for irreducible characters on $\mathrm{SL}(2, F_v)$, by Harish-Chandra for archimedean $v$ and by Sally and Shalika [197][27] for nonarchimedean $v$. The proof was given in Section 1 of [158], the first half of the paper, by answering Questions A and B from Section 1.1. We note that there is a recent extension to $\mathrm{SL}(n)$ by Daniel Johnstone [109], using some interesting new methods.

In the second half of the paper, Langlands studied two functions

$$\theta_f(a) = \left( \prod_{v \in S} f_v^G(a_v) \right) \mathbf{1}^S(a^S)$$

and

$$\phi_f(a') = \left( \prod_{v \in S} f_v'(a_v') \right) \mathbf{1}^S((a')^S),$$

for regular adelic variables $a, a' \in \mathbb{A}_F$ in $\mathscr{A}(G)$, and a large (variable) finite set $S$ of places. The second function is really a composite of the transfers from $G$ to all of the dihedral groups $G'$, with their local Steinberg–Hitchin bases $\mathscr{A}(G_v')$ embedded in $\mathscr{A}(G_v)$. The regular set in $\mathscr{A}(G_v)$ becomes a disjoint union of (the regular points in) their domains, and $f_v'(a_v')$ is then the value of the transfer mapping for the group $G_v'$ attached to the domain that contains $a_v'$. Langlands was interested in the singularities of these functions. For archimedean $v$, these were obtained from special cases of the general results of Harish-Chandra, mentioned for $\mathrm{GL}(2)$ as Harish-Chandra families in Section 7 of this volume, and reviewed briefly for general groups in the last section. For nonarchimedean $v$, the singularities are more complex for general groups. However, their qualitative behavior is well understood in terms of the Shalika germs introduced in [216].

Langlands' interest was in the specialization to $G = \mathrm{SL}(2)$ of the results established for general (simply connected) $G$ in the second half of [74]. These consist of the truncated form of Poisson summation, and its application to the contribution of the trivial automorphic representation of $G$ to the geometric side of the trace formula. Combining various arguments, Langlands was able to establish a more tractable form of Poisson summation for both of the functions $\theta_f$ and $\phi_f$. The singularities of the local components $f_v^G(a_v)$ and $f_v'(a_v')$ were of course an important part of this. Once Poisson summation was in place, the singularities could then be used to analyze the asymptotic behavior of the Fourier transforms

$$\widehat{\theta}_f(a), \widehat{\phi}_f(a), \qquad a \in \mathbb{A}.$$

With the comparison of trace formulas for $G$ and $\{G'\}$ in mind, Langlands then considered the difference

$$S_{\mathrm{reg,ell}}(f) - \sum_{G'} S'_{\mathrm{reg,ell}}(f') \sim \sum_{b \in \mathscr{A}(G,F)} \theta_f(b) - \phi_f(b), \tag{154}$$

---

[27] As Langlands remarks, the proofs for this announcement are not published. They were originally circulated as a long preprint by Sally and Shalika, part of which was later published as the fundamental article [216] on Shalika germs. However, it was the complementary part that contained the $p$-adic characters, and that was unfortunately never published.

a linear form in $f$ that becomes

$$\sum_{b \in \mathscr{A}(G,F)} \left( \widehat{\theta}_f(b) - \widehat{\phi}_f(b) \right)$$

after the application of Poisson summation [158, (5.5) and (5.12)]. In the rest of the paper, he added some remarks on possible next steps, which would include the removal of the contribution $\widehat{\theta}_f(0)$ of the trivial automorphic representation of SL(2) from (154). His main concern would be to gain information about the proposed limits (149) and their variants.

My discussion of the paper [158] has been superficial, as is no doubt clear from my use of the symbol $\sim$ in (154). In particular, I have said nothing of the truncated Poisson summation formula from [74]. However, it does seem that (154) is a natural identity, with important implications for the general Beyond Endoscopic comparison of stable trace formulas.

A different approach to Poisson summation was introduced by Ali Altug in his 2003 thesis, written under the supervision of Langlands and Peter Sarnak. It was published later, with applications, in the three papers [3], [4], [5]. It will be instructive for us to discuss each of these papers.

In contrast to [74], where Poisson summation was established for a general class of test functions $f$ on a general group $G$, Altug's methods apply to a restricted class of functions $f$ on the particular group $G = \mathrm{GL}(2)$. However, they lead to much sharper estimates. Moreover, they would seem in principle to be applicable to more general groups, even if the technical problems in extending them will be formidable. We shall describe the first paper [3], which is devoted to Poisson summation, in some detail.

Following [3], we take $G$ to be the group $\mathrm{GL}(2)$, $F = \mathbb{Q}$, and

$$f = f_\infty \cdot f^\infty = f_\infty \cdot f^{\infty,p} \cdot f_p^k, \quad p \text{ prime}, k \in \mathbb{Z}^{\geq 0},$$

to be a test function in the space

$$C_c^\infty(Z_+ \backslash G(\mathbb{A})), \qquad Z_+ = A_G(\mathbb{R})^0 = \left\{ \begin{pmatrix} u & 0 \\ 0 & u \end{pmatrix} : u > 0 \right\},$$

as follows. The archimedean component $f_\infty$ is any function in $C_c^\infty(Z_+ \backslash G(\mathbb{R}))$, while $f^{\infty,p}$ is the characteristic function of the standard open, maximal compact subgroup $K^{\infty,p}$ of $G(\mathbb{A}^{\infty,p})$, and $f_p^k$ is the product of $p^{-k/2}$ with the characteristic function of the open compact subset

$$\{X \in M_{2 \times 2}(\mathbb{Z}_p) : |\det X|_p = p^{-k}\}$$

of $G(\mathbb{Q}_p)$. This is of course a very special case, but it was well suited to illustrating Altug's techniques.

For the general regular elliptic part $I_{\mathrm{ell,reg}}(f)$ of the trace formula at $f$, we are interested in the value of the normalized orbital integral $f_G(\gamma)$ of $f$ at a point

$\gamma \in \Gamma_{\mathrm{ell,reg}}(\mathbb{Q})$. It is a product

$$f_G(a) = f_{\infty,G}(a) \cdot |D(\gamma)|^{-\frac{1}{2}} \mathrm{Orb}(\gamma, f^\infty),$$

where

$$|D(\gamma)| = |D(\gamma)|_\infty = (|D(\gamma)|^\infty)^{-1}$$

while

$$a = (b, a_2), \qquad b \in \mathbb{Q}, a_2 \in \mathbb{Q}^*,$$

is the bijective image of $\gamma$ in $\mathscr{A}_{\mathrm{irred}}(1, \mathbb{Q})$, and $f_{\infty,G}(a) = f_{\infty,G}(\gamma)$. One sees from the properties of the unramified orbital integral

$$\mathrm{Orb}(\gamma, f^\infty) = \mathrm{Orb}(\gamma, f_p^k) \mathrm{Orb}(\gamma, f^{\infty,p}),$$

and the definition of the function $f_p^k$, that $f_G(a)$ vanishes unless the irreducible monic, quadratic polynomial $p_a(\lambda)$ has integral coefficients, with constant term $a_2$ equal to $\pm p^k$. It then follows that $I_{\mathrm{ell,reg}}(f)$ is equal to the expression

$$\sum_{\varepsilon \in \{\pm 1\}} \sum_{b \in \mathbb{Z}_{\mathrm{irred}}^\varepsilon} f_{\infty,G}^\varepsilon(b) \cdot |D(\gamma)|^{-\frac{1}{2}} \mathrm{vol}(\gamma) \mathrm{Orb}(\gamma, f^\infty), \tag{155}$$

where we have written $\mathbb{Z}_{\mathrm{irred}}^\varepsilon = \mathscr{B}_{\mathrm{irred}}^\varepsilon(1, \mathbb{Z})$ for the set of integers $b \in \mathbb{Z}$ such that the pair $a = (b, \varepsilon p^k)$ lies in $\mathscr{A}_{\mathrm{irred}}(1, \mathbb{Z})$, and

$$f_{\infty,G}^\varepsilon(b) = f_{\infty,G}(b, \varepsilon p^k) = f_{\infty,G}(a).$$

We are following the discussion in §3 of the paper [3]. As was noted there, the expansion (155) of $I_{\mathrm{ell,reg}}(f)$ is where any discussion of Poisson summation would begin. The inner sum is of course over a set of integers $b \in \mathbb{Z}$. Is there a natural extension of the values of the summands to a Schwartz function on $\mathbb{R}$ to which Poisson summation could be applied? The answer is clearly negative. For a start, we must not forget that Poisson summation applies to the sum of the values of a suitable function on $\mathbb{R}$ over the lattice $\mathbb{Z}$, not a linear combination.

There are in fact a number of obstacles. The most immediately daunting is perhaps the volume factor $\mathrm{vol}(\gamma)$ in (155). It depends very much on $b$ as an integer, and in particular, on the splitting field $E_a$ over $\mathbb{Q}$ of the quadratic polynomial

$$p_a(\lambda) = p_b^\varepsilon(\lambda), \qquad a = (b, \varepsilon p^k),$$

over $\mathbb{Z}$. The same goes for the factor $\mathrm{Orb}(\gamma, f^\infty)$. It is a product of $q$-adic orbital integrals of the global spherical function

$$\left( \prod_{q \neq p} f_q \right) f_p^k$$

at $\gamma$. It too depends on the splitting field of $p_b^\varepsilon(\lambda)$. The archimedean factor $f_{\infty,G}^\varepsilon(b)$ is more amenable. It does extend to the natural function of $b \in \mathbb{R}$ given by the normalized archimedean orbital integrals of $f_\infty$. This should in principle be the function to which one would try to apply Poisson summation. However, it comes with its own problems as a function of $\mathbb{R}$, namely the singularities given by Harish-Chandra's jump conditions. These occur at points $b \in \mathbb{R}$ at which the characteristic polynomials $p_b^\varepsilon(\lambda)$ have repeated factors over $\mathbb{R}$, which is to say that the corresponding discriminants vanish. And finally, there is the problem that the sum over $b$ is not over the full lattice $\mathbb{Z}$ in $\mathbb{R}$, but the subset $\mathbb{Z}_{\text{irred}}^\varepsilon$ of $\mathbb{Z}$. Altug dealt with all of these problems.

For the volume coefficient $\text{vol}(\gamma)$ in (155), the first step was to apply the Dirichlet class number formula for the quadratic extension $E_\gamma = E_a$ of $\mathbb{Q}$. Since the volume is essentially the regulator of $E_\gamma/\mathbb{Q}$, the formula can be written

$$\text{vol}(\gamma) = |D_\gamma|^{\frac{1}{2}} L\Big(1, \Big(\frac{D_\gamma}{\cdot}\Big)\Big),$$

where $D_\gamma$ is the discriminant of the quadratic extension $E_\gamma$, $\Big(\frac{D_\gamma}{\cdot}\Big)$ is the Kronecker symbol, and $L(\cdot, \Big(\frac{D_\gamma}{\cdot}\Big))$ is the Dirichlet $L$-function of $E_\gamma/\mathbb{Q}$ [3, §2.2.2]. At first glance, this seems to replace one problematic coefficient with two new ones!

However, as the discriminant of the field $E_\gamma$, $D_\gamma$ is closely related to the discriminant $D(\gamma)$ of the irreducible quadratic polynomial

$$p_\gamma(\lambda) = p_a(\lambda) = p_b^\varepsilon(\lambda).$$

This occurs in its own right as the factor $|D(\gamma)|^{-\frac{1}{2}}$ in (155). Their joint contribution equals

$$|D(\gamma)|^{-\frac{1}{2}} |D_\gamma|^{\frac{1}{2}} = |s_\gamma^2|^{-\frac{1}{2}} = s_\gamma^{-1},$$

for a positive integer $s_\gamma$. This last integer is in turn closely related to the product $\text{Orb}(\gamma, f^\infty)$ of $q$-adic orbital integrals in (155). Langlands had already computed the product in his initial paper on the subject [155, Lemma 1 and §2.5] as

$$\text{Orb}(\gamma, f^\infty) = p^{-\frac{k}{2}} \prod_{f \mid s_\gamma} f \cdot \Big(\prod_{q \mid f} \Big(1 - \Big(\frac{D_\gamma}{q}\Big) q^{-1}\Big)\Big). \tag{156}$$

Its product with $s_\gamma^{-1}$ and the other factor $L(1, \Big(\frac{D_\gamma}{\cdot}\Big))$ above equals

$$p^{-\frac{k}{2}} \sum_{f \mid s_\gamma} (f/s_\gamma) \Big(\prod_{q \mid f} \Big(1 - \Big(\frac{D_\gamma}{q}\Big)\Big)\Big) L\Big(1, \Big(\frac{D_\gamma}{\cdot}\Big)\Big),$$

which becomes

$$p^{-\frac{k}{2}} \sum_{f \mid s_\gamma} (1/f) L\Big(1, \Big(\frac{D(\gamma)/f^2}{\cdot}\Big)\Big)$$

with a change of variable in the sum over $f$, and the definition

$$L\Big(1,\Big(\frac{D(\gamma)/f^2}{\cdot}\Big)\Big) = \prod_{q|f}\Big(L_q\Big(1,\Big(\frac{D_\gamma}{q}\Big)\Big)^{-1}\Big)\cdot L\Big(1,\Big(\frac{D_\gamma}{\cdot}\Big)\Big)$$

of the Dirichlet $L$-function for the discriminant $D(\gamma)/f^2$. The expression (155) for $I_{\mathrm{reg,ell}}(f)$ then becomes

$$p^{-\frac{k}{2}}\sum_\varepsilon\sum_b f^\varepsilon_{\infty,G}(b)\Big(\sum_{f|s_\gamma}(1/f)L\Big(1,\Big(\frac{D(\gamma)/f^2}{\cdot}\Big)\Big)\Big).\qquad(157)$$

(See [3, (2) and (3)].)

Altug was then able to treat the Dirichlet $L$-values at 1 with what is known as the approximate functional equation. To this end, he attached the more general Dirichlet series

$$L(s,\delta) = \sum_{f^2|\delta}'\frac{1}{f^{2s-1}}L\Big(s,\Big(\frac{\delta/f^2}{\cdot}\Big)\Big)$$

to any discriminant $\delta$. Thus, $\delta$ is an integer congruent to 0 or 1 modulo 4, while $\Big(\frac{\delta/f^2}{\cdot}\Big)$ is the Kronecker symbol again, and $\sum_{f^2|\delta}'$ stands for the sum over the integers $f$ such that $\delta/f^2$ is also congruent to 0 or 1 modulo 4 [3, §3.1]. Its value at $\delta = D(\gamma)$ and $s = 1$ equals the expression in (157) in brackets. In particular, it is constructed as a product of the Dirichlet $L$-value $L\Big(1,\Big(\frac{D_\gamma}{\cdot}\Big)\Big)$ with the values of finitely many $q$-adic orbital integrals (with $k = 0$). This enhanced $L$-function was introduced by Zagier in 1977 [255], who established a functional equation linking its values at $s$ and $1-s$. He would probably not have been aware of its interpretation in terms of $q$-adic orbital integrals on $\mathrm{GL}(2,\mathbb{Q}_q)$. However, it is interesting to think that these objects were implicit in his quite different setting not long after they had also become a part of Langlands' study of base change for $\mathrm{GL}(2)$, with its implications for the Artin conjecture and ultimately Fermat's Last Theorem.

The approximate functional equation applies to any reasonable Dirichlet series with functional equation [99, §10.6]. Altug used it to express the sum $L(1,D(\gamma))$ in the brackets of (157) by an ungainly but more tractable expression. The process takes the value $L(s,\delta)$ of the general $L$-function above back to a Dirichlet series, but one whose summands over $n$ are weighted so as to converge absolutely, in which the original coefficients are averaged against a well behaved test function $F$. There is a second term, a contour integral that would give trouble on its own. However, a change of contour, together with an application of the actual functional equation for $L(s,\delta)$ transforms the troublesome integral to a weighted Dirichlet series for $L(1-s,\delta)$, in which the original coefficients are averaged over a second well behaved function $H$. (See the discussion of [3], and in particular, Proposition 3.4 and Corollary 3.5). With this, the seemingly insurmountable problem of the arithmetic dependence on $\gamma$ of the original coefficients $\mathrm{vol}(\gamma)$ in (155) can be controlled.

We note in passing that that the product $\zeta(s)L(s,\delta)$ is called the *Dedekind zeta function of the (monogenic) order*

$$R_\delta = \mathbb{Z}[\lambda]/(p_\delta(\lambda))$$

in the quadratic field

$$E_\delta = \mathbb{Q}[\lambda]/(p_\delta(\lambda)),$$

where $p_\delta(\lambda)$ is the irreducible characteristic polynomial attached to $\delta$. (We are assuming here that $\delta = D(\gamma)$ as above.) Its analogue for $\mathrm{GL}(n+1)$ over a number field $F$ is the topic of the paper [254] of Z. Yun. He established its functional equation and class number formula as Theorem 1 of the paper. To do so, he defined the local factors as functions of $s$ that satisfy their own functional equations [254, Theorem 2.5]. Yun then expressed the values at $s = 1$ of these local factors in terms of local orbital integrals at $\gamma$ [254, Corollary 4.6]. This paper is likely to be important in any attempt to generalize the methods of Altug. However, it does not provide an analogue for $\mathrm{GL}(n+1)$ of the explicit formula (156) of Langlands for the orbital integrals on $\mathrm{GL}(2)$. Something of the sort appears to be essential for higher rank, but there are some qualitative differences even between the cases of $\mathrm{GL}(2)$ and $\mathrm{GL}(3)$. (See [117].) I have been told that the perverse sheaves attached to unramified orbital integrals by Ngô acquire singularities in generalizing from $\mathrm{GL}(2)$ to $\mathrm{GL}(3)$.

Returning to the regular elliptic expansion (155) for $\mathrm{GL}(2)$, we recall the other two obstacles mentioned above. The first concerned the singularities of the function $f^\varepsilon_{\infty,G}(b) = f_{\infty,G}(a)$ at the zero set of the discriminant $D(\gamma) = D(a)$. To deal with it, Altug observed that for any Schwartz function $\Phi$ in $\mathbb{R}$ and any $\alpha > 0$, the product

$$\Phi(|D(a)|^{-\alpha})f^\varepsilon_{\infty,G}(b), \qquad b \in \mathbb{R},$$

is also a Schwartz function of $\mathbb{R}$ ([3, Proposition 4.1]). Altug built this into the Schwartz functions $F$ and $H$ he obtained from the approximate functional equation. He had taken care to include a supplementary parameter in each of the arguments in [3, (4')] ($A^{-1}$ for $F$ and $A$ for $H$), which he then simply set equal to $|D(\gamma)|^\alpha$, for any number $\alpha$ with $0 < \alpha < 1$. The resulting expression in the brackets of (157) was then sufficient to make the complementary factor $f^\varepsilon_{\infty,G}(b)$ in (157) the restriction to $\mathbb{Z}^\varepsilon_{\mathrm{irred}}$ of a well defined Schwartz function on $\mathbb{R}$.

The second and last obstruction was that the sum over $b$ in (157) was only over the subset $\mathbb{Z}^\varepsilon_{\mathrm{irred}}$ of the lattice $\mathbb{Z}$ of $\mathbb{R}$ (corresponding to *irreducible* characteristic polynomials). Altug simply added the complementary set of orbital integrals $f^\varepsilon_{\infty,G}(b)$ to the sum in (157). The trace formula actually calls for *weighted* orbital integrals (the term (iv) on page 516 of [103]) to be taken here, but we have been dealing with approximations ever since we identified $I_{\mathrm{reg,ell}}(f)$ with the geometric side of the trace formula. We write $\bar{I}_{\mathrm{reg,ell}}(f)$ for the expression (157) with $b$ summed over $\mathbb{Z}$ rather than the subset $\mathbb{Z}^\varepsilon_{\mathrm{irred}}$.

For the record, we write out the expression for (157) obtained by Altug for $\bar{I}_{\mathrm{reg,ell}}(f)$ (with the extra summands for $b \in \mathbb{Z}$), even though our discussion is still lacking some of the finer details. It is

$$p^{-\frac{k}{2}} \sum_{\varepsilon} \sum_{b \in \mathbb{Z}} f_G^{\varepsilon}(b) \sum_{f^2 \mid D^{\varepsilon}(b)} \frac{1}{f} \sum_{\ell=1}^{\infty} \frac{1}{\ell} \left( \frac{D^{\varepsilon}(b)/f^2}{\ell} \right) \cdot$$

$$\left[ F\left( \frac{\ell f^2}{|D^{\varepsilon}(b)|^{\alpha}} \right) + \ell f^2 |D^{\varepsilon}(b)|^{-\frac{1}{2}} H\left( \frac{\ell f^2}{|D^{\varepsilon}(b)|^{1-\alpha}} \right) \right].$$

(This is essentially the first expression given in the proof of Theorem 4.2 from [3].) Since everything converges absolutely, the sum over $b \in \mathbb{Z}$ can be taken inside the sum over $f$ and $\ell$. It is still not yet quite possible to apply Poisson summation to the sum over $b$. For while the last expression in the square brackets extends to a Schwartz function of $b \in \mathbb{R}$, there are also the coefficients $\left( \frac{D^{\varepsilon}(b)/f^2}{\ell} \right)$ that depend arithmetically on $b \in \mathbb{Z}$. Altug's solution was to break this sum over $b$ into a double sum over the finite subset

$$C(\ell, f) = \{b \,(\mathrm{mod}\ 4\ell f^2) : D^{\varepsilon}(b) \equiv 0 \,(\mathrm{mod}\ f^2), D^{\varepsilon}(b)/f \equiv 0, 1 \,(\mathrm{mod}\ 4)\}$$

of congruence classes modulo $4\ell f^2$, and the infinite affine lattice

$$\{m \in \mathbb{Z} : m \equiv b \,(\mathrm{mod}\ 4\ell f^2)\}.$$

Since

$$\left( \frac{D^{\varepsilon}(b+m)/f^2}{\ell} \right) = \left( \frac{D^{\varepsilon}(b)/f^2}{\ell} \right),$$

the coefficients could then be taken outside the sum over $m$, allowing him then to apply Poisson summation to the sum.

The final result is then

$$\bar{I}_{\mathrm{ell,reg}}(f) = \sum_{\xi \in \mathbb{Z}} \widehat{\bar{I}}_{\mathrm{ell,reg}}(\xi, f), \tag{158}$$

where $\widehat{\bar{I}}_{\mathrm{ell,reg}}(\xi, f)$ equals an expression

$$\sum_{\varepsilon} \sum_{f=1}^{\infty} \frac{1}{f^3} \sum_{\ell=1}^{\infty} \frac{1}{\ell^2} \widehat{f}_{G,\ell,f}^{\varepsilon}(\xi) K_{\ell,f}(\xi, \varepsilon p^k),$$

for a Fourier transform

$$\widehat{f}_{G,\ell,f}^{\varepsilon}(\xi) = \int_{\mathbb{R}} (f_G^{\varepsilon}(x) \phi_{\ell,f}^{F,H}(x) \exp\left( \frac{-x\xi}{2\ell f^2 p^{-\frac{k}{2}}} \right) dx,$$

with the function

$$\phi_{\ell,f}^{F,H}(x) = F(\ell f^2 |D^{\varepsilon}(x)|^{-\alpha}) + (\ell f^2)|D^{\varepsilon}(x)|^{-\frac{1}{2}} H(\ell f^2 |D^{\varepsilon}(x)|^{\alpha-1}),$$

and for a finite exponential, Kloosterman-like sum

$$\mathrm{Kl}_{\ell,f}(\xi, \varepsilon p^k) = \sum_{b \in C(\ell,f)} \left( \frac{D^{\varepsilon}(b)/f^2}{\ell} \right) \exp\left( \frac{b\xi}{4\ell f^2} \right).$$

This is Theorem 4.2 of [3], one of the two main results of the paper.

The second main result of Altug's paper [3] is Theorem 6.1. It is a reworking of the formula for the constant term $\widehat{I}_{\mathrm{ell,reg}}(0,\xi)$ in the expansion (158), obtained from the functional equations for the $L$-functions $L(s,\delta)$, some elementary but quite elaborate identities established in §5 of [3] (Lemmas 5.1–5.3, Corollary 5.4), and various changes in the contour integral over $x$ in the original formula. (Altug writes the formula for $\widehat{I}(0,f)$, and the integral over $x$ that it contains, as the expression $(13)_{\xi=0}$ displayed prior to the statement of Theorem 6.1.)

The formula in the assertion of Altug's Theorem 6.1, which we will not quote, seems to make his formula Theorem 4.2 (distilled in (158) above) look simple by comparison! However, it has a compelling logic. It expresses $\widehat{I}_{\mathrm{ell,reg}}(0,f)$ as a sum of two simple integrals over $x \in \mathbb{R}$, together with a third complicated integral over $x$. In his last result [3, Lemma 6.2], Altug identifies the two simple integrals as the character $\mathrm{tr}(\mathbf{1}(f))$ of the trivial one-dimensional representation of $G(\mathbb{A})$ at $f$, and the supplementary term[28] on the spectral side given by (vi) on p. 647 of [103]. The hope is that the contribution of the third integral vanishes in the putative limits (149) attached to (irreducible) representations $r \neq 1$ of $\mathrm{GL}(2,\mathbb{C})$. Altug puts Theorem 4.2 and Theorem 6.1 together as a complex formula for $I_{\mathrm{reg,ell}}(f)$, which he states as Theorem $1.128^0$ in the introduction to [3].

This completes our discussion of (158), Altug's version of Poisson summation for $\mathrm{GL}(2)$. It is surprisingly complex, even as it applies only to a special case within $\mathrm{GL}(2)$. We have emphasised it in our discussion because such detail seems to be a necessary prelude for the kind of estimates that will ultimately be needed. The role of the more elementary version of Poisson summation in [74] and [158], which applies to any group $G$, is different. It was offered simply as evidence for the basic idea, that of using some form of Poisson summation on the Steinberg–Hitching base to recognize the ultimate contribution of the noncuspidal discrete spectrum to the geometric side of the trace formula.

In his second paper [4], Altug wrote $f^p$ for the function $f = f_\infty \cdot f^{\infty,p} \cdot f_p^k$ from [3] with $k = 1$. He then showed that its value at the third (complicated) integral in the formula for $\widehat{I}(0,f^p)$ from Theorem 6.1 of [3] (adjusted$28^0$ as in the statement of Theorem 4.1 of [4]) is bounded in $p$ [4, formula (00) from Corollary 4.3]. This was actually relatively simple. Much deeper were his estimates of the terms $\widehat{I}_{\mathrm{ell,reg}}(\xi,f)$ with $\xi \neq 0$ in (158). In Theorem 4.4, he established an estimate

---

[28] What I have stated here is not strictly correct. The second simple integral in Theorem 6.1 actually equals *twice* the supplementary term (vi). However, the excess really belongs naturally in the complicated third integral, as Altug observes in a restatement of Theorem 1.1 for his subsequent paper. (See Theorem 4.1 in [4].)
.

$$| \sum_{\xi \neq 0} \widehat{I}(\xi, f^p)| \leq c_\infty p^{\frac{1}{4}},$$

where $c_\infty$ is a constant that depends only on $f_\infty$, but not $p$. This was a consequence of analytic results on asymptotic properties of Fourier transforms in Appendix A [4, Theorems A.14 and A.15, and Corollary 11.16], as well as arithmetic properties of the character sums $\mathrm{Kl}_{\ell,f}(\xi, \pm p)$ in Appendix B [4, Corollary B.8]. Combined with a straightforward analysis in §3 in [4] of the remaining terms in the full trace formula for $\mathrm{tr}(R_{\mathrm{cusp}}(f^p))$, he arrived at his main result, the estimate

$$\mathrm{tr}(R_{\mathrm{cusp}}(f^p)) = O(p^{\frac{1}{4}}) \tag{159}$$

of [4, Theorem 1.1].

The estimate (159) represents a partial bound towards the Ramanujan conjecture for Hecke eigenvalues of Maass forms. It is the same estimate that had been established in 1980 by Kuznetsov [126] by what is now regarded as a special case of the relative trace formula. As Altug remarks, the importance of (159) is in its method of proof by the Arthur–Selberg trace formula, which has more structure, and has been established for general groups. The full Ramanujan conjecture would be a consequence of functoriality, but is still far from known. It would amount to a bound

$$\mathrm{tr}(R_{\mathrm{cusp}}(f^p)) = O(p^\varepsilon)$$

for every $\varepsilon > 0$. The estimate (159) is thus intermediate between the full Ramanujan conjecture and the elementary bound

$$\mathrm{tr}(R_{\mathrm{cusp}}(f^p)) = O(p^{\frac{1}{2}})$$

represented by the 1-dimensional representation of $\mathrm{GL}(2, \mathbb{A})$. We recall from our discussion of the Langlands–Shahidi method that in establishing functoriality for the irreducible 4 and 5 dimensional representations of $\widehat{G} = \mathrm{GL}(2, \mathbb{C})$, Shahidi and Kim obtained bounds that are sharper than (159). However the full Ramanujan conjecture (for Hecke operators of Maass forms) would require functoriality for all irreducible representations of $\widehat{G}$.

The third paper [5] of Altug established an $r$-trace formula (148) in the special case of $G = \mathrm{GL}(2)$ and $F = \mathbb{Q}$ at hand. For this, he restricted the basic function $f = f_\infty f^\infty$ further, by taking $f^\infty$ to be the characteristic function of the standard maximal compact subgroup $K^\infty$ of $G(\mathbb{A})$, and $f_\infty$ to be a *cuspidal* function on $G(\mathbb{R})$, with tempered characters supported on the discrete series representation parameterized by a fixed integer $k \geq 3$. With these data, Altug considered the proposed limit (149) when $r$ is the standard, two dimensional representation of $\widehat{G} = \mathrm{GL}(2, \mathbb{C})$. Combining the intricate analytic results of [4], he established that the limit exists, and equals 0. In other words, the $r$-trace formula in this case is

$$I^r_{\mathrm{cusp}}(f) = 0,$$

according to the definition (148). This is what is expected. For the only cuspidal automorphic representations $\pi$ of $G = \mathrm{GL}(2)$ that are *proper* functorial images (which is to say, not *primitive* in the language of §9), should be of "CM" or "Galois type", attached to two dimensional representations of the Weil group $W_\mathbb{Q}$. This implies that $\pi(f)$ actually vanishes[29], given the choice of $f^\infty$ and the fact that the class number of $\mathbb{Q}$ is 1. Therefore $\pi$ should contribute nothing to $I^r_{\mathrm{cusp}}(f)$. It is a general fact that a primitive representation $\pi$ would also contribute nothing to $I^r_{\mathrm{cusp}}(f)$, for any irreducible nontrivial representation $r$.

We should include a historical remark at this point before concluding our discussion of Altug's work. When Langlands first introduced his ideas in the precursor [154] of his published article [155], Sarnak had reservations about the form of the proposed limit (149) obtained from the order of poles at $s = 1$ of the $L$-functions $L(s, \pi, r)$. He was concerned that proving and calculating a limit (149), difficult under any circumstances, might be even more intractable with the form of the right hand side. He suggested in the letter [199] to Langlands replacing the sum over $p$ on the right hand side of (149) by a sum over $n \in \mathbb{N}$. This amounts essentially to a change of the weighting coefficient $m_\pi(r)$ in (143), the order of the pole of $L(s, \pi, r)$ at $s = 1$, to another coefficient $n_\pi(r)$, the *residue* of $L(s, \pi, r)$ at $s = 1$. It would then open the possibility of applying Poisson summation again, this time to the new sum over $n$. The Tauberian theorem we quoted from [208] would also hold in this context, applied to the Dirichlet series for $L(s, \pi, r)$ rather than its logarithmic derivative. Langlands himself appears to have been ambivalent about this suggestion. For the new coefficients would no longer be additive in $r$, complicating the anticipated spectral (primitive) Beyond Endoscopic decomposition of what we are calling the $r$-trace formula.

The suggestion was first taken up by A. Venkatesh, a student of Sarnak at the time. He actually worked with the Kuznetsov formula for $\mathrm{GL}(2)$, a special case of Jacquet's proposed "relative trace formula", rather than the trace formula itself. The technical difficulties simplify in this setting. Using the Poisson summation formula proposed by Sarnak, Venkatesh considered the more complex case of the three dimensional, symmetric square representation $r_2$ of $\mathrm{GL}(2)$. In his 2002 thesis [241] and a subsequent paper [242], he was able to separate the contribution of the relevant forms of CM-type, which is to say the automorphic representations $\pi$ of $\mathrm{GL}(2)$ attached to two-dimensional representations of the Weil group $W_F$, by using the residues at $s = 1$ of the $L$-functions $L(s, \pi, r_2)$. His results were also more general in that they applied to many number fields $F$ in place of $\mathbb{Q}$, and automorphic forms that were ramified, which is to say of level greater than 1. The work of Venkatesh was a breakthrough. It would be very interesting to make a careful comparison of

---

[29] This property should remain valid if $r$ is replaced by an $(m+1)$-dimensional symmetric power $r_m$, for $m > 1$. We would therefore again expect that $I^{r_m}_{\mathrm{cusp}}(f) = 0$, even though this would be much more difficult to prove. On the other hand, if $f^\infty$ is a more general function, or $\mathbb{Q}$ is replaced by a more general field $F$, we would not expect $I^{r_m}_{\mathrm{cusp}}(f)$ to vanish. For a table of multiplicities

$$-m_\pi(r_m) = [1_{L_G} : r_m \circ \phi],$$

where $\pi$ corresponds to a 2-dimensional Galois representation $\phi$ of $\Gamma_\mathbb{Q}$, see page 6 of [154].

his techniques with those of Altug. In this report, we will settle for a few general remarks at the end on the possible future roles on the two kinds of "trace" formulas.

Altug took up Sarnak's suggestion in his third paper [5]. The Tauberian theorem for the new weighting coefficient attached to any $G$, $\pi$ and $r$ would take the general form

$$n_\pi(r) = \lim_{X \to \infty} |X|^{-1} \sum_{n < X} \operatorname{tr}(T(n, \pi, r))$$

where[30].

$$T(n, \pi, r) = \pi^S(h_n^r), \qquad S = \{\infty\},$$

is the (unramified) Hecke operator for $\pi^S$, $r$ and $n$ such that

$$L^\infty(s, \pi, r) = L^S(s, \pi, r) = \sum_{n=0}^\infty \operatorname{tr}(T(n, \pi, r))n^{-s}.$$

Altug again confined himself to the case that $r = r_1$ is the standard two-dimensional representation of the group $^LG = \widehat{G} = \operatorname{GL}(2, \mathbb{C})$. It is likely that his methods could be extended to $r = r_2$, and to the more general fields $F$ and functions $f$ treated by Venkatesh, but he did not attempt to deal with the increased complexity that would arise from the classical trace formula for GL(2) in this paper. In fact, as we noted earlier, he restricted $f = f_\infty f^\infty$ further so as to be supported on the unramified automorphic representations $\pi = \pi_\infty \pi^\infty$ such that $\pi_\infty$ corresponds to an automorphic form in the space $S_k$ of (holomorphic) cusp forms of weight $k$ (and level 1), for fixed $k \geq 1$. Then

$$T_k(n) = \bigoplus_\pi T(n, \pi, r)$$

is the $n$th Hecke operator acting on the (finite-dimensional) complex vector space $S_k$.

Altug's main result, Theorem 1.1 of [5], is an estimate

---

[30] Motivated by notation from the beginning of the section, we have written

$$h_n^r = \bigotimes_p h_{n,p}^r, \qquad n = \prod_{p \notin S} p^{n_p},$$

here for the function in the unramified Hecke algebra

$$C_c(K^S \backslash G(\mathbb{A}^S)/K^S) = \widetilde{\bigotimes_{p \notin S}} C_c(\mathbf{1}_p \backslash G(\mathbb{Q}_p)/\mathbf{1}_p)$$

whose Satake transform $\widehat{h}_n^r$ at a family of semisimple classes $c = \{c_p : p \notin S\}$ in $^LG$ equals the product

$$\widehat{h}_n^r(c^s) = \prod_p \widehat{h}_{n,p}(c_p) = \prod_p \operatorname{tr}((S_{n_p}r)(c_p))$$

of characters of symmetric powers $S_{n_p}r$ of $r$ at the points of $c_p$. (See [44])

$$\sum_{n<X} \operatorname{tr}(T_k(n)) = O_{k,\varepsilon}(X^{\frac{31}{32}+\varepsilon}),$$

for any $\varepsilon > 0$. In particular, for any such $\pi$, the function

$$|X|^{-1} \sum_{n<X} \operatorname{tr}(T(n,\pi,r))$$

converges strongly to 0 as $X$ approaches $\infty$. In Corollary 1.2 of [5], Altug specialized this to an estimate for the Ramanujan $\Delta$-function (of weight 12 and level 1), or rather a property of the corresponding $L$-function $L(s,\Delta)$. In the remarks following the statement of the corollary [5, p. 3–5], he discussed other interesting consequences, including the original $r$-trace formula.

The rest of Altug's third paper [5] is devoted to the proof of Theorem 1.1. As he observed prior to the statement of the theorem, it puts together all of the work of his previous two papers. In Section 2 he gave a broad outline of the proof of the theorem, and the proof of its Corollary 1.2. In Section 3, he dealt with the supplementary (nonelliptic, noncuspidal) terms in the trace formula for GL(2). Section 4 concerns the elliptic terms, which remain the basic objects. It is the heart of the paper. Section 4.1 is a short review and further analysis of the results of [3], notably the Poisson summation formula (158) for the sums over $m$ introduced prior to its statement. Section 4.2.1 begins with a heuristic discussion of estimates provided by subsequent Theorems 4.9, 4.11 and 4.13, and then proceeds with their proof. Section 4.2.2 contains the final estimates given by the last Theorem 4.15 and its Corollary 4.16. It is in their proof that the critical second application of Poisson summation comes, the one for the sum over $n$. (The character sum $\operatorname{Kl}_{\ell,f}(\xi,n)$ is observed here to be periodic in $n$ modulo $4\ell f^2$, allowing among other things, a change from the sum over $n > 0$ to other sums over $n \neq 0$.) The final Section 5 contains further real and $p$-adic analysis. This was used in the proofs of estimates (Proposition 5.2 and Corollary 5.9) postponed from Section 4.

This completes our discussion of the work of Altug. As a last word, it might be worthwhile to recapitulate each of his three papers in a sentence or two. The first one [3] establishes Poisson summation, and shows how the nontempered, trivial representation occurs naturally in the Fourier transform term with $\xi = 0$. The second paper [4] is an estimate that asserts that the complementary part of the $\xi = 0$ term as well as the remaining $\xi \neq 0$ terms are exponentially smaller than the trivial (one-dimensional) character. However, the error bound is still exponentially larger than the tempered singular term ((vi) in [103]) that was also removed from the summand with $\xi = 0$. The first two papers apply to a single formula, for a function $f^p$ that depends on $p$ and $k = n_p$. The third paper [5] applies to a weighted average of such formulas, parameterized by positive integers $n = \prod_p p^{n_p}$. The averaging process decreases the error term considerably further, to more than what it would be for any isolated tempered automorphic representation on the spectral side, and in particular, for the singular tempered representations from term (vi) of [103] that were removed from the original summand with $\xi = 0$.

The amount of space we have devoted to Altug's work might appear surprising. However, it is likely to serve as a concrete foundation for the future study of Beyond Endoscopy, at least as it follows Langlands' ideas for using the (stable) trace formula. Our discussion has also given us opportunities to illustrate aspects of the basic strategy, as it might apply in practice. To be sure, Altug's papers are rather imposing, particularly when one begins with the statement of the main result [3, Theorem 1.1] of his first paper, and its refinement [4, Theorem 4.1] for the second. But despite their complexity, the terms in the stated formulas are elementary. Their derivation and later application depend more than anything on basic analysis, despite the fact that they are leading to new techniques. There is also a suggestive unity to the three papers. Each one solves a specific problem in the simplest of cases, using only the methods of the trace formula. Taken in succession, they match sequential steps laid out by Langlands in his general strategy for Beyond Endoscopy.

An interesting problem would be to relate Altug's work with the second half of Langlands' paper [158]. In general, one might eventually want to establish a primitive (stable) trace formula, whose spectral side contains only the cuspidal, tempered representations that are not proper functorial images. One could imagine an inductive definition for any $G$ obtained by subtracting from its stable trace formula the primitive trace formulas attached to all proper beyond endoscopic data.[31]

For $G = \mathrm{GL}(2)$ and $G'$ a proper beyond endoscopic datum, the derived group $G'_{\mathrm{der}}$ would then be the trivial group $1_G$, so $(G', \mathscr{G}', \xi')$ would therefore correspond to an irreducible two dimensional representation of the Weil group. As we have noted, every such object ramifies over $\mathbb{Q}$, and therefore contributes nothing. To carry out the proposal, one would therefore want to extend Altug's results to more general functions $f = f_\infty f^\infty$, or more general number fields $F$, or ideally, both. From the resulting sum over $\xi \neq 0$, one could then subtract the expression obtained by Langlands by Poisson summation on the sum of the two dimensional characters on $W_F$ (informed perhaps by the more transparent formula described on p. 1 of [109]). Langlands' expression does not actually include characters on $W_F$ with finite image, those given by the irreducible characters of dihedral, tetrahedral, octahedral and icosahedral type on the Galois group $\Gamma_F$. This might be the real point. Is it too much to hope that we might recognize something of the sum of these characters in the more concrete difference of the original two sums over $\xi$? We are now talking about the deepest aspect of Beyond Endoscopy, the one that Langlands regarded as the true essence of the problem.

A part of Langlands' later paper [157] is pertinent to this last point. (The paper is described as a prologue, but the article on Functoriality and Reciprocity to which it refers has not been written.) First of all, Section 4 contains his reflections on the letter of Sarnak [199], and in particular, on the relative merits of working with the residues of either the $L$-functions $L(s, \pi, \rho)$ themselves or of their logarithmic derivatives. Altug's proof of an $r$-trace formula for $\mathrm{GL}(2)$, described in the remarks

---

[31] I have yet to think carefully about this. As suggested in the discussion of stable transfer earlier in the section, it would presumably follow the proposed construction of $L_F$ from Section 9. In particular, the beyond endoscopic groups $G'$ should perhaps have $G'_{\mathrm{der}}$ simply connected and $\dim G'_{\mathrm{der}} < \dim G_{\mathrm{der}}$.

in §1.1 of [5], is encouraging. It suggests that we might be able to have the best of both worlds.

It is the larger Section 5 of [157] that concerns Langlands' thoughts on complex Galois representations. He describes it as the arithmetic side of Beyond Endoscopy, as opposed perhaps to the analytic side later studied by Altug. Langlands discusses in some detail a paper [57] of Dedekind on quaternionic extensions $F$ of $\mathbb{Q}$, the extensions with Galois group the quaternionic group $Q_8$ of order 8 that the reader will recall (with $K/L$ in place of $F/\mathbb{Q}$) from the diagram at the end of §7. His thoughts on this might be described as "hard arithmetic", as opposed to "hard analysis". They can be taken as a reflection of Langlands' view, expressed in several places, that functoriality for complex Galois representations will be heart of Beyond Endoscopy. They are where the subject began, as the attempt by Langlands to extend the Artin reciprocity law, and thereby create a nonabelian class field theory. The Principle of Functoriality was of course his answer. However, the precepts of Beyond Endoscopy do not so far include arithmetic techniques with the power to make further progress on Artin's conjecture.

Another part of the paper [157] is devoted to an entirely different topic. It represents the beginnings of what Langlands called the *geometric theory*, as distinct from what is often referred to as the geometric Langlands program. Very roughly speaking, if the original (arithmetic) Langlands program concerns finite extensions of a global field $F$, a number field or the field of rational functions on a nonsingular projective curve over a finite field, these geometric programs are over the complex numbers. They concern the finite extensions of the field $F$ of meromorphic functions on a compact nonsingular Riemann surface $X$. Their starting point is the set $\mathrm{Bun}_X(G)$ of (equivalence classes of) holomorphic principal $G$-bundles on $X$, for a complex reductive group $G$. For any $x \in X$, $F_x$ is used to denote the field of formal complex Laurent series at $x$, and $\mathscr{O}_x$ is the subring of formal power series at $x$. There is then a natural bijection

$$\mathrm{Bun}_X(G) \cong G(F) \backslash G(\mathbb{A}_F)/K_F, \tag{160}$$

where

$$G(\mathbb{A}_F) = \prod_{x \in X}^{\sim} G(F_x)$$

and

$$K_F = \prod_{x \in X} G(\mathscr{O}_x).$$

It is clear from this that the theory should bear at least some formal resemblance to automorphic representation theory, even as its content ought to be primarily geometric.

The original geometric Langlands program was not initiated by Langlands. It is founded on notions from abstract algebraic geometry, accompanied by various sophisticated techniques from category theory. In particular, $\mathrm{Bun}_X(G)$ is treated as an algebraic stack, while the analogue of an automorphic representation takes the form of a perverse sheaf on the stack. Rather than try to comment further on these notions, let me refer the reader to the clearly written Bourbaki lecture [76] by Gaitsgory, who credits ideas of Beilinson, Deligne, Drinfeld and Laumon for the origins of the subject, and who is himself responsible for more recent progress. We note also that "Geometric Langlands" has had a significant influence in string theory. (See the Bourbaki lecture of Frenkel [73], and references there.)

Langlands wanted to understand the geometric theory in more concrete terms. In particular, he wanted to apply the methods of differential geometry and harmonic analysis in place of abstract algebraic geometry. His goal was to attach explicit (unramified) Hecke operators to the moduli space (160), together with corresponding explicit Hecke eigenvalues. This was not done in the prologue [157]. Langlands worked six additional years, at length posting the long paper [160] with his results in Russian, intending it especially for the Russian-American mathematicians in the field. He then worked further to convert it to the paper [162] in English, which he posted finally in 2020. I have read only the short note [161] describing his aims and results in quite general terms, but I shall look forward to reading the long paper.

The background for the paper is interesting. Langlands was motivated by the early paper [28] of Atiyah, written before his much better known work on the index and fixed point theorems. In it, Atiyah presented a concrete classification of the set of complex principal bundles of dimension $n$ over a complex torus, which is to say, the set (160) with $G = \mathrm{GL}(n)$ and $X$ an elliptic curve [28, Theorem 7]. This would be an *exact* fundamental domain for the space (160), rather than an approximate fundamental domain of the kind attached to a Siegel set in the arithmetic theory. Informed by this paper (which I am told is hard going), Langlands turned next to Atiyah's widely read paper [29] with Bott. He then set about studying the space of complex plane bundles over an elliptic curve, the special case of Atiyah's classification with $n = 2$. Langlands has reported that the desired Hecke operators and their eigenvalues appeared, in what seems to have been rather dramatic fashion, only at the very end.

In his note [161], Langlands has proposed a number of ways to attempt to extend these results. The first step might be to establish them for $\mathrm{GL}(n)$-bundles over $X$, the setting of Atiyah's classification. One could then attempt to replace $\mathrm{GL}(n)$ by an arbitrary complex group $G$, and $X$ by an arbitrary (compact, nonsingular) Riemann surface.[32] Langlands suggests that these last extensions will be much more difficult, since they would require among other things, major extensions of Atiyah's classification. One would of course also have to be armed with a firm command of the Langlands paper [162], a serious prerequisite to be sure. Assuming all of this could be established, it would then be interesting to compare the results with the geometric Langlands program based on abstract algebraic geometry. Finally, one could consider "ramification", or rather what would be ramification in the arithmetic theory. That is, one would consider $G$-bundles with extra structure, just as one takes elliptic curves with level structure in the classical theory of modular forms. It would entail replacing the group $K_F$ by a (normal) subgroup of finite co-dimension, defined by the analogue of congruence conditions ([157, p. 55]). A solution would have particular interest, since the question was initially not widely investigated through abstract algebraic geometry. (See however [8] for remarkable recent progress in this subject.)

---

[32] The case that $X$ is the Riemann sphere is presumably easy. As Atiyah noted at the beginning of his article, Grothendieck has shown that any vector bundle over a rational curve is a direct sum of line bundles.

Langlands suggests that each of these extensions will be difficult (if not as difficult as the many arithmetic problems in Beyond Endoscopy). However, they would probably have wide appeal. They seem to be revealing a new hidden structure in the objects commonly studied by differential geometers and topologists, which in turn is suggestive of the arithmetic structure in the original Langlands program.

We return to the discussion of Beyond Endoscopy with some final comments on the different possible approaches to functoriality. I have not mentioned the results of Braverman and Kazhdan [39], and the subsequent contributions of Ngô and others [187], [48], [38]. The idea is to try to establish the analytic continuation and functional equation for various automorphic $L$-functions $L(s, \pi, r)$ directly, and then perhaps use converse theorems to establish functoriality. This has direct roots in the work of Hecke and his original converse theorem, which we discussed in Section 6 in our review of Jacquet–Langlands [103]. It is the opposite of Langlands' fundamental idea of using functoriality to establish the general analytic continuation and functional equation from the basic case of principal $L$-functions for $\mathrm{GL}(n)$. I have not read these papers, and can make no further comment, even though I have enjoyed lectures on the subject. Braverman and Kazhdan's results also motivated L. Lafforgue. He has conjectured a nonlinear Poisson summation formula, which would be a consequence of functoriality, but which conversely would imply functoriality [128], [129].

We shall finish the report with some speculative remarks on the strategic differences between the (Arthur–Selberg) trace formula and the (Jacquet) relative "trace" formula[33]. Langlands' strategy applied to the former, while the work of Venkatesh for $\mathrm{GL}(2)$ on the latter faced fewer technical difficulties. (See also [194], [195], [196].) We are of course a long way from realizing any of the general goals of Beyond Endoscopy, but it is reasonable to try to plan how best to move forward. I believe strongly in the trace formula (perhaps not surprisingly), and I will try to express some of my reasons for this. However, I have limited experience with the period formula, so the reader should keep an open mind. Besides, there is a different point of view that I will express at the very end.

The trace formula is more highly developed. It has a clear general structure, and each of its terms has a natural, well defined source. It is reasonable to think that each of the terms will also have well defined role in Beyond Endoscopy. It is actually the stable trace formula that would be applied to groups other than $GL_n$. This is more sophisticated, but it is completely general, and has an equally rigid structure.

---

[33] This is not really a trace formula. In general, it would actually be a formula for periods of automorphic forms, rather than for characters of automorphic representations. Let me take the liberty of calling it (and any one of its variants) the *period formula* in what follows. We should note that the periods here differ from the Grothendieck periods attached to motives discussed in §9. They are defined as integrals of automorphic forms on $G$ against cycles defined by subgroups $H$ of $G$, which can sometimes represent parings between cohomology classes and homology classes. However, they make sense also for automorphic forms that are not motivic. We can call them *automorphic periods*, as opposed to the *motivic periods* from §9.

The period formula, at least in the Kuznetsov formula that was applied by Venkatesh to $GL(2)$, depends on the fact that a cuspidal automorphic representation of $GL(2)$ has a Whittaker model. For a quasisplit classical group, it is now known that every tempered, cuspidal $L$-packet contains a representation with a Whittaker model [23, §8.3]. To exploit this, however, one would want to work with local and global $L$-packets, which leads us back to the stable trace formula. For exceptional groups, the property is not known, even though it is widely expected. However, a proof would seem to require a theory of endoscopy for general groups. It may be that general endoscopy would have to be a consequence of Beyond Endoscopy but this is speculative. Keep in mind, however, that one really would want a full theory of Beyond Endoscopy for *all* groups, including exceptional groups. For it would be needed just to be able to classify the automorphic representations of, say, general linear groups, in terms of functorial images of primitive representations of smaller groups. This is closely related to the question of an explicit construction of the automorphic Galois group $L_F$ as described in Section 9.

For some more serious speculation (!), consider the following. Beyond Endoscopy is a proposal for attacking the Principle of Functoriality. What about its companion, Langlands' Reciprocity Conjecture? Is it possible that some extension of Beyond Endoscopy might have to be used to establish the two conjectures together? Langlands suggests something of the sort in his interesting and suggestive article [159]. It seems quite plausible to me. If so, the theory of general Shimura varieties would assume a central role beyond what it already holds. As a highly developed and mature theory, ultimately based on the arithmetic theory of reductive groups, it would presumably become a foundation for Reciprocity rather than just a source of interesting examples, much as the theory of Shimura varieties attached to $GL(2)$ is a foundation of the STW-conjecture for elliptic curves over $\mathbb{Q}$. Of course it would also demand fundamental new techniques, of which Wiles' proof, and the more recent work of Taylor and collaborators, would represent a beginning.

Take for example the role of automorphic representations that are *finite*[34]. These are the automorphic representations that ought to give nonabelian class field theory, surely the ultimate motivation. As we have already observed, they are to be considered the objects at the heart of Beyond Endoscopy, and for which the proof of Functoriality would presumably be the deepest. But they also represent Artin motives, the most basic of motives. The Reciprocity Conjecture in this case becomes the corresponding case of Functoriality. Since Beyond Endoscopy is supposed to lead to a general proof of Functoriality, it would be surprising if it did not have some major role in the proof of Reciprocity.

It is not hard to think of other analogies between Functoriality and Reciprocity that might be hinting at a common proof. My point is simply this. If the intuition

---

[34] By this I mean the analogue for any $G$ over $F$ of classical modular forms of weight 0 or 1, sometimes called modular forms of type $A_{00}$. In general, they correspond to homomorphisms of the global Weil group $W_F$ to the $L$-group ${}^L G_F$ that factor through the quotient $\Gamma_F$, or equivalently, whose image in $\widehat{G}$ is finite.

we are taking is valid, we would be well advised to formulate arguments in terms
of the stable trace formula. For we would have to come to the problems with a deep
understanding of the automorphic properties of general Shimura varieties. As we
have seen, these are almost as closely tied to the stable trace formula as are the
automorphic representations themselves.

As a final thought, let me consider a different way of viewing the original ques-
tion. It is a philosophical query, related to the duality between automorphic forms
and automorphic representations. Automorphic forms go back to the end of the
nineteenth century with Poincaré, followed by the successive generalizations of
Hilbert, Siegel and Harish-Chandra. The notion of an automorphic representation
came much later, after the adelic language had become common. As we have noted,
the term itself seems to have first appeared in Borel's 1976 Bourbaki lecture. But
it was Langlands who emphasized the dichotomy between the two notions. As we
have seen, this is expressed in the Corvallis papers [35], [145]. The two formulas
we are considering reflect this dichotomy. Their parallel origins are clearly viewed
in the following simple diagram

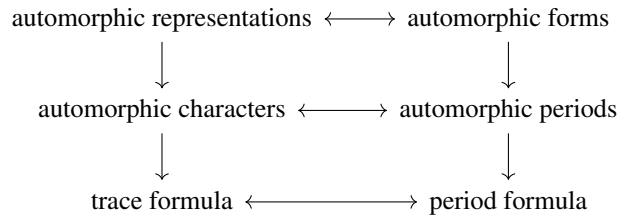

Do the two formulas (or classes of formulas) play dual roles? Or is the diagram
a red herring? Do they give information that is sometimes complementary, or will
they ultimately reduce to the same identities, whatever the circumstances? The two
formulas are in any event sufficiently different that they could both be applied sep-
arately to Beyond Endoscopy, and then compared. It would not matter if the results
overlap. In fact, it would be very useful to have a clear understanding of their com-
mon properties. This would give us a broader perspective for when we run into
difficult problems that demand new ideas, as we most surely will!